%% file: main.tex
\documentclass[12pt]{article}
\usepackage[utf8]{inputenc}
\usepackage{amsthm,amsmath,amssymb,graphicx,enumerate,bbm,dsfont,subcaption,authblk}
\usepackage{natbib}
\usepackage[colorlinks=true,linkcolor=violet,citecolor=violet]{hyperref}\usepackage{tikz,graphicx,multicol}
\usepackage[ruled,vlined]{algorithm2e}
\usepackage{geometry}
\allowdisplaybreaks

\title{Hypothesis testing for eigenspaces of covariance matrix\thanks{Research supported by the NSF grants DMS-1662139 and DMS-1712591,  and the NIH grant 2R01-GM072611-14.}}

\date{}

\input Definitions.tex

\author[]{Igor Silin\thanks{Corresponding author. E-mail: {isilin@princeton.edu}} }
\author[]{Jianqing Fan}
\affil[]{Department of Operations Research and Financial Engineering,\\ Princeton University}

\begin{document}

\maketitle
\begin{abstract}
    Eigenspaces of covariance matrices play an important role in statistical machine learning, arising in variety of modern algorithms. Quantitatively, it is convenient to describe the eigenspaces in terms of spectral projectors. This work focuses on hypothesis testing for the spectral projectors, both in one- and two-sample scenario.
    We present new tests, based on a specific matrix norm developed in order to utilize the structure of the spectral projectors.
    A new resampling technique of independent interest is introduced and analyzed: it serves as an alternative to the well-known multiplier bootstrap, significantly reducing computational complexity of bootstrap-based methods.
    We provide theoretical guarantees for the type-I error of our procedures, which remarkably improve the previously obtained results in the field. Moreover, we analyze power of our tests.
    Numerical experiments illustrate good performance of the proposed methods compared to previously developed ones.
\end{abstract}


\input Content-arxiv.tex

\bibliographystyle{apalike}

\end{document}

%% file: Definitions.tex
\newtheorem{theorem}{Theorem}[section]
\newtheorem{corollary}[theorem]{Corollary}
\newtheorem{lemma}[theorem]{Lemma}
\newtheorem{proposition}[theorem]{Proposition}
\newtheorem{assumption}{Assumption}[section]
\newtheorem{definition}{Definition}[section]
\newtheorem{remark}[]{Remark}
\newtheorem{example}[]{Example}
\numberwithin{equation}{section}      
\numberwithin{remark}{section}
\numberwithin{example}{section}

\DeclareMathOperator{\eqdef}{\;\stackrel{\operatorname{def}}{=}\;}

\DeclareMathOperator{\Prob}{\mathds{P}}
\DeclareMathOperator{\E}{\mathds{E}}
\DeclareMathOperator{\Var}{\mathds{V}ar}
\DeclareMathOperator{\Cov}{\mathds{C}ov}

\DeclareMathOperator{\R}{\mathds{R}}

\newcommand{\Id}{\mathbf{I}}
\newcommand{\Oo}{\mathbf{O}}
\DeclareMathOperator{\Tr}{\mathsf{Tr}}
\DeclareMathOperator{\rank}{\mathsf{rank}}
\DeclareMathOperator{\supp}{\mathsf{supp}}
\DeclareMathOperator{\cspan}{\mathsf{span}}
\DeclareMathOperator{\Fr}{\mathsf{F}}
\DeclareMathOperator{\mymax}{\mathsf{max}}
\DeclareMathOperator{\CONST}{\mathsf{C}}
\newcommand{\T}{\top}

\newcommand{\Sph}{\mathcal{S}}

\newcommand{\eps}{\varepsilon}

\newcommand{\Xdata}{X_1, \ldots, X_n}
\newcommand{\data}{\pmb{\mathbb{X}}}
\newcommand{\Xdataa}{X_1^a, \ldots, X_{2n_a}^a}
\newcommand{\dataa}{\pmb{\mathbb{X}^a}}
\newcommand{\Xdatab}{X_1^b, \ldots, X_{2n_b}^b}
\newcommand{\datab}{\pmb{\mathbb{X}^b}}

\DeclareMathOperator{\relr}{\mathbf{r}}
\DeclareMathOperator{\rr}{\mathsf{d}}
\DeclareMathOperator{\rre}{\overline{\mathsf{d}}}

\DeclareMathOperator{\Eb}{\bar{E}}

\newcommand{\St}{\mathbf{\Sigma}}
\newcommand{\Se}{\mathbf{\widehat{\Sigma}}}
\newcommand{\Su}{\mathbf{\Sigma}}
\newcommand{\So}{\mathbf{\overline{\Sigma}}}
\newcommand{\Sb}{\mathbf{\Sigma}^B}
\newcommand{\Sba}{\mathbf{\Sigma}^B_a}
\newcommand{\Sbb}{\mathbf{\Sigma}^B_b}
\newcommand{\Sf}{\mathbf{\Sigma}^F}
\newcommand{\Sfa}{\mathbf{\Sigma}^F_a}
\newcommand{\Sfb}{\mathbf{\Sigma}^F_b}
\newcommand{\Sp}{\mathbf{\widetilde{\Sigma}}}
\newcommand{\J}{\mathcal{J}}
\newcommand{\Pp}{\mathbf{P}}
\newcommand{\Pt}{\mathbf{P_\J}}
\newcommand{\Ptc}{\mathbf{P_{\J^c}}}
\newcommand{\Pta}{\mathbf{P}_a}
\newcommand{\Ptb}{\mathbf{P}_b}
\newcommand{\Pb}{\mathbf{P}_\J^B}
\newcommand{\Pf}{\mathbf{P}_\J^F}
\newcommand{\Ph}{\mathbf{P^{\circ}}}
\newcommand{\Ptr}{\mathbf{P_r}}
\newcommand{\Pts}{\mathbf{P_s}}
\newcommand{\Ptt}{\mathbf{P_t}}
\newcommand{\Ptrr}{\mathbf{P_{\geq r_0}}}
\newcommand{\Ptrp}{\mathbf{P_{r^{\prime}}}}
\newcommand{\Ptsp}{\mathbf{P_{s^{\prime}}}}
\newcommand{\Pe}{\mathbf{\widehat{P}_\J}}
\newcommand{\Pea}{\mathbf{\widehat{P}}_a}
\newcommand{\Peb}{\mathbf{\widehat{P}}_b}
\newcommand{\Poa}{\mathbf{\overline{P}}_a}
\newcommand{\Pob}{\mathbf{\overline{P}}_b}
\newcommand{\Pba}{\mathbf{P}_{a}^B}
\newcommand{\Pbb}{\mathbf{P}_{b}^B}
\newcommand{\Pfa}{\mathbf{P}_{a}^F}
\newcommand{\Pfb}{\mathbf{P}_{b}^F}
\newcommand{\Per}{\mathbf{\widehat{P}_r}}

\newcommand{\Ppp}{\mathbf{\widetilde{P}_\J}}
\newcommand{\Ppr}{\mathbf{\widetilde{P}_r}}
\newcommand{\Pstar}{\mathbf{P^*}}
\newcommand{\Pbar}{\mathbf{\overline{P}}}

\newcommand{\Resr}{\mathbf{R_r}}
\newcommand{\Pert}{\mathbf{E}}

\newcommand{\Qt}{\mathcal{Q}^{(1)}}
\newcommand{\Qb}{\mathcal{Q}^{(1)}_B}
\newcommand{\Qf}{\mathcal{Q}^{(1)}_F}
\newcommand{\Qr}{\mathcal{Q}^{(1)}_R}

\newcommand{\Qtab}{\mathcal{Q}^{(2)}}
\newcommand{\Qbab}{\mathcal{Q}^{(2)}_B}
\newcommand{\Qfab}{\mathcal{Q}^{(2)}_F}
\newcommand{\Qrab}{\mathcal{Q}^{(2)}_R}

\newcommand{\Qtsp}{\widetilde{\mathcal{Q}}^{(1)}}
\newcommand{\Qbsp}{\widetilde{\mathcal{Q}}^{(1)}_B}
\newcommand{\Qfsp}{\widetilde{\mathcal{Q}}^{(1)}_F}

\newcommand{\Qtspab}{\widetilde{\mathcal{Q}}^{(2)}}
\newcommand{\Qbspab}{\widetilde{\mathcal{Q}}^{(2)}_B}
\newcommand{\Qfspab}{\widetilde{\mathcal{Q}}^{(2)}_F}

\newcommand{\clow}{\underline{\kappa}}
\newcommand{\chigh}{\overline{\kappa}}
\newcommand{\cond}{\kappa}
\newcommand{\condab}{\kappa_{a,b}}
\newcommand{\clowab}{\underline{\kappa}_{a,b}}
\newcommand{\chighab}{\overline{\kappa}_{a,b}}

\newcommand{\scc}{\psi_{n}} 
\newcommand{\scca}{\psi_{n_a}} 
\newcommand{\sccb}{\psi_{n_b}} 
\newcommand{\sccab}{\psi_{n_a\land n_b}} 
\newcommand{\bcc}{\widetilde{\psi}_{n}} 
\newcommand{\bccab}{\widetilde{\psi}_{n_a\land n_b}} 
\newcommand{\fcc}{\widetilde{\psi}_{n}} 
\newcommand{\fccab}{\widetilde{\psi}_{n_a\land n_b}} 

\newcommand{\bcd}{\Delta_B} 
\newcommand{\bcdab}{\Delta_B^{a,b}} 
\newcommand{\fcd}{\Delta_F} 
\newcommand{\fcdab}{\Delta_F^{a,b}} 

\newcommand{\myset}{\mathcal{D}}

\newcommand{\AAA}{\mathbf{A}}
\newcommand{\BB}{\mathbf{B}}
\newcommand{\FF}{\mathbf{F}}

\newcommand{\ff}{\mathbf{f}}
\newcommand{\ixi}{\pmb{\xi}}
\newcommand{\iXi}{\pmb{\Xi}}
\newcommand{\HH}{\mathbf{H}}

%% file: Content-arxiv.tex
\section{Introduction}
    \input source/Introduction.tex

\input source/Structure-arxiv.tex

\section{Setup and statistical problem} \label{Setup}
    \input source/Setup.tex

\section{Testing procedure} \label{Test}
    \input source/Test.tex

\section{Theoretical properties} \label{Theory}
    \input source/Theory.tex

\section{Application to Factor Models} \label{S:FM}
    \input source/FM.tex

\section{Simulation studies} \label{Numerical}
    \input source/Numerical.tex

\section{Discussion}  \label{Discussion}

\input source/Discussion.tex

\newpage
\section{Main proofs} \label{Proofs}
    \input source/Proof.tex

\newpage
\appendix

    \section{Auxiliary results from literature} \label{App:A}
        \input source/AuxLiterature.tex

    \section{Auxiliary proofs} \label{S:AuxProofs}
        \input source/proofs/Proof_TechLemmas.tex

%% file: source/Introduction.tex
\subsection{Background}
We consider a traditional statistical scenario, where we observe $n$ i.i.d. zero-mean random vectors $X_1, \ldots, X_n$ in dimension $d$. Let $X$ be a generic random vector with the same distribution. The geometric structure of the data is described by the covariance matrix
\begin{equation}
    \begin{aligned}
        \St = \E\left[ XX^\T\right].
    \end{aligned}
    \nonumber
\end{equation}
The simplest estimator of $\St$ is the sample covariance matrix
\begin{equation}
    \begin{aligned}
        \Se = \frac{1}{n}\sum\limits_{i=1}^n X_i X_i^\T.
    \end{aligned}
    \nonumber
\end{equation}
The covariance matrix estimation is one of the fundamental problems in statistics: it extends far beyond the sample covariance matrix and has been very well studied under various structural assumptions and different robustification techniques. Some representative works over the past decade include \citet{C1, C2, C3, C4, C5, Koltchinskii_CIAMBFSCO, C7, C6}, among many others. Problem of hypothesis testing for covariance matrix was considered in \cite{Cai}.

However, in order to develop successful methods for modern machine learning problems, one has to go further then just covariance matrices. In particular, eigenstructure of the covariance matrix contains a lot of meaningful information:
\begin{itemize}
    \item In dimension reduction, Principal Component Analysis (PCA) (\cite{Pearson}) projects given high-dimensional observations onto low-dimensional subspace spanned by some number of the leading eigenvectors.
    \item Factor Models (e.g. \cite{FM1,FM2,FM3,FM4}), surprisingly closely related to PCA (see \cite{POET}), also make use of the eigenstructure of the covariance matrix to estimate underlying factors and loadings.
    \item Spectral methods in clustering and community detection (\cite{vonLuxburg}) rely on the eigenvectors of specifically constructed Laplacian matrix (which in some cases can be modelled as covariance matrix).
\end{itemize}
(See \cite{Fan_PCA} for the exposition of problems that can be approached with Spectral/PCA-based techniques.)
To that end, a careful statistical analysis is required for the eigenvectors, or, more generally, for the spectral projector of the covariance matrix $\St$:
\begin{equation}
    \begin{aligned}
        \Pt = \sum\limits_{k \in \mathcal{I}_\J} u_k u_k^\T,
    \nonumber
    \end{aligned}
\end{equation}
where $\{ u_k \}_{k=1}^d$ is an orthonormal basis of ordered eigenvectors of $\St$, $\J$ specifies the set of eigenspaces of interest and the set $\mathcal{I}_\J$ consists of the indices of the respective eigenvectors. Its empirical version $\Pe$ is computed from $\Se$.
The reason why we focus on the spectral projectors rather than working directly with the eigenvectors is that there is always an ambiguity in eigenvectors, while spectral projectors are in one-to-one correspondence with the subspace spanned by the eigenvectors, which is really what plays a role.
Together with the mentioned progress on the covariance matrix estimation, the prominent Davis-Kahan inequality (\cite{DK}) makes the question of statistical estimation of the true spectral projector relatively easy. In contrast, statistical inference (uncertainty quantification, hypothesis testing and confidence sets) for eigenspaces, or in particular for principal components, is significantly less studied but longstanding problem.

\cite{Anderson} was the first paper to study asymptotic distribution of an eigenvector of the sample covariance matrix $\Se$. It proposes the following asymptotically $\chi^2_{d-1}$-distributed statistic to test whether the $k$-th eigenvector $u_k$ of $\St$ is equal (up to a sign) to some specified unit vector $u^\circ$:
\begin{equation}
    \begin{aligned}
        n \left(\lambda_k(\Se) \,{u^\circ}^\T \Se^{-1} u^\circ + {u^\circ}^\T \Se u^\circ/\lambda_k(\Se)  -2\right),
    \nonumber
    \end{aligned}
\end{equation}
where $\lambda_k(\Se)$ denotes the $k$-th eigenvalue of $\Se$. Le Cam's asymptotic theory was utilized in \cite{LAN1} to derive a test for the same problem in case of elliptical distributions, while \cite{LAN2} studied the test from \cite{LAN1} even further in the regime where the spectral gap vanishes. Some other asymptotic results for subspaces spanned by eigenvectors are derived in \cite{Tyler1, Tyler2}.

The two-sample problem also has a long history dating several decades back.
A descriptive technique for comparison of principal components of two or more groups was discussed in \cite{Krzanowski1}. Accompanying empirical results were presented in \cite{Krzanowski2}. A more theoretically justified approach was suggested by \cite{Schott_1}, which considers the test statistic
\begin{equation}
    \begin{aligned}
        \sum\limits_{k=1}^m \left[ \lambda_k(\Se_a) + \lambda_k(\Se_a) - \lambda_k(\Se_a +\Se_b) \right],
    \nonumber
    \end{aligned}
\end{equation}
where $\Se_a$ and $\Se_b$ are the sample covariance matrices of two samples $X_1^a, \ldots, X_{n_a}^a$ and $X_1^a, \ldots, X_{n_b}^b$, respectively.
It is proven that the limiting distribution of this test statistic under null is generalized $\chi^2$. A more sophisticated, again asymptotically $\chi^2$, test statistic was developed in \cite{Schott_2}. Furthermore, \cite{Fujioka} proposed a method, based on the trace of the specific matrix:
\begin{equation}
    \begin{aligned}
        \Tr\left[ {U_2^{(a)}}^\T {U_1^{(b)}} {U_1^{(b)}}^\T {U_2^{(a)}}\right],
    \nonumber
    \end{aligned}
\end{equation}
where $U_1^{(b)} = [u_1^{(b)}, \ldots, u_m^{(b)}]$ consists of the leading $m$ eigenvectors of $\Se_b$ and $U_2^{(a)} = [u_{m+1}^{(a)}, \ldots, u_d^{(a)}]$ consists of the last $(d-m)$ eigenvectors of $\Se_a$.

The above methods are asymptotic and are valid only for a fixed dimension $d$ and a sample size $n$ growing to infinity.
A new line of research in this area was initiated by
\cite{Koltchinskii_NAACOSPOSC}, which obtained the normal approximation for the squared Frobenius distance $\| \Pe - \Pt\|_{\Fr}^2$, providing finite sample error bounds for Kolmogorov distance.
However, this result could not be used directly for the statistical inference as the mean and the variance of the normal distribution approximating $\| \Pe - \Pt \|_{\Fr}^2$ depend on the true unknown $\St$; the idea of splitting the sample into three parts to estimate mean and variance was just mentioned.
A follow-up paper of \cite{Koltchinskii_NARPCA} formalized this sample splitting idea, and derived completely data-driven test statistic with known approximating distribution.
Another approximation was proposed in \cite{Naumov}. The focus of that paper is on constructing confidence sets for the true spectral projector, so the multiplier bootstrap was employed to deal with unknown parameters of the limiting distribution.
\cite{Silin_1} proposed to use Bayesian inference instead of bootstrap, at the same time extending the results from \cite{Naumov} to non-Gaussian data.
Even though these works did not pose the hypothesis testing problem, it is straightforward to develop one-sample tests based on their results.

\subsection{Contributions}
The aim of this work is to develop statistical procedures for testing one- and two-sample hypotheses about underlying eigenspaces of covariance matrices.

We try to address the following challenges:
\begin{itemize}
    \item High-dimensionality. Rates obtained in the previous works require $d^3 \ll n$ (except the cases of Gaussian data with small effective rank), which significantly restricts the applicability of the proposed methods.
    \item Heavy-tailed data. While most of the discussed literature focuses only on the Gaussian data, the ability to move beyond Gaussian or sub-Gaussian distributions is crucial for modern applications, especially in finance.
    \item Computational complexity. Over the past decades, the multiplier bootstrap (also called wild bootstrap) has been one of the main tools for statistical inference. However, to generate one bootstrap sample, a statistician needs to perform $O(n)$ operations (to generate $n$ bootstrap weights), which can lead to intractable running time of a bootstrap-based procedure.
\end{itemize}
Our main contributions can be summarized as follows:
\begin{itemize}
    \item We develop new statistical procedures for one- and two-sample hypothesis testing for eigenspaces of covariance matrix. The tests are based on newly developed matrix norm, which is designed by taking the structure of spectral projectors into account. In fact, we develop a family of tests that provides flexibility in controlling the trade-off between the closeness to the desired type-I error and power.
    \item We propose a new resampling technique of independent interest, which can be considered as an alternative to the multiplier bootstrap. While possessing the same statistical properties, this technique reduces computational complexity by $n$ times compared to the bootstrap.
    \item Our theoretical results for the presented procedures include both validity guarantees (type-I error close to the desired level) as well as power analysis (probability of rejection of null hypothesis goes to one under alternative). The results do not rely on the Gaussianity or sub-Gaussianity of the data. Moreover, we demonstrate a significant progress in obtaining dimension-free bounds for this problem. In some setups (e.g. Factor Models) the dependence on $d$ is remarkably improved compared to the previous works.
    \item The numerical study confirms good properties of our algorithms. The proposed procedures outperform the variety of previously developed methods in a wide diversity of settings.
\end{itemize}

%% file: source/Structure-arxiv.tex
\subsection{Structure of the paper}
The paper is organized as follows.
We conclude the introduction with defining necessary notations in Subsection~\ref{Notations}.
The general framework and problem formulation are presented in Section~\ref{Setup}. The proposed testing procedures are described in Section~\ref{Test}. Their theoretical properties are analyzed in Section~\ref{Theory}.
In Section~\ref{S:FM} we apply the developed methods to Factor Models. Section~\ref{Numerical} provides some numerical simulations.
The comparison with other works is presented in Section~\ref{Discussion}. Section~\ref{Proofs} is devoted to the main proofs. Finally, Appendix~\ref{App:A} and Appendix~\ref{S:AuxProofs} gather auxiliary results and proofs, respectively.

\subsection{Notations} \label{Notations}

The following notations are used throughout the work.
For positive integers $k$ and $l$, we write $[k]$ as shorthand for the set $\{ 1, 2, \ldots, k\}$ and $[k:l]$ for $\{ k, k+1, \ldots, l\}$.
The space of $k$-dimensional real-valued vectors is denoted by $\R^k$.
The space of real-valued matrices of size $k \times l$ is denoted by $\R^{k\times l}$.
We use $0_k$ for the zero vector in $\R^k$, $\Oo_{k\times l}$ for $k\times l$ matrix of zeros and $\Id_k$ for the identity matrix of size $k \times k$.
For a matrix $A$, we denote by $A_{[i:j],[k:l]}$ its submatrix formed by intersection of rows $\{i, i+1, \ldots, j\}$ and columns $\{ k, k+1, \ldots, l \}$. Let $\supp[x]$ denote the support of a vector $x$.

For a vector $x\in\R^k$, $\| x\|$ denotes its Euclidean norm.
By $S^{k-1}$ we denote unit sphere in $\R^k$.
For a matrix $A\in\R^{k\times l}$, notations $\| A \|$, $\| A \|_{\Fr}$ and $\| A \|_*$ mean spectral norm (largest singular value), Frobenius norm (square root of sum of squared singular values) and nuclear norm (sum of singular values), respectively, while $\| A \|_{\mymax}$ denotes maximal absolute elementwise norm.
$\Tr[\cdot]$ and $\rank[\cdot]$ stand for trace and rank.

For two real numbers $a$ and $b$, by $a\lor b$ and $a\land b$ we mean their maximum and minimum, respectively.
The relation $ a \lesssim b $ means that there exists an absolute constant $ C $, different from place to place, such that $ a \leq Cb $, while
$ a \asymp b $ means that $ a \lesssim b $ and $ b \lesssim a $. When this constant has a subscript or argument, i.e. $C_\gamma$ or $C(\gamma)$, it specifies that this constant may be different for different values of variable $\gamma$, but does not depend on anything else. 

%% file: source/Setup.tex
\subsection{Setup}
Let $\Xdata$ be i.i.d. mean zero random vectors in $\R^d$ and $X$ be a generic random vector from the same distribution. We store the observed data in a matrix
\begin{equation} \begin{aligned}
    \data = [\Xdata] \in \R^{d \times n}.
\nonumber \end{aligned} \end{equation}
The covariance matrix of the data is
\begin{equation} \begin{aligned}
    \St = \Cov[X] = \E\left[XX^{\T}\right] \in \R^{d\times d}.
\nonumber \end{aligned} \end{equation}
Typically, $\St$ is unknown, and one estimates it using its sample version $\Se$:
\begin{equation} \begin{aligned}
    \Se = \frac{1}{n} \sum\limits_{i=1}^n X_i X_i^{\T} \in \R^{d\times d}.
\nonumber \end{aligned} \end{equation}

Let us introduce some notations. Let $ \sigma_{1} \geq \ldots \geq \sigma_{d} $ be the ordered eigenvalues of $ \St $ (assume all eigenvalues are strictly positive).
Suppose that among them there are $ q $ distinct eigenvalues $  \mu_{1} > \ldots > \mu_q  $.
Introduce groups of indices $ \mathcal{I}_{r} = \{ j\in[d]: \mu_{r} = \sigma_{j}\} $ and denote by $ m_{r} $ the multiplicity factor $ |\mathcal{I}_{r}| $  for all $ r \in [q] $.
The corresponding eigenvectors are denoted as $ u_{1}, \ldots, u_{d} $.  Define projector on $r$-th eigenspace as $\Ptr = \sum\limits_{k \in \mathcal{I}_r} u_k u_k^\T$ for $r\in[q]$.
Similarly, suppose that $ \Se $ has $ d $ eigenvalues $ \widehat{\sigma}_{1} > \ldots > \widehat{\sigma}_{d} $ (distinct with probability one).
The corresponding eigenvectors are $ \widehat{u}_{1},\ldots, \widehat{u}_{d} $.

Suppose we are interested in the sum of some of the $q$ eigenspaces of $\St$.
In particular, let
\begin{equation} \begin{aligned}[c]
        \J = \{ r_1,\; r_1 + 1, \; \ldots,\; r_2\}
\nonumber \end{aligned} \end{equation}
be a set of consecutive indices of eigenspaces of interest. Define also
\begin{equation}
    \begin{aligned}[c]
        \mathcal{I}_{\J} = \bigcup\limits_{r\in\J} \mathcal{I}_r.
    \nonumber
    \end{aligned}
\end{equation}
Quantitatively, sum of the $\J$ eigenspaces of $\St$ is described by the projector onto this subspace, defined as
\begin{equation} \begin{aligned}
    \Pt \eqdef \sum_{r \in \J} \Ptr = \sum_{r \in \J}
    \sum_{k \in \mathcal{I}_r} u_k u_k^\T = \sum_{k \in \mathcal{I}_\J} u_k u_k^\T \in \R^{d\times d}.
\nonumber \end{aligned} \end{equation}
Its empirical counterpart is given by
\begin{equation} \begin{aligned}
    \Pe \eqdef \sum_{r \in \J} \Per = \sum_{r \in \J}
    \sum_{k \in \mathcal{I}_r} \widehat{u}_k \widehat{u}_k^{\T} = \sum_{k \in \mathcal{I}_\J} \widehat{u}_k \widehat{u}_k^{\T} \in \R^{d\times d}.
\nonumber \end{aligned} \end{equation}
The rank of these projectors is $m \eqdef |\mathcal{I}_\J| = \sum\limits_{r\in\J} m_r$. As an example, when $ \mathcal{I}_\J = \{ 1, \ldots, m\} $, then $\Pt$ consists of the projector onto the eigenspace spanned by the eigenvectors of the top $m$ distinguished eigenvalues, while $ \Pe $ is its empirical counterpart. For brevity, we occasionally will be using the notation $\Ptc \eqdef \Id_d - \Pt$.

\subsection{Statistical problem}
One may be interested in testing hypothesis about $\Pt$.
The hypothesis testing problem
\begin{equation}
    \begin{aligned}
        H_0^{(1)}:\;\;\Pt = \Ph\hspace{1.3cm}\text{vs}\hspace{1.3cm}H_1^{(1)}:\;\;\Pt \neq \Ph
    \nonumber
    \end{aligned}
\end{equation}
for a given projector $\Ph$ of rank $m$ is the main focus of our work.

Two-sample problem is also of great interest.
Suppose we have two i.i.d. samples: $X_1^{a}, \ldots, X_{2n_a}^{a}$ and $X_1^{b}, \ldots, X_{2n_b}^{b}$ (it will be clear later why we denote the sizes of the samples as $2n_a$ and $2n_b$; assume for simplicity they are even numbers). As previously, we store them as $\dataa$ and $\datab$. Let the true covariance matrix of the first sample be $\St_{a}$ and the true covariance matrix of the second sample be $\St_{b}$. Let $\J_a$ be a set of consecutive indices of eigenspaces of $\St_a$ and $\J_b$ be a set of consecutive indices of eigenspaces of $\St_b$. Sets $\mathcal{I}_{\J_a}$ and $\mathcal{I}_{\J_b}$ contain the indices of the associated ordered eigenvectors; it makes sense to require $|\mathcal{I}_{\J_a}| = |\mathcal{I}_{\J_b}| = m$ (so in both one- and two-sample problems $m$ denotes the dimension of the subspace being tested). We denote by $\Pta$ and $\Ptb$ the corresponding projectors of rank $m$:
\begin{equation}
    \begin{aligned}
        &\Pta = \sum\limits_{k\in \mathcal{I}_{\J_a}} u_k^a {u_k^a}^\T,\\
        &\Ptb = \sum\limits_{k\in \mathcal{I}_{\J_b}} u_k^b {u_k^b}^\T,
    \nonumber
    \end{aligned}
\end{equation}
where $\{ u_k^a \}_{k=1}^d$ and $\{ u_k^b \}_{k=1}^d$ are the sets of ordered (w.r.t. the associated eigenvalues) eigenvectors of $\St_a$ and $\St_b$.
Here, in order to avoid excessive sub- and superscripts, we slightly abuse the notation: $\Pta$ and $\Ptb$ should not be confused with $\Ptr$ for $r\in[q]$ from one-sample case.
In addition to the one-sample problem stated above, we will propose a method for the following hypothesis testing problem
\begin{equation}
    \begin{aligned}
        H_0^{(2)}:\;\;\Pta=\Ptb\hspace{1.3cm}\text{vs}\hspace{1.3cm}H_1^{(2)}:\;\;\Pta \neq \Ptb\,.
    \nonumber
    \end{aligned}
\end{equation}

In both of these problems, a statistician is often given a desired level of the test $\alpha$.
However, in most of the situations (including our setting), creating a reasonable test with type-I error exactly $\alpha$ is difficult or impossible.
Our goal is to develop tests, whose type-I errors will be close to the level $\alpha$, and provide finite sample guarantees for the discrepancy between them.

%% file: source/Test.tex
Previous works of \cite{Koltchinskii_NAACOSPOSC,Koltchinskii_NARPCA,Naumov,Silin_1} considered the Frobenius norm $\sqrt{n} \| \Pe - \Pt \|_{\Fr}$ and would suggest this object as a basis for one-sample testing procedure. Another interesting random quantity to analyze would be the spectral norm
    \begin{equation}
        \begin{aligned}
            \Qtsp \eqdef \sqrt{n} \| \Pe - \Pt \|.
        \nonumber
        \end{aligned}
    \end{equation}
(Here and further the superscript specifies whether we are in context of one-sample or two-sample problem.)
However, as we will see, current techniques doesn't allow us to obtain an approximation to the distribution of $\Qtsp$ that is accurate in high dimensions. This prevents us from developing a test based on this random quantity, and forces us to construct a new, less conventional and more problem-specific, matrix norm that will have better theoretical properties.

The matrix norm, which our test statistic will be based on, is introduced in the following definition.

\begin{definition} \label{D:norm}
Let $\Pp \in \R^{d\times d}$ be a projector of rank $m$.
Fix $\Gamma = [\Gamma_1 \;\Gamma_2] \in \R^{d\times d}$ with $\Gamma_1\in\R^{d\times m} , \Gamma_2 \in \R^{d\times(d-m)}$ satisfying
\begin{equation}
    \begin{aligned}
        &\Gamma_1 \Gamma_1^\T = \Pp, \;\;\Gamma_1^\T \Gamma_1 = \Id_m,\\
        &\Gamma_2 \Gamma_2^\T = \Id_d - \Pp, \;\; \Gamma_2^\T \Gamma_2 = \Id_{d-m}.
    \label{Gamma_Prop}
    \end{aligned}
\end{equation}
Let also $s_1 \in[m] $ and $s_2 \in [d-m]$.
Then, for any symmetric matrix $A \in \R^{d\times d}$ define
\begin{equation}
    \begin{aligned}
        \| A \|_{(\Pp, \Gamma, s_1, s_2)} &\eqdef \frac{1}{2}\| \Gamma_1^\T A \Gamma_1 \| + \frac{1}{2}\| \Gamma_2^\T A \Gamma_2 \| +\\
        & \qquad\qquad + \max\limits_{\substack{k \in [m-s_1+1]\\ l \in [d-m-s_2+1]}}
             \left\| [\Gamma_1^\T A \Gamma_2]_{[k:(k+s_1-1)], [l:(l+s_2-1)]} \right\|.
    \nonumber
    \end{aligned}
\end{equation}
\end{definition}

Let us briefly describe the role of $\Gamma_1, \Gamma_2$ and $s_1, s_2$ in the above definition. As can be seen from \eqref{Gamma_Prop}, the columns of $\Gamma_1$ form an orthonormal basis in the subspace associated with $\Pp$, while the columns of $\Gamma_2$ form an orthonormal basis in the orthogonal complement. Thus, $\Gamma_1$ can be found as the set of eigenvectors of $\Pp$ corresponding to the eigenvalue $1$ of multiplicity $m$, and $\Gamma_2$ is the set of eigenvectors of $\Pp$ corresponding to the eigenvalue $0$ of multiplicity $(d-m)$, i.e. the eigendecomposition of $\Pp$  looks like
\begin{equation}
\begin{aligned}
	\Pp = [\Gamma_1 \; \Gamma_2] \begin{bmatrix} \Id_m & \Oo_{m\times(d-m)} \\ \Oo_{(d-m)\times m} & \Oo_{(d-m)\times(d-m)} \end{bmatrix} \begin{bmatrix} \Gamma_1^\T \\ \Gamma_2^\T \end{bmatrix}.
\nonumber
\end{aligned}
\end{equation}
 This rotation is necessary for our future theoretical analysis. Clearly, $\Gamma_1$ and $\Gamma_2$ satisfying \eqref{Gamma_Prop} are not unique, but a specific choice will not play any role in the sequel. 
 The first two terms will be negligible under null hypothesis while allowing us to improve the power of the test (``power enhancement"); the third term is the main term that will give us the desired approximation.
The integers $s_1$ and $s_2$ parametrize the family of norms and give flexibility in the test that we will develop: as we will see, the test based on the norm with $s_1=s_2=1$ will have better guarantees under null hypothesis and weaker power (less omnibus), while taking largest possible values $s_1 = m$, $s_2 = d-m$ yields the test with potentially unstable behaviour under $H_0$ but omnibus.
Figure~\ref{draw1} further explains Definition~\ref{D:norm}.
\input source/draws/draw1.tex

We state some useful properties of this operator in the next proposition.
\begin{proposition}[Properties of $\| \cdot \|_{(\Pp, \Gamma, s_1, s_2)}$] \label{P:Properties0}
Fix arbitrary $\Pp, \Gamma = [\Gamma_1\; \Gamma_2], s_1, s_2$ as in Definition~\ref{D:norm}. Then, the following holds:
\begin{enumerate}[(i)]
 \item $\| \cdot \|_{(\Pp, \Gamma, s_1, s_2)}$ is indeed a norm on the space of symmetric matrices.
 \item $\| \cdot \|_{(\Pp, \Gamma, s_1, s_2)}$ is equivalent to the spectral norm: for any symmetric $A\in\R^{d\times d}$
 \begin{equation}
    \begin{aligned}
        \frac{1}{2} \sqrt{\frac{s_1}{m} \cdot\frac{s_2}{d-m}}\cdot\| A \| \leq \| A \|_{(\Pp, \Gamma, s_1, s_2)} \leq 2\| A \|.
    \nonumber
    \end{aligned}
\end{equation}
\end{enumerate}
\end{proposition}

\subsection{One-sample test}
Our one-sample test is based on the following random quantity
\begin{equation}
    \begin{aligned}
        \Qt \eqdef \sqrt{n} \| \Pe - \Ph \|_{(\Ph, \Gamma^\circ, s_1, s_2)},
    \nonumber
    \end{aligned}
\end{equation}
where $\Gamma^\circ$ satisfying properties \eqref{Gamma_Prop} for $\Ph$ (as in Definition~\ref{D:norm}) is chosen arbitrarily.

\begin{remark}[Link between $\Qt$ and $\Qtsp$]
    Under $H_0^{(1)}$ it holds $\Ph = \Pt$, and the random quantity of interest becomes
    $\Qt = \sqrt{n} \| \Pe - \Pt \|_{(\Pt, \Gamma^\circ, s_1, s_2)}$. Take $s_1 = m$, $s_2 = d-m$. Note that even in this case,
    \begin{equation}
        \begin{aligned}
             &\| \Pe - \Pt \|_{(\Pt, \Gamma^\circ, m, d-m)} \neq \| \Pe - \Pt \|,
        \nonumber
        \end{aligned}
    \end{equation}
    though we have bounds as in Proposition~\ref{P:Properties0}.
    However, as will be seen in the proofs, due to a specific structure of spectral projectors, it holds
    \begin{equation}
        \begin{aligned}\| \Pe - \Pt \|_{(\Pt, \Gamma^\circ, m, d-m)} \approx \| \Pe-\Pt\|
        \nonumber
        \end{aligned}
    \end{equation}
    up to higher-order terms with high probability,
    and, moreover, $\Qtsp = \sqrt{n}\| \Pe - \Ph \|$ can also be used as the test statistics with the same theoretical guarantees under null hypothesis as for $\sqrt{n}\| \Pe - \Ph \|_{(\Ph, \Gamma^\circ, m, d-m)}$.
\end{remark}

If we knew the quantiles $q^{(1)}(\alpha)$ of the distribution of $\Qt$ under $H_0^{(1)}$, we would use the following test
\begin{equation}
    \begin{aligned}
        \phi_\alpha(\data) = \mathbbm{1}\left\{\sqrt{n} \| \Pe - \Ph \|_{(\Ph, \Gamma^\circ, s_1, s_2)} \geq q^{(1)}(\alpha)\right\},
    \nonumber
    \end{aligned}
\end{equation}
which has type-I error exactly $\alpha$. However, in practice the distribution of $\Qt$ is unavailable to us, since even if we could obtain closed-from approximation to it, it would depend heavily on the underlying unknown covariance $\St$. Hence, we suggest two approaches to approximate $q^{(1)}(\alpha)$.

\subsubsection*{Approach 1: Bootstrap-based test}
Let us apply the idea of multiplier bootstrap to approximate the unknown distribution of $\Qt$ under $H_0^{(1)}$.
Consider $\eta_1, \ldots, \eta_n \stackrel{i.i.d.}{\sim} \mathcal{N}(1,1)$.
Define $\Sb \eqdef \frac{1}{n}\sum\limits_{i=1}^n \eta_i X_i X_i^{\top}$ and the corresponding projector $\Pb$ from $\Sb$. Consider the random quantity
\begin{equation}
    \begin{aligned}
        \Qb \eqdef \sqrt{n} \| \Pb - \Pe \|_{(\Ph, \Gamma^\circ, s_1, s_2)}.
    \nonumber
    \end{aligned}
\end{equation}
The hope is that $(\Qb \,|\, \data) \stackrel{d}{\approx} \Qt$ under $H_0^{(1)}$ with high probability. At the same time, the distribution of $(\Qb\,|\, \data)$ is available to us and can be sampled to find its $\alpha$-quantile $q^{(1)}_B(\alpha)$.

\subsubsection*{Approach 2: Frequentist-Bayes related test}
We also propose another resampling technique to approximate the unknown distribution of $\Qt$ under $H_0^{(1)}$.
Consider $Z_1, \ldots, Z_n \stackrel{i.i.d.}{\sim} \mathcal{N}(0,\Se)$.
Define $\Sf \eqdef \frac{1}{n}\sum\limits_{i=1}^n Z_i Z_i^{\top}$ and the corresponding projector $\Pf$ from $\Sf$. Consider the random quantity
\begin{equation}
    \begin{aligned}
        \Qf \eqdef \sqrt{n} \| \Pf - \Pe \|_{(\Ph, \Gamma^\circ, s_1, s_2)}.
    \nonumber
    \end{aligned}
\end{equation}
Similarly to Approach 1, we expect that $(\Qf \,|\, \data) \stackrel{d}{\approx} \Qt$ under $H_0^{(1)}$ with high probability. Again, the distribution of $(\Qf \,|\, \data)$ is available to us and can be sampled in order to find its $\alpha$-quantile $q^{(1)}_F(\alpha)$. Note that instead of sampling $Z_1, \ldots, Z_n \stackrel{i.i.d.}{\sim} \mathcal{N}(0,\Se)$, we can directly generate $\Sf\sim\frac{1}{n}\cdot \mbox{Wishart}(n, \Se)$, which is more computationally efficient.

\begin{remark}[Relation to Frequentist Bayes] \label{FBConnection}
    One may be curious why we call Approach 2 ``Frequentist-Bayes related''. It turns out, that this resampling method somehow arises from the Bayesian inference conducted in \cite{Silin_1}. Due to space limitations, we do not elaborate on this connection in our work.
\end{remark}

Based on one of the presented resampling strategies, we summarize our test method as in Algorithm~\ref{Algorithm}.

\begin{algorithm}[]
\SetAlgoLined
\vspace{0.2cm}
\textbf{Input: }Data $\data = [\Xdata]$, set $\mathcal{I}_\J$, null hypothesis projector $\Ph$ of rank $|\mathcal{I}_\J|$, \\\hspace{1.55cm}desired level $\alpha$.\\
\textbf{Hyperparameters:} $s_1, s_2$, number of resampling iterations $N$.\\
    Set $m := |\mathcal{I}_\J|$;\\
    Compute $\Se := \frac{1}{n}\sum\limits_{i=1}^n X_i X_i^{\top}$;\\
    Compute the corresponding projector $\Pe$ from $\Se$;\\
    Fix $\Gamma^\circ_1 \in\R^{d\times m}$ such that ${\Gamma^\circ_1}^\T \Gamma^\circ_1 = \Id_m$ and $\Gamma^\circ_1 {\Gamma^\circ_1}^\T = \Ph$;\\
    Fix $\Gamma_2^\circ \in \R^{d\times (d-m)}$ such that ${\Gamma^\circ_2}^\T \Gamma^\circ_2 = \Id_{d-m}$ and $\Gamma^\circ_2 {\Gamma^\circ_2}^\T = \Id_d - \Ph$;\\
    Apply Bootstrap-based resampling:\\
    \For{$k = 1, \ldots, N$}{
        Sample $\eta_1, \ldots, \eta_n \stackrel{i.i.d.}{\sim} \mathcal{N}(1,1)$;\\
        Compute $\Sb := \frac{1}{n}\sum\limits_{i=1}^n \eta_i X_i X_i^{\top}$;\\
        Compute the corresponding projector $\Pb$ from $\Sb$;\\
        Compute $k$-th realization $\Qr(k) := \sqrt{n} \| \Pb - \Pe \|_{(\Ph, \Gamma^\circ, s_1, s_2)}$;
     }
     or apply Frequentist-Bayes related resampling:\\
    \For{$k = 1, \ldots, N$}{
        Sample $\Sf := \frac{1}{n}\cdot Wishart(n, \Se)$;\\
        Compute the corresponding projector $\Pf$ from $\Sf$;\\
        Compute $k$-th realization $\Qr(k) := \sqrt{n} \| \Pf - \Pe \|_{(\Ph, \Gamma^\circ, s_1, s_2)}$;
     }
     Compute $q^{(1)}_R(\alpha) := \alpha$-quantile of $\{\Qr(k) \}_{k=1}^{N}$;\\
 \KwResult{$\phi^R_\alpha(\data) := \mathbbm{1}\{\sqrt{n} \| \Pe - \Ph \|_{(\Ph, \Gamma^\circ, s_1, s_2)} \geq q^{(1)}_R(\alpha)\}$,\\$\hspace{1.62cm}\text{p-value}(\data) := \frac{1}{N}\sum\limits_{k=1}^{N} \mathbbm{1}\{ \sqrt{n} \| \Pe - \Ph \|_{(\Ph, \Gamma^\circ, s_1, s_2)} \geq \Qr(k)\} $. }
 \caption{One-sample testing procedure}
 \label{Algorithm}
\end{algorithm}

\subsection{Two-sample test}
In one-sample problem we have a null hypothesis projector $\Ph$ given to us, and can use it in our test statistic. Specifically, we use $\| \cdot \|_{(\Ph, \Gamma^\circ, s_1, s_2)}$-norm. In contrast, in two-sample problem we have only two samples, while no $\Ph$ is provided, so the one-sample procedure cannot be straightforwardly extended, as it is not clear what norm to use.

To overcome this difficulty, we split each of the samples $\Xdataa$ and $\Xdatab$ into two equal parts.  The second part will be used to learn $\Ph$ and $\Gamma^\circ$ and the first part will be used to construct a test.  More specifically, define
\begin{equation}
    \begin{aligned}
        &\Se_a \eqdef \frac{1}{n_a} \sum\limits_{i=1}^{n_a} X_i^a {X_i^a}^\T,\;\;\;
        \So_a \eqdef \frac{1}{n_a} \sum\limits_{i=n_a+1}^{2n_a} X_i^a {X_i^a}^\T,\\
        &\Se_b \eqdef \frac{1}{n_b} \sum\limits_{i=1}^{n_b} X_i^b {X_i^b}^\T,\;\;\;
        \So_b \eqdef \frac{1}{n_a} \sum\limits_{i=n_b+1}^{2n_b} X_i^b {X_i^b}^\T.
    \nonumber
    \end{aligned}
\end{equation}
Denote by $\Pea, \Poa$ the corresponding projectors of $\Se_a, \So_a$ associated with $\mathcal{I}_{\J_a}$, and by $\Peb, \Pob$ the corresponding projectors of $\Se_b, \So_b$ associated with $\mathcal{I}_{\J_b}$.
Introduce
\begin{equation}
    \begin{aligned}
        \Pbar \eqdef \arg\min\limits_{\substack{\mathbf{P}:\text{ projector,} \\ \rank(\mathbf{P}) = m}} \left\{ \| \mathbf{P} - \Poa \|_{\Fr}^2 + \| \mathbf{P} - \Pob \|_{\Fr}^2\right\}.
    \label{Def:Pbar}
    \end{aligned}
\end{equation}
One can show that it can be easily computed: $\Pbar = \Psi \Psi^\T$, where $\Psi \in\R^{d\times m}$ consists of the eigenvectors of $(\Poa + \Pob)$ associated with $m$ largest eigenvalues. Fix $\overline{\Gamma}$ satisfying properties \eqref{Gamma_Prop} for $\Pbar$ in an arbitrary way.  Define the symmetric (w.r.t. change of sample $a$ and sample $b$) test statistic
\begin{equation}
    \begin{aligned}
        \Qtab \eqdef \sqrt{\frac{n_an_b}{n_a+n_b}} \| \Pea - \Peb \|_{(\Pbar, \overline{\Gamma},s_1, s_2)}.
    \nonumber
    \end{aligned}
\end{equation}

To estimate its distribution, we again employ one of the presented approaches: Bootstrap-based or Frequentist-Bayes related. Both of them are straightforwardly extended from one-sample case and lead to the following random quantities:
\begin{equation}
    \begin{aligned}
        \Qbab = \sqrt{\frac{n_an_b}{n_a+n_b}}\|(\Pba - \Pea)-(\Pbb-\Peb)\|_{(\Pbar, \overline{\Gamma}, s_1, s_2)}
    \nonumber
    \end{aligned}
\end{equation}
and
\begin{equation}
    \begin{aligned}
        \Qfab = \sqrt{\frac{n_an_b}{n_a+n_b}}\|(\Pfa-\Pea)-(\Pfb-\Peb)\|_{(\Pbar,\overline{\Gamma}, s_1, s_2)}.
    \nonumber
    \end{aligned}
\end{equation}
Note that here $\Pba, \Pbb, \Pfa, \Pfb$ correspond only to the first halves of the samples and has nothing to do with the second halves.
Now the hope is that $(\Qbab \,|\, \dataa, \datab) \stackrel{d}{\approx} (\Qtab\,|\,\overline{\Gamma})$ and $(\Qfab \,|\, \dataa,\datab) \stackrel{d}{\approx} (\Qtab\,|\,\overline{\Gamma})$ with high probability.
This brings us directly to Algorithm~\ref{Algorithm2}.

\begin{algorithm}[]
\SetAlgoLined
\vspace{0.2cm}
\textbf{Input: }Data $\dataa = [\Xdataa]$ and $\datab = [\Xdatab]$, \\\hspace{1.55cm}sets $\mathcal{I}^a$ and $\mathcal{I}^b$ of the same size, desired level $\alpha$.\\
\textbf{Hyperparameters:} $s_1, s_2$, number of resampling iterations $N$.\\
    Set $m := |\mathcal{I}^a| = |\mathcal{I}^b|$;\\
    Compute $\Se_a := \frac{1}{n_a}\sum\limits_{i=1}^{n_a} X_i^a {X_i^a}^{\top}$ and $\So_a := \frac{1}{n_a}\sum\limits_{i=n_a+1}^{2n_a} X_i^a {X_i^a}^{\top}$;\\
    Compute $\Se_b := \frac{1}{n_b}\sum\limits_{i=1}^{n_b} X_i^b {X_i^b}^{\top}$ and $\So_b := \frac{1}{n_b}\sum\limits_{i=n_b+1}^{2n_b} X_i^b {X_i^b}^{\top}$;\\
    Compute the corresponding projectors $\Pea$, $\Poa$, $\Peb$ and $\Pob$;\\
    Compute $\Pbar := \arg\min\limits_{\substack{\mathbf{P}:\text{ projector,} \\ \rank(\mathbf{P}) = m}} \left\{ \| \mathbf{P} - \Poa \|_{\Fr}^2 + \| \mathbf{P} - \Pob \|_{\Fr}^2\right\}$ using eigendecomposition;\\
    Fix $\overline{\Gamma}_1 \in\R^{d\times m}$ such that $\overline{\Gamma}_1^\T \overline{\Gamma}_1 = \Id_m$ and $\overline{\Gamma}_1 \overline{\Gamma}_1^\T = \Pbar$;\\
    Fix $\overline{\Gamma}_2 \in \R^{d\times (d-m)}$ such that $\overline{\Gamma}_2^\T \overline{\Gamma}_2 = \Id_{d-m}$ and $\overline{\Gamma}_2 \overline{\Gamma}_2^\T = \Id_d - \Pbar$;\\
    Apply Bootstrap-based resampling:\\
    \For{$k = 1, \ldots, N$}{
        Sample $\eta_1^a, \ldots, \eta_{n_a}^a \stackrel{i.i.d.}{\sim} \mathcal{N}(1,1)$ and $\eta_1^b, \ldots, \eta_{n_b}^b \stackrel{i.i.d.}{\sim} \mathcal{N}(1,1)$;\\
        Compute $\Sba := \frac{1}{n_a}\sum\limits_{i=1}^{n_a} \eta_i^a X_i^a {X_i^a}^{\top}$
        and $\Sbb := \frac{1}{n_b}\sum\limits_{i=1}^{n_b} \eta_i^b X_i^b {X_i^b}^{\top}$;\\
        Compute the corresponding projectors $\Pba$ from $\Sba$ and $\Pbb$ from $\Sbb$;\\
        Compute $k$-th realization $\Qrab(k) := \sqrt{\frac{n_an_b}{n_a+n_b}} \| (\Pba-\Pea) - (\Pbb-\Peb) \|_{(\Pbar, \overline{\Gamma}, s_1, s_2)}$;
     }
     or apply Frequentist-Bayes related resampling:\\
    \For{$k = 1, \ldots, N$}{
        Sample $\Sfa := \frac{1}{n_a}\cdot Wishart(n_a, \Se_a)$ and $\Sfb := \frac{1}{n_b}\cdot Wishart(n_b, \Se_b)$;\\
        Compute the corresponding projectors $\Pfa$ from $\Sfa$
        and $\Pfb$ from $\Sfb$;\\
        Compute $k$-th realization $\Qrab(k) := \sqrt{\frac{n_an_b}{n_a+n_b}} \| (\Pfa-\Pea) - (\Pfb-\Peb) \|_{(\Pbar, \overline{\Gamma}, s_1, s_2)}$;
     }
     Compute $q^{(2)}_R(\alpha) := \alpha$-quantile of $\{\Qrab(k) \}_{k=1}^{N}$;\\
 \KwResult{$\phi^R_\alpha(\dataa; \datab) := \mathbbm{1}\left\{\sqrt{\frac{n_an_b}{n_a+n_b}}\| \Pea - \Peb \|_{(\Pbar, \overline{\Gamma}, s_1, s_2)}  \geq q^{(2)}_R(\alpha)\right\}$,\\
 \hspace{1.62cm}$\text{p-value}(\dataa; \datab) := \frac{1}{N} \cdot\sum\limits_{k=1}^{N}  \mathbbm{1}\left\{  \sqrt{\frac{n_an_b}{n_a+n_b}}\| \Pea - \Peb \|_{(\Pbar, \overline{\Gamma}, s_1, s_2)} \geq \Qrab(k)\right\} $. }
 \caption{Two-sample testing procedure}
 \label{Algorithm2}
\end{algorithm}

%% file: source/draws/draw1.tex
\newcommand{\Rect}[5]{
    \draw[#1] (#2,#3) rectangle(#2+#4,#3-#5);
}

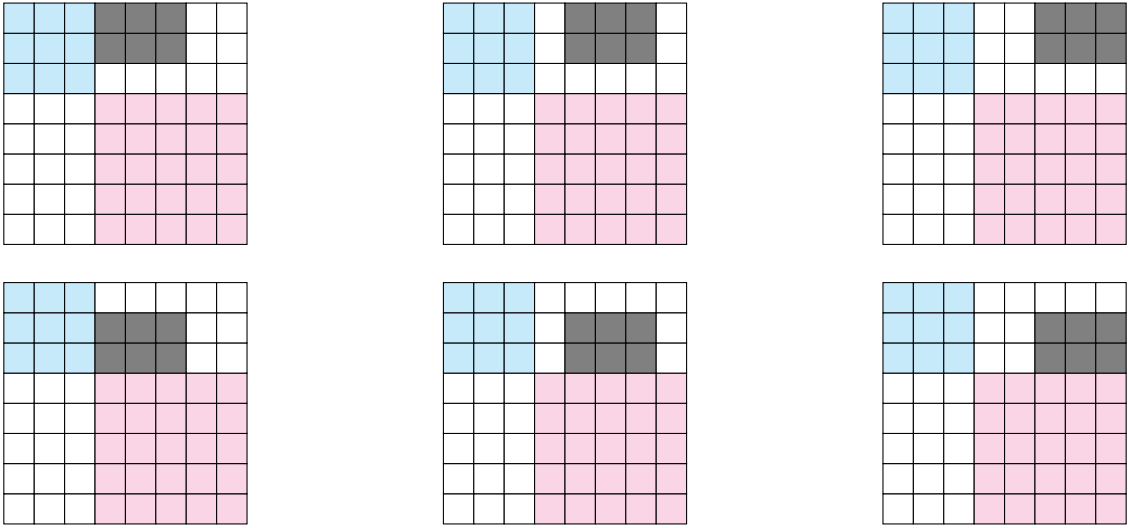
\begin{figure}[!h]
\begin{multicols}{3}

\begin{center}
\begin{tikzpicture}[scale=0.4]

\Rect{fill=white}{0}{0}{1}{1}
\Rect{fill=white}{1}{0}{1}{1}
\Rect{fill=white}{2}{0}{1}{1}
\Rect{fill=magenta!20}{3}{0}{1}{1}
\Rect{fill=magenta!20}{4}{0}{1}{1}
\Rect{fill=magenta!20}{5}{0}{1}{1}
\Rect{fill=magenta!20}{6}{0}{1}{1}
\Rect{fill=magenta!20}{7}{0}{1}{1}

\Rect{fill=white}{0}{1}{1}{1}
\Rect{fill=white}{1}{1}{1}{1}
\Rect{fill=white}{2}{1}{1}{1}
\Rect{fill=magenta!20}{3}{1}{1}{1}
\Rect{fill=magenta!20}{4}{1}{1}{1}
\Rect{fill=magenta!20}{5}{1}{1}{1}
\Rect{fill=magenta!20}{6}{1}{1}{1}
\Rect{fill=magenta!20}{7}{1}{1}{1}

\Rect{fill=white}{0}{2}{1}{1}
\Rect{fill=white}{1}{2}{1}{1}
\Rect{fill=white}{2}{2}{1}{1}
\Rect{fill=magenta!20}{3}{2}{1}{1}
\Rect{fill=magenta!20}{4}{2}{1}{1}
\Rect{fill=magenta!20}{5}{2}{1}{1}
\Rect{fill=magenta!20}{6}{2}{1}{1}
\Rect{fill=magenta!20}{7}{2}{1}{1}

\Rect{fill=white}{0}{3}{1}{1}
\Rect{fill=white}{1}{3}{1}{1}
\Rect{fill=white}{2}{3}{1}{1}
\Rect{fill=magenta!20}{3}{3}{1}{1}
\Rect{fill=magenta!20}{4}{3}{1}{1}
\Rect{fill=magenta!20}{5}{3}{1}{1}
\Rect{fill=magenta!20}{6}{3}{1}{1}
\Rect{fill=magenta!20}{7}{3}{1}{1}

\Rect{fill=white}{0}{4}{1}{1}
\Rect{fill=white}{1}{4}{1}{1}
\Rect{fill=white}{2}{4}{1}{1}
\Rect{fill=magenta!20}{3}{4}{1}{1}
\Rect{fill=magenta!20}{4}{4}{1}{1}
\Rect{fill=magenta!20}{5}{4}{1}{1}
\Rect{fill=magenta!20}{6}{4}{1}{1}
\Rect{fill=magenta!20}{7}{4}{1}{1}

\Rect{fill=cyan!20}{0}{5}{1}{1}
\Rect{fill=cyan!20}{1}{5}{1}{1}
\Rect{fill=cyan!20}{2}{5}{1}{1}
\Rect{fill=white}{3}{5}{1}{1}
\Rect{fill=white}{4}{5}{1}{1}
\Rect{fill=white}{5}{5}{1}{1}
\Rect{fill=white}{6}{5}{1}{1}
\Rect{fill=white}{7}{5}{1}{1}

\Rect{fill=cyan!20}{0}{6}{1}{1}
\Rect{fill=cyan!20}{1}{6}{1}{1}
\Rect{fill=cyan!20}{2}{6}{1}{1}
\Rect{fill=gray}{3}{6}{1}{1}
\Rect{fill=gray}{4}{6}{1}{1}
\Rect{fill=gray}{5}{6}{1}{1}
\Rect{fill=white}{6}{6}{1}{1}
\Rect{fill=white}{7}{6}{1}{1}

\Rect{fill=cyan!20}{0}{7}{1}{1}
\Rect{fill=cyan!20}{1}{7}{1}{1}
\Rect{fill=cyan!20}{2}{7}{1}{1}
\Rect{fill=gray}{3}{7}{1}{1}
\Rect{fill=gray}{4}{7}{1}{1}
\Rect{fill=gray}{5}{7}{1}{1}
\Rect{fill=white}{6}{7}{1}{1}
\Rect{fill=white}{7}{7}{1}{1}

\end{tikzpicture}
\end{center}

\begin{center}
\begin{tikzpicture}[scale=0.4]

\Rect{fill=white}{0}{0}{1}{1}
\Rect{fill=white}{1}{0}{1}{1}
\Rect{fill=white}{2}{0}{1}{1}
\Rect{fill=magenta!20}{3}{0}{1}{1}
\Rect{fill=magenta!20}{4}{0}{1}{1}
\Rect{fill=magenta!20}{5}{0}{1}{1}
\Rect{fill=magenta!20}{6}{0}{1}{1}
\Rect{fill=magenta!20}{7}{0}{1}{1}

\Rect{fill=white}{0}{1}{1}{1}
\Rect{fill=white}{1}{1}{1}{1}
\Rect{fill=white}{2}{1}{1}{1}
\Rect{fill=magenta!20}{3}{1}{1}{1}
\Rect{fill=magenta!20}{4}{1}{1}{1}
\Rect{fill=magenta!20}{5}{1}{1}{1}
\Rect{fill=magenta!20}{6}{1}{1}{1}
\Rect{fill=magenta!20}{7}{1}{1}{1}

\Rect{fill=white}{0}{2}{1}{1}
\Rect{fill=white}{1}{2}{1}{1}
\Rect{fill=white}{2}{2}{1}{1}
\Rect{fill=magenta!20}{3}{2}{1}{1}
\Rect{fill=magenta!20}{4}{2}{1}{1}
\Rect{fill=magenta!20}{5}{2}{1}{1}
\Rect{fill=magenta!20}{6}{2}{1}{1}
\Rect{fill=magenta!20}{7}{2}{1}{1}

\Rect{fill=white}{0}{3}{1}{1}
\Rect{fill=white}{1}{3}{1}{1}
\Rect{fill=white}{2}{3}{1}{1}
\Rect{fill=magenta!20}{3}{3}{1}{1}
\Rect{fill=magenta!20}{4}{3}{1}{1}
\Rect{fill=magenta!20}{5}{3}{1}{1}
\Rect{fill=magenta!20}{6}{3}{1}{1}
\Rect{fill=magenta!20}{7}{3}{1}{1}

\Rect{fill=white}{0}{4}{1}{1}
\Rect{fill=white}{1}{4}{1}{1}
\Rect{fill=white}{2}{4}{1}{1}
\Rect{fill=magenta!20}{3}{4}{1}{1}
\Rect{fill=magenta!20}{4}{4}{1}{1}
\Rect{fill=magenta!20}{5}{4}{1}{1}
\Rect{fill=magenta!20}{6}{4}{1}{1}
\Rect{fill=magenta!20}{7}{4}{1}{1}

\Rect{fill=cyan!20}{0}{5}{1}{1}
\Rect{fill=cyan!20}{1}{5}{1}{1}
\Rect{fill=cyan!20}{2}{5}{1}{1}
\Rect{fill=gray}{3}{5}{1}{1}
\Rect{fill=gray}{4}{5}{1}{1}
\Rect{fill=gray}{5}{5}{1}{1}
\Rect{fill=white}{6}{5}{1}{1}
\Rect{fill=white}{7}{5}{1}{1}

\Rect{fill=cyan!20}{0}{6}{1}{1}
\Rect{fill=cyan!20}{1}{6}{1}{1}
\Rect{fill=cyan!20}{2}{6}{1}{1}
\Rect{fill=gray}{3}{6}{1}{1}
\Rect{fill=gray}{4}{6}{1}{1}
\Rect{fill=gray}{5}{6}{1}{1}
\Rect{fill=white}{6}{6}{1}{1}
\Rect{fill=white}{7}{6}{1}{1}

\Rect{fill=cyan!20}{0}{7}{1}{1}
\Rect{fill=cyan!20}{1}{7}{1}{1}
\Rect{fill=cyan!20}{2}{7}{1}{1}
\Rect{fill=white}{3}{7}{1}{1}
\Rect{fill=white}{4}{7}{1}{1}
\Rect{fill=white}{5}{7}{1}{1}
\Rect{fill=white}{6}{7}{1}{1}
\Rect{fill=white}{7}{7}{1}{1}

\end{tikzpicture}
\end{center}

\begin{center}
\begin{tikzpicture}[scale=0.4]

\Rect{fill=white}{0}{0}{1}{1}
\Rect{fill=white}{1}{0}{1}{1}
\Rect{fill=white}{2}{0}{1}{1}
\Rect{fill=magenta!20}{3}{0}{1}{1}
\Rect{fill=magenta!20}{4}{0}{1}{1}
\Rect{fill=magenta!20}{5}{0}{1}{1}
\Rect{fill=magenta!20}{6}{0}{1}{1}
\Rect{fill=magenta!20}{7}{0}{1}{1}

\Rect{fill=white}{0}{1}{1}{1}
\Rect{fill=white}{1}{1}{1}{1}
\Rect{fill=white}{2}{1}{1}{1}
\Rect{fill=magenta!20}{3}{1}{1}{1}
\Rect{fill=magenta!20}{4}{1}{1}{1}
\Rect{fill=magenta!20}{5}{1}{1}{1}
\Rect{fill=magenta!20}{6}{1}{1}{1}
\Rect{fill=magenta!20}{7}{1}{1}{1}

\Rect{fill=white}{0}{2}{1}{1}
\Rect{fill=white}{1}{2}{1}{1}
\Rect{fill=white}{2}{2}{1}{1}
\Rect{fill=magenta!20}{3}{2}{1}{1}
\Rect{fill=magenta!20}{4}{2}{1}{1}
\Rect{fill=magenta!20}{5}{2}{1}{1}
\Rect{fill=magenta!20}{6}{2}{1}{1}
\Rect{fill=magenta!20}{7}{2}{1}{1}

\Rect{fill=white}{0}{3}{1}{1}
\Rect{fill=white}{1}{3}{1}{1}
\Rect{fill=white}{2}{3}{1}{1}
\Rect{fill=magenta!20}{3}{3}{1}{1}
\Rect{fill=magenta!20}{4}{3}{1}{1}
\Rect{fill=magenta!20}{5}{3}{1}{1}
\Rect{fill=magenta!20}{6}{3}{1}{1}
\Rect{fill=magenta!20}{7}{3}{1}{1}

\Rect{fill=white}{0}{4}{1}{1}
\Rect{fill=white}{1}{4}{1}{1}
\Rect{fill=white}{2}{4}{1}{1}
\Rect{fill=magenta!20}{3}{4}{1}{1}
\Rect{fill=magenta!20}{4}{4}{1}{1}
\Rect{fill=magenta!20}{5}{4}{1}{1}
\Rect{fill=magenta!20}{6}{4}{1}{1}
\Rect{fill=magenta!20}{7}{4}{1}{1}

\Rect{fill=cyan!20}{0}{5}{1}{1}
\Rect{fill=cyan!20}{1}{5}{1}{1}
\Rect{fill=cyan!20}{2}{5}{1}{1}
\Rect{fill=white}{3}{5}{1}{1}
\Rect{fill=white}{4}{5}{1}{1}
\Rect{fill=white}{5}{5}{1}{1}
\Rect{fill=white}{6}{5}{1}{1}
\Rect{fill=white}{7}{5}{1}{1}

\Rect{fill=cyan!20}{0}{6}{1}{1}
\Rect{fill=cyan!20}{1}{6}{1}{1}
\Rect{fill=cyan!20}{2}{6}{1}{1}
\Rect{fill=white}{3}{6}{1}{1}
\Rect{fill=gray}{4}{6}{1}{1}
\Rect{fill=gray}{5}{6}{1}{1}
\Rect{fill=gray}{6}{6}{1}{1}
\Rect{fill=white}{7}{6}{1}{1}

\Rect{fill=cyan!20}{0}{7}{1}{1}
\Rect{fill=cyan!20}{1}{7}{1}{1}
\Rect{fill=cyan!20}{2}{7}{1}{1}
\Rect{fill=white}{3}{7}{1}{1}
\Rect{fill=gray}{4}{7}{1}{1}
\Rect{fill=gray}{5}{7}{1}{1}
\Rect{fill=gray}{6}{7}{1}{1}
\Rect{fill=white}{7}{7}{1}{1}

\end{tikzpicture}
\end{center}

\begin{center}
\begin{tikzpicture}[scale=0.4]

\Rect{fill=white}{0}{0}{1}{1}
\Rect{fill=white}{1}{0}{1}{1}
\Rect{fill=white}{2}{0}{1}{1}
\Rect{fill=magenta!20}{3}{0}{1}{1}
\Rect{fill=magenta!20}{4}{0}{1}{1}
\Rect{fill=magenta!20}{5}{0}{1}{1}
\Rect{fill=magenta!20}{6}{0}{1}{1}
\Rect{fill=magenta!20}{7}{0}{1}{1}

\Rect{fill=white}{0}{1}{1}{1}
\Rect{fill=white}{1}{1}{1}{1}
\Rect{fill=white}{2}{1}{1}{1}
\Rect{fill=magenta!20}{3}{1}{1}{1}
\Rect{fill=magenta!20}{4}{1}{1}{1}
\Rect{fill=magenta!20}{5}{1}{1}{1}
\Rect{fill=magenta!20}{6}{1}{1}{1}
\Rect{fill=magenta!20}{7}{1}{1}{1}

\Rect{fill=white}{0}{2}{1}{1}
\Rect{fill=white}{1}{2}{1}{1}
\Rect{fill=white}{2}{2}{1}{1}
\Rect{fill=magenta!20}{3}{2}{1}{1}
\Rect{fill=magenta!20}{4}{2}{1}{1}
\Rect{fill=magenta!20}{5}{2}{1}{1}
\Rect{fill=magenta!20}{6}{2}{1}{1}
\Rect{fill=magenta!20}{7}{2}{1}{1}

\Rect{fill=white}{0}{3}{1}{1}
\Rect{fill=white}{1}{3}{1}{1}
\Rect{fill=white}{2}{3}{1}{1}
\Rect{fill=magenta!20}{3}{3}{1}{1}
\Rect{fill=magenta!20}{4}{3}{1}{1}
\Rect{fill=magenta!20}{5}{3}{1}{1}
\Rect{fill=magenta!20}{6}{3}{1}{1}
\Rect{fill=magenta!20}{7}{3}{1}{1}

\Rect{fill=white}{0}{4}{1}{1}
\Rect{fill=white}{1}{4}{1}{1}
\Rect{fill=white}{2}{4}{1}{1}
\Rect{fill=magenta!20}{3}{4}{1}{1}
\Rect{fill=magenta!20}{4}{4}{1}{1}
\Rect{fill=magenta!20}{5}{4}{1}{1}
\Rect{fill=magenta!20}{6}{4}{1}{1}
\Rect{fill=magenta!20}{7}{4}{1}{1}

\Rect{fill=cyan!20}{0}{5}{1}{1}
\Rect{fill=cyan!20}{1}{5}{1}{1}
\Rect{fill=cyan!20}{2}{5}{1}{1}
\Rect{fill=white}{3}{5}{1}{1}
\Rect{fill=gray}{4}{5}{1}{1}
\Rect{fill=gray}{5}{5}{1}{1}
\Rect{fill=gray}{6}{5}{1}{1}
\Rect{fill=white}{7}{5}{1}{1}

\Rect{fill=cyan!20}{0}{6}{1}{1}
\Rect{fill=cyan!20}{1}{6}{1}{1}
\Rect{fill=cyan!20}{2}{6}{1}{1}
\Rect{fill=white}{3}{6}{1}{1}
\Rect{fill=gray}{4}{6}{1}{1}
\Rect{fill=gray}{5}{6}{1}{1}
\Rect{fill=gray}{6}{6}{1}{1}
\Rect{fill=white}{7}{6}{1}{1}

\Rect{fill=cyan!20}{0}{7}{1}{1}
\Rect{fill=cyan!20}{1}{7}{1}{1}
\Rect{fill=cyan!20}{2}{7}{1}{1}
\Rect{fill=white}{3}{7}{1}{1}
\Rect{fill=white}{4}{7}{1}{1}
\Rect{fill=white}{5}{7}{1}{1}
\Rect{fill=white}{6}{7}{1}{1}
\Rect{fill=white}{7}{7}{1}{1}

\end{tikzpicture}
\end{center}

\begin{center}
\begin{tikzpicture}[scale=0.4]

\Rect{fill=white}{0}{0}{1}{1}
\Rect{fill=white}{1}{0}{1}{1}
\Rect{fill=white}{2}{0}{1}{1}
\Rect{fill=magenta!20}{3}{0}{1}{1}
\Rect{fill=magenta!20}{4}{0}{1}{1}
\Rect{fill=magenta!20}{5}{0}{1}{1}
\Rect{fill=magenta!20}{6}{0}{1}{1}
\Rect{fill=magenta!20}{7}{0}{1}{1}

\Rect{fill=white}{0}{1}{1}{1}
\Rect{fill=white}{1}{1}{1}{1}
\Rect{fill=white}{2}{1}{1}{1}
\Rect{fill=magenta!20}{3}{1}{1}{1}
\Rect{fill=magenta!20}{4}{1}{1}{1}
\Rect{fill=magenta!20}{5}{1}{1}{1}
\Rect{fill=magenta!20}{6}{1}{1}{1}
\Rect{fill=magenta!20}{7}{1}{1}{1}

\Rect{fill=white}{0}{2}{1}{1}
\Rect{fill=white}{1}{2}{1}{1}
\Rect{fill=white}{2}{2}{1}{1}
\Rect{fill=magenta!20}{3}{2}{1}{1}
\Rect{fill=magenta!20}{4}{2}{1}{1}
\Rect{fill=magenta!20}{5}{2}{1}{1}
\Rect{fill=magenta!20}{6}{2}{1}{1}
\Rect{fill=magenta!20}{7}{2}{1}{1}

\Rect{fill=white}{0}{3}{1}{1}
\Rect{fill=white}{1}{3}{1}{1}
\Rect{fill=white}{2}{3}{1}{1}
\Rect{fill=magenta!20}{3}{3}{1}{1}
\Rect{fill=magenta!20}{4}{3}{1}{1}
\Rect{fill=magenta!20}{5}{3}{1}{1}
\Rect{fill=magenta!20}{6}{3}{1}{1}
\Rect{fill=magenta!20}{7}{3}{1}{1}

\Rect{fill=white}{0}{4}{1}{1}
\Rect{fill=white}{1}{4}{1}{1}
\Rect{fill=white}{2}{4}{1}{1}
\Rect{fill=magenta!20}{3}{4}{1}{1}
\Rect{fill=magenta!20}{4}{4}{1}{1}
\Rect{fill=magenta!20}{5}{4}{1}{1}
\Rect{fill=magenta!20}{6}{4}{1}{1}
\Rect{fill=magenta!20}{7}{4}{1}{1}

\Rect{fill=cyan!20}{0}{5}{1}{1}
\Rect{fill=cyan!20}{1}{5}{1}{1}
\Rect{fill=cyan!20}{2}{5}{1}{1}
\Rect{fill=white}{3}{5}{1}{1}
\Rect{fill=white}{4}{5}{1}{1}
\Rect{fill=white}{5}{5}{1}{1}
\Rect{fill=white}{6}{5}{1}{1}
\Rect{fill=white}{7}{5}{1}{1}

\Rect{fill=cyan!20}{0}{6}{1}{1}
\Rect{fill=cyan!20}{1}{6}{1}{1}
\Rect{fill=cyan!20}{2}{6}{1}{1}
\Rect{fill=white}{3}{6}{1}{1}
\Rect{fill=white}{4}{6}{1}{1}
\Rect{fill=gray}{5}{6}{1}{1}
\Rect{fill=gray}{6}{6}{1}{1}
\Rect{fill=gray}{7}{6}{1}{1}

\Rect{fill=cyan!20}{0}{7}{1}{1}
\Rect{fill=cyan!20}{1}{7}{1}{1}
\Rect{fill=cyan!20}{2}{7}{1}{1}
\Rect{fill=white}{3}{7}{1}{1}
\Rect{fill=white}{4}{7}{1}{1}
\Rect{fill=gray}{5}{7}{1}{1}
\Rect{fill=gray}{6}{7}{1}{1}
\Rect{fill=gray}{7}{7}{1}{1}

\end{tikzpicture}
\end{center}

\begin{center}
\begin{tikzpicture}[scale=0.4]

\Rect{fill=white}{0}{0}{1}{1}
\Rect{fill=white}{1}{0}{1}{1}
\Rect{fill=white}{2}{0}{1}{1}
\Rect{fill=magenta!20}{3}{0}{1}{1}
\Rect{fill=magenta!20}{4}{0}{1}{1}
\Rect{fill=magenta!20}{5}{0}{1}{1}
\Rect{fill=magenta!20}{6}{0}{1}{1}
\Rect{fill=magenta!20}{7}{0}{1}{1}

\Rect{fill=white}{0}{1}{1}{1}
\Rect{fill=white}{1}{1}{1}{1}
\Rect{fill=white}{2}{1}{1}{1}
\Rect{fill=magenta!20}{3}{1}{1}{1}
\Rect{fill=magenta!20}{4}{1}{1}{1}
\Rect{fill=magenta!20}{5}{1}{1}{1}
\Rect{fill=magenta!20}{6}{1}{1}{1}
\Rect{fill=magenta!20}{7}{1}{1}{1}

\Rect{fill=white}{0}{2}{1}{1}
\Rect{fill=white}{1}{2}{1}{1}
\Rect{fill=white}{2}{2}{1}{1}
\Rect{fill=magenta!20}{3}{2}{1}{1}
\Rect{fill=magenta!20}{4}{2}{1}{1}
\Rect{fill=magenta!20}{5}{2}{1}{1}
\Rect{fill=magenta!20}{6}{2}{1}{1}
\Rect{fill=magenta!20}{7}{2}{1}{1}

\Rect{fill=white}{0}{3}{1}{1}
\Rect{fill=white}{1}{3}{1}{1}
\Rect{fill=white}{2}{3}{1}{1}
\Rect{fill=magenta!20}{3}{3}{1}{1}
\Rect{fill=magenta!20}{4}{3}{1}{1}
\Rect{fill=magenta!20}{5}{3}{1}{1}
\Rect{fill=magenta!20}{6}{3}{1}{1}
\Rect{fill=magenta!20}{7}{3}{1}{1}

\Rect{fill=white}{0}{4}{1}{1}
\Rect{fill=white}{1}{4}{1}{1}
\Rect{fill=white}{2}{4}{1}{1}
\Rect{fill=magenta!20}{3}{4}{1}{1}
\Rect{fill=magenta!20}{4}{4}{1}{1}
\Rect{fill=magenta!20}{5}{4}{1}{1}
\Rect{fill=magenta!20}{6}{4}{1}{1}
\Rect{fill=magenta!20}{7}{4}{1}{1}

\Rect{fill=cyan!20}{0}{5}{1}{1}
\Rect{fill=cyan!20}{1}{5}{1}{1}
\Rect{fill=cyan!20}{2}{5}{1}{1}
\Rect{fill=white}{3}{5}{1}{1}
\Rect{fill=white}{4}{5}{1}{1}
\Rect{fill=gray}{5}{5}{1}{1}
\Rect{fill=gray}{6}{5}{1}{1}
\Rect{fill=gray}{7}{5}{1}{1}

\Rect{fill=cyan!20}{0}{6}{1}{1}
\Rect{fill=cyan!20}{1}{6}{1}{1}
\Rect{fill=cyan!20}{2}{6}{1}{1}
\Rect{fill=white}{3}{6}{1}{1}
\Rect{fill=white}{4}{6}{1}{1}
\Rect{fill=gray}{5}{6}{1}{1}
\Rect{fill=gray}{6}{6}{1}{1}
\Rect{fill=gray}{7}{6}{1}{1}

\Rect{fill=cyan!20}{0}{7}{1}{1}
\Rect{fill=cyan!20}{1}{7}{1}{1}
\Rect{fill=cyan!20}{2}{7}{1}{1}
\Rect{fill=white}{3}{7}{1}{1}
\Rect{fill=white}{4}{7}{1}{1}
\Rect{fill=white}{5}{7}{1}{1}
\Rect{fill=white}{6}{7}{1}{1}
\Rect{fill=white}{7}{7}{1}{1}

\end{tikzpicture}
\end{center}

\end{multicols}
\caption{
    Graphical illustration of how $\| A \|_{(\Pp, \Gamma, s_1, s_2)}$ is computed. 
    In this example, we take $d = 8$, $m = 3$, $s_1 = 2$, $s_2 = 3$.
    Consider the rotated matrix $\widetilde{A} = \Gamma^\T A \Gamma$ and split it into four blocks: $m\times m$ top left block (blue), $(d-m)\times(d-m)$ bottom right block (pink), bottom left $(d-m)\times m$ block (white) and top right $m\times(d-m)$ block.
    Then $\| A \|_{(\Pp, \Gamma, s_1, s_2)}$ is computed as half of the sum of spectral norms of blue and pink blocks, plus the largest spectral norm of gray submatrices, for which we have $(m-s_1+1)\cdot(d-m-s_2+1) = 6$ options.
    }
\label{draw1}
\end{figure}

%% file: source/Theory.tex
Before stating our assumptions and theoretical results, we introduce some important characteristics of the true covariance $\St$ that will appear in the error bounds.
In particular, the relative rank of $\St$ (see \cite{Wahl}) is
\begin{equation}
    \begin{aligned}
        \relr_r(\St) \eqdef \sum\limits_{s\neq r} \frac{m_s \mu_s}{|\mu_r-\mu_s|} + \frac{m_r \mu_r}{\min(\mu_{r-1} - \mu_r, \mu_r - \mu_{r+1})} \;\;\;\text{ for all }\;r\in [q].
    \nonumber
    \end{aligned}
\end{equation}
It turns out that the following quantity will play role of effective dimension:
\begin{equation}
    \begin{aligned}
        \rr_\J(\St) \eqdef \left( \sum\limits_{r\in\J} \left( \relr_r(\St) \sqrt{\sum\limits_{s\neq r} \frac{m_r\mu_r m_s \mu_s}{(\mu_r-\mu_s)^2}} \right) \right)^{2/3}.
    \nonumber
    \end{aligned}
\end{equation}
Other important quantities appearing in the theorems are
\begin{equation}
    \begin{aligned}
    \clow_\J(\St) \eqdef \min\limits_{r\in\J, s\notin\J} \frac{\sqrt{\mu_r \mu_s}}{|\mu_r - \mu_s|},\;\;\;
    \chigh_\J(\St) \eqdef \max\limits_{r\in\J, s\notin\J} \frac{\sqrt{\mu_r \mu_s}}{|\mu_r - \mu_s|},\;\;\;
    \cond_\J(\St) \eqdef \chigh_\J(\St)/\clow_\J(\St).
    \nonumber
    \end{aligned}
\end{equation}
The last quantity $\cond_\J(\St)$ can be interpreted as a kind of condition number, but with respect to splitting the eigenvalues into two groups associated with $\J$ and $\J^c$.
Throughout the paper, when it does not cause ambiguity (in context of one-sample problem) we write $\rr, \clow, \chigh, \cond$ instead of $\rr_\J(\St), \clow_\J(\St), \chigh_\J(\St), \cond_\J(\St)$, respectively, to keep the notation light.
\begin{remark}
	Later in Section~\ref{S:FM} we will focus on factor models. Under the standard assumptions imposed in that field, we will see that that above quantities in this case are of the following order:
	\begin{equation}
	\begin{aligned}
		\rr \asymp m^{5/3},\;\;\;\clow \asymp \chigh \asymp \frac{1}{\sqrt{d}},\;\;\;\cond \asymp 1,
	\nonumber
	\end{aligned}
	\end{equation}
where $m$ plays role of the number of common factors in the model. The fact that the effective dimension $\rr$ in this situation does not depend on the full dimension $d$ (can even be finite) will help us to significantly weaken the relation between $d$ and $n$ required for the validity of the tests, compared the the previous works.
\end{remark}


\subsection{Assumptions}
We start by specifying the assumptions which will be required in our theorems.
\begin{assumption}[Uncorrelatedness] \label{Uncorrelatedness}
    \label{A: independence}
    $v^\T \Pt X X^\T \Pt \widetilde{v}$ and $w^\T \Ptc X X^\T \Ptc \widetilde{w}$ are uncorrelated for all $v, \widetilde{v}, w, \widetilde{w} \in \R^d$.
\end{assumption}
\begin{remark}
    Any of the following conditions is sufficient for Assumption~\ref{Uncorrelatedness}:
    \begin{enumerate}[(i)]
        \item $\Pt X$ and $\Ptc X$ are independent (these random vectors are always orthogonal, and consequently uncorrelated; this condition is somewhat stronger);
        \item The components of $\St^{-1/2} X$ are independent;
        \item $X$ is Gaussian random vector.
    \end{enumerate}
    Additionally note, that (iii) implies (ii), (ii) implies (i).
\end{remark}

\begin{assumption}[Tail bound]
    \label{A: tails}
    $\St^{-1/2} X$ is jointly sub-Weibull random vector with parameter $0 < \beta \leq 2$ (see \cite{Weibull}). That is, there exists a constant $c > 0$ such that
    \begin{equation}
            \begin{aligned}
            	\| \St^{-1/2}X \|_{J, \psi_\beta} \eqdef \sup\limits_{u\in S^{d-1}} \| u^{\top} \St^{-1/2}X\|_{\psi_\beta} \leq c < \infty,
                \nonumber
            \end{aligned}
        \end{equation}
    where $\|\cdot\|_{\psi_\beta}$ is the Orlicz norm for $\psi_\beta = e^{x^\beta} - 1$.
    The following tail bound takes place:
    \begin{equation}
            \begin{aligned}
            	& \Prob\left[ |u^{\top} \St^{-1/2}X| \geq t \right] \leq 2\,\exp\left(-\left(t/c\right)^\beta\right),
                \nonumber
            \end{aligned}
        \end{equation}
    for all $u\in S^{d-1}$ and $t > 0$.
\end{assumption}
\begin{remark}
    Case $\beta = 2$ corresponds to sub-Gaussian distribution of $X$. We restrict ourselves to the case $\beta \leq 2$, since it is unreasonable to expect tails lighter than Gaussian in applications. Nevertheless, our results extend easily to $\beta > 2$ by replacing $\beta$ with $(\beta \land 2)$ in all of the further error bounds.
\end{remark}

Let us introduce some auxiliary quantities and rates, which will appear in our bounds:
    \begin{equation}
        \begin{aligned}
        p = p_{d,n,s_1,s_2}&\eqdef \exp\left( (s_1+s_2)\log(3n) + 2\log(d) \right),\\
	    \scc &\eqdef C_\beta c^2  \left( \sqrt{\frac{\log(n)+\log(d)}{n}} + \frac{(\log(n))^{1/\beta} (\log(n)+\log(d))^{2/\beta}}{n}\right),\\
        \bcc &\eqdef C c^2 \frac{ (\log(n) + \log(2d^2))^{\frac{2}{\beta}+\frac{1}{2}} }{\sqrt{n}},\\
        \zeta[\delta] &\eqdef \delta \left( \log\left(\frac{ep}{\delta}\right)\right)^{1/2} \;\;\;\text{ for all }\;\delta > 0,\\
         \vartheta[\delta] &\eqdef \delta^{1/3} \left( \log\left(\frac{ep}{\delta}\right)\right)^{2/3} \;\;\;\text{ for all }\;\delta > 0.
            \nonumber
        \end{aligned}
    \end{equation}
The constants $C_\beta$ and $C$ are properly chosen and come from the proofs of the theorems in the sequel. The functions $\zeta[\cdot]$ and $\vartheta[\cdot]$ are introduced just for convenience to avoid long expressions with logarithmic factors.
Now we state an additional assumption.
\begin{assumption} \label{A:additional}
    The following holds:
\begin{enumerate}[(i)]
    \item $\scc \max\limits_{r\in\J} \relr_r(\St) \leq 1/12$.
    \item $\bcc \max\limits_{r\in\J} \relr_r(\St) \leq 1/12$.
\end{enumerate}
\end{assumption}

\subsection{Validity}
\subsubsection{One-sample test}
In this subsection we work under $H_0^{(1)}$, so that $\Pt = \Ph$ and
\begin{equation}
            \begin{aligned}
            	\Qt = \sqrt{n}\| \Pe - \Pt \|_{(\Pt, \Gamma^\circ, s_1, s_2)}.
            \nonumber
            \end{aligned}
        \end{equation}
Our first result provides approximation for the distribution of $\Qt$ under $H_0^{(1)}$.
\begin{theorem}[One-sample test; test statistics approximation]
    \label{Th:GA1}
    Let the data $\data$ satisfy Assumptions \ref{A: independence}, \ref{A: tails}, \ref{A:additional}(i).
    Then there exists a Gaussian vector $Y \in \R^p$, with
    specific covariance structure (presented in the proof) that depends on $\St$, such that under $H_0^{(1)}$ holds
     \begin{equation}
            \begin{aligned}
            	&\sup\limits_{z\in\R} \left| \Prob\left[ \Qt \leq z \right] - \Prob\left[ \max\limits_{j\in[p]} Y_j \leq z \right] \right| \leq \Diamond^{(1)},
            \nonumber
            \end{aligned}
        \end{equation}
    where
            \begin{align}
                \Diamond^{(1)} &= C_\cond \left\{\Diamond^{GA} + \zeta\left[\sqrt{n}\scc^2\rr^{3/2}/\clow\right] \right\},\nonumber\\
            	\Diamond^{GA} &=
            	8^{3/(2\beta)} \left( \frac{( \log(pn))^7}{n} \right)^{1/8}
            	+ c^2 \left( \frac{(\log(2pn^2))^{3+4/\beta}}{n} \right)^{1/2}. \label{Eq:Diamond}
            \end{align}

    Moreover, the same result holds for spectral norm test statistics $\Qtsp = \sqrt{n}\| \Pe - \Ph\|$, if we take $s_1 = m$, $s_2 = d-m$.
\end{theorem}
This theorem gives understanding of how to prove the next two validity results. For the validity of Approach 1 we need to define an additional quantity and an assumption on it.
\begin{assumption} \label{A:additionalB}
    Define
    \begin{equation}
        \begin{aligned}
        \bcd &\eqdef
        	C_{\beta} \,c^4 \, \cond^2 \left( \sqrt{\frac{\log(pn)}{n}} + \frac{(\log(n)^{2/\beta} ( \log(pn))^{4/\beta}}{n}\right).
            \nonumber
        \end{aligned}
    \end{equation}
    Here again $C_\beta$ comes from the corresponding proof. Suppose $\bcd \leq 1/2$.
\end{assumption}

\begin{theorem}[One-sample test; validity of Approach 1]
    \label{Th:B1}
    Let the data $\data$ satisfy Assumptions \ref{A: independence}, \ref{A: tails}, \ref{A:additional}. Also suppose Assumption~\ref{A:additionalB} is fulfilled.
    Then under $H_0^{(1)}$ with probability $1-1/n$
     \begin{equation}
            \begin{aligned}
            	&\sup\limits_{z\in\R} \left| \Prob\left[ \Qt \leq z \right] - \Prob\left[ \Qb \leq z \,|\,\data\right] \right| \leq \Diamond^{(1)}_B,
            \nonumber
            \end{aligned}
        \end{equation}
    where
    \begin{equation}
            \begin{aligned}
            	\Diamond^{(1)}_B \eqdef
            	C_\cond \left\{ \Diamond^{GA} + \zeta\left[\sqrt{n}(\bcc + \scc)^2\rr^{3/2}/\clow\right] + \vartheta[\bcd] \right\}
            \nonumber
            \end{aligned}
        \end{equation}
    with $\Diamond^{GA}$ from \eqref{Eq:Diamond}.

    Moreover, the same result holds for spectral norm test statistics $\Qtsp = \sqrt{n}\| \Pe - \Ph\|$ and $\Qbsp = \sqrt{n}\|\Pb-\Pe\|$, if we take $s_1 = m$, $s_2 = d-m$.
\end{theorem}

Similarly, the validity of Approach 2 requires the following assumption.
\begin{assumption} \label{A:additionalF}
Define
    \begin{equation}
        \begin{aligned}
        \fcd &\eqdef
            |\J|\, C_\beta c^2 \cond^2 \left( \sqrt{\frac{\log(pn)}{n}} + \frac{(\log(n))^{1/\beta} (\log(pn))^{2/\beta}}{n} \right),
            \nonumber
        \end{aligned}
    \end{equation}
    As above, $C_\beta$ comes from the corresponding proof. Suppose $\fcd \leq 1/2$.
\end{assumption}

\begin{theorem}[One-sample test; validity of Approach 2]
    \label{Th:F1}
    Let the data $\data$ satisfy Assumptions \ref{A: independence}, \ref{A: tails}, \ref{A:additional}. Also suppose Assumption~\ref{A:additionalF} is fulfilled.
    Then under $H_0^{(1)}$ with probability $1-1/n$
     \begin{equation}
            \begin{aligned}
            	&\sup\limits_{z\in\R} \left| \Prob\left[ \Qt \leq z \right] - \Prob\left[ \Qf \leq z \,|\,\data\right] \right| \leq \Diamond^{(1)}_F,
            \nonumber
            \end{aligned}
        \end{equation}
    where
    \begin{equation}
            \begin{aligned}
            	 \Diamond^{(1)}_F \eqdef
            	C_\cond \left\{ \Diamond^{GA} + \zeta\left[\sqrt{n}(\fcc + \scc)^2\rr^{3/2}/\clow\right] + \vartheta[\fcd] \right\},
            \nonumber
            \end{aligned}
        \end{equation}
    with $\Diamond^{GA}$ from \eqref{Eq:Diamond}.

    Moreover, the same result holds for spectral norm test statistics $\Qtsp = \sqrt{n}\| \Pe - \Ph\|$ and $\Qfsp = \sqrt{n}\| \Pf - \Pe\|$, if we take $s_1 = m$, $s_2 = d-m$.
\end{theorem}

The previous two results imply that both Approach 1 and Approach 2 have type-I error close to the desired level $\alpha$. This is formalized in the following Corollary.

\begin{corollary}[One-sample test; type-I error] \label{Corollary1}
    (i) Assume the conditions of Theorem~\ref{Th:B1} are fulfilled. Define
        \begin{equation}
            \begin{aligned}
            q^{(1)}_B(\alpha) \eqdef
                \inf \left\{ \gamma > 0: \Prob\left[ \Qb > \gamma \,\big|\,\data\right] \leq \alpha \right\}.
            \nonumber
            \end{aligned}
        \end{equation}
        Then
        \begin{equation}
            \begin{aligned}
                \sup\limits_{\alpha\in(0; 1)} \left| \Prob\left[\Qt > q^{(1)}_B(\alpha)\right] - \alpha\right| \leq \Diamond_B^{(1)} + \frac{1}{n},
            \nonumber
            \end{aligned}
        \end{equation}
        where $\Diamond_B^{(1)}$ is the total error term from Theorem~\ref{Th:B1}.
    \\
    (ii) Assume the conditions of Theorem~\ref{Th:F1} are fulfilled. Define
        \begin{equation}
            \begin{aligned}
            q^{(1)}_F(\alpha) \eqdef
                \inf \left\{ \gamma > 0: \Prob\left[ \Qf > \gamma \,\big|\,\data\right] \leq \alpha \right\}.
            \nonumber
            \end{aligned}
        \end{equation}
        Then
        \begin{equation}
            \begin{aligned}
                \sup\limits_{\alpha\in(0; 1)} \left| \Prob\left[\Qt > q^{(1)}_F(\alpha)\right] - \alpha\right| \leq \Diamond_F^{(1)} + \frac{1}{n},
            \nonumber
            \end{aligned}
        \end{equation}
        where $\Diamond_F^{(1)}$ is the total error term from Theorem~\ref{Th:F1}.
\end{corollary}

\begin{remark} \label{rate}
    For the sake of illustration, let us treat $\beta$ as fixed and omit the logarithmic terms.
    Then the error bounds on the Kolmogorov distance in the previous theorems become more transparent and can be bounded by:
        \begin{equation}
            \begin{aligned}
            	C_\cond \left\{ \left( \frac{(s_1+s_2)^7}{n} \right)^{1/8} + \left( \frac{(s_1+s_2)^{4/\beta+2}}{n} \right)^{1/3}
            	+ \frac{1}{\clow} \left( \frac{(s_1+s_2)\,\rr^3}{n}\right)^{1/2} \right\},
            \nonumber
            \end{aligned}
        \end{equation}
    which in case of the spectral norm reduces to
    \begin{equation}
            \begin{aligned}
            	C_\cond \left\{ \left( \frac{d^7}{n} \right)^{1/8} + \left( \frac{d^{4/\beta+2}}{n} \right)^{1/3} \right\}.
            \nonumber
            \end{aligned}
        \end{equation}
    Note additionally, that with slightly different technique used in previous works of  \cite{Koltchinskii_NAACOSPOSC,Naumov,Silin_1}, $\rr^3$ can be replaced by $d^2$. This will improve our bound in case when $\rr \asymp d$, however will be worse if $\rr \ll d$. Since the main motivation behind our work is Factor Models, where $\rr \asymp m^{5/3} \ll d$, we choose to present the result with $\rr^3$. We preview the bound that will be obtained in case of Factor Models (take $s_1 = s_2 = 1$ for simplicity):
    \begin{equation}
            \begin{aligned}
            	C \left\{ \frac{1}{n^{1/8}}  + m^{5/2} \sqrt{\frac{d}{n}} \right\}.
            \nonumber
            \end{aligned}
        \end{equation}
\end{remark}

\subsubsection{Two-sample test}
Similar theoretical properties are obtained for the two-sample problem. Before we state them, we introduce one more version of effective dimension that will show up:
\begin{equation}
    \begin{aligned}
        \rre_\J(\St) \eqdef \sum\limits_{r\in\J} \sum\limits_{s\notin\J} \frac{m_r\mu_r m_s \mu_s}{(\mu_r-\mu_s)^2}.
    \nonumber
    \end{aligned}
\end{equation}
Also, we define ``total effective dimensions'' for two samples as
\begin{equation}
    \begin{aligned}
        \rr_{a,b} &\eqdef \left( \rr_{\J_a}(\St_a)^{3/2} + \rr_{\J_b}(\St_b)^{3/2} \right)^{2/3},\\
        \rre_{a,b} &\eqdef \left( \rre_{\J_a}(\St_a)^{1/2} + \rre_{\J_b}(\St_b)^{1/2} \right)^{2},
    \nonumber
    \end{aligned}
\end{equation}
and
\begin{equation}
    \begin{aligned}
        \chighab \eqdef \chigh_{\J_a}(\St_a) \lor \chigh_{\J_b}(\St_b),\;\;
        \clowab \eqdef \clow_{\J_a}(\St_a) \land \clow_{\J_b}(\St_b),\;\;
        \condab \eqdef \frac{\chighab}{\clowab}.
    \nonumber
    \end{aligned}
\end{equation}
Define $p_{a,b}$ in a similar fashion as $p$, but with $n$ replaced by $n_a + n_b$.

\begin{theorem}[Two-sample test; test statistic approximation] \label{Th:GA2}
Let the data $\dataa$ and $\datab$ satisfy Assumptions \ref{A: independence}, \ref{A: tails}, \ref{A:additional}(i) (with $n$ replaced by $n_a \land n_b$).
Additionally, assume $\rr_{a,b}^{3/2} \rre_{a,b}^{1/2} \sccab \leq \rr_{a,b}^{3/2} + \rre_{a,b}$.
    Then there exists a Gaussian vector $Y^{a,b} \in \R^{p_{a,b}}$, with
    specific covariance structure (presented in the proof) that depends on $\St_a$ and $\St_b$, such that under $H_0^{(2)}$ holds
     \begin{equation}
            \begin{aligned}
            	&\sup\limits_{z\in\R} \left| \Prob\left[ \Qtab \leq z \,\big|\,\overline{\Gamma}\right] - \Prob\left[ \max\limits_{j\in[p]} Y_j^{a,b} \leq z \,\big|\,\overline{\Gamma}\right] \right| \leq \Diamond^{(2)},
            \nonumber
            \end{aligned}
        \end{equation}
    with probability $1-1/n_a-1/n_b$, where
    \begin{equation}
            \begin{aligned}
                \Diamond^{(2)} &\eqdef
                C_{\condab} \left\{\Diamond^{GA} + \zeta\left[\sqrt{\frac{n_an_b}{n_a+n_b}} \sccab^2 (\rr_{a,b}^{3/2} + \rre_{a,b})/\clowab \right] \right\},\\
            	\Diamond^{GA} &\eqdef
            	8^{3/(2\beta)} \left( \frac{( \log(p_{a,b}(n_a+n_b)))^7}{n_a+n_b} \right)^{1/8}
            	+ c^2 \left( \frac{(\log(2p_{a,b}(n_a+n_b)^2))^{3+4/\beta}}{n_a+n_b} \right)^{1/2} +\\
            	&\qquad \qquad + \frac{1}{n_a} + \frac{1}{n_b}.
                \label{Eq:Diamond_ab}
            \end{aligned}
        \end{equation}
\end{theorem}

To state validity of Approach 1 and Approach 2 in two-sample problem, let $\bcdab$ be defined as $\bcd$ and $\fcdab$ be defined as $\fcd$ with $n$, $\cond$, $|\J|$ replaced with $n_a\land n_b$, $\condab$, $|\J_a|\lor|\J_b|$, respectively. Then we have the following theorems.

\begin{theorem}[Two-sample test; validity of Approach 1] \label{Th:B2}
Let the data $\dataa$ and $\datab$ satisfy Assumptions \ref{A: independence}, \ref{A: tails}, \ref{A:additional} (with $n$ replaced by $n_a \land n_b$). Additionally, assume $\bcdab \leq 1/2$ and $\rr_{a,b}^{3/2} \rre_{a,b}^{1/2} \bccab \leq \rr_{a,b}^{3/2} + \rre_{a,b}$.
    Then under $H_0^{(2)}$ with probability $1-1/n_a-1/n_b$
     \begin{equation}
            \begin{aligned}
            	&\sup\limits_{z\in\R} \left| \Prob\left[ \Qtab \leq z \,\big|\,\overline{\Gamma} \right] - \Prob\left[ \Qbab \leq z \,\big|\,\dataa,\datab\right] \right| \leq \Diamond^{(2)}_B,
            \nonumber
            \end{aligned}
        \end{equation}
        where
        \begin{equation}
            \begin{aligned}
            	\Diamond^{(2)}_B &\eqdef
            	C_{\condab} \Bigg\{ \Diamond^{GA} + \zeta\left[\sqrt{\frac{n_an_b}{n_a+n_b}}(\sccab+\bccab)^2 (\rr_{a,b}^{3/2} + \rre_{a,b})/\clowab \right] + \vartheta\left[\bcdab\right] \Bigg\},
            \nonumber
            \end{aligned}
        \end{equation}
    with $\Diamond^{GA}$ from \eqref{Eq:Diamond_ab} .
\end{theorem}

\begin{theorem}[Two-sample test; validity of Approach 2] \label{Th:F2}
Let the data $\dataa$ and $\datab$ satisfy Assumptions \ref{A: independence}, \ref{A: tails}, \ref{A:additional} (with $n$ replaced by $n_a \land n_b$). Additionally, assume $\fcdab \leq 1/2$ and $\rr_{a,b}^{3/2} \rre_{a,b}^{1/2} \fccab \leq \rr_{a,b}^{3/2} + \rre_{a,b}$.
    Then under $H_0^{(2)}$ with probability $1-1/n_a-1/n_b$
     \begin{equation}
            \begin{aligned}
            	&\sup\limits_{z\in\R} \left| \Prob\left[ \Qtab \leq z \,\big|\,\overline{\Gamma} \right] - \Prob\left[ \Qfab \leq z \,\big|\,\dataa,\datab\right] \right| \leq
            	\Diamond^{(2)}_F,
            \nonumber
            \end{aligned}
        \end{equation}
    where
    \begin{equation}
            \begin{aligned}
            	\Diamond^{(2)}_F
            	\eqdef
            	C_{\condab}\Bigg\{ \Diamond^{GA} + \zeta\left[\sqrt{\frac{n_an_b}{n_a+n_b}} (\sccab+\fccab)^2 (\rr_{a,b}^{3/2} + \rre_{a,b})/\clowab\right] +
            	\vartheta\left[\fcdab\right] \Bigg\},
            \nonumber
            \end{aligned}
        \end{equation}
    with $\Diamond^{GA}$ from \eqref{Eq:Diamond_ab}.
\end{theorem}
\begin{remark}
    The condition $\rr_{a,b}^{3/2} \rre_{a,b}^{1/2} \fccab \leq \rr_{a,b}^{3/2} + \rre_{a,b}$ in the theorems above is technical and is imposed just to slightly simplify the bounds.
\end{remark}

Similarly to one-sample case, we have the following guarantees for type-I error.
\begin{corollary}[Two-sample test; type-I error] \label{Corollary2}
    (i) Assume the conditions of Theorem~\ref{Th:B2} are fulfilled. Define
        \begin{equation}
            \begin{aligned}
            q^{(2)}_B(\alpha) \eqdef
                \inf \left\{ \gamma > 0: \Prob\left[ \Qbab > \gamma \,\big|\,\data\right] \leq \alpha \right\}.
            \nonumber
            \end{aligned}
        \end{equation}
        Then
        \begin{equation}
            \begin{aligned}
                \sup\limits_{\alpha\in(0; 1)} \left| \Prob\left[\Qtab > q^{(2)}_B(\alpha)\right] - \alpha\right| \leq \Diamond_B^{(2)} + \frac{1}{n_a} + \frac{1}{n_b},
            \nonumber
            \end{aligned}
        \end{equation}
        where $\Diamond_B^{(2)}$ is the complete error term from Theorem~\ref{Th:B2}.
    \\
    (ii) Assume the conditions of Theorem~\ref{Th:F2} are fulfilled. Define
        \begin{equation}
            \begin{aligned}
            q^{(2)}_F(\alpha) \eqdef
                \inf \left\{ \gamma > 0: \Prob\left[ \Qf > \gamma \,\big|\,\data\right] \leq \alpha \right\}.
            \nonumber
            \end{aligned}
        \end{equation}
        Then
        \begin{equation}
            \begin{aligned}
                \sup\limits_{\alpha\in(0; 1)} \left| \Prob\left[\Qtab > q^{(2)}_F(\alpha)\right] - \alpha\right| \leq \Diamond_F^{(2)} + \frac{1}{n_a} + \frac{1}{n_b},
            \nonumber
            \end{aligned}
        \end{equation}
        where $\Diamond_F^{(2)}$ is the complete error term from Theorem~\ref{Th:F2}.
\end{corollary}

As in one-sample case, more transparent expression for the error bound can be seen in Remark~\ref{rate}, with $n$ replaced by $n_a\land n_b$ and $\rr$ replaced by $\rr_{a,b} + \rre_{a,b}^{2/3}$.

\begin{remark}
    One can also consider two-sample tests based on the spectral norm. In this case, there is no need to split the sample and condition on $\overline{\Gamma}$, as there are no unknown rotations involved.
    Similar results holds true for spectral norm test statistics 
    \begin{equation}
            \begin{aligned}
            &\Qtspab = \sqrt{\frac{n_an_b}{n_a+n_b}}\| \Pea - \Peb\|,\\ 
            &\Qbspab = \sqrt{\frac{n_an_b}{n_a+n_b}}\|(\Pba-\Pea)-(\Pbb-\Peb)\|,\\
            &\Qfspab = \sqrt{\frac{n_an_b}{n_a+n_b}}\|(\Pfa-\Pea)-(\Pfb-\Peb)\|,
            \nonumber
            \end{aligned}
        \end{equation}
   if we put $s_1 = m$, $s_2 = d-m$ in the error bounds (here $n_a$ and $n_b$ are the sizes of the whole samples $a$ and $b$).
\end{remark}

\subsection{Power analysis}

After we understood the behavior of our procedures under the null hypothesis, we are also interested in the behavior under the alternatives. In particular, the question is whether the power of our procedure goes to $1$ and under which conditions.
We first answer this question for the one-sample test.

\begin{theorem}[One-sample test; power] \label{Th:Power1}
    Under $H_1^{(1)}$ assume $\|\Pt - \Ph\|_{(\Ph, \Gamma^\circ, s_1, s_2)} \geq \lambda_n/\sqrt{n}$, where $\lambda_n$ satisfies
    \begin{equation}
            \begin{aligned}
            \liminf\limits_{n\to\infty}\; \frac{\lambda_n}{\sqrt{n}\left((\scc+\bcc)\rre^{1/2} + (\scc+\bcc)^2\rr^{3/2} \right)} \geq \CONST
            \label{Cond: lambda}
            \end{aligned}
        \end{equation}
     for some absolute constant $\CONST > 0$.
    Then
    \begin{enumerate}[(i)]
        \item the power of Approach~1
        \begin{equation}
            \begin{aligned}
                 \Prob\left[\Qt > \gamma^{(1)}_B(\alpha)\right] \to 1 \;\;\text{ as }n\to \infty;
            \nonumber
            \end{aligned}
        \end{equation}
        \item the power of Approach~2
        \begin{equation}
            \begin{aligned}
                 \Prob\left[\Qt > \gamma^{(1)}_F(\alpha)\right] \to 1 \;\;\text{ as }n\to \infty.
            \nonumber
            \end{aligned}
        \end{equation}
    \end{enumerate}
\end{theorem}
An important question is how restrictive the assumption $\|\Pt - \Ph\|_{(\Ph, \Gamma^\circ, s_1, s_2)} \geq \lambda_n/\sqrt{n}$ is. It would be more natural to assume $\|\Pt - \Ph\| \geq \lambda_n/\sqrt{n}$, which is significantly weaker condition in the worst case, when the bound $\|\Pt - \Ph\| \leq 2\sqrt{\frac{m(d-m)}{s_1 s_2}} \|\Pt - \Ph\|_{(\Ph, \Gamma^\circ, s_1, s_2)}$ is close to being tight. However, due to the first two ``power enhancement'' terms in the definition of $\|\cdot\|_{(\Pp, \Gamma, s_1, s_2)}$, the bound is tight only in very specific cases, while if $\Pt - \Ph$ is random (not adversarially chosen), we can expect $\|\Pt - \Ph\|_{(\Ph, \Gamma^\circ, s_1, s_2)} \asymp\|\Pt - \Ph\|$ (from our numerical experiences). This makes the assumptions of Theorem~\ref{Th:Power1} reasonable.

Now we provide guarantees for power of the two-sample tests.
\begin{theorem}[Two-sample test; power] \label{Th:Power2}
    Under $H_1^{(2)}$ assume
    \begin{equation}
            \begin{aligned}
            \|\Pta - \Ptb\| \geq \lambda_{n_a,n_b} \cdot 2\sqrt{\frac{m(d-m)}{s_1 s_2}} \sqrt{\frac{n_a+n_b}{n_an_b}},
            \nonumber
            \end{aligned}
        \end{equation}
        where $\lambda_{n_a,n_b}$ satisfies
        \begin{equation}
            \begin{aligned}
            \liminf\limits_{n_a,n_b\to\infty} \; \frac{\lambda_{n_a,n_b}}{\sqrt{n_an_b/(n_a+n_b)}\left((\sccab+\bccab)\rre_{a,b}^{1/2} + (\sccab+\bccab)^2\rr_{a,b}^{3/2} \right)} \geq \CONST
            \nonumber
            \end{aligned}
        \end{equation}
     for some absolute constant $\CONST > 0$.
        Then
    \begin{enumerate}[(i)]
        \item the power of Approach~1
        \begin{equation}
            \begin{aligned}
                 \Prob\left[\Qtab > \gamma^{(2)}_B(\alpha)\right] \to 1 \;\;\text{ as }n_a, n_b\to\infty;
            \nonumber
            \end{aligned}
        \end{equation}
        \item the power of Approach~2
        \begin{equation}
            \begin{aligned}
                 \Prob\left[\Qtab > \gamma^{(2)}_F(\alpha)\right] \to 1 \;\;\text{ as }n_a, n_b\to\infty.
            \nonumber
            \end{aligned}
        \end{equation}
    \end{enumerate}
\end{theorem}
One can notice that here, unlike the one-sample case, we make the assumption for the spectral norm
        \begin{equation}
            \begin{aligned}
            \|\Pta - \Ptb\| \geq \lambda_{n_a,n_b} \cdot 2\sqrt{\frac{m(d-m)}{s_1 s_2}} \sqrt{\frac{n_a+n_b}{n_an_b}},
            \nonumber
            \end{aligned}
        \end{equation}
because it is more convenient for the proof. However, we pay the factor $2\sqrt{m(d-m)/(s_1s_2)}$ for avoiding spectral norm as our test statistic. Again, this factor corresponds to worst-case scenario for very specific $\Pta - \Ptb$, while in most cases this condition can be much weaker, i.e.
        \begin{equation}
            \begin{aligned}
            \|\Pta - \Ptb\| \gtrsim \lambda_{n_a,n_b} \sqrt{\frac{n_a+n_b}{n_an_b}}
            \nonumber
            \end{aligned}
        \end{equation}
for the most of non-adversarial choices of pairs of $\Pta$ and $\Ptb$.

Due to the space limitations, the optimality analysis of the presented tests is left for the future work.

%% file: source/FM.tex
Factor model (FM) specifies the data generating process to be
\begin{equation}
    \begin{aligned}
        X_i = \BB \ff_i + \ixi_i\;\;\;\text{for } i\in[n],
     \label{fan1}
    \end{aligned}
\end{equation}
where
\begin{equation}
    \begin{aligned}
        &\BB \in \R^{d\times m} \text{ is deterministic loading matrix},\\
        &\ff_i \in \R^m \text{ is a random vector of $m$ common factors},\\
        &\ixi_i \in \R^d \text{ is a random idiosyncratic component}.
    \nonumber
    \end{aligned}
\end{equation}
If we put $\FF \eqdef [\ff_1, \ldots, \ff_n]^\T \in \R^{n\times m}$ and $\iXi \eqdef [\ixi_1, \ldots, \ixi_n] \in \R^{d\times n}$, the factor model can be rewritten in matrix form
\begin{equation}
    \begin{aligned}
        \data = \BB \FF^\T + \iXi.
    \nonumber
    \end{aligned}
\end{equation}
It is natural to assume that $\{ \ixi_i \}_{i=1}^n$ are uncorrelated with $\{ \ff_i \}_{i=1}^n$.
In addition, for simplicity, we assume $\{ \ff_i \}_{i=1}^n$ are i.i.d. and so are $\{ \ixi_i \}_{i=1}^n$. In literature this assumption is often relaxed to strong mixing condition, allowing weak dependence between pairs of consecutive factors and idiosyncratic components; we believe our general results can be extended to weakly dependent $\Xdata$ as well, however we stick to original i.i.d. framework to avoid technical details.

As in the literature, it is important to mention, that FM is not identifiable. In particular, for any invertible $\HH \in \R^{m \times m}$ it holds $\BB \FF^\T = (\BB \HH)(\HH^{-1} \FF^\T)$, so that the loading matrix $\BB\HH$ and the factors $\FF(\HH^{-1})^\T$ are as good in explaining $\data$ as $\BB$ and $\FF$. However, $\cspan[\BB]$ and $\cspan[\FF]$ are identifiable; indeed, for any $\HH$ as above holds $\cspan[\BB] = \cspan[\BB \HH]$ and $\cspan[\FF] = \cspan\left[\FF (\HH^{-1})^\T\right]$. Our work exploits terminology of spectral projectors rather than subspaces, so we remind that there is one-to-one correspondence between subspace $\cspan[\AAA]$ and the projector $\AAA(\AAA^\T \AAA)^{-1}\AAA^\T$ for any matrix $\AAA$ with linearly independent columns.

\begin{remark}
We also note that usually it is assumed that $\Cov[\ff_1] = \Id_m$ and $\BB^\T \BB$ is diagonal, in order to bring some concreteness in derivations.  This reduces ambiguity of parametrization but does not solve completely the identifiability issue. For our purposes this assumption does not play any role.
\end{remark}

The covariance matrix under this model looks like
\begin{equation}
    \begin{aligned}
        \St = \BB \Cov[\ff_1]\, \BB^\T + \Cov[\ixi_1].
    \label{fan2}
    \end{aligned}
\end{equation}
We will be interested in $m$ principal eigenvectors of $\St$. Define $\J$ so that $\mathcal{I}_\J = \{ 1,\ldots, m \}$, implying that $\Pt$ is the projector onto subspace spanned by $m$ principal eigenvectors of $\St$.

Before formulating specific hypothesis and applying our general scheme, let us recall standard assumptions from FM literature and their implications on the rate in our general framework.
\begin{assumption} \label{FM_assumption}
    The eigenvalues of $\St$ are distinct and there exist absolute constants $L_1, L_2, L_3 > 0$ such that
    \begin{itemize}
        \item $L_1d \geq \mu_1 > \ldots > \mu_m \geq L_2d$.
        \item $L_2 \geq \mu_{m+1} > \ldots > \mu_d \geq L_3$.
    \end{itemize}
\end{assumption}
\noindent This assumption readily implies the following:
\begin{equation}
	\begin{aligned}
		\rr \asymp m^{5/3},\;\;\;\clow \asymp \chigh \asymp \frac{1}{\sqrt{d}},\;\;\;\cond \asymp 1,
	\nonumber
	\end{aligned}
	\end{equation}
so if we were to apply out testing procedure for the $m$-dimensional principal eigenspace, the error rate would be
        \begin{equation}
            \begin{aligned}
            	C \left\{ \left( \frac{(s_1+s_2)^7}{n} \right)^{1/8} + \left( \frac{(s_1+s_2)^{4/\beta+2}}{n} \right)^{1/3}
            	+ m^{5/2} \left( \frac{(s_1+s_2)\,d}{n}\right)^{1/2} \right\},
            \nonumber
            \end{aligned}
        \end{equation}
    or, if $s_1=s_2=1$, simply
      \begin{equation}
            \begin{aligned}
            	C \left\{ \frac{1}{n^{1/8}}  + m^{5/2} \sqrt{\frac{d}{n}} \right\}.
            \nonumber
            \end{aligned}
        \end{equation}

Now we demonstrate how our general framework reduces to testing loading matrices.  From \eqref{fan2}, it is clear that $\BB$ is closely related to the space spanned by the eigenvectors of top $m$ eigenvalues; see
Proposition~\ref{FM1} below.  At the same time, by multiplying $\BB^\T$ on both sides of equation \eqref{fan1}, assuming the that noise is smoothed out, we have $\ff_i \approx (\BB^\T \BB)^{-1} \BB^\T X_i$.  Since the matrix $\BB^\T \BB$ plays only the normalization role, $\BB^\T X_i$ is really the estimate of the latent factors.  Thus, testing $\BB$ lying in a specific space amounts to testing whether the latent factors are the known factors such as the famous Fama-French 3-factor or 5-factor models, see \cite{FF1, FF2}.

Suppose we have some guess $\BB^\circ$ for the unknown underlying loading matrix $\BB$. Define $\Ph = \BB^\circ({\BB^\circ}^\T \BB^\circ)^{-1}{\BB^\circ}^\T$, i.e. $\Ph$ is projector onto $\cspan[\BB^\circ]$. Another scenario could be that instead of $\BB^\circ$ we are given projector $\Ph$ from the very beginning.
The corresponding testing problem writes as  
\begin{equation}
    \begin{aligned}
        H_0:\;\; \cspan[\BB]= \cspan[\BB^\circ]\hspace{1.1cm}
        \text{vs}\hspace{1.1cm}H_1:\;\; \cspan[\BB] \not = \cspan[\BB^\circ]
    \nonumber
    \end{aligned}
\end{equation}
or equivalently
\begin{equation}
    \begin{aligned}
        H_0:\;\;\Ph = \BB(\BB^\T \BB)^{-1} \BB^\T \hspace{1.1cm}\text{vs}\hspace{1.1cm}H_1:\;\;\Ph \neq \BB(\BB^\T \BB)^{-1} \BB^\T
    \nonumber
    \end{aligned}
\end{equation}
The following proposition attempts to bridge the gap between this testing problem and our framework.
\begin{proposition} \label{FM1}
    Under Assumption~\ref{FM_assumption}, it holds
    \begin{equation}
    \begin{aligned}
        \| \Pt - \BB(\BB^\T \BB)^{-1} \BB^\T \| = O\left( \frac{1}{d} \right).
    \nonumber
    \end{aligned}
\end{equation}
\end{proposition}
We notice that the true projector $\Pt$ of $\St$ onto the $m$ principal directions is not exactly corresponds to $\cspan[\BB]$. This is not satisfactory for us, because if we push this additional error term through the proof, we will get another $\sqrt{n/d}$ term in the final bound, while we already have $\sqrt{d/n}$ term. This will make our results meaningless.
However, we can artificially remove the contribution of idiosyncratic components to our underlying eigenspaces by assuming additional conditions on the interaction between factors and idiosyncratic components.
\begin{proposition} \label{FM2}
    Under condition $\Cov[\ixi_1] \, \BB = \Oo_{d\times m}$, it holds
    \begin{equation}
    \begin{aligned}
        \Pt = \BB(\BB^\T \BB)^{-1} \BB^\T.
    \nonumber
    \end{aligned}
\end{equation}
\end{proposition}

\noindent This allows to rewrite the hypothesis in familiar form:
\begin{equation}
    \begin{aligned}
        H_0:\;\;\Pt = \Ph\hspace{1.3cm}\text{vs}\hspace{1.3cm}H_1:\;\;\Pt \neq \Ph
    \nonumber
    \end{aligned}
\end{equation}
Now we can directly apply Algorithm~\ref{Algorithm}.
Our procedure uses $\Pe$, which naturally arises in FM context, since in POET (see \cite{POET}) $\BB$ is estimated by
\begin{equation}
    \begin{aligned}
        \widehat{\BB} \eqdef [\widehat{\sigma}_1 \widehat{u}_1, \ldots, \widehat{\sigma}_m \widehat{u}_m],
    \nonumber
    \end{aligned}
\end{equation}
leading exactly to $\widehat{\BB}(\widehat{\BB}^\T \widehat{\BB})^{-1} \widehat{\BB}^\T = \Pe$.  As before, one can use the test statistic  $\Qt = \sqrt{n}\| \Pe - \Ph\|_{(\Ph, \Gamma^\circ, s_1,s_2)}$, whose distribution can be approximated by one of two approaches.

Likewise, the two-sample problem can be solved. As explained above, this amounts to test whether the latent factors are the same between two groups (e.g. treatment vs. control, pre-financial crisis vs. post financial crisis).  If we have two samples $\Xdataa$ and $\Xdatab$ generated from two FMs, e.g. with $\BB^a, \FF^a$ and $\BB^b, \FF^b$, respectively, and we want to understand whether their loading matrices span the same subspace, this is equivalent to testing
\begin{equation}
    \begin{aligned}
        H_0:\;\;\Pta=\Ptb\hspace{1.3cm}\text{vs}\hspace{1.3cm}H_1:\;\;\Pta \neq \Ptb\,
    \nonumber
    \end{aligned}
\end{equation}
where $\Pta$ and $\Ptb$ are projectors onto $m$ principal eigenspaces of underlying covariance matrices of the samples $a$ and $b$, respectively. Algorithm~\ref{Algorithm2} can be employed to deal with such a problem. 

%% file: source/Numerical.tex
\subsection{Construction of the covariance matrix for null and alternative hypothesis}
In order to clearly describe the setup of our experiments, we start with a toy example on a plane, i.e. $d=2$.
Suppose we are interested in testing whether the first principal direction of our $2$-dimensional data is aligned with the first axis of our coordinate system, i.e. $\J = \{ 1\}$, $m=1$ and the hypothetical leading eigenvector is $u^\circ = [1, 0]^\T$. In this case, the spectral projector is
\begin{equation}
    \begin{aligned}
        \Ph = \begin{bmatrix}
        1 & 0 \\ 0 & 0
        \end{bmatrix}.
    \nonumber
    \end{aligned}
\end{equation}
To empirically study the power of our method, we generate the data in such a way that its leading eigenvector is obtained by rotating $u^\circ$ by angle $\varphi$, i.e. $u_{\varphi} = [\cos\varphi, \sin\varphi]^\T$, while the orthogonal direction is given by $v_{\varphi} = [-\sin\varphi, \cos\varphi]^\T$. The associated true projectors are
\begin{equation}
    \begin{aligned}
        \Pp_1 = u_{\varphi} u_{\varphi}^\T = \begin{bmatrix}
        \cos^2\varphi & \cos\varphi\,\sin\varphi \\ \cos\varphi\,\sin\varphi & \sin^2\varphi
        \end{bmatrix},\;\;
        \Pp_2 = v_{\varphi} v_{\varphi}^\T = \begin{bmatrix}
        \sin^2\varphi & -\cos\varphi\,\sin\varphi \\ -\cos\varphi\,\sin\varphi & \cos^2\varphi
        \end{bmatrix}.
    \nonumber
    \end{aligned}
\end{equation}
If the corresponding variances along these directions (eigenvalues of the covariance matrix) are $\mu_1$ and $\mu_2$ ($\mu_1 > \mu_2$), then the true covariance matrix of the data is formed as
\begin{equation}
    \begin{aligned}
        \St^{(\varphi)} = \mu_1 \Pp_1 + \mu_2 \Pp_2 =
        \begin{bmatrix}
        \mu_1\cos^2\varphi + \mu_2 \sin^2\varphi & (\mu_1-\mu_2)\,\cos\varphi\,\sin\varphi \\ (\mu_1-\mu_2)\,\cos\varphi\,\sin\varphi & \mu_1\sin^2\varphi + \mu_2 \cos^2\varphi
        \end{bmatrix}.
    \nonumber
    \end{aligned}
\end{equation}
Thereby, $\varphi = 0$ corresponds to null hypothesis $\Ph = \Pp_1$, while under the alternative the larger deviations of angle $\varphi$ from $0$ (until one point) mean that $\Pp_1$ is further from $\Ph$.
So, the suitable data for our experiment can be generated from some distribution with the covariance matrix $\St^{(\varphi)}$ with varying $\varphi$.
Our goal is to verify that with growing $\varphi$ our methods reject the null hypothesis more often and eventually this probability approaches $1$ and check the size of the test when $\varphi = 0$.
\begin{figure}
\begin{center}
\includegraphics[scale=0.3]{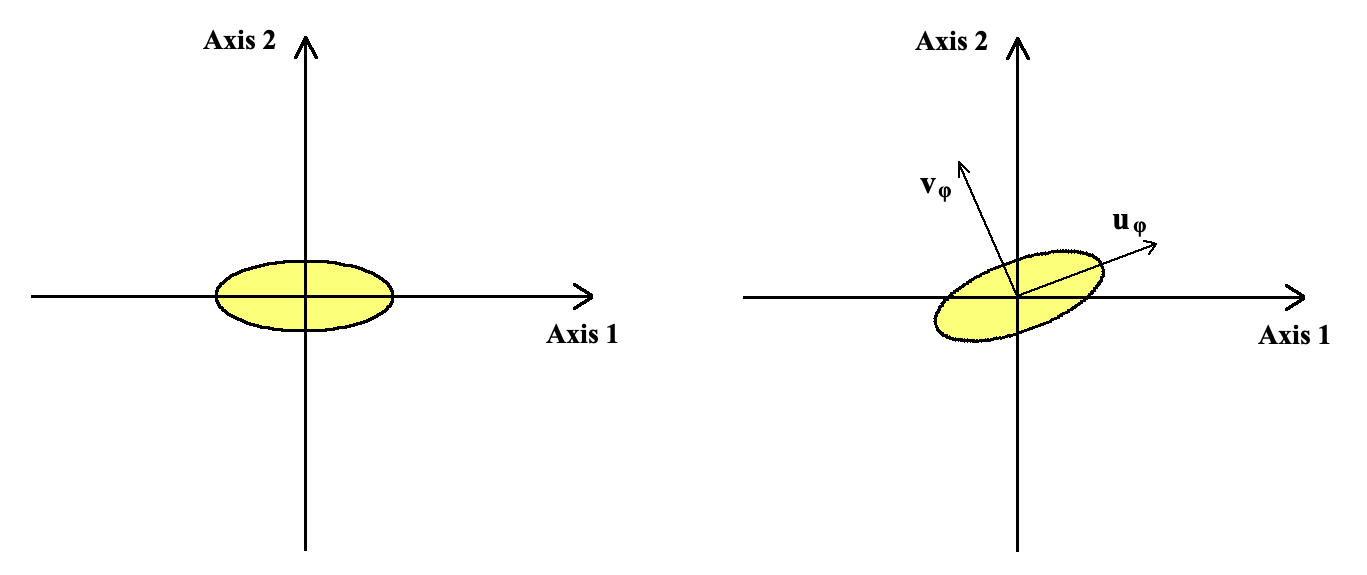}
\caption{Construction of alternative hypothesis data. Here in yellow we depict sub-level sets of Gaussian density with mean zero and covariance $\St^{(0)}$ and $\St^{(\varphi)}$ in order to understand how the clouds of data look like in each case.}
\label{cloud}
\end{center}
\end{figure}
See Figure~\ref{cloud} for visualization of the construction.

Now we extend this setting to higher dimensions. In dimension $d$ we are interested in the subspace spanned by $m$ leading eigenvectors, i.e. $\J = \{ 1, \ldots, m\}$ (for simplicity we explain the procedure assuming all eigenvalues are distinct; it will be clear how to extend the construction to the case of multiplicities within the first $m$ eigenvalues and within the last $(d-m)$ eigenvalues). Without loss of generality, we assume that under null hypothesis the eigenvectors are aligned with the axes of the coordinate system, so that the hypothetical projector is
\begin{equation}
    \begin{aligned}
        \Ph = \begin{bmatrix}
        \Id_m & \Oo_{m\times(d-m)}\\
        \Oo_{(d-m)\times m} & \Oo_{(d-m)\times(d-m)}
        \end{bmatrix},
    \nonumber
    \end{aligned}
\end{equation}
and the default covariance matrix is diagonal with descending entries and characterized only by eigenvalues $\mu_1 > \ldots > \mu_m > \mu_{m+1} > \ldots > \mu_d$:
\begin{equation}
    \begin{aligned}
        \St^{(0)} = \begin{bmatrix}
        \mu_1 & & & & &\\
         & \ddots & & & &\\
         & & \mu_m & & &\\
         & & & \mu_{m+1} & &\\
         & & & & \ddots &\\
         & & & & & \mu_d\\
        \end{bmatrix},
    \nonumber
    \end{aligned}
\end{equation}
To generate the data under alternative, we rotate the plane containing the first and $(m+1)$-th axes by the angle $\varphi$, i.e. the leading eigenvector becomes
\begin{equation}
    \begin{aligned}
        u_\varphi = [ \underbrace{\cos\varphi, 0, \ldots, 0}_{m}, \sin\varphi, 0, \ldots, 0]^\T,
    \nonumber
    \end{aligned}
\end{equation}
while $(m+1)$-th eigenvector turns into
\begin{equation}
    \begin{aligned}
        v_\varphi = [ \underbrace{-\sin\varphi, 0, \ldots, 0}_{m}, \cos\varphi, 0, \ldots, 0]^\T.
    \nonumber
    \end{aligned}
\end{equation}
The covariance matrix is then
\begin{equation}
    \begin{aligned}
        \St^{(\varphi)} = \mu_1 u_\varphi u_\varphi^\T + \sum\limits_{r=2}^m \mu_r \Pp_r + \mu_{m+1} v_\varphi v_\varphi^\T + \sum\limits_{r=m+2}^d \mu_r \Pp_r,
    \nonumber
    \end{aligned}
\end{equation}
or explicitly
\begin{equation}
    \begin{aligned}
        \St^{(\varphi)} = \begin{bmatrix}
         \mu_1\cos^2\varphi + \mu_{m+1} \sin^2\varphi & 0 & \ldots & 0 & (\mu_1-\mu_{m+1})\,\cos\varphi\,\sin\varphi & 0 & \ldots & 0\\
         0& \mu_2 & & & 0 & & &\\
         \vdots& & \ddots & & \vdots & & &\\
         0& & & \mu_m& 0& & &\\
         (\mu_1-\mu_{m+1})\,\cos\varphi\,\sin\varphi & 0 & \ldots  & 0 & \mu_1\sin^2\varphi + \mu_{m+1} \cos^2\varphi  & & &\\
         0& & & & & \mu_{m+2} & &\\
         \vdots& & & & & & \ddots &\\
         0& & & & & & & \mu_d\\
        \end{bmatrix}.
    \nonumber
    \end{aligned}
\end{equation}

\subsection{Description of regimes and data distributions}
In our experiments we focus on three regimes:
\begin{itemize}
    \item Factor Model regime: we take $m=8$, and $\mu_1 = 5 d$, $\mu_2 = 4 d$, $\mu_3 = 3.5 d$, $\mu_4 = 3 d$, $\mu_5 = 2.5 d$, $\mu_6 = 2 d$, $\mu_7 = 1.5 d$, $\mu_8 = d$ and $\mu_{9}, \ldots, \mu_d$ uniformly distributed in $[0.5; 1.5]$ and sorted.
    \item Spiked regime: we take $m=1$ and $\mu_1 = 10$, $\mu_2 = 6$, $\mu_3 = 3$, $\mu_{4} = 1$ (of multiplicity $d-3$).
    \item Decay regime: we take $m=5$ and $\mu_1 = 10$, $\mu_2 = 9$, $\mu_3 = 8$, $\mu_4 = 7$, $\mu_5 = 6$, $\mu_k = 2^{-(k-6)}$ for $k = 6,\ldots, d$.
\end{itemize}
We consider two types of data distributions:
\begin{itemize}
    \item Gaussian distribution with mean zero and proper covariance $\St$.
    \item Laplace distribution: we generate components of vector $\widetilde{X}$ as independent Laplace r.v.'s with scale parameter $1/\sqrt{2}$ (so that each component has unit variance), and then put our observation $X = \St^{1/2} \widetilde{X}$ (so that $X$ has covariance matrix $\St$).
\end{itemize}
Once we fix a regime, a data distribution and methods that we want to compare, we conduct the simulations for the sample size in range $n\in \{ 500, 1500, 5000 \}$ and the dimension in range $d\in\{ 50, 150\}$.  Significance level is fixed to be $\alpha=0.05$.
In one-sample problem, for each $n$, $d$ we perform the following:
\begin{itemize}
    \item We try a number of angles $\varphi$ (including $\varphi = 0$) --- they are chosen differently in different settings in order to illustrate the transition of the power from $\alpha$ to $1$.
    \item For each nonzero angle $\varphi$ we generate $100$ samples $\data$ of size $n$ in dimension $d$ with the covariance matrix $\St^{(\varphi)}$ specified by the formula above and regime. For angle $\varphi = 0$ we generate $1000$ samples, since it is important to estimate type-I error accurately.
    \item For each sample we apply each method to test hypothesis $\Pt = \Ph$ vs $\Pt \neq \Ph$. Since some of the methods are resampling-based, we fix the number of resampling $N = 2000$.
    \item We estimate the power (for non-zero angles) and type-I error (for $\varphi=0$) of the test simply as the fraction of samples, for which the null hypothesis was rejected.
\end{itemize}
The steps for two-sample problem are similar, but $\dataa$ generated from distribution with covariance matrix $\St^{(\varphi)}$ and $\datab$ generated from distribution with covariance matrix $\St^{(-\varphi)}$.
The testing problem is changed accordingly.

\subsection{Experimental results}
Now we describe three scenarios. In each scenario we compare our methods with the methods from previous literature, suitable for each particular situation.
\\
\noindent\textbf{Scenario 1}\\
In this scenario we consider one-sample problem in Factor Model regime with Laplace distribution.
We compare the following methods already discussed above:
\begin{itemize}
    \item ``\textsf{Fr-Bootstrap}'': Frobenius norm test statistic + Bootstrap (\cite{Naumov}).
    \item ``\textsf{Fr-Bayes}'': Frobenius norm test statistic + Frequentist Bayes (\cite{Silin_1}).
    \item ``\textsf{Spectral-1}'': Spectral norm test statistic $\Qtsp$ + Approach 1.
    \item ``\textsf{Spectral-2}'': Spectral norm test statistic $\Qtsp$ + Approach 2.
    \item ``\textsf{New-1}'': $\| \cdot \|_{(\Pp, \Gamma, s_1, s_2)}$-norm test statistic $\Qt$ (with $s_1 = s_2 = 1$) + Approach 1.
    \item ``\textsf{New-2}'': $\| \cdot \|_{(\Pp, \Gamma, s_1, s_2)}$-norm test statistic $\Qt$ (with $s_1 = s_2 = 1$) + Approach 2.
\end{itemize}
In some of the settings (with relatively large $n$ and $d$) we are not able to run bootstrap based methods (``\textsf{Fr-Bootstrap}'', ``\textsf{Spectral-1}'', ``\textsf{New-1}''), since their computational time is too intensive.

\begin{figure}
     \begin{subfigure}{0.32\textwidth}
         \centering
         \includegraphics[width=\textwidth]{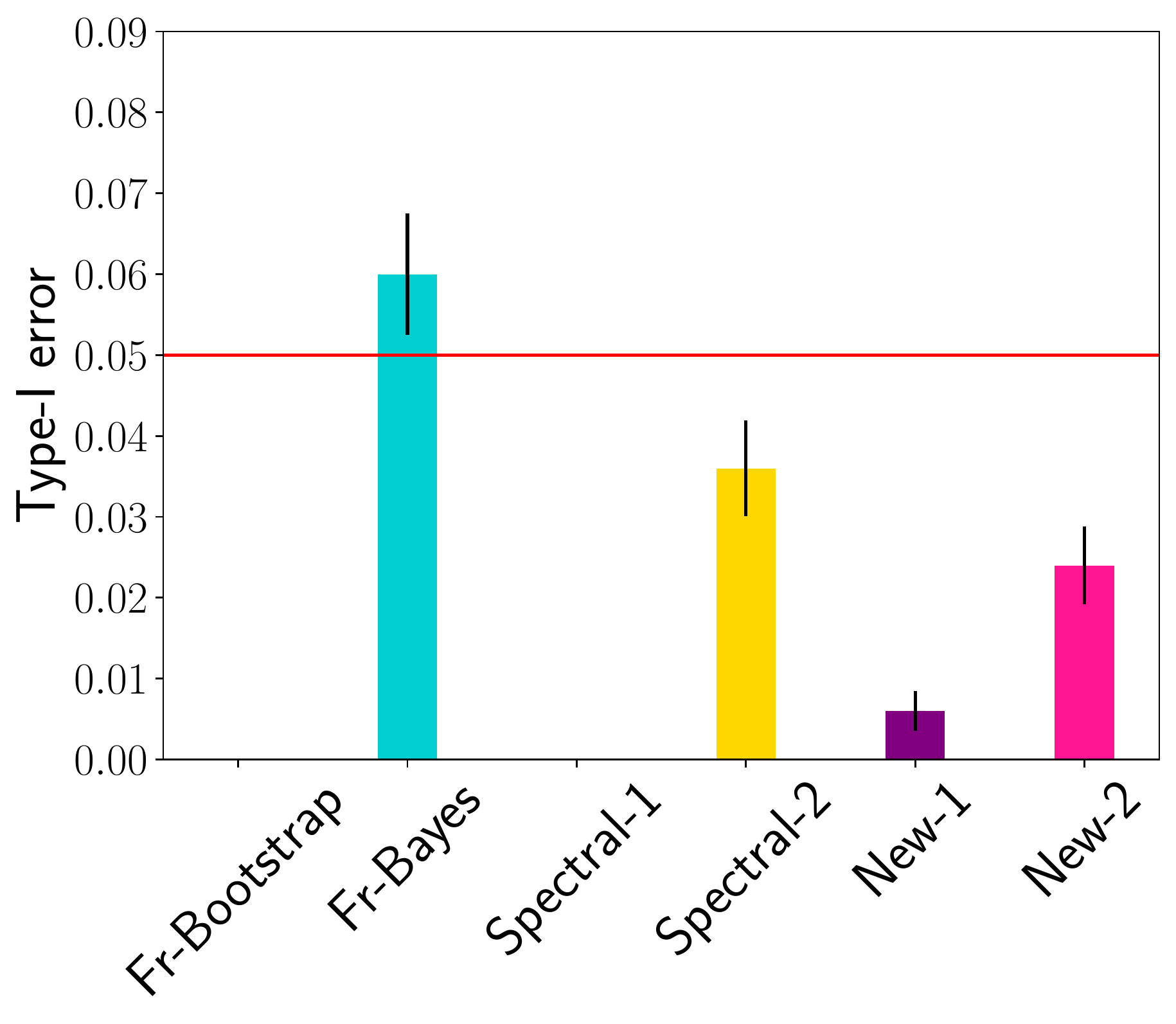}
         \vfill
         \includegraphics[width=\textwidth]{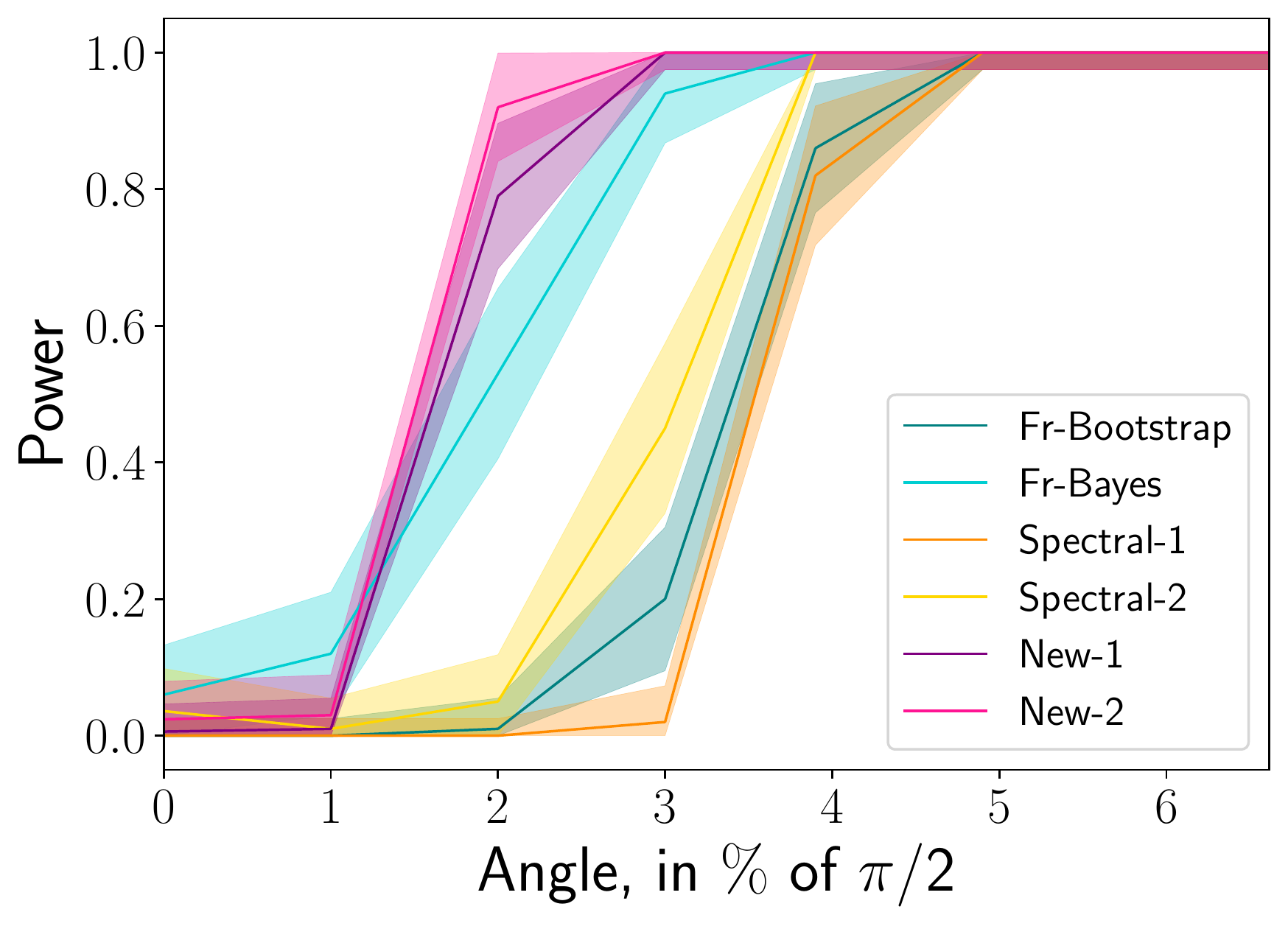}
         \caption{$n=500, \;d=50$}
     \end{subfigure}
     \hfill
     \begin{subfigure}{0.32\textwidth}
         \centering
         \includegraphics[width=\textwidth]{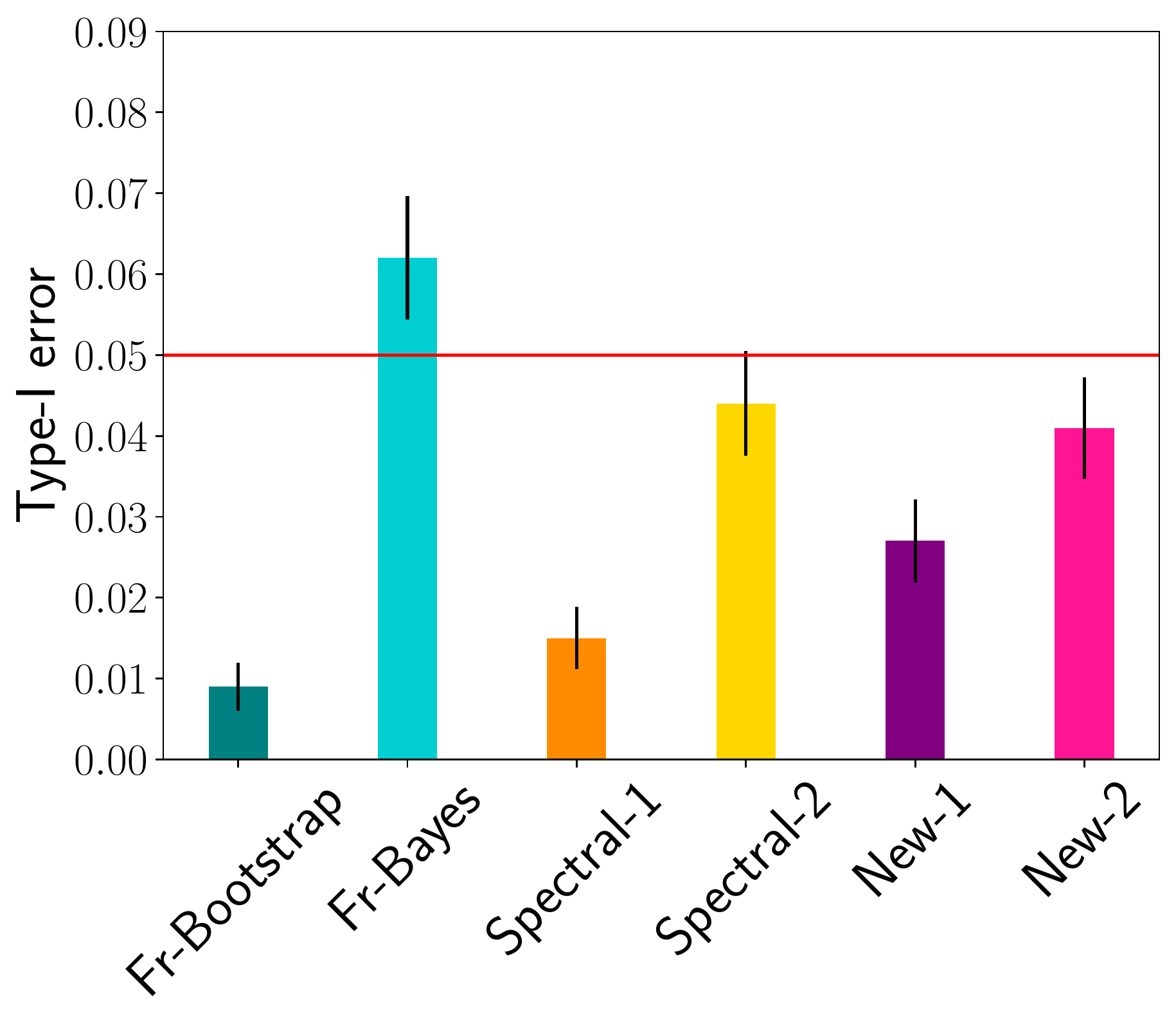}
         \vfill
         \includegraphics[width=\textwidth]{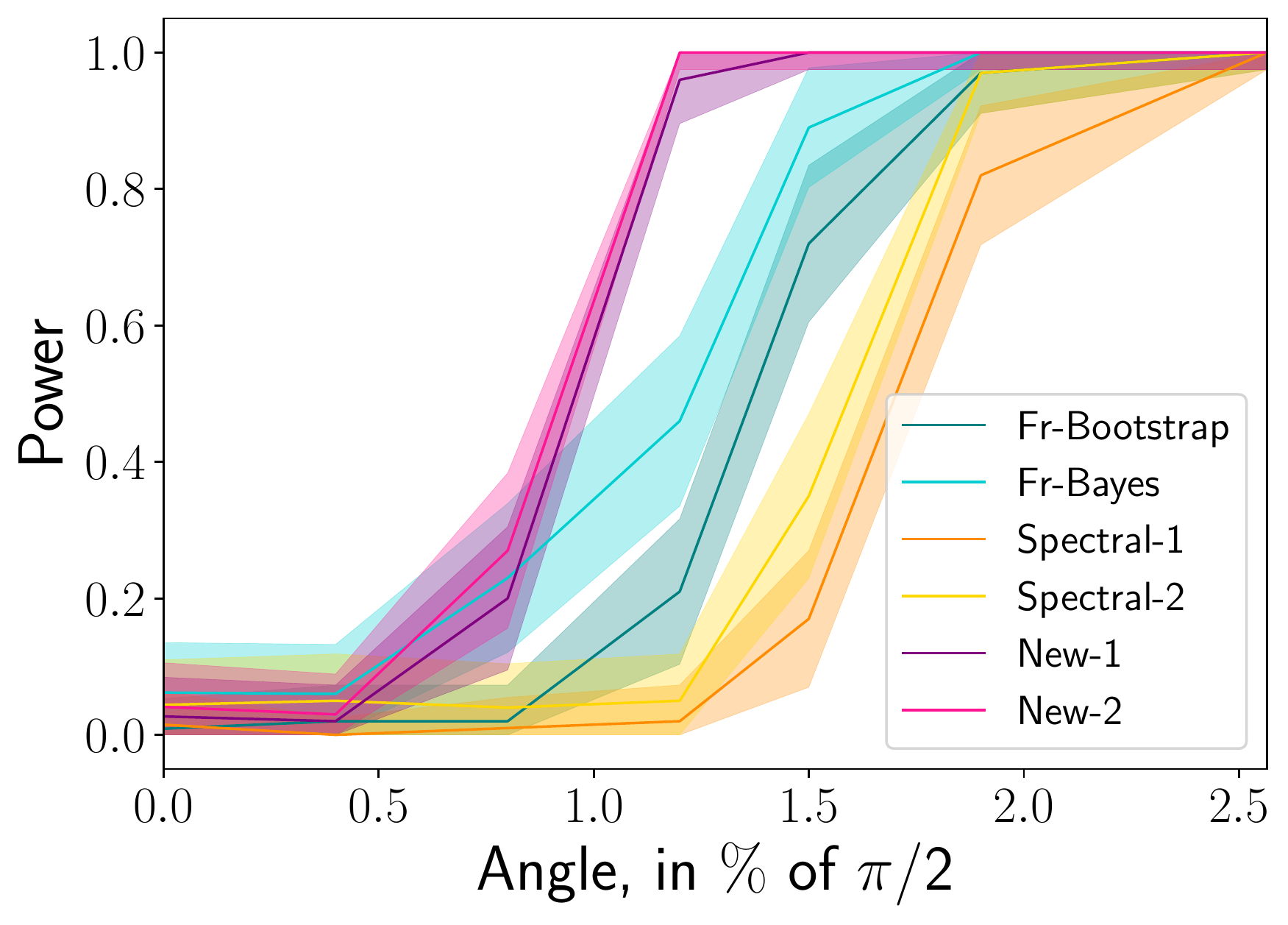}
         \caption{$n=1500, \;d=50$}
     \end{subfigure}
     \hfill
     \begin{subfigure}{0.32\textwidth}
         \centering
         \includegraphics[width=\textwidth]{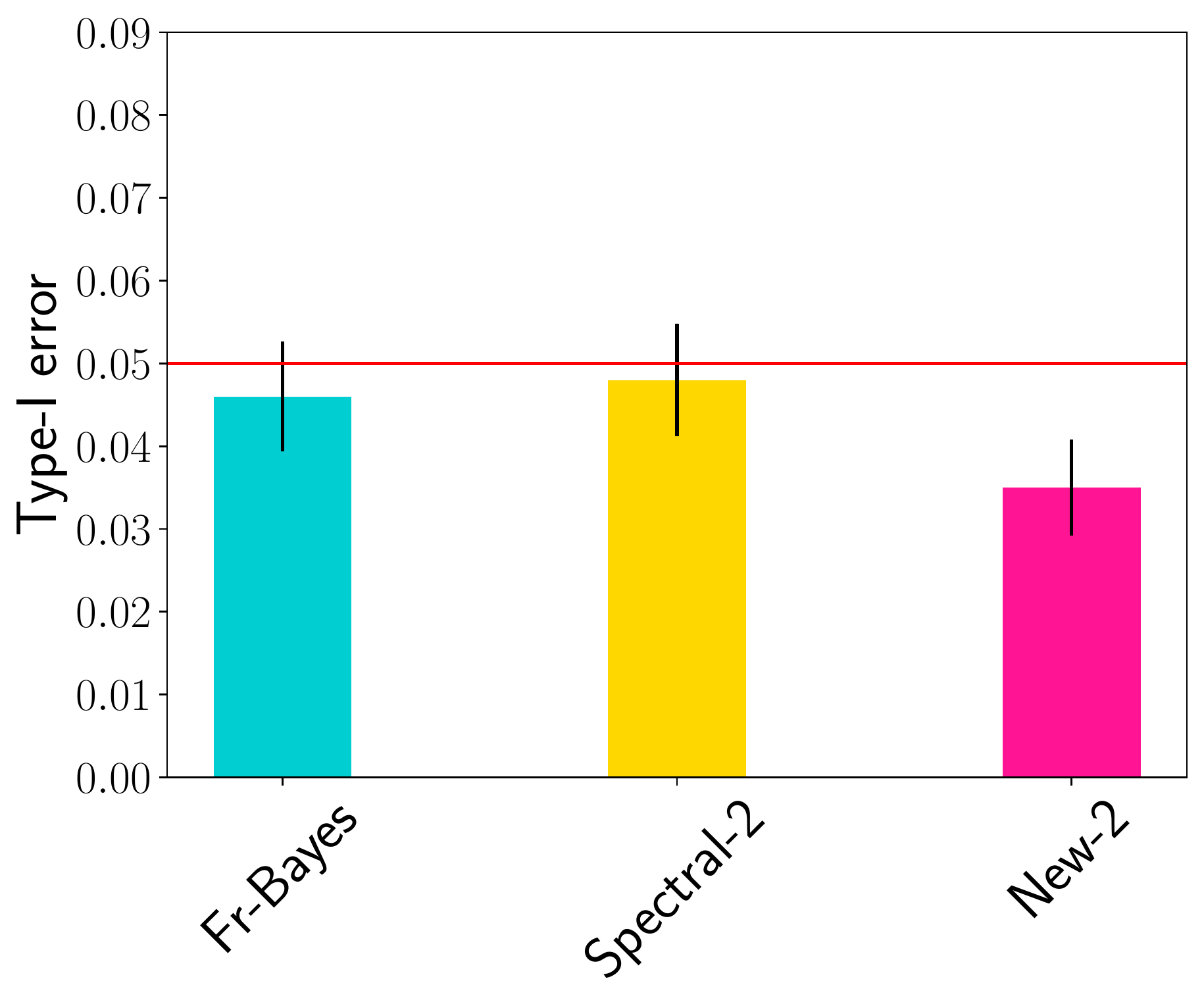}
         \vfill
         \includegraphics[width=\textwidth]{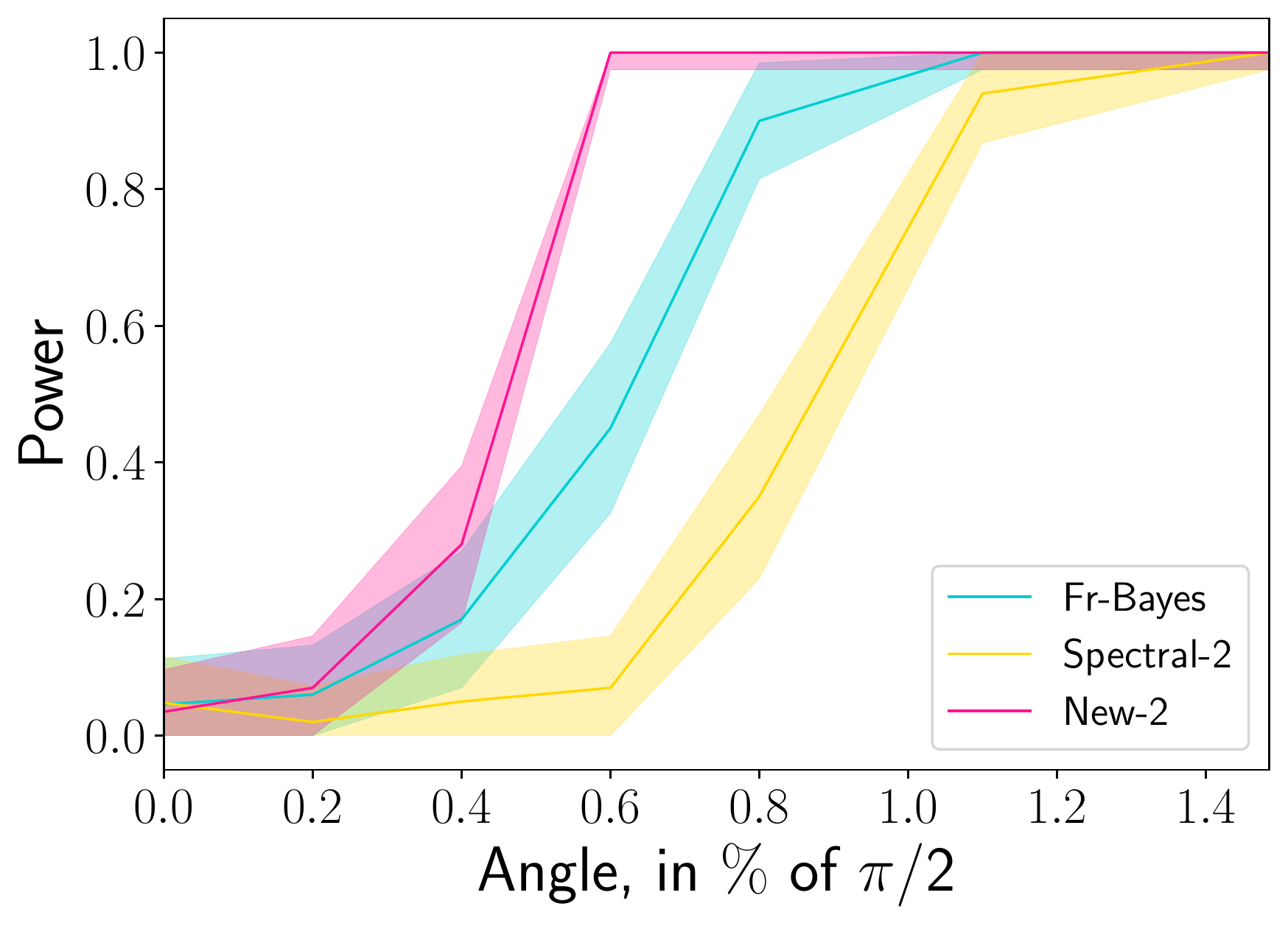}
         \caption{$n=5000, \;d=50$}
     \end{subfigure}
     \vfill
     \vspace{1.5cm}
     \begin{subfigure}{0.32\textwidth}
         \centering
         \includegraphics[width=\textwidth]{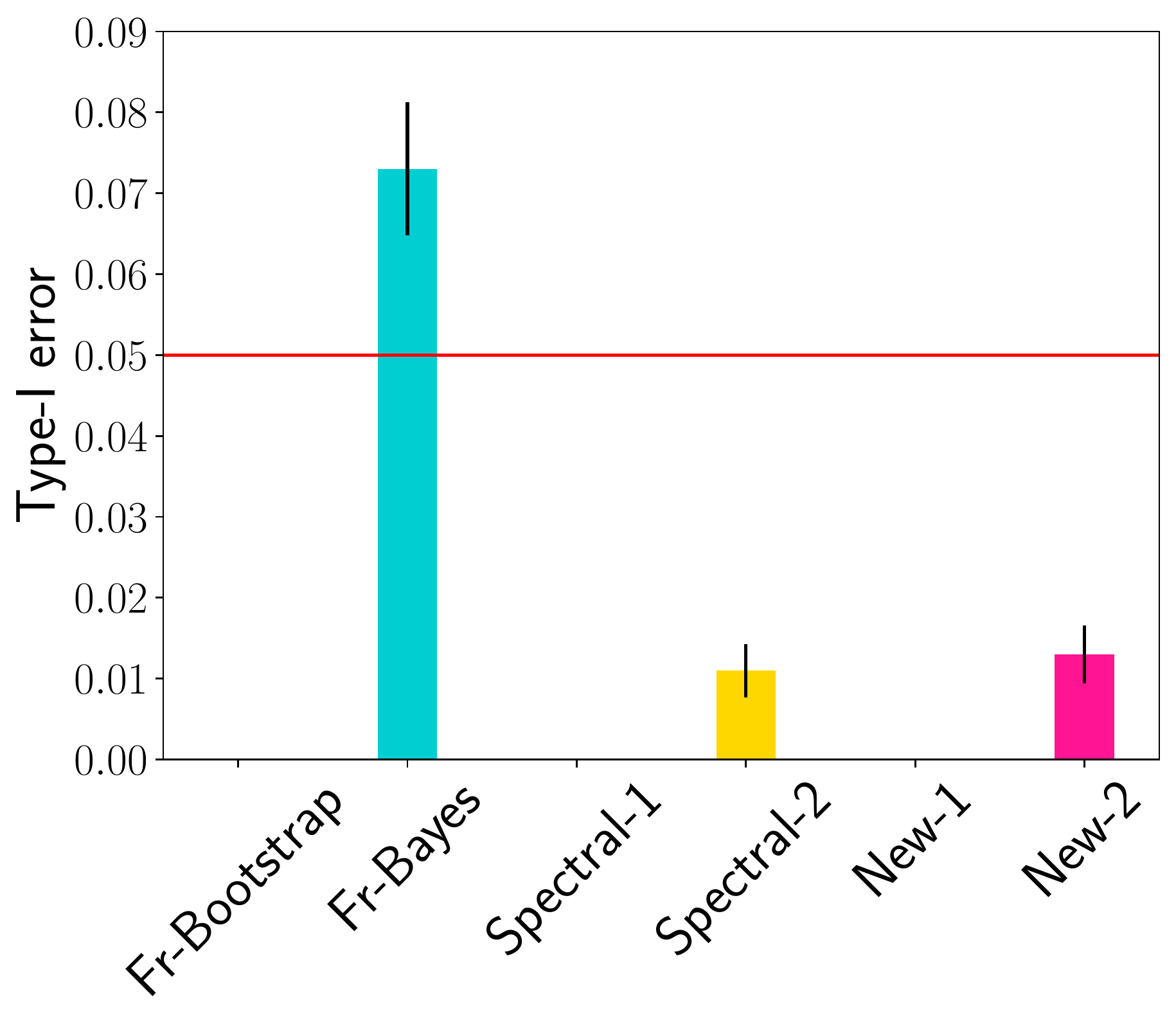}
         \vfill
         \includegraphics[width=\textwidth]{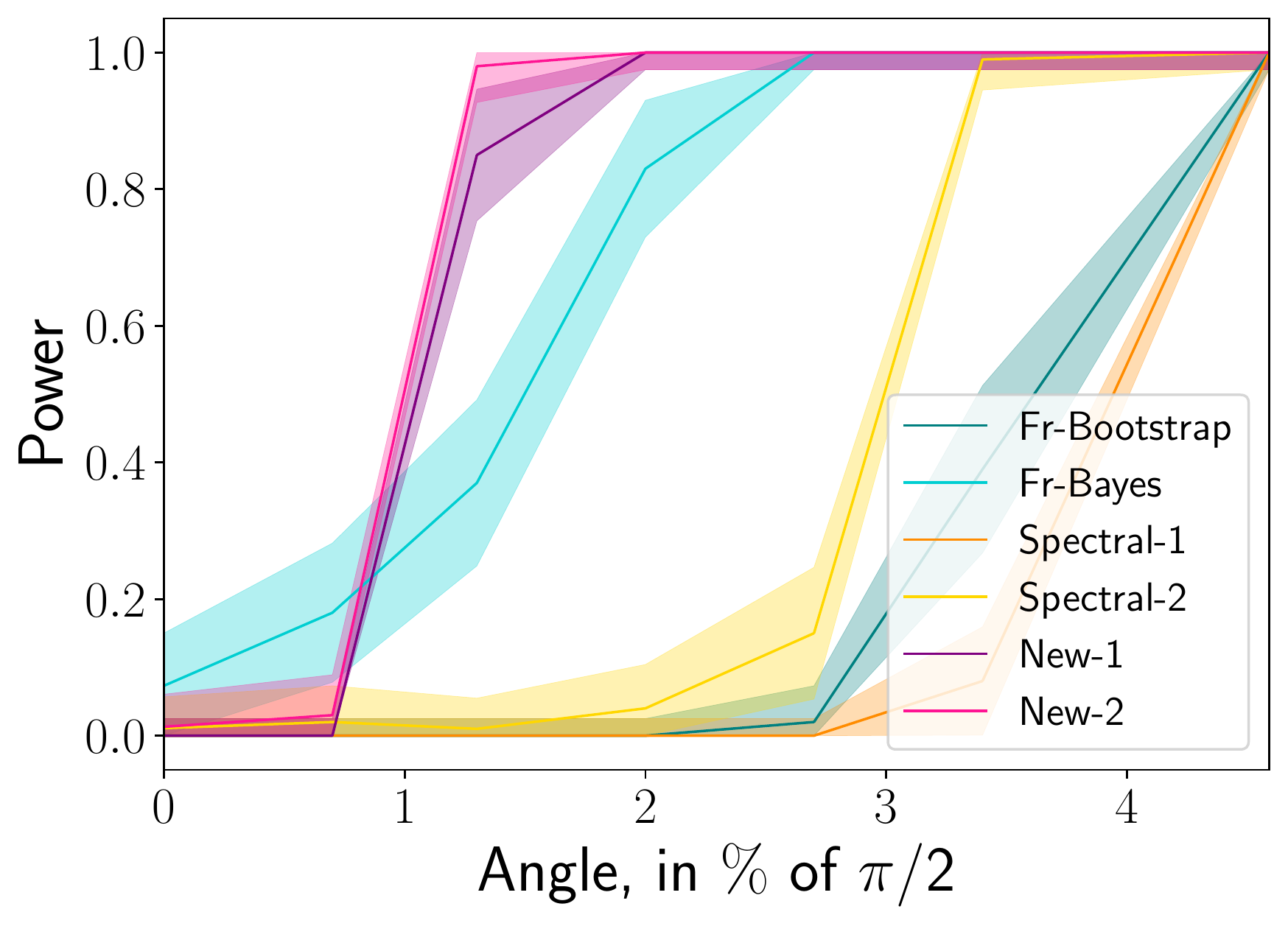}
         \caption{$n=500, \;d=150$}
     \end{subfigure}
     \hfill
     \begin{subfigure}{0.32\textwidth}
         \centering
         \includegraphics[width=\textwidth]{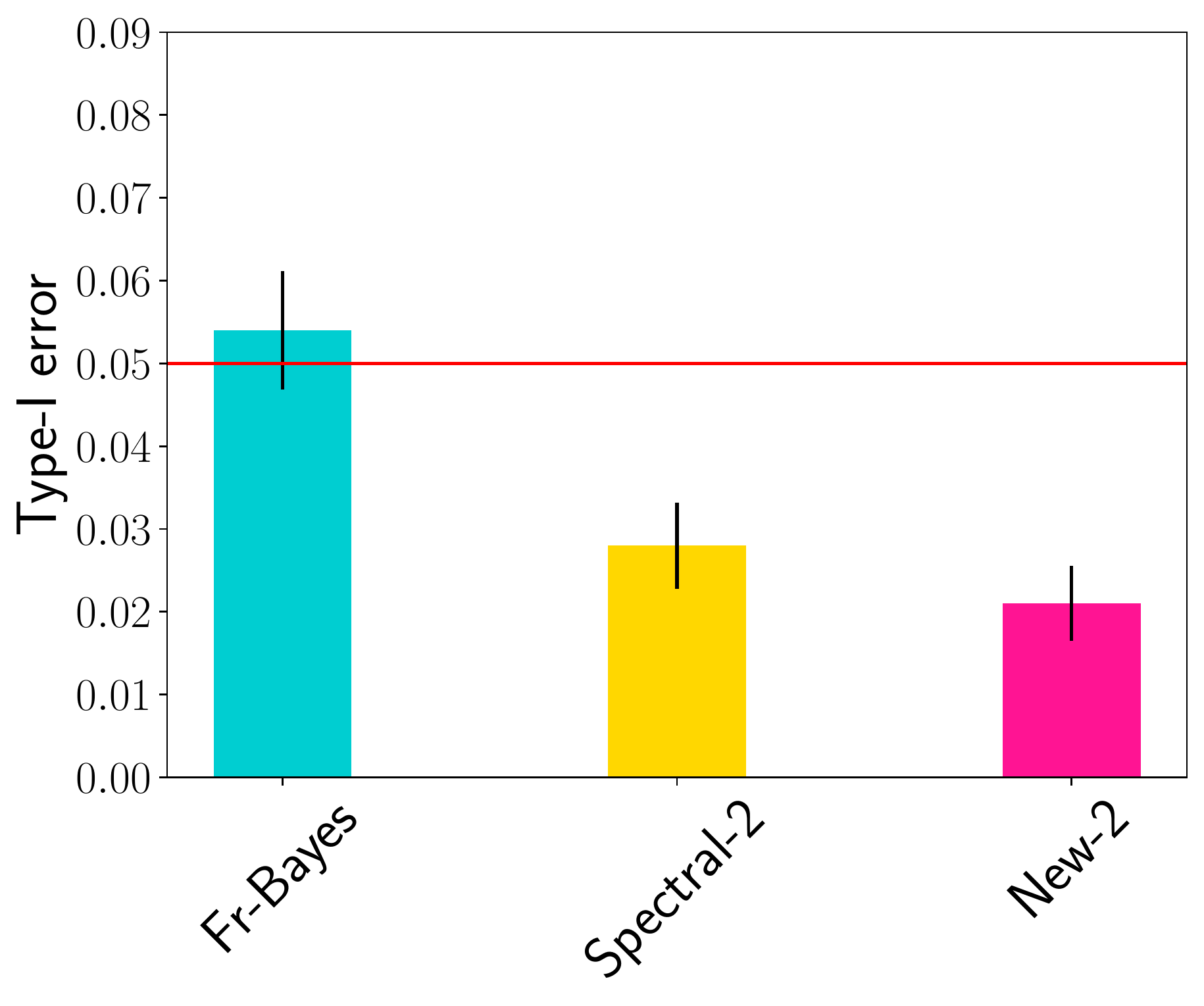}
         \vfill
         \includegraphics[width=\textwidth]{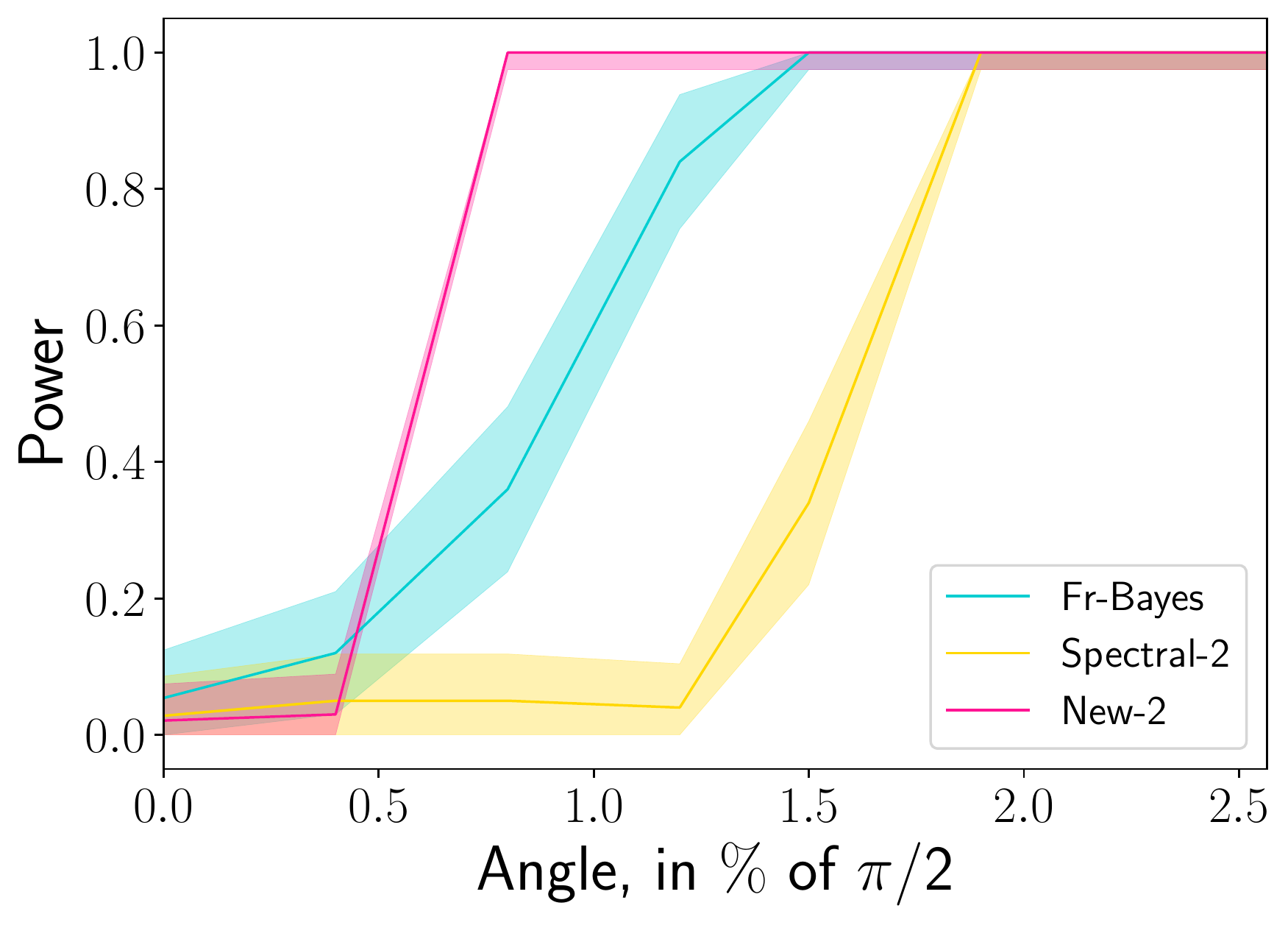}
         \caption{$n=1500, \;d=150$}
     \end{subfigure}
     \hfill
     \begin{subfigure}{0.32\textwidth}
         \centering
         \includegraphics[width=\textwidth]{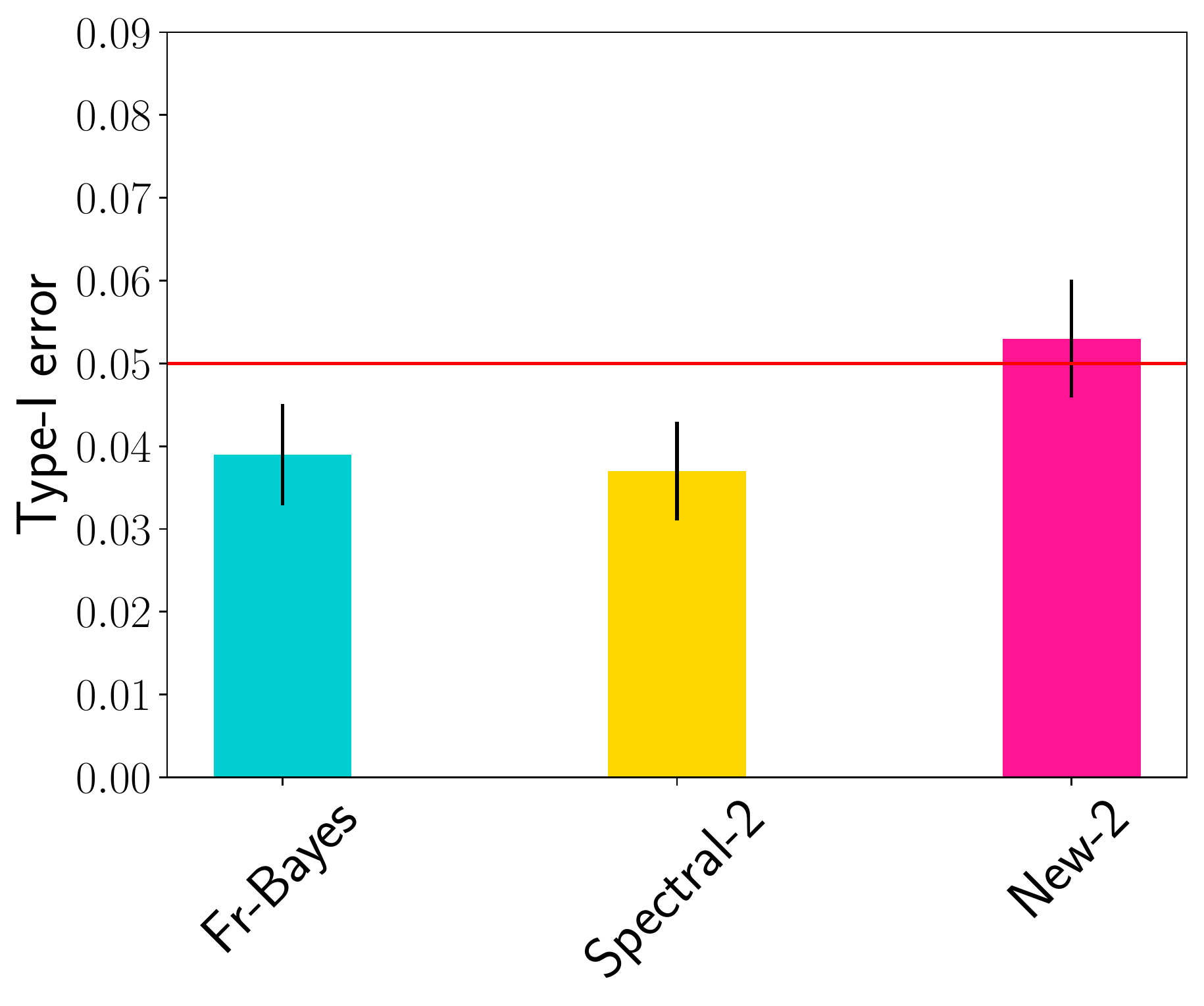}
         \vfill
         \includegraphics[width=\textwidth]{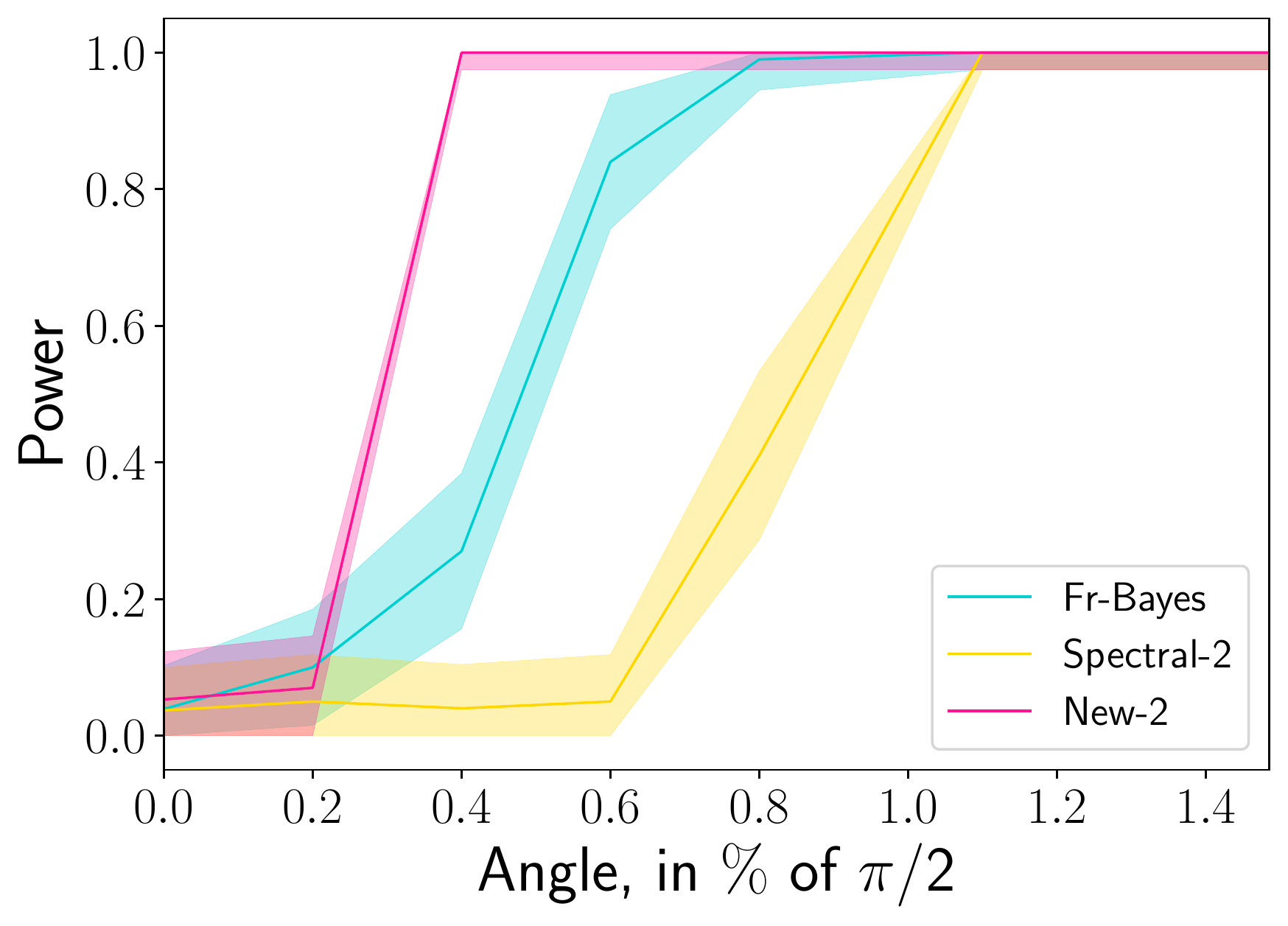}
         \caption{$n=5000, \;d=150$}
     \end{subfigure}
        \caption{Experiments for Scenario 1: One-sample problem, FM regime with $m=8$, Laplace distribution. }
        \label{Exp1}
\end{figure}

\begin{figure}
     \begin{subfigure}{0.32\textwidth}
         \centering
         \includegraphics[width=\textwidth]{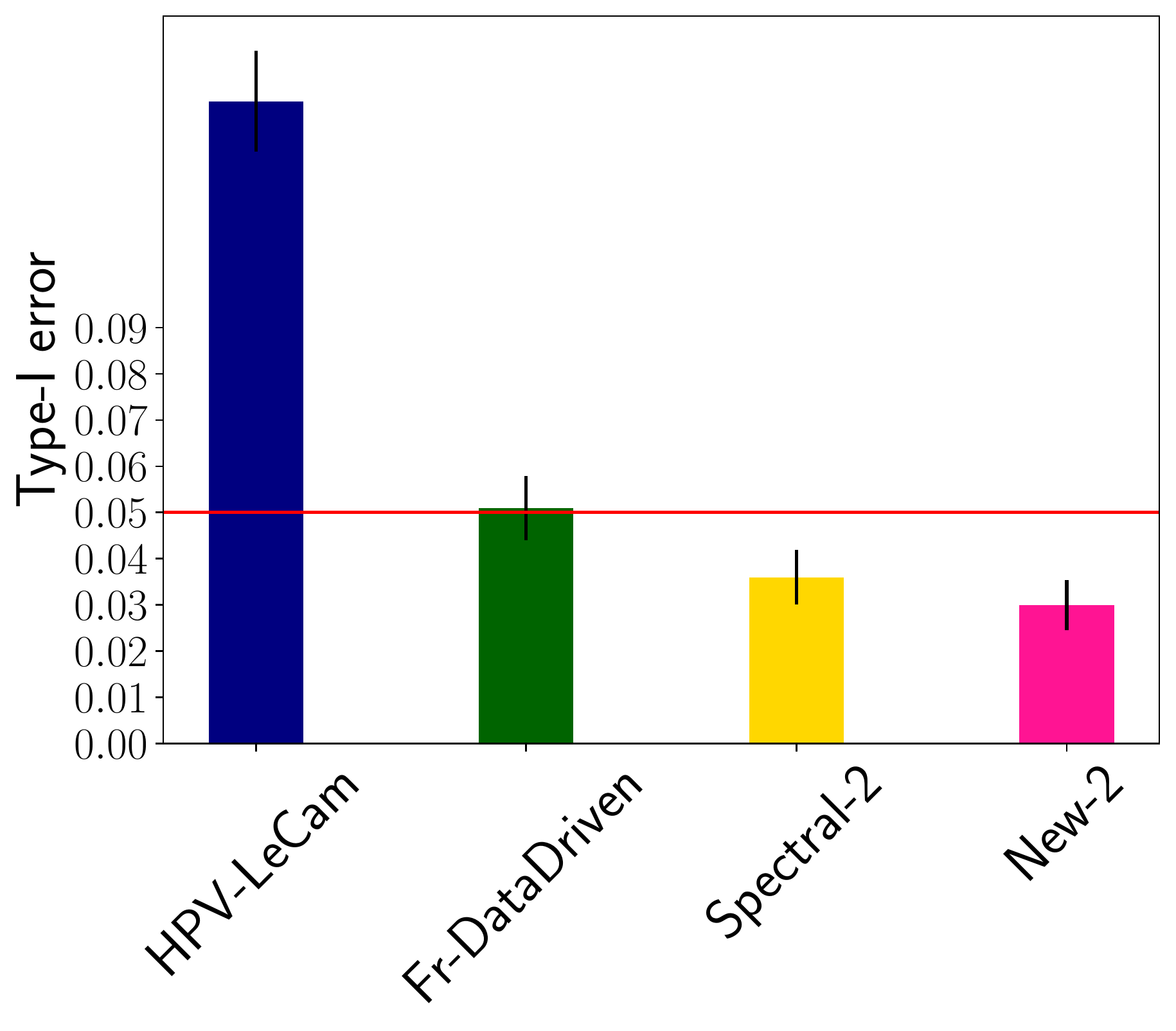}
         \vfill
         \includegraphics[width=\textwidth]{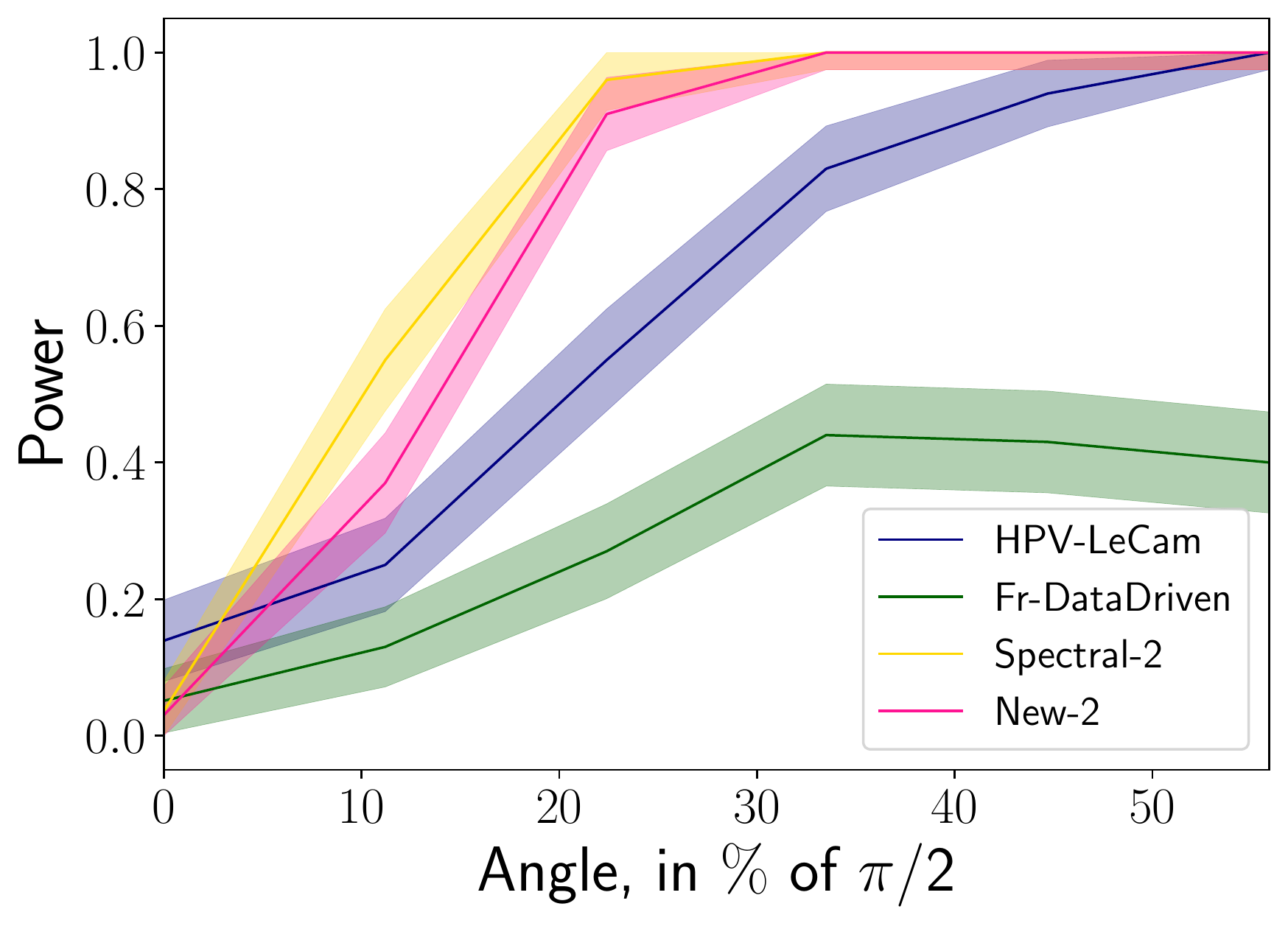}
         \caption{$n=500, \;d=50$}
     \end{subfigure}
     \hfill
     \begin{subfigure}{0.32\textwidth}
         \centering
         \includegraphics[width=\textwidth]{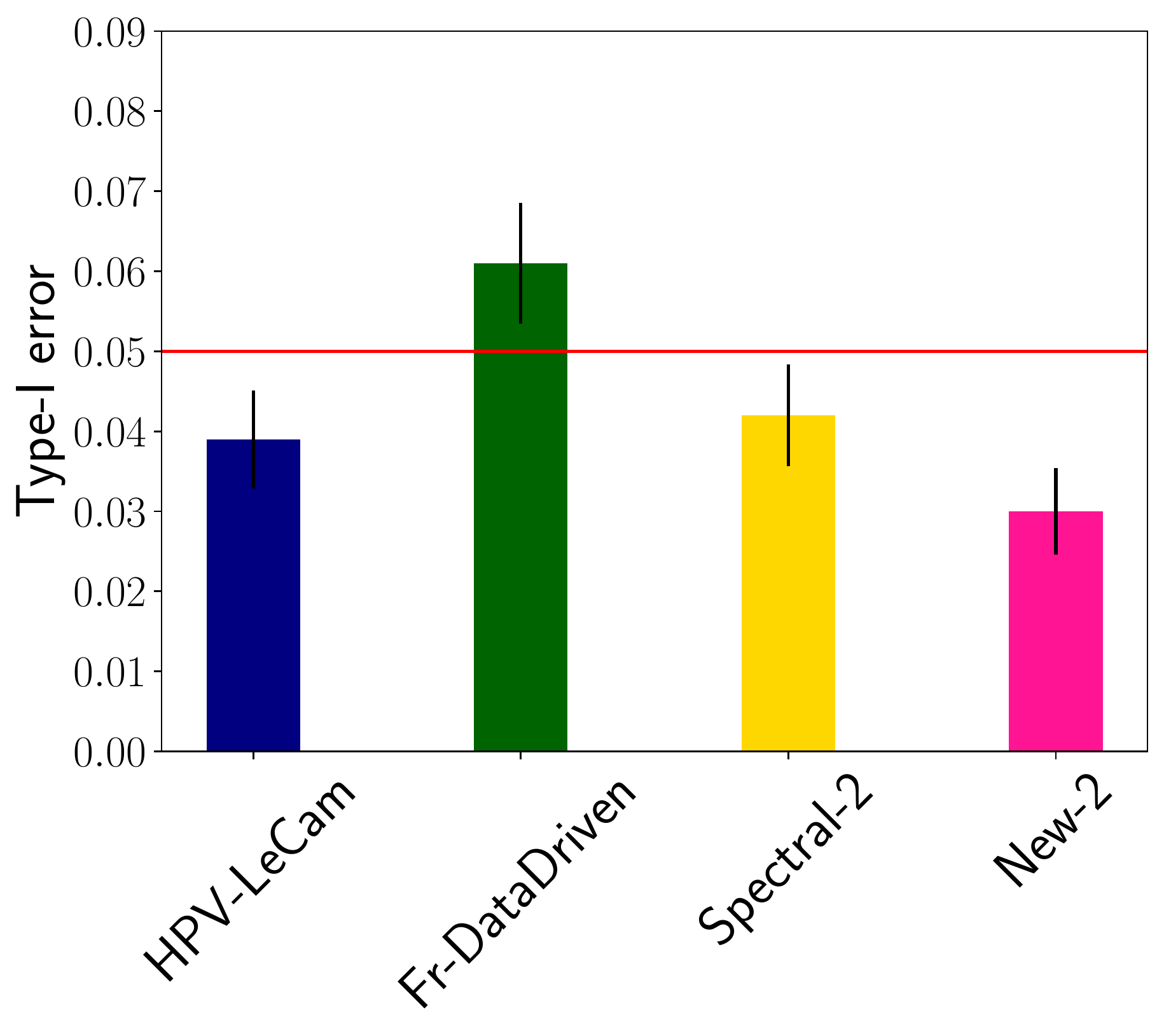}
         \vfill
         \includegraphics[width=\textwidth]{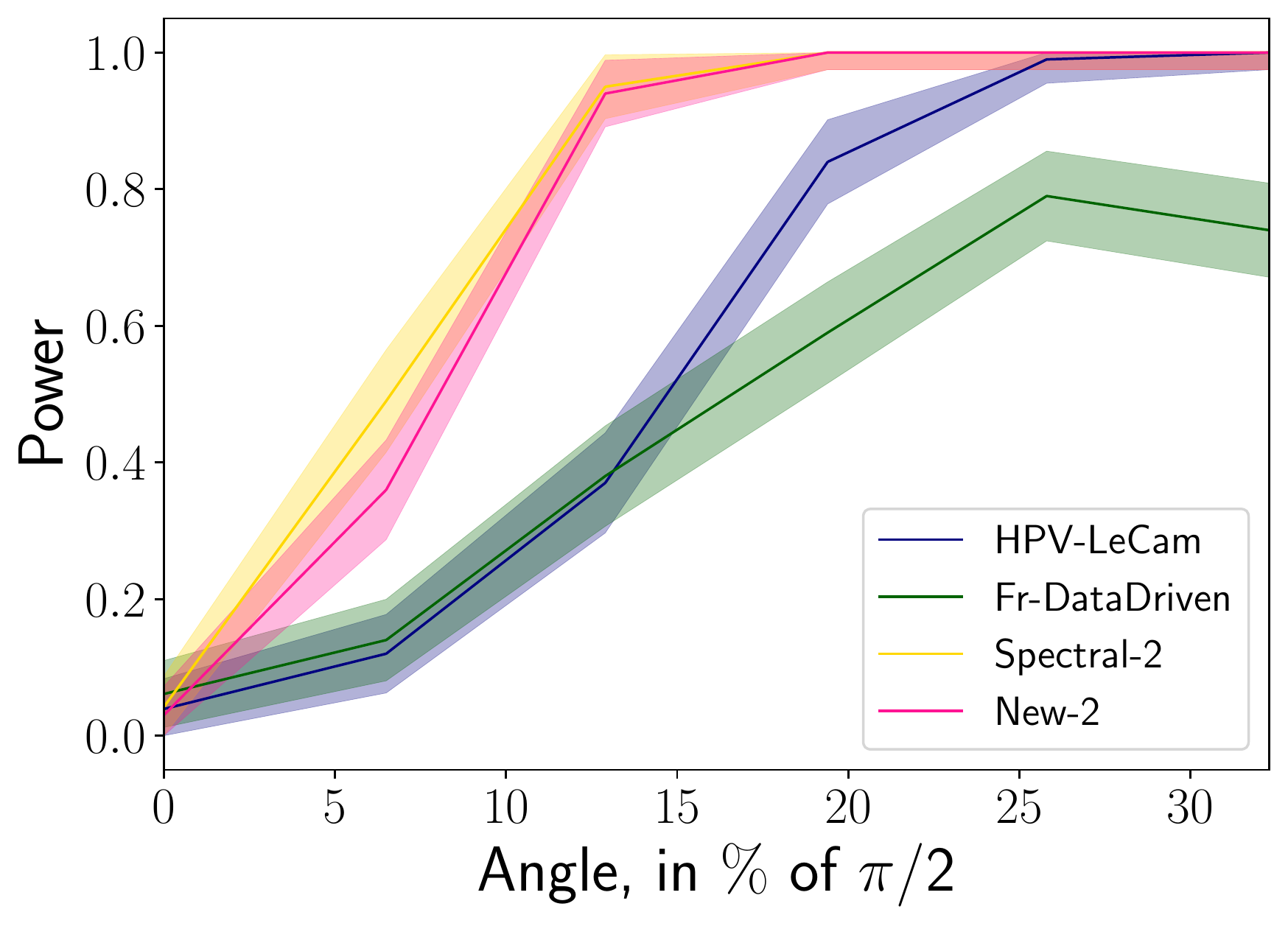}
         \caption{$n=1500, \;d=50$}
     \end{subfigure}
     \hfill
     \begin{subfigure}{0.32\textwidth}
         \centering
         \includegraphics[width=\textwidth]{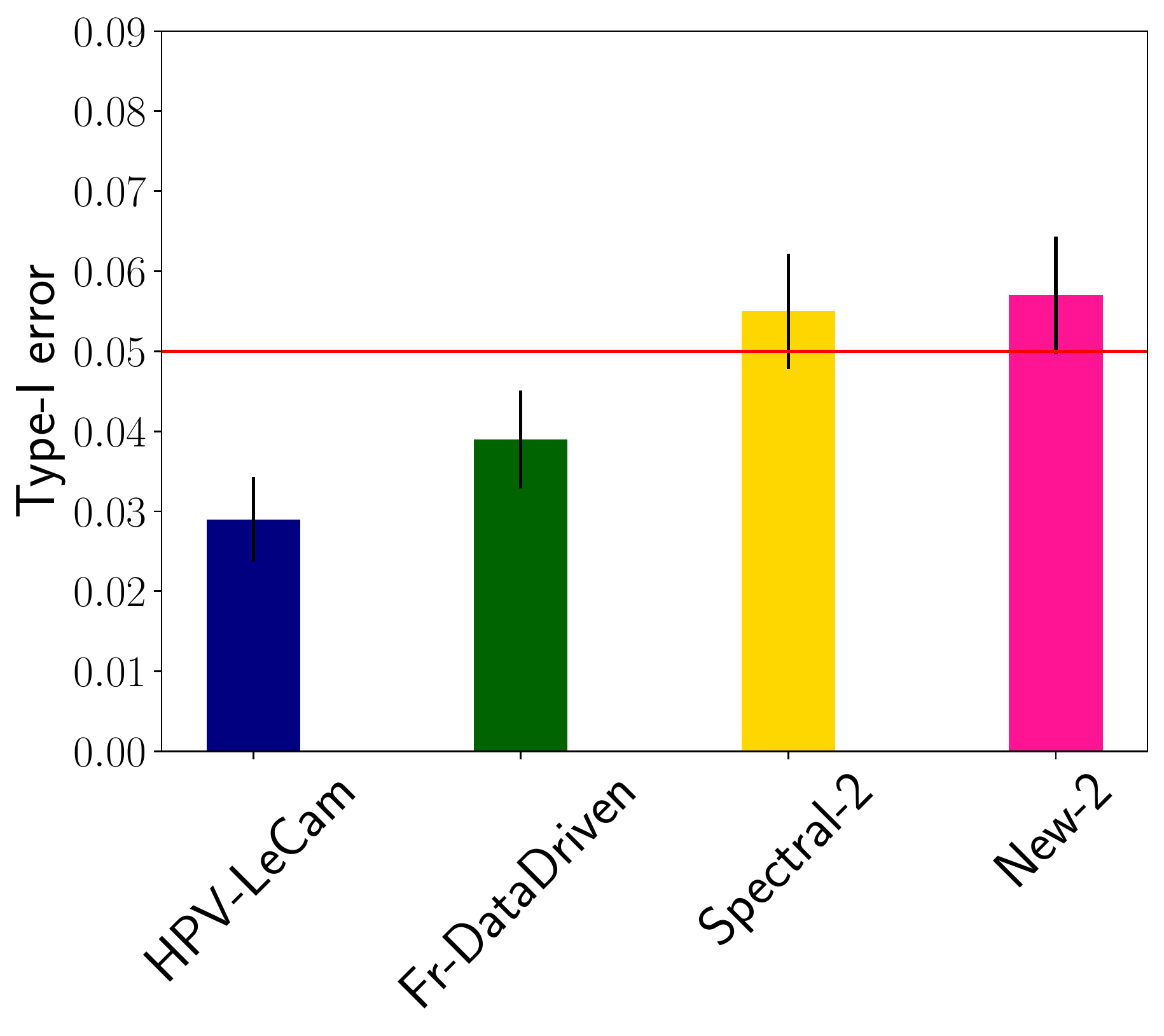}
         \vfill
         \includegraphics[width=\textwidth]{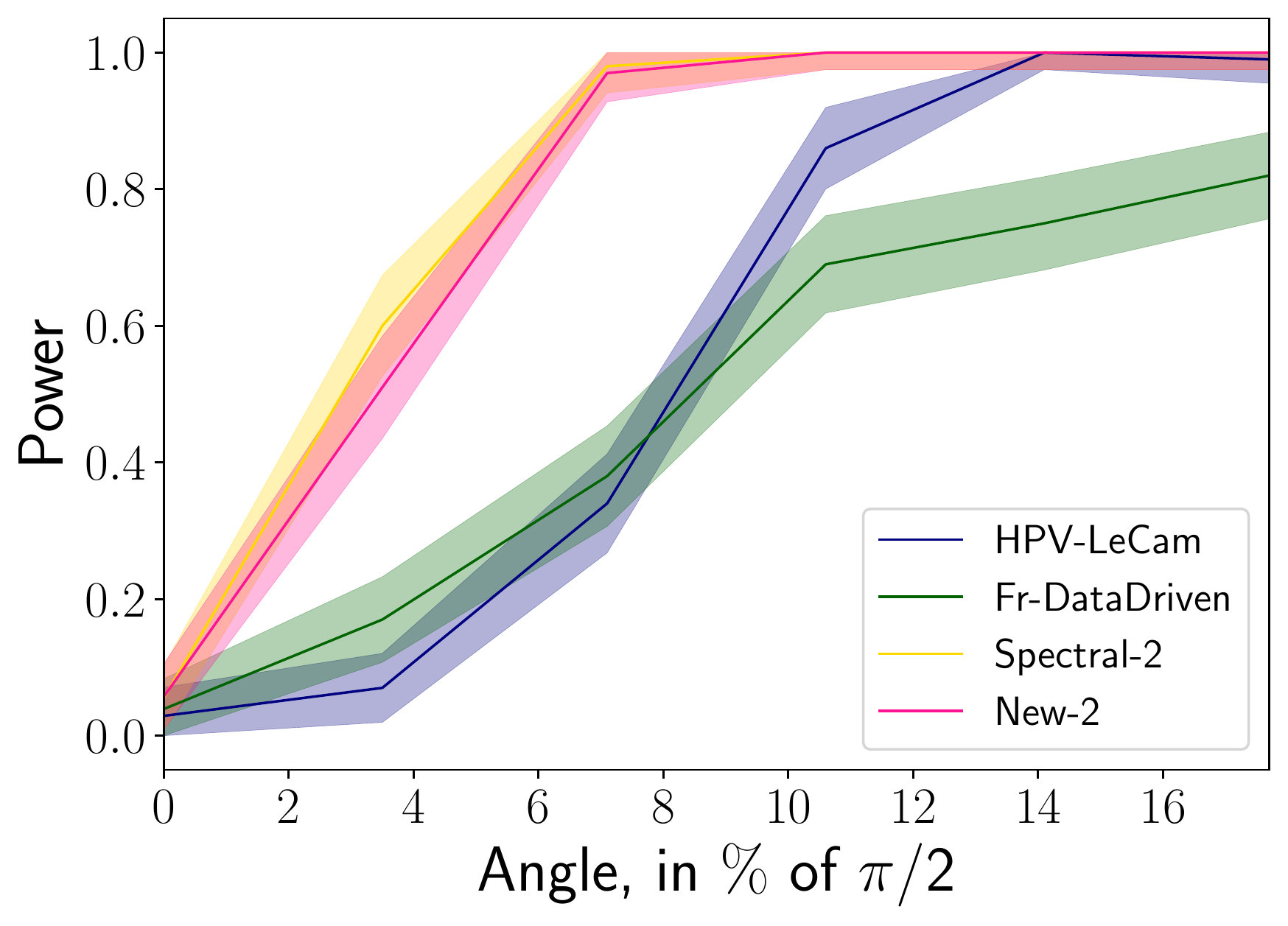}
         \caption{$n=5000, \;d=50$}
     \end{subfigure}
     \vfill
     \vspace{1.5cm}
     \begin{subfigure}{0.32\textwidth}
         \centering
         \includegraphics[width=\textwidth]{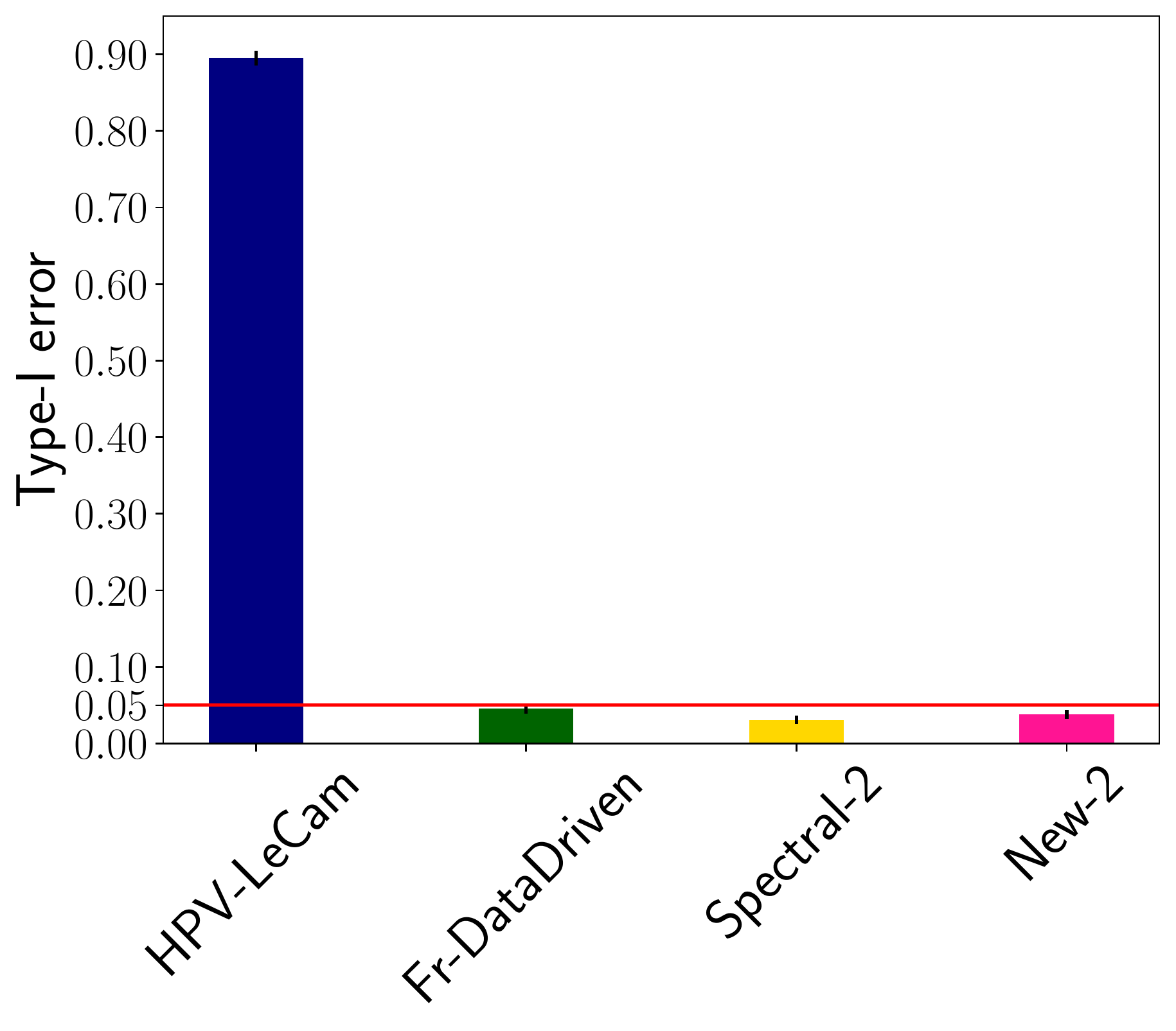}
         \vfill
         \includegraphics[width=\textwidth]{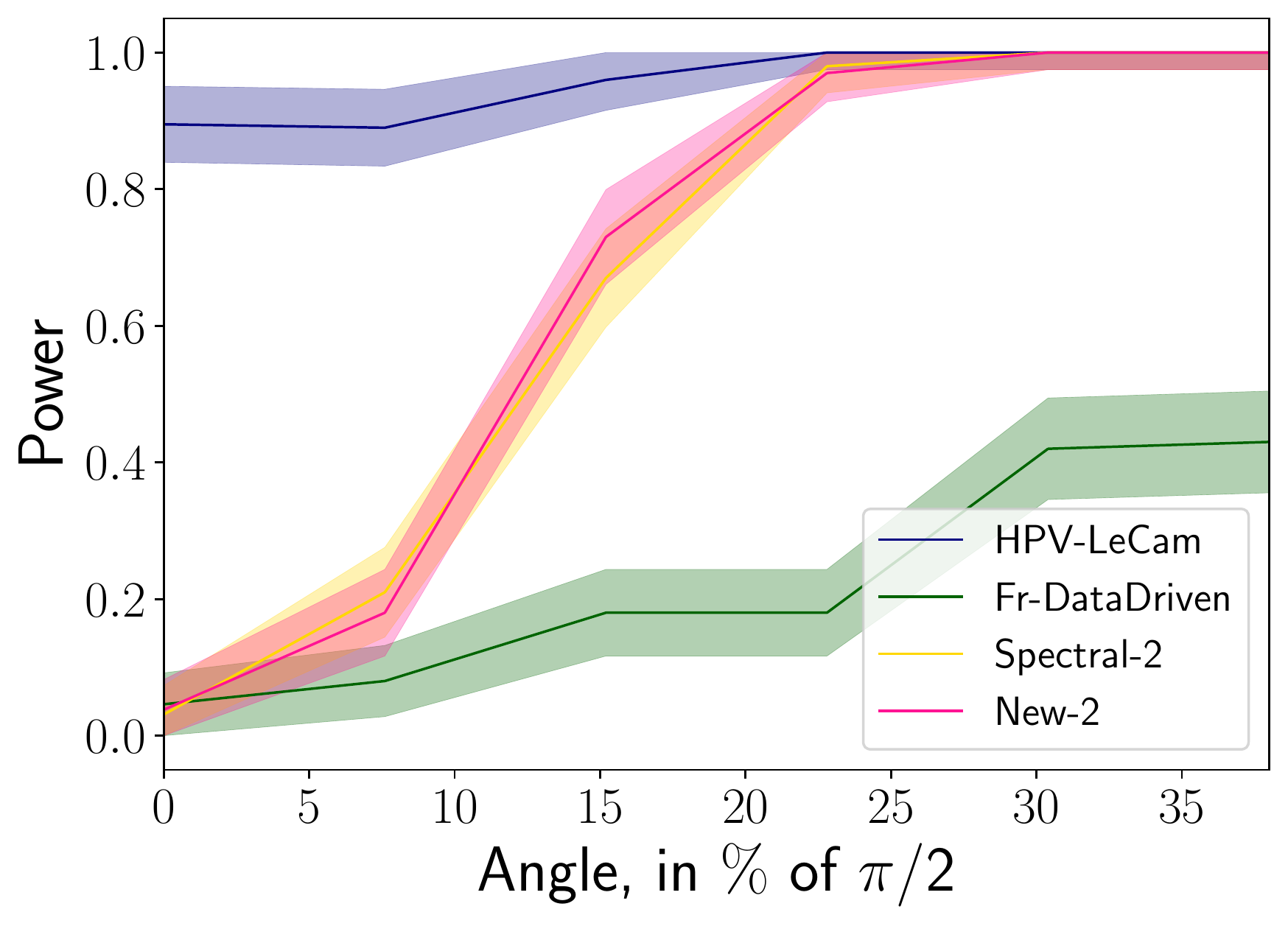}
         \caption{$n=500, \;d=150$}
     \end{subfigure}
     \hfill
     \begin{subfigure}{0.32\textwidth}
         \centering
         \includegraphics[width=\textwidth]{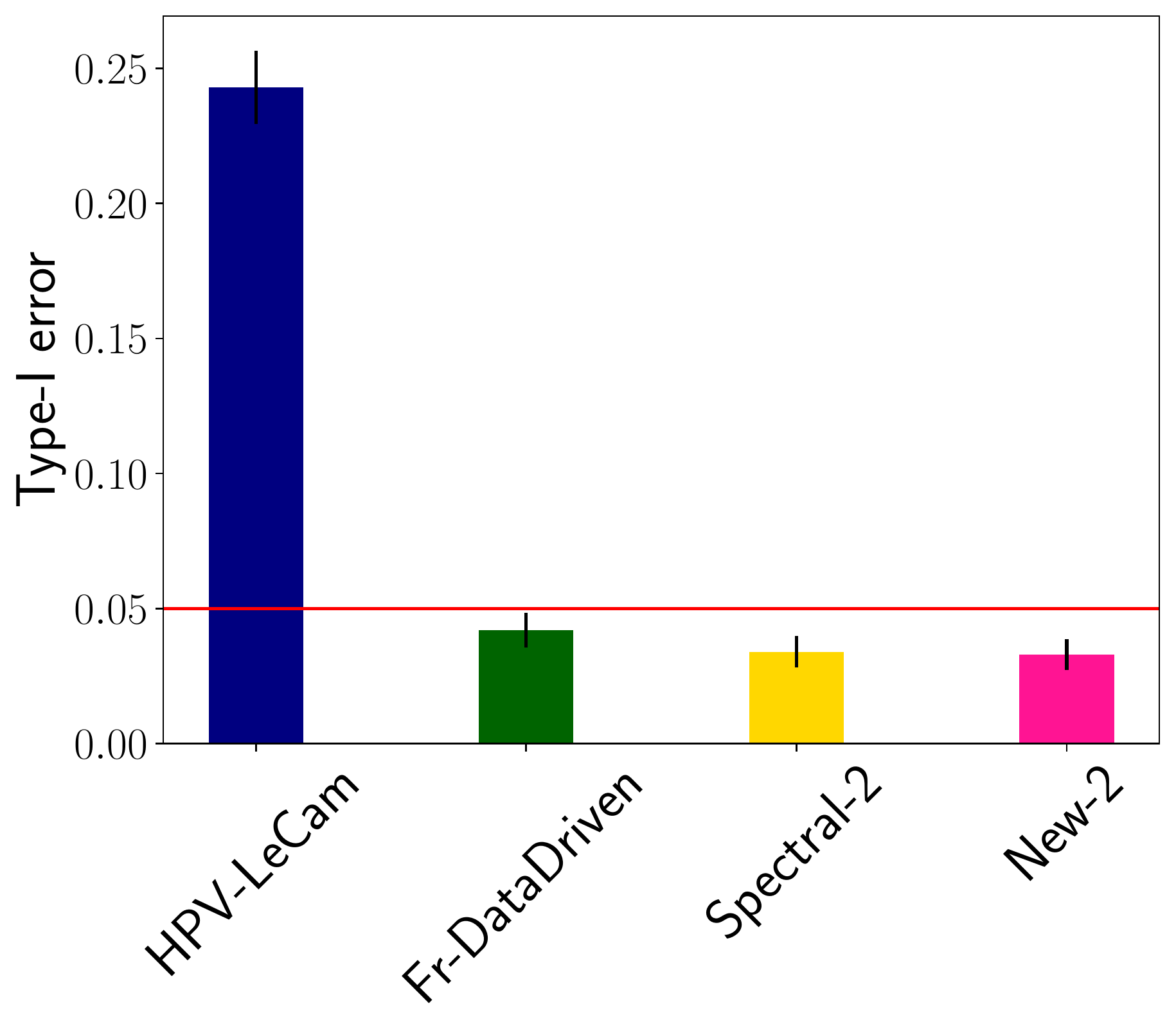}
         \vfill
         \includegraphics[width=\textwidth]{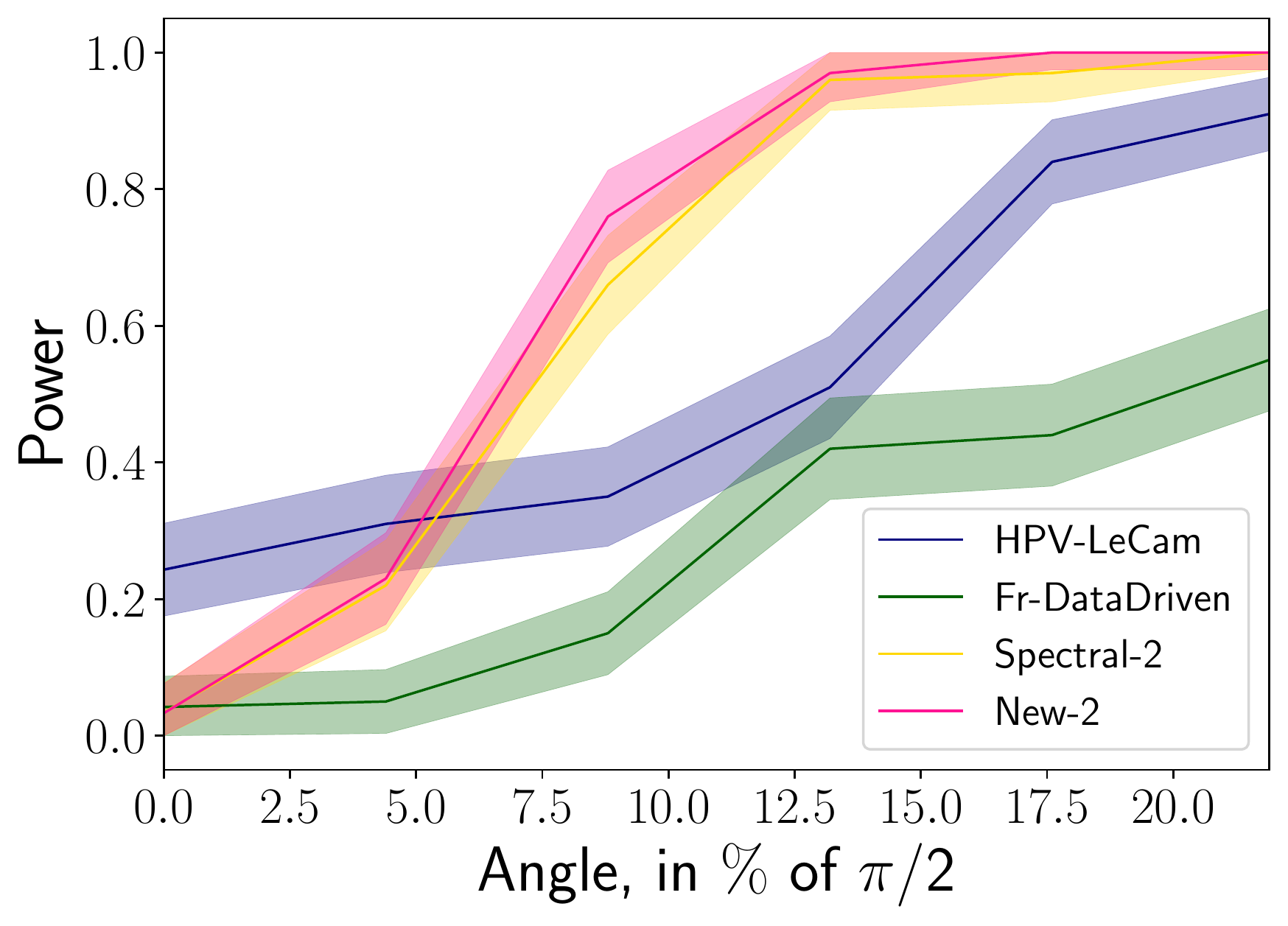}
         \caption{$n=1500, \;d=150$}
     \end{subfigure}
     \hfill
     \begin{subfigure}{0.32\textwidth}
         \centering
         \includegraphics[width=\textwidth]{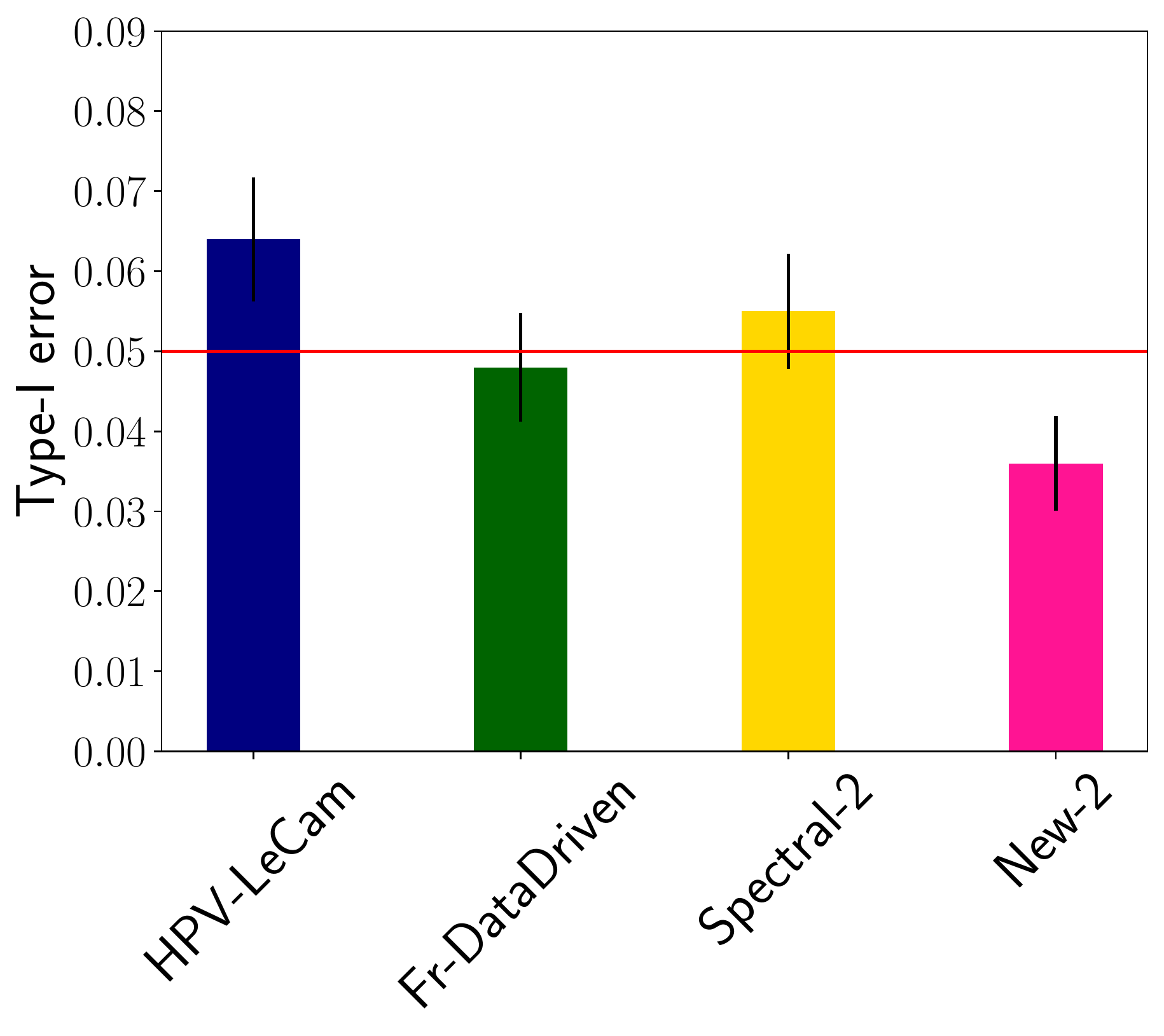}
         \vfill
         \includegraphics[width=\textwidth]{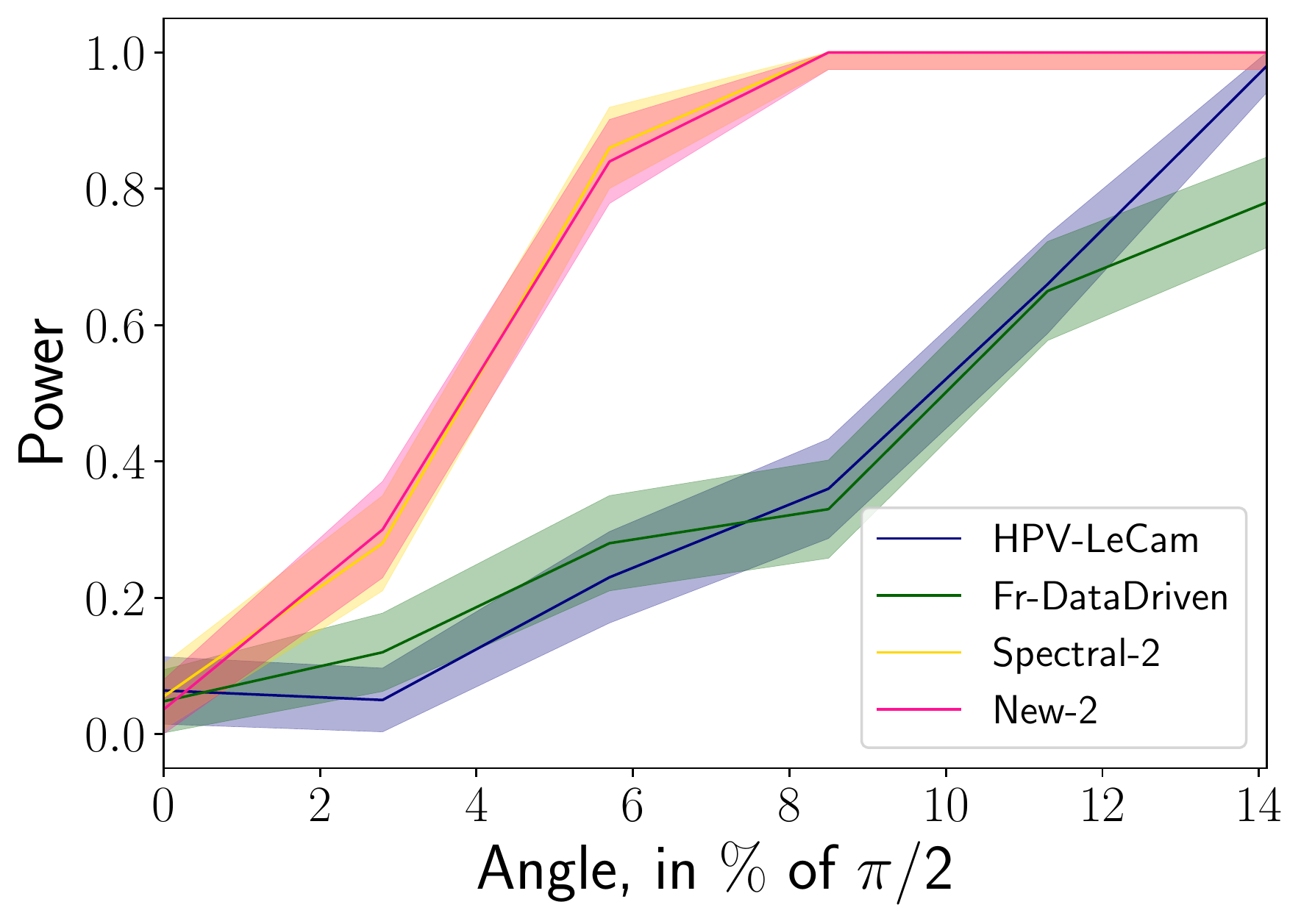}
         \caption{$n=5000, \;d=150$}
     \end{subfigure}
        \caption{Experiments for Scenario 2: One-sample problem, spiked regime with $m=1$, Laplace distribution.}
        \label{Exp2}
\end{figure}

\begin{figure}
     \begin{subfigure}{0.32\textwidth}
         \centering
         \includegraphics[width=\textwidth]{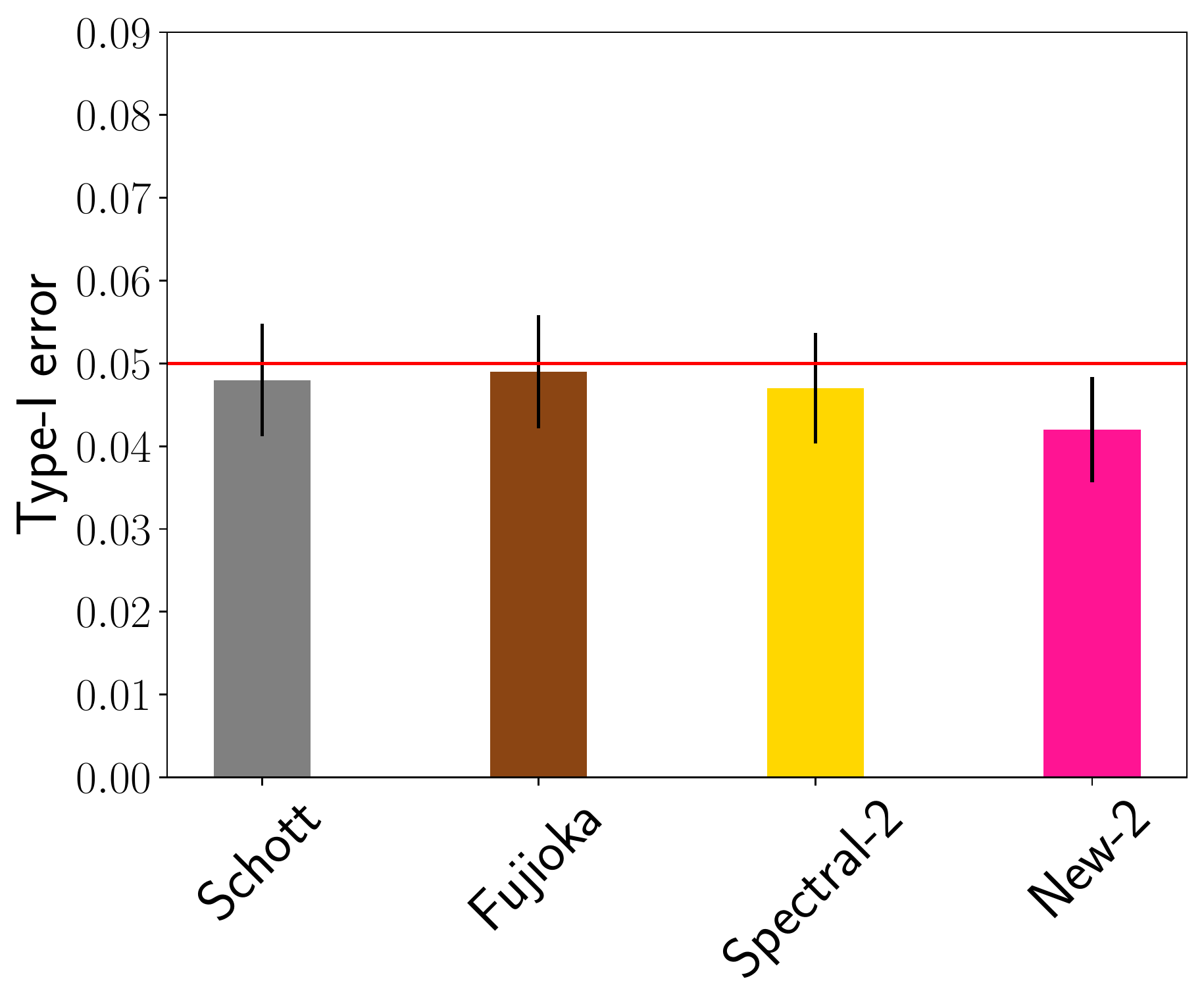}
         \vfill
         \includegraphics[width=\textwidth]{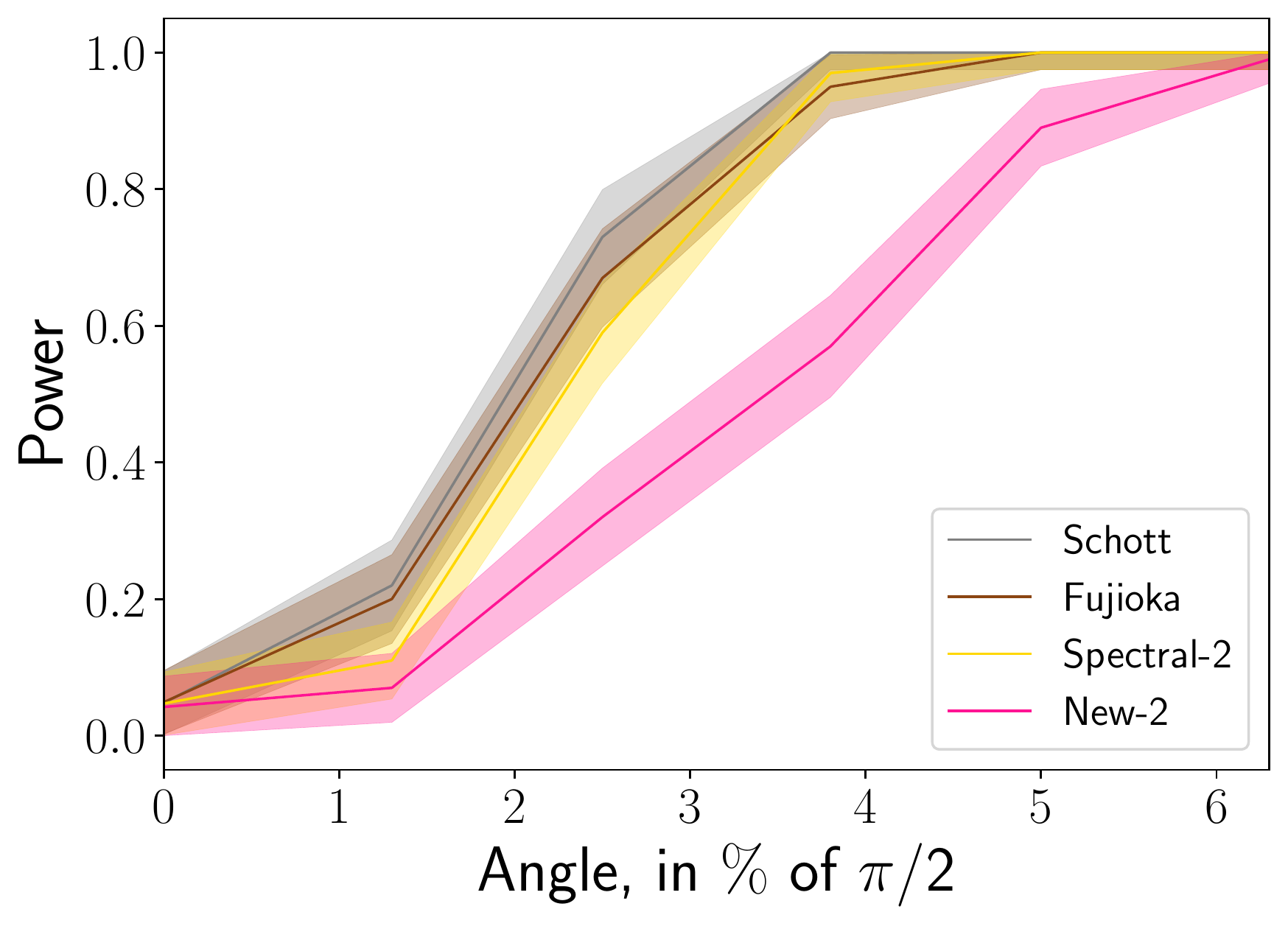}
         \caption{$n=500, \;d=50$}
     \end{subfigure}
     \hfill
     \begin{subfigure}{0.32\textwidth}
         \centering
         \includegraphics[width=\textwidth]{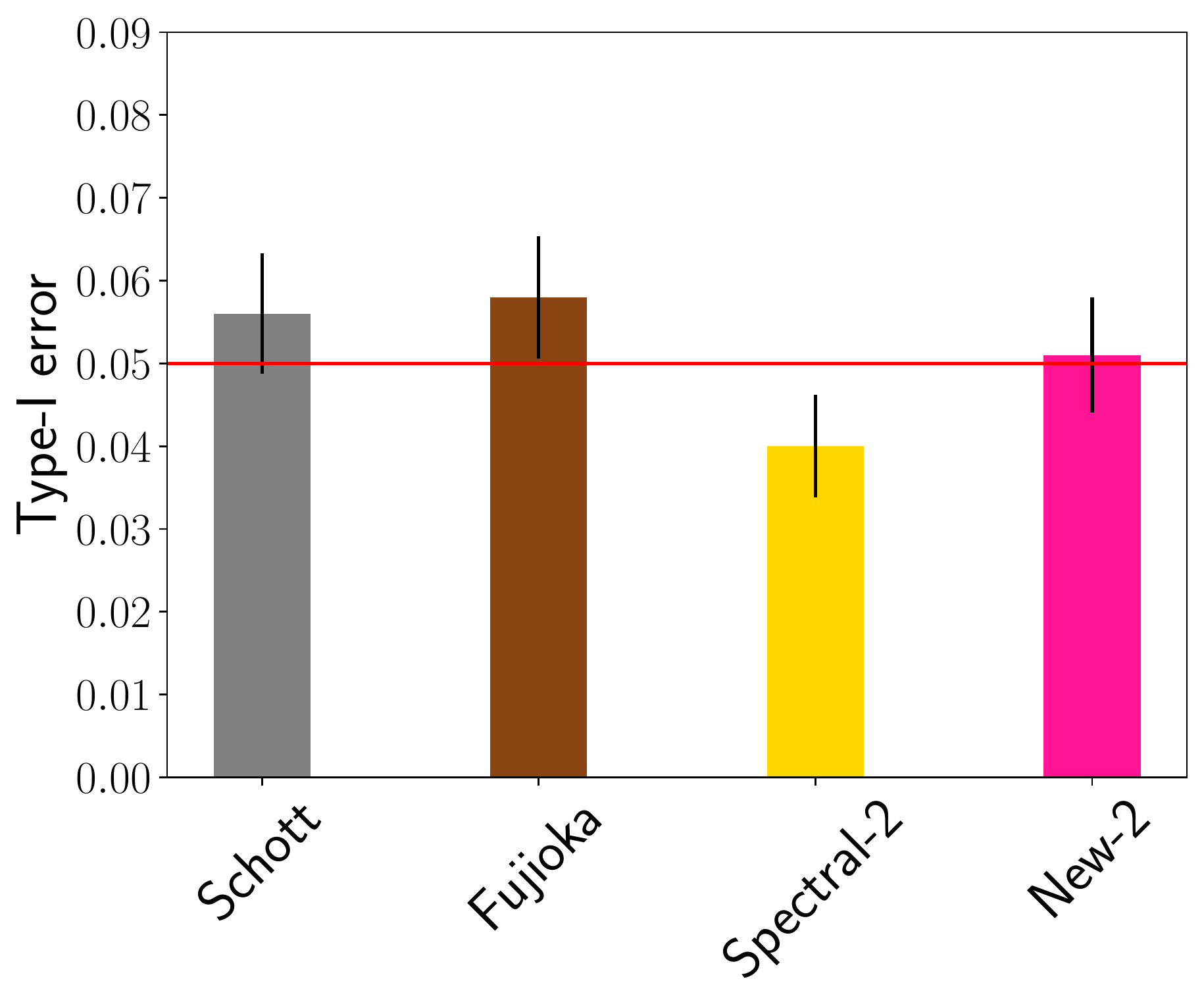}
         \vfill
         \includegraphics[width=\textwidth]{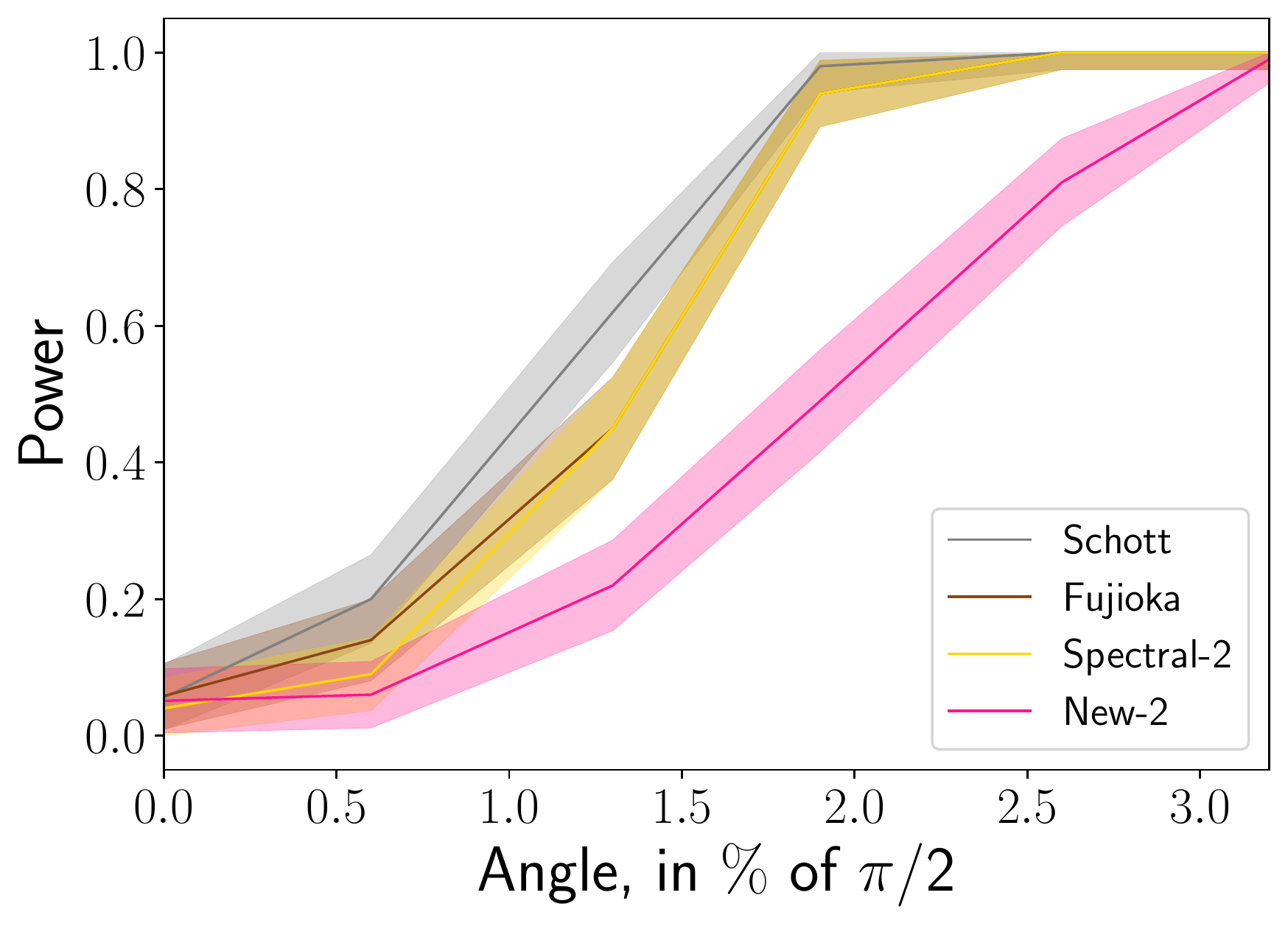}
         \caption{$n=1500, \;d=50$}
     \end{subfigure}
     \hfill
     \begin{subfigure}{0.32\textwidth}
         \centering
         \includegraphics[width=\textwidth]{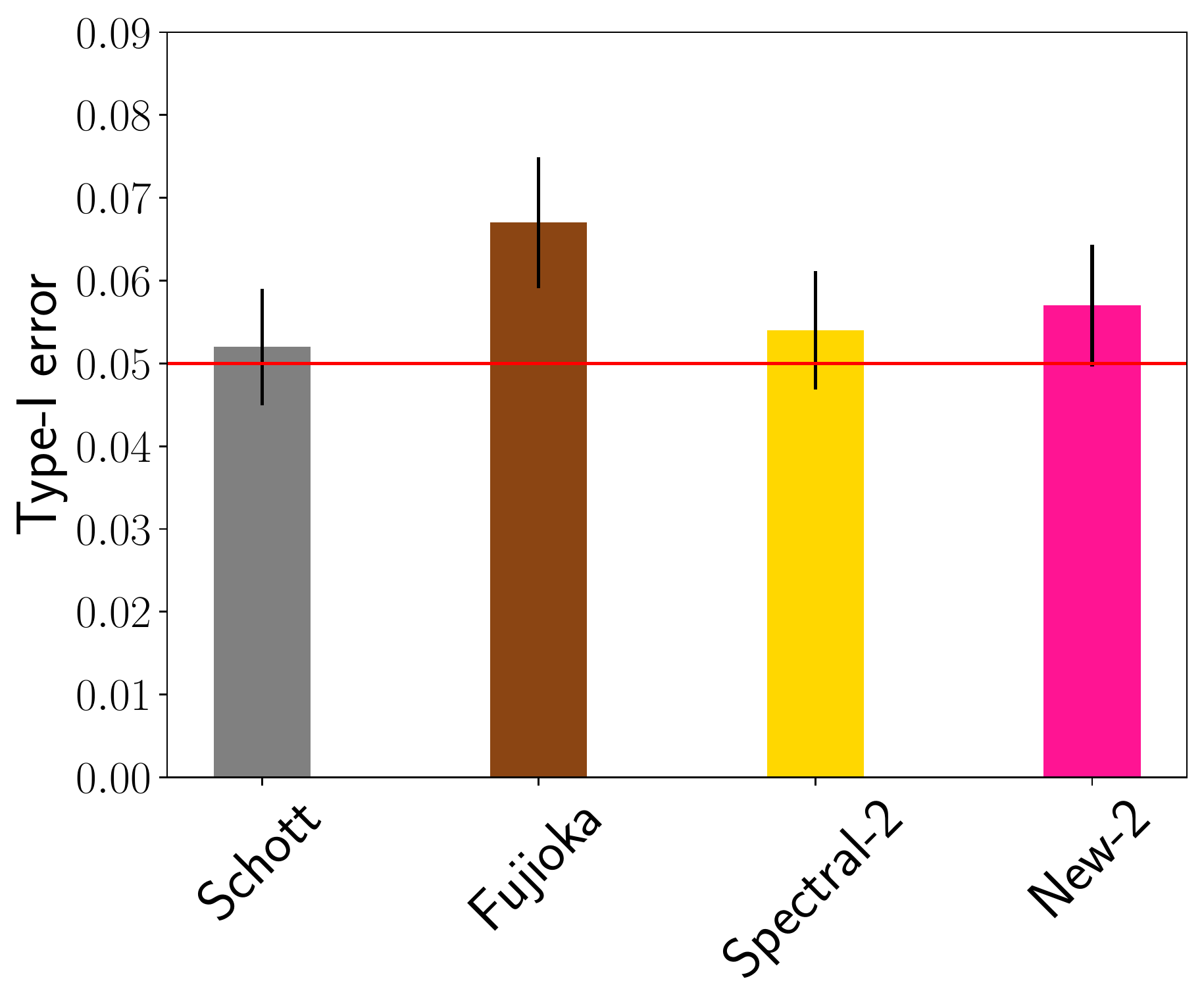}
         \vfill
         \includegraphics[width=\textwidth]{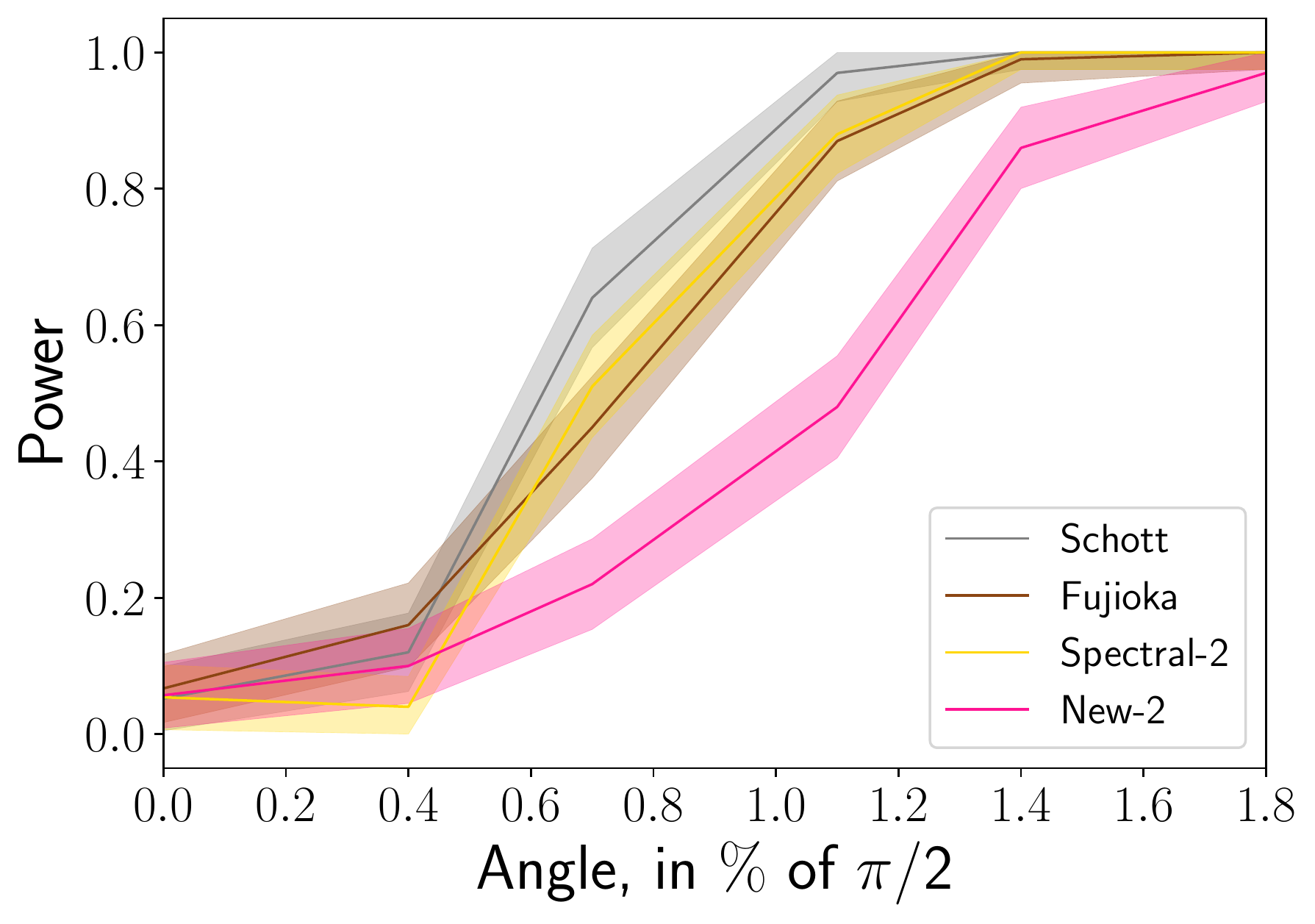}
         \caption{$n=5000, \;d=50$}
     \end{subfigure}
     \vfill
     \vspace{1.5cm}
     \begin{subfigure}{0.32\textwidth}
         \centering
         \includegraphics[width=\textwidth]{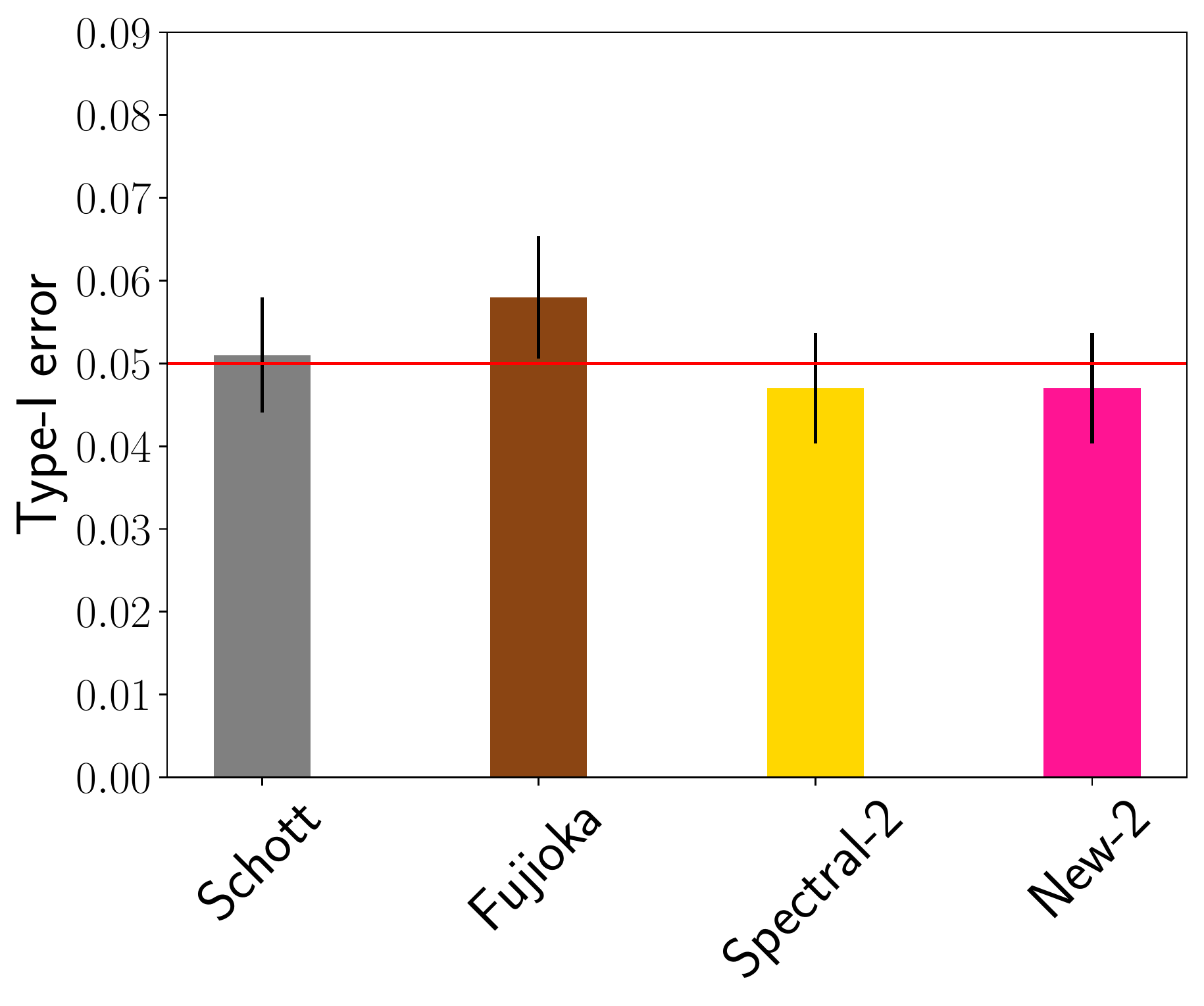}
         \vfill
         \includegraphics[width=\textwidth]{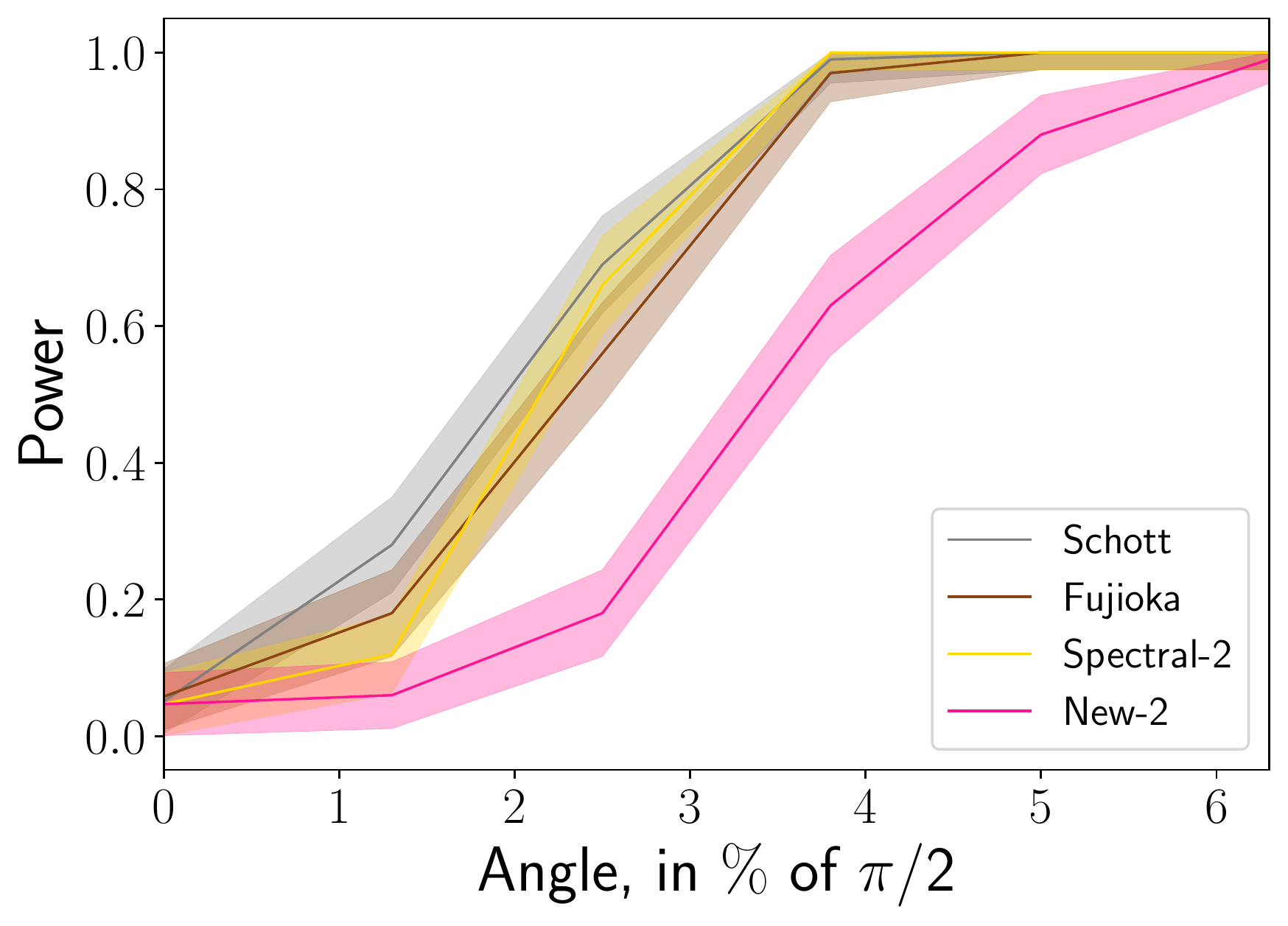}
         \caption{$n=500, \;d=150$}
     \end{subfigure}
     \hfill
     \begin{subfigure}{0.32\textwidth}
         \centering
         \includegraphics[width=\textwidth]{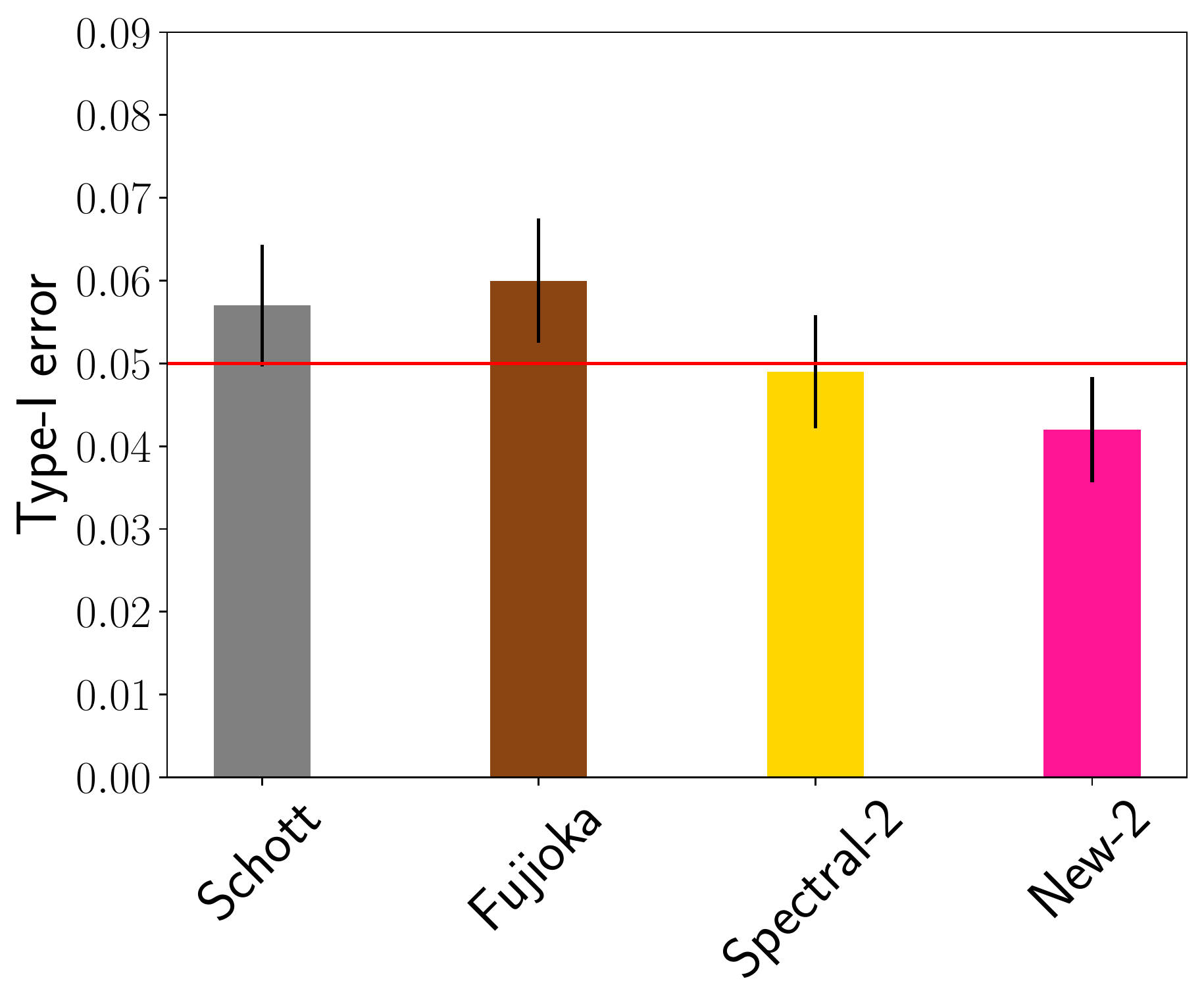}
         \vfill
         \includegraphics[width=\textwidth]{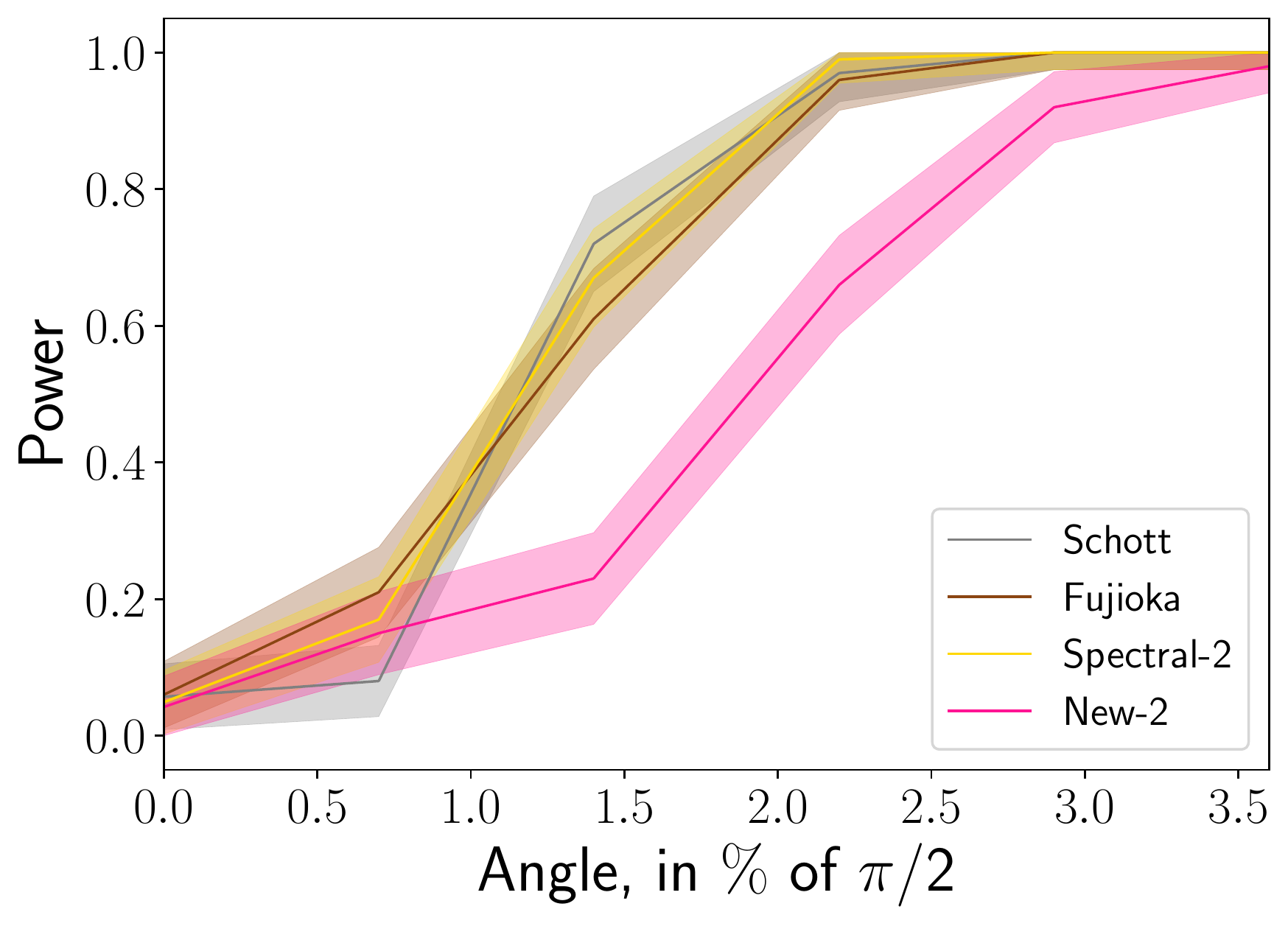}
         \caption{$n=1500, \;d=150$}
     \end{subfigure}
     \hfill
     \begin{subfigure}{0.32\textwidth}
         \centering
         \includegraphics[width=\textwidth]{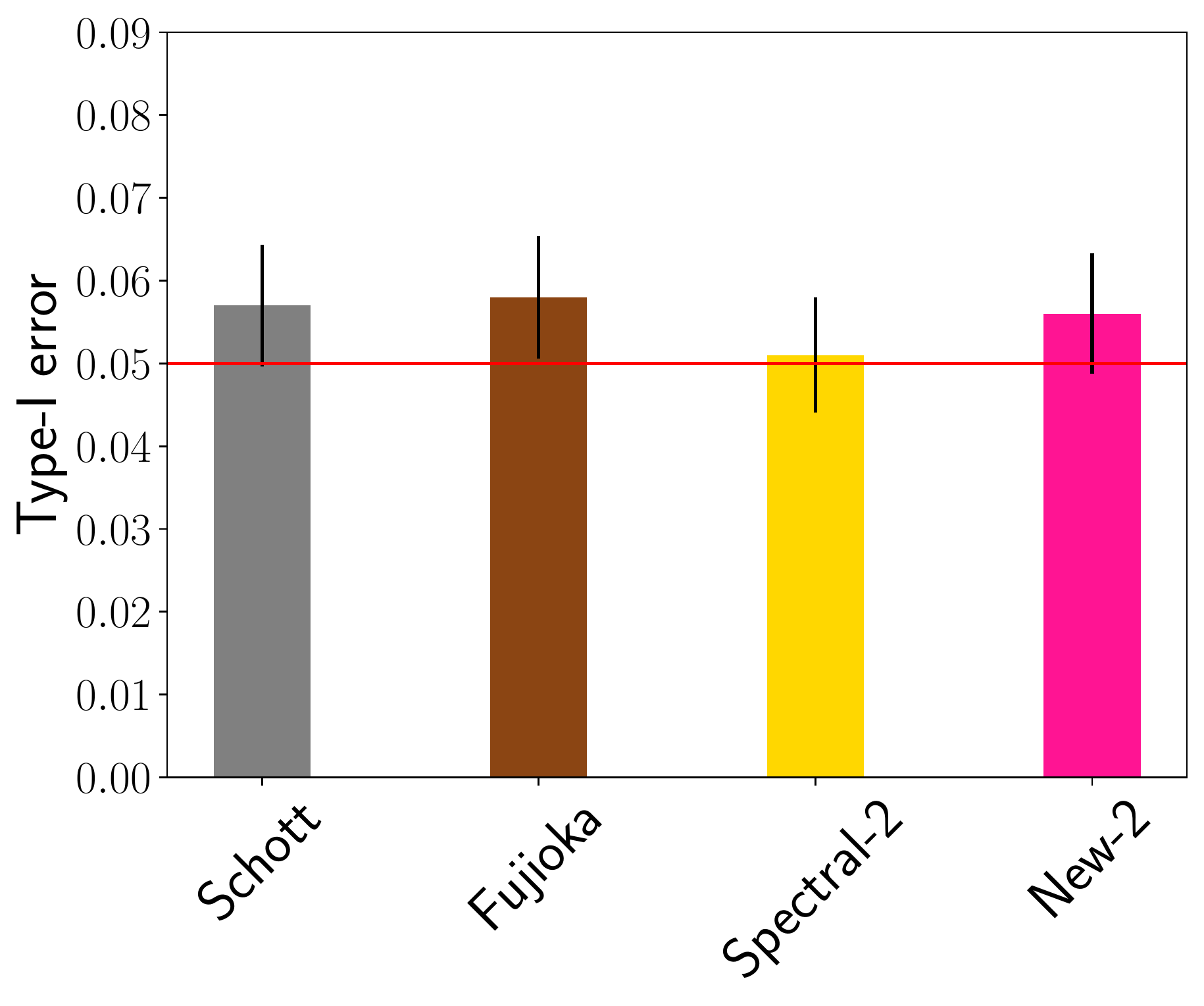}
         \vfill
         \includegraphics[width=\textwidth]{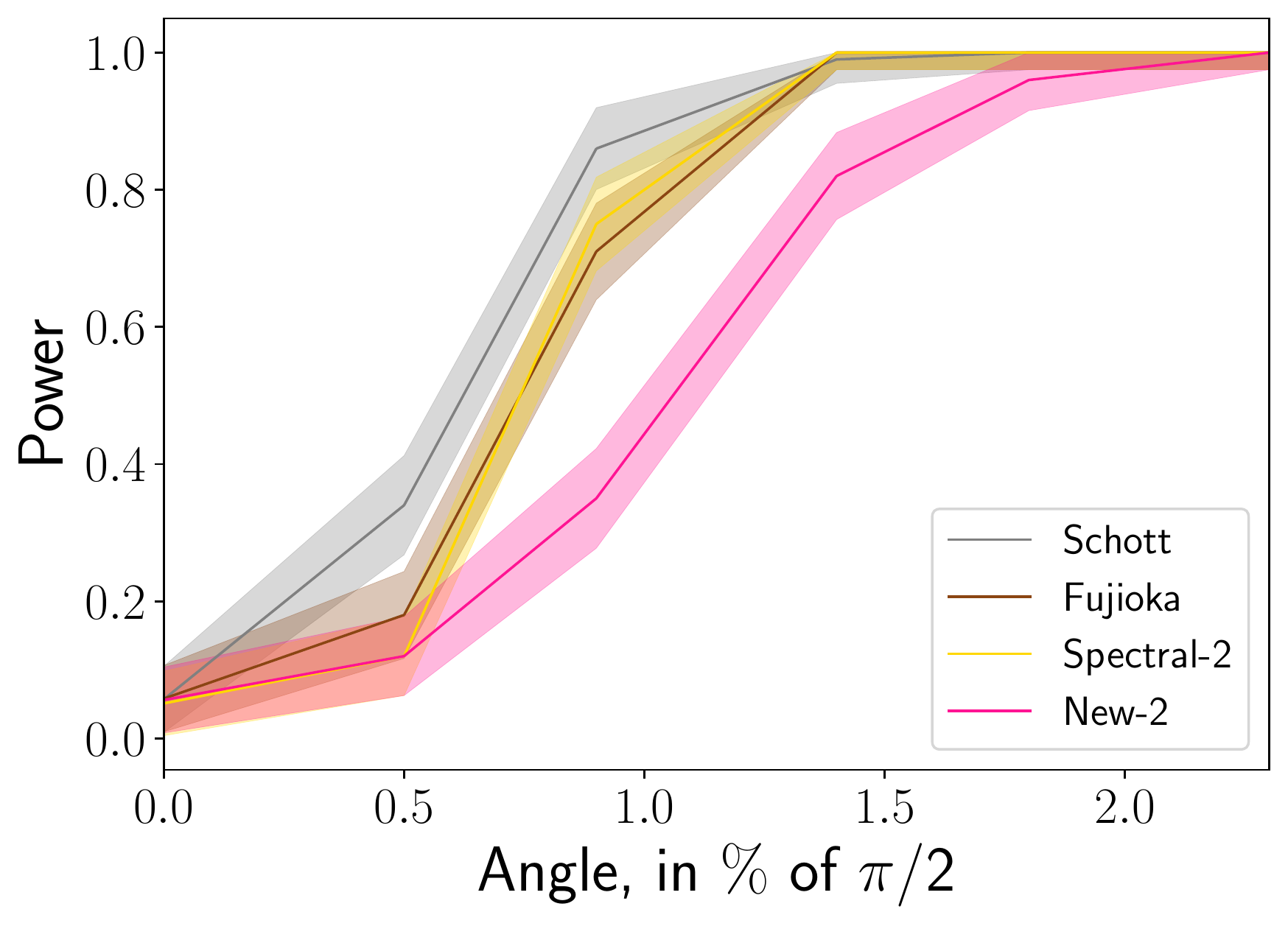}
         \caption{$n=5000, \;d=150$}
     \end{subfigure}
        \caption{Experiments for Scenario 3: Two-sample problem, decay regime with $m=5$, Laplace distribution.}
        \label{Exp3}
\end{figure}

The results are presented in Figure~\ref{Exp1}. We observe that the bootstrap-based methods are extremely conservative, almost never rejecting null hypothesis. This also implies very weak power of such methods. Our procedures are also slightly conservative, but the power of ``\textsf{New-1}'' and ``\textsf{New-2}'' significantly outperforms the methods based on Spectral and Frobenius norms. In further scenarios, we exclude bootstrap-based approaches.
\\
\\
\noindent\textbf{Scenario 2}\\
The next scenario considers one-sample problem in spiked regime with Gaussian distribution.
We compare the following methods:
\begin{itemize}
    \item ``\textsf{HPV-LeCam}'': Le Cam optimal test (\cite{LAN1})
    \item ``\textsf{Fr-DataDriven}'': Frobenius norm test statistic + Sample splitting strategy (\cite{Koltchinskii_NARPCA})
    \item ``\textsf{Spectral-2}'': same as in the previous scenario.
    \item ``\textsf{New-2}'': same as in the previous scenario.
\end{itemize}
The results are presented in Figure~\ref{Exp2}. Main conclusion here is that ``\textsf{HPV-LeCam}'' dramatically fails when the sample size is not significantly larger than the dimension, and even in the opposite case its power is quite weak. Quite unexpectedly, ``\textsf{Fr-DataDriven}'' behaves well under null (even though this method requires the dimension going to infinity). However, its power is inferior to ``\textsf{Spectral-2}'' and ``\textsf{New-2}'', which perform very similar in this setting, though again slightly conservative under null.
\\
\\
\noindent\textbf{Scenario 3}\\
In the last scenario we focus on two-sample problem in decay regime with Laplace distribution.
We compare the following methods:
\begin{itemize}
    \item ``\textsf{Schott}'': the procedure proposed in \cite{Schott_1}
    \item ``\textsf{Fujioka}'': the procedure proposed in \cite{Fujioka}.
    \item ``\textsf{Spectral-2}'': Spectral norm test statistic $\Qtspab$ + Approach 2.
    \item ``\textsf{New-2}'': $\| \cdot \|_{(\Pp, \Gamma, s_1, s_2)}$-norm test statistic $\Qtab$ (with $s_1 = s_2 = 1$) + Approach 2.
\end{itemize}
The results are presented in Figure~\ref{Exp3}. The quality of the four compared methods in this scenario is approximately similar, with ``\textsf{New-2}'' being slightly weaker then the others. The explanation is that ``\textsf{New-2}'' requires to split the sample into two halves, and effectively only half of the data is used to measure the discrepancy. Another important observation is that in this scenario with decay regime, the behaviour under null is very stable, and doesn't really seem to depend on the dimension. This promises that the type-I error bounds can be made dependant on some notion of effective rank, rather then the full dimension.

%% file: source/Discussion.tex
\subsection{Building confidence sets for $\Pt$}
Even though our work focuses on hypothesis testing, the idea can be used for constructing the confidence sets for the true spectral projector $\Pt$ from the data $[X_1, \ldots, X_{2n}] = \data$. Split the sample into two equal parts, compute the sample covariance matrices $\Se$ and $\So$ based on the first and second halves of the sample, respectively. Let $\Pe$ and $\overline{\Pp}_\J$ be the corresponding projectors. Fix $\overline{\Gamma}$ satisfying \eqref{Gamma_Prop} for $\overline{\Pp}_\J$. For a given confidence level $(1-\alpha)$, consider sets
            \begin{equation}
                \begin{aligned}
                &\mathcal{CS}^{1-\alpha}_B(\data) \eqdef \left\{ \Pp \in \R^{d\times d} \text{projector of rank }m: \sqrt{n} \| \Pp - \Pe \|_{(\overline{\Pp}_\J, \overline{\Gamma}, s_1,s_2)} \leq q^B_{\alpha} \right\},\\
                &\mathcal{CS}^{1-\alpha}_F(\data) \eqdef \left\{ \Pp \in \R^{d\times d} \text{projector of rank }m: \sqrt{n} \| \Pp - \Pe \|_{(\overline{\Pp}_\J, \overline{\Gamma}, s_1,s_2)} \leq q^F_{\alpha} \right\},
                \nonumber
                \end{aligned}
            \end{equation}
            where $q^B_\alpha$ and $q^F_\alpha$ are the $\alpha$-quantiles of $(\sqrt{n} \| \Pb - \Pe \|_{(\overline{\Pp}_\J, \overline{\Gamma}, s_1,s_2)} \;|\; \data)$ and \\$(\sqrt{n} \| \Pf - \Pe \|_{(\overline{\Pp}_\J, \overline{\Gamma}, s_1,s_2)} \;|\; \data)$, respectively ($\Pb$ and $\Pf$ depend on the first half of the sample only). One can show that the theoretical properties of the coverage probabilities of these sets are similar to the theoretical properties of type-I error in one-sample testing problem.

\subsection{Why do we use this test statistic?}
The reasons behind very non-trivial construction of our test statistic are partially similar to \cite{WXZhou}, which addresses similar testing problem, but for covariance matrices, rather than spectral projectors. Specifically, \cite{WXZhou} try to approximate the distribution of $\|\Se-\St\|$ and apply the bootstrap inference for it. So, their original idea is to work with spectral norm.
However, to approximate the distribution of $\| \Se-\St\|$ they require $d^9 \ll n$, and to approximate the distribution of $\| \Pe-\Pt\|$ we need $d^7 \ll n$. To avoid this restrictive assumption, \cite{WXZhou} introduce a parameter $s$, which helps to connect $\| \Se-\St\|$ with $\| \Se-\St\|_{\mymax}$ in a ``smooth'' way: in particular, they consider the $s$-sparse largest eigenvalue
\begin{equation}
    \begin{aligned}
        \sup\limits_{v\in\mathbbm{V}(s,d)} v^\T (\Se-\St) v,
    \nonumber
    \end{aligned}
\end{equation}
where $\mathbbm{V}(s,d) \eqdef \{ v\in\Sph^{d-1}:\;|\supp(v)|\leq s\}$. The quality of approximation of distribution of this quantity is expressed in terms of $s^9/n$ (omitting logarithmic factors), as we take supremum over smaller set. We could try to follow the same logic and work with the test statistic (here we focus one one-sample framework for simplicity)
\begin{equation}
    \begin{aligned}
        \sup\limits_{v\in\mathbbm{V}(s,d)} v^\T (\Pe-\Pt) v,
    \nonumber
    \end{aligned}
\end{equation}
however, in this case it follows from our proof that under $H_0^{(1)}$
\begin{equation}
    \begin{aligned}
        \Var[v^\T (\Pe-\Pt) v] \asymp C_\St\cdot v^\T\Pt v\cdot v^\T \Ptc v,
    \nonumber
    \end{aligned}
\end{equation}
which prevents us from applying the main tool in our analysis, Gaussian approximation (see \cite{GAR+Max}, Theorem 2.2), as we cannot guarantee the lower bound 
\begin{equation}
    \begin{aligned}
    \Var[v^\T (\Pe-\Pt) v] \geq c_1 > 0
    \nonumber
    \end{aligned}
\end{equation}
    uniformly over $v$. In the case of covariance matrices from \cite{WXZhou} this problem is solved either by assuming $\lambda_{min}(\St) \geq c_1 > 0$, or by considering the normalized version
\begin{equation}
    \begin{aligned}
        \sup\limits_{v\in\mathbbm{V}(s,d)} \frac{v^\T (\Se-\St) v}{v^\T\St v}.
    \nonumber
    \end{aligned}
\end{equation}
In our situation this normalization would lead to
\begin{equation}
    \begin{aligned}
        \sup\limits_{v\in\mathbbm{V}(s,d)} \frac{v^\T (\Pe-\Pt) v}{\sqrt{v^\T\Pt v\cdot v^\T (\Id_d-\Pt) v}},
    \nonumber
    \end{aligned}
\end{equation}
which is a reasonable object in theory. However, from practical prospective such a normalization leads to computational issues (in addition to intractability of combinatorial optimization over $\mathbbm{V}(s,d)$).
It turns out that the rotation $\Gamma$ that we introduce in Definition~\ref{D:norm} (here $\Gamma$ corresponds to $\Pt$) plays role of the normalization and can be used instead: specifically, we could consider
\begin{equation}
    \begin{aligned}
        \sup\limits_{\substack{v\in\mathbbm{V}(s_1,m)\\w\in\mathbbm{V}(s_2,d-m)}} [v^\T w^\T] \, \Gamma^\T (\Pe-\Pt) \Gamma \begin{bmatrix} v\\w \end{bmatrix}.
    \nonumber
    \end{aligned}
\end{equation}
One may check that in this case under $H_0^{(1)}$
\begin{equation}
    \begin{aligned}
        \Var\left[[v^\T w^\T] \Gamma^\T (\Pe-\Pt) \Gamma \begin{bmatrix} v\\w \end{bmatrix} \right] \asymp C_\St,
    \nonumber
    \end{aligned}
\end{equation}
which makes Gaussian approximation applicable. Moreover, to avoid computation intractability caused by optimization over $\mathbbm{V}(s,d)$, or more specifically in our case $\mathbbm{V}(s_1,m)$ and $\mathbbm{V}(s_2,d-m)$, we replace them by sets $\myset^m_{s_1}$ and $\myset^{d-m}_{s_2}$, where $\myset^d_s$ consist of unit vectors in $\R^d$, whose support consists of $s$ consecutive coordinates. So, it is another way to provide smooth connection between extreme cases using the parameters $s_1$ and $s_2$, and it would lead to
\begin{equation}
    \begin{aligned}
        \sup\limits_{\substack{v\in\myset^m_{s_1}\\w\in\myset^{d-m}_{s_2}}} [v^\T w^\T] \, \Gamma^\T (\Pe-\Pt) \Gamma \begin{bmatrix} v\\w \end{bmatrix},
    \nonumber
    \end{aligned}
\end{equation}
which can be written also as
\begin{equation}
    \begin{aligned}
        \sup\limits_{\substack{v\in\myset^m_{s_1}\\w\in\myset^{d-m}_{s_2}}} \left( v^\T \Gamma_1^\T (\Pe-\Pt) \Gamma_1  v + w^\T \Gamma_2^\T (\Pe-\Pt) \Gamma_2  w + 2v^\T \Gamma_1^\T (\Pe-\Pt) \Gamma_2  w \right).
    \nonumber
    \end{aligned}
\end{equation}
We go even further, and the last step of explaining the reasons behind specific construction of our test statistic is the observation that under $H_0^{(1)}$, due to specific structure of spectral projectors, the first two terms in the above display become negligible, and in fact we can replace them with $\| \Gamma_1^\T (\Pe-\Pt) \Gamma_1 \|$ and $\|\Gamma_2^\T (\Pe-\Pt) \Gamma_2\|$, which still will be negligible. The reason behind this change is that while the properties under null hypothesis are not spoiled, the discrimination power under alternative hypothesis of the sum of these spectral norms is better rather that of $\sup\limits_{\substack{v\in\myset^m_{s_1}\\w\in\myset^{d-m}_{s_2}}} \left( v^\T \Gamma_1^\T (\Pe-\Pt) \Gamma_1  v + w^\T \Gamma_2^\T (\Pe-\Pt) \Gamma_2  w \right)$. In other words, this allows to gain in power for free (``power enhancement''). This leads to the final version of our test statistic 
\begin{equation}
    \begin{aligned}
        \| \Gamma_1^\T (\Pe-\Pt) \Gamma_1 \| + \| \Gamma_2^\T (\Pe-\Pt) \Gamma_2 \| + 2 \sup\limits_{\substack{v\in\myset^m_{s_1}\\w\in\myset^{d-m}_{s_2}}} v^\T \Gamma_1^\T (\Pe-\Pt) \Gamma_2  w,
    \nonumber
    \end{aligned}
\end{equation}
which gives an equivalent (up to a factor $2$) definition of our norm.

\subsection{Comparison with covariance matrix testing in \cite{WXZhou}}

        As discussed above, \cite{WXZhou} focuses on a problem on bootstrap inference for $s$-sparse largest eigenvalue of $(\Se-\St)$, and, consequently, applies the results to hypotheses testing setting for covariance matrices. Now we want to compare different aspects of our work and \cite{WXZhou}.

        While we deal with a different problem of hypothesis testing for spectral projectors, the results rely on similar idea of Gaussian approximation for maxima of sums random vectors after $\varepsilon$-net argument for supremum. Also, both works try to replace spectral norm to get better rates: \cite{WXZhou} works with $s$-sparse largest eigenvalue, and we consider $\| \cdot \|_{(\Pp, \Gamma, s_1,s_2)}$-norm. Here we highlight what differs our work from \cite{WXZhou}, apart from the fact that the objects of interest for us are spectral projectors, rather than covariance matrices.
        \begin{itemize}
            \item As can be seen from the previous subsection, generalization of $s$-sparse largest eigenvalue norm is not straightforward.
            \item New norm brings computational tractability for the test statistic. \cite{WXZhou} claims that they compute $s$-sparse largest eigenvalue using truncated power method, but in fact this method doesn't apply to their framework, see \cite{TruncPM}.
            \item Proofs of \cite{WXZhou} and ours are based on different results. While proof of \cite{WXZhou} uses coupling inequality for maxima of sums of random vectors (see Corollary~4.1 of \cite{GAR+Sup}),
            we employ Gaussian approximation technique (see Theorem~2.2 of \cite{GAR+Max}).
            Though these results are closely tied and derived by the same authors, it turns out that the latter allows to obtain slightly better rate: for instance, if we consider test statistic based on spectral norm,
            the results of \cite{WXZhou} require (omitting logarithmic terms) $d^9/n \ll 1$, while we assume a bit weaker $d^7/n \ll 1$.

            \item Bootstrap inference has been very popular for this kind of statistical problems where the limiting distribution of the test statistic depends on unknown parameters of the model and/or, in addition, is very unfriendly to work with, as in our case. However, as we already mentioned, multiplier bootstrap suffers from one computational issue, since to generate one bootstrap sample, one needs to generate $n$ random weights $\eta_1, \ldots, \eta_n$. Hence, in our work, along with the standard bootstrap method (Approach 1) we suggested the method, linked to Frequentist Bayes.
            The computational complexity of Approach 2 (specifically, of its ``resampling'' stage) does not depend on $n$, hence, it is significantly more efficient than Approach 1.
            \item Continuing the previous point, we note that the implementation of the bootstrap procedure in \cite{WXZhou} does not allow to build confidence sets, since their bootstrap test statistic $\widetilde{B}_{max}$ (see equation (2.3) from \cite{WXZhou}) depends on $\St$, thus is known only under $H_0$. In general, testing hypotheses and constructing confidence sets are known to be dual problems; however, constructing confidence sets is more difficult in a sense that we never know $\St$ in this case, while for testing hypotheses we have a hypothetical covariance $\mathbf{\Sigma^\circ}$, which can be used in test statistic and coincides with the true one under $H_0$.

            In contrast, we provide the procedure for building confidence sets.

            \item Finally, the results of \cite{WXZhou} assume sub-Gaussian data distribution, even though they mention that it can be relaxed. In modern applications the data are often heavy-tailed, and the extension beyond sub-Gaussianity becomes crucial. We make use of results from recent paper of \cite{Weibull} that provides user-friendly framework for dealing with sub-Weibull distributions considered in our work.

        \end{itemize}

\subsection{Comparison with previous works on inference for projectors}
Previous works on the topic are
\cite{Koltchinskii_NAACOSPOSC, Koltchinskii_NARPCA, Naumov, Silin_1}. Here we discuss how our work differs from these papers.

\begin{itemize}

    \item All of the mentioned works do not state the problem as hypothesis testing, and hence, they do not analyze power of the test that can be proposed based on their results. Furthermore, two-sample case was not considered as well.
    \item \cite{Koltchinskii_NAACOSPOSC, Koltchinskii_NARPCA, Naumov} rely on Gaussianity of the data distributions, which is, undoubtedly, extremely restrictive. \cite{Silin_1} work under significantly weaker ``covariance concentration condition'' (cov. conc.) of the form $\| \Se-\St\| \leq \delta_{n,d} \|\St\|$ with high probabiliy, and the rate $\delta_{n,d}$ appears in their bounds. In our work, no parametric assumption is imposed as well. Moreover, as was already mentioned, not only our results apply to sub-Gaussian case, but extend to distributions with heavier tails.
    \item The quantity of interest in all of the mentioned works is squared Frobenius distance between projectors $\| \Pe - \Pt \|_{\Fr}^2$, and the limiting distributions in \cite{Koltchinskii_NAACOSPOSC, Naumov, Silin_1} depends on the unknown true covariance $\St$, hence, statistical inference requires some special treatment.  \cite{Koltchinskii_NAACOSPOSC} suggest splitting the sample into three parts to make statistical inference. In \cite{Koltchinskii_NARPCA} the developed limiting distribution doesn't depend on $\St$, which makes statistical inference straightforward.
    \cite{Naumov} uses multiplier bootstrap, which is computationally intensive, as we pointed out previously. \cite{Silin_1} suggests Bayesian inference, that actually serves as a basis for our Approach 2.

    Unlike these works, we consider completely different test statistic, and though the limiting distribution does depend on the underlying true covariance as well, we present both the multiplier approach and the approach emerging from Bayesian perspective to make a valid calibration for our test.
    \item All of the mentioned works use linear approximation for spectral projectors and bound the remainder term as in Lemma~2 of \cite{Koltchinskii_AACBFBFOSPOSC}. While for Gaussian data distribution this result is sufficient to state dimension-free bounds (thanks to the sample covariance concentration in terms of effective rank, see \cite{Koltchinskii_CIAMBFSCO}, Corollary 2), without Gaussianity the appearance of the term $\sqrt{d^2/n}$ in the error bounds seems to be inevitable. In contrast, to bound the remainder term in linear approximation for projectors, we use new tight results from \cite{Wahl}, which allow to state the bounds in terms of relative rank. As a result, the dependence on the dimension is much better for example in Factor Model setting.
\end{itemize}

We summarize the comparison of the methods in the following table. The column ``Complexity'' specifies how many times we need to compute the corresponding norm in the procedures. In the last two columns we compare the required relations between the dimension $d$ and the sample size $n$ in two important regimes: Factor Model (FM) regime and Spiked ($\Tr[\St] \asymp d$) regime.
\begin{center}
\begin{tabular}{|l||c|c|c|c|c|}
\hline
 Method & Idea &  Data & Complexity & FM rate & Spiked rate \\ \hline \hline
 Koltchinskii, & Data-driven & Gaussian & $O(d^2)$  & not appl. & $1 \ll d \ll n$   \\
 Lounici (2017b,c) & &  &  &  &  \\ \hline
 Naumov et al. (2019) & Bootstrap & Gaussian & $O(Nnd^2)$  & $d^2 \ll n$ & $d^6 \ll n$   \\ \hline
 Silin, Spokoiny (2018) & Bayes & Cov. conc. & $O(Nd^2)$  & $d^{3.5} \ll n$ & $d^{3.5} \ll n$   \\ \hline
 Our Approach 1 & Bootstrap & Sub-Weibull & $O(Nnd^2)$  & $d \ll n$ & $d^3 \ll n$   \\ \hline
 Our Approach 2 & $\approx$ Bayes & Sub-Weibull & $O(Nd^2)$ & $d \ll n$ & $d^3 \ll n$   \\ \hline\hline
\end{tabular}
\end{center}
\vspace{0.4cm}
We again mention that in the Spiked regime the condition $d^3 \ll n$ can be improved to $d^2 \ll n$, but we do not pursue this goal in current work but rather focus on FM regime. 

%% file: source/Proof.tex
\subsection{Key ingredients and outline of the proof}
We start by presenting the key ingredients that our main theorems relies on.
All the lemmas stated in this subsection are either borrowed from the literature or proved below in the end of the paper.

\subsubsection{Concentration of sample covariance for sub-Weibull distributions}
\input source/lemmas/Lemmas_KC.tex

\subsubsection{Linear approximation of projectors}
\input source/lemmas/Lemmas_JW.tex

\subsubsection{Alternative representation of $\| \cdot \|_{(\Pp, \Gamma, s_1, s_2)}$}
\input source/lemmas/Lemma_prop.tex

\subsubsection{$\eps$-net argument} \label{eps-net-arg}
    \input source/lemmas/Description_epsnet.tex
    
\subsubsection{Gaussian approximation, anti-concentration and comparison for maxima}
\input source/lemmas/Lemmas_CCK.tex

\subsection{Proof of Theorem \ref{Th:GA1}}
    \input source/proofs/Proof_GA1.tex

\subsection{Proof of Theorem \ref{Th:B1}}
    \input source/proofs/Proof_B1.tex

\subsection{Proof of Theorem \ref{Th:F1}}
    \input source/proofs/Proof_F1.tex
    
\subsection{Proof of Theorem \ref{Th:GA2}}
    \input source/proofs/Proof_GA2.tex
    
\subsection{Proofs of Theorem \ref{Th:B2} and Theorem \ref{Th:F2}}
    \input source/proofs/Proof_BF2.tex
    
\subsection{Proofs of Corollary \ref{Corollary1} and Corollary \ref{Corollary2}}
    \input source/proofs/Proof_C12.tex

\subsection{Proof of Theorem \ref{Th:Power1} and \ref{Th:Power2}}
    \input source/proofs/Proof_Power1.tex


%% file: source/lemmas/Lemmas_KC.tex
\noindent The first lemma describes how the sample covariance concentrates under Assumption \ref{A: tails}. Our concentration for covariance is written in somewhat specific form; the reason for that will be justified in the next subsection.
\begin{lemma} \label{L:KC_SCconc}
    Let the data satisfy Assumption \ref{A: tails}.
    Then the following bound holds with probability $1-1/n$:
    \begin{equation}
        \begin{aligned}
	    \max\limits_{s,t\in[q]}\frac{\| \Pts (\Se - \St) \Ptt \|_{\Fr}}{\sqrt{m_s \mu_s m_t \mu_t}} \leq \scc,
        \nonumber
        \end{aligned}
    \end{equation}
    where 
    \begin{equation}
        \begin{aligned}
	    \scc \eqdef C_\beta c^2 \left( \sqrt{\frac{\log(n)+\log(d)}{n}} + \frac{(\log(n))^{1/\beta} (\log(n)+\log(d))^{2/\beta}}{n}\right).
        \nonumber
        \end{aligned}
    \end{equation}
\end{lemma}

%% file: source/lemmas/Lemmas_JW.tex
The projector of a covariance matrix is a complicated nonlinear operator. We use machinery from \cite{Wahl}, unlike previous works (e.g. \cite{Naumov, Silin_1}) which were based on \cite{Koltchinskii_AACBFBFOSPOSC}, Lemma 2. Novel result from \cite{Wahl} allows to obtain linear approximation for spectral projectors with remainder term bounded by dimension-free rate even for non-Gaussian distributions.
We slightly modify it to prepare it for our framework.
\begin{lemma} \label{L:KL}
    Let $\Sp$ be perturbed covariance matrix.
    Take 
    \begin{equation}
        \begin{aligned}
        	x = \max\limits_{s,t\in[q]} \frac{\| \Pts (\Sp-\St) \Ptt \|_{\Fr}}{\sqrt{m_s\mu_s m_t\mu_t}}
        \nonumber
        \end{aligned}
    \end{equation}
    and assume
    \begin{equation}
        \begin{aligned}
        	x \max\limits_{r\in\J} \relr_r(\St) \leq \frac{1}{6}.
        \label{A:JW}
        \end{aligned}
    \end{equation}
    Then following decomposition holds:
    \begin{equation}
        \begin{aligned}
    	\Ppp - \Pt 
    	= 
    	L_{\J}(\Sp-\St) + R_{\J}(\Sp-\St),
        \label{L:PD}
        \end{aligned}
    \end{equation}
    	where the linear part is 
    \begin{equation}
        \begin{aligned}
    	L_{\J}(\Sp-\St)
    	& \eqdef &
    	\sum\limits_{r \in {\J}} \sum\limits_{s \notin {\J}} 
    			\frac{\Ptr (\Sp-\St) \Pts + \Pts (\Sp-\St) \Ptr}{\mu_r - \mu_s}
        \nonumber
        \end{aligned}
    \end{equation}
    and the remainder term satisfies    
    \begin{equation}
        \begin{aligned}
    	\| R_\J(\Sp-\St)\| \leq
    	\| R_\J(\Sp-\St)\|_{\Fr} \leq Cx^2 \sum\limits_{r\in\J} \left( \relr_r(\St) \sqrt{\sum\limits_{s\neq r} \frac{m_r\mu_r m_s \mu_s}{(\mu_r-\mu_s)^2}} \right) = Cx^2 \rr_\J(\St)^{3/2}.
        \label{L:LA}
        \end{aligned}
    \end{equation}
    
\end{lemma}

The following lemma deals with the remainder term using the previous lemma.
\begin{lemma} \label{L:KL_gen}
    Let $\Sp$ be perturbed covariance matrix (potentially depending on the data and additional source of randomness; e.g. $\Sb$ or $\Sf$) and $\Ppp$ is the corresponding projector.
    Assume for some rate $\widetilde{\psi}_n$ holds
    \begin{equation}
        \begin{aligned}
        	\Prob\left[ \max\limits_{s,t\in[q]} \frac{\| \Pts (\Sp-\Se) \Ptt \|_{\Fr}}{\sqrt{m_s\mu_s m_t\mu_t}} \leq \widetilde{\psi}_n \,\Big|\, \data \right] \geq 1- \frac{1}{n}
        \label{prob_assumption}
        \end{aligned}
    \end{equation}
    with probability $1-1/n$,
    and 
    \begin{equation}
        \begin{aligned}
        	&(\scc \lor \widetilde{\psi}_n) \max\limits_{r\in\J} \relr_r(\St) \leq \frac{1}{12},
        \nonumber
        \end{aligned}
    \end{equation}
Then the following approximation holds with probability $1-2/n$
    \begin{equation}
        \begin{aligned}
    	&\Prob\left[ \left| \| \Ppp - \Pe \|_{(\Pt,\Gamma^\circ,s_1,s_2)} - \| L_\J(\Sp-\Se) \|_{(\Pt,\Gamma^\circ,s_1,s_2)} \right|
    	\leq
    	C (\scc+\widetilde{\psi}_n)^2 \rr_\J(\St)^{3/2} \,\Big|\,\data\right] \geq 
    	\\&\qquad\geq 1-\frac{1}{n}.
        \label{LA_bound}
        \end{aligned}
    \end{equation}
\end{lemma}
As one can see, the setting of the lemma is pretty general, and we are going to apply it in the sequel with different $\Sp$.

%% file: source/lemmas/Lemma_prop.tex
For the clarity of presentation, in Section~\ref{Test} we introduced user-friendly definition of \\$\| \cdot \|_{(\Pp, \Gamma, s_1, s_2)}$ that doesn't require any extra definitions. However, in our proofs it will be more convenient to work with $\| \cdot \|_{(\Pp, \Gamma, s_1, s_2)}$ expressed in a slightly different way.

\begin{definition} \label{myset}
    Let $k$ be an integer and $s\in[k]$. 
    Define
    \begin{equation} 
        \begin{aligned}
            \myset^k_s \eqdef 
            \left\{ x\in\Sph^{k-1}\;\big|\; \max\limits_{x_j\neq 0} (j) - \min\limits_{x_j\neq 0} (j) \leq s-1 \right\}
            = \bigcup\limits_{l=0}^{k-s} \left\{ [0_l^\T, y^\T, 0_{k-l-s}^\T ]^\T\;\big| \;\;y\in\Sph^{s-1} \right\}.
        \nonumber 
        \end{aligned} 
    \end{equation}
\end{definition}

\begin{example}
    Consider, first, extreme cases:
    \begin{itemize}
        \item $s = 1$: we have $\myset^k_1 = \{ \pm e_j\}_{j=1}^k$, where $e_j \in \R^k$ is $j$-th standard basis vector.
        \item $s = k$: we have $\myset^k_k = \Sph^{k-1}$.
    \end{itemize}
    To illustrate Definition~\ref{myset} in a less trivial case, take $k=7$ and $s = 3$. Then $\myset^k_s$ consists of $k$-dimensional unit vectors with support contained in the shadow area of one of the following $k-s+1=5$ variants, depicted on Figure~\ref{draw2}.
\end{example}

\input source/draws/draw2.tex

\begin{lemma}[Additional properties of $\| \cdot \|_{(\Pp, \Gamma, s_1, s_2)}$] \label{P:Properties1}
Fix arbitrary $\Pp$ of rank $m$, $\Gamma = [\Gamma_1\; \Gamma_2]$, $s_1$, $s_2$ as in Definition~\ref{D:norm}. Then, the following holds:
\begin{enumerate}[(i)]
 \item $\| \cdot \|_{(\Pp,\Gamma,s_1,s_2)}$ can be alternatively represented as
\begin{equation} 
    \begin{aligned}
        \| A \|_{(\Pp, \Gamma, s_1, s_2)} \eqdef \frac{1}{2}\| \Gamma_1^\T A \Gamma_1 \| + \frac{1}{2}\| \Gamma_2^\T A \Gamma_2 \| + 
        \sup\limits_{\substack{v \in \myset^m_{s_1} \\ w \in \myset^{d-m}_{s_2}}} 
        v^\T \Gamma_1^\T A \,\Gamma_2 w.
    \nonumber 
    \end{aligned} 
\end{equation}
 \item For a symmetric $A\in\R^{d\times d}$  of the form $A = \Pp A (\Id_d - \Pp) + (\Id_d - \Pp) A \Pp$ we have
 \begin{equation} 
    \begin{aligned}
        \| A \|_{(\Pp,\Gamma,s_1,s_2)} = \sup\limits_{\substack{v \in \myset^m_{s_1}\\ w \in \myset^{d-m}_{s_2}}}  v^\T \Gamma_1^\T  A \Gamma_2 w .
    \nonumber 
    \end{aligned} 
\end{equation}
If, additionally, $s_1 = m$ and $s_2 = d-m$, then
 \begin{equation} 
    \begin{aligned}
        \| A \|_{(\Pp,\Gamma,m,d-m)} = \| A \|,
    \nonumber 
    \end{aligned} 
\end{equation}
while in case $s_1 = 1$ and $s_2 = 1$ we get
 \begin{equation} 
    \begin{aligned}
        \| A \|_{(\Pp,\Gamma, 1, 1)} = \| \Gamma_1^\T A \Gamma_2 \|_{\mymax}.
    \nonumber 
    \end{aligned} 
\end{equation}
\end{enumerate}
\end{lemma}

%% file: source/draws/draw2.tex
\newcommand{\myRect}[5]{
    \draw[#1] (#2,#3) rectangle(#2+#4,#3-#5);
}

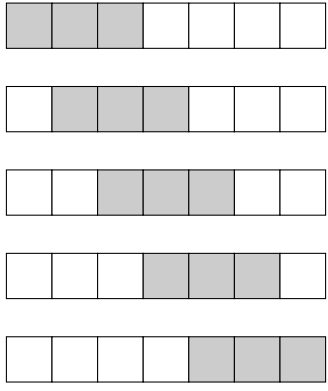
\begin{figure}[!h]

\begin{center}
\begin{tikzpicture}[scale=0.6]

\myRect{fill=gray!40}{0}{0}{1}{1}
\myRect{fill=gray!40}{1}{0}{1}{1}
\myRect{fill=gray!40}{2}{0}{1}{1}
\myRect{fill=white}{3}{0}{1}{1}
\myRect{fill=white}{4}{0}{1}{1}
\myRect{fill=white}{5}{0}{1}{1}
\myRect{fill=white}{6}{0}{1}{1}

\end{tikzpicture}
\end{center}
\begin{center}
\begin{tikzpicture}[scale=0.6]

\myRect{fill=white}{0}{0}{1}{1}
\myRect{fill=gray!40}{1}{0}{1}{1}
\myRect{fill=gray!40}{2}{0}{1}{1}
\myRect{fill=gray!40}{3}{0}{1}{1}
\myRect{fill=white}{4}{0}{1}{1}
\myRect{fill=white}{5}{0}{1}{1}
\myRect{fill=white}{6}{0}{1}{1}

\end{tikzpicture}
\end{center}
\begin{center}
\begin{tikzpicture}[scale=0.6]

\myRect{fill=white}{0}{0}{1}{1}
\myRect{fill=white}{1}{0}{1}{1}
\myRect{fill=gray!40}{2}{0}{1}{1}
\myRect{fill=gray!40}{3}{0}{1}{1}
\myRect{fill=gray!40}{4}{0}{1}{1}
\myRect{fill=white}{5}{0}{1}{1}
\myRect{fill=white}{6}{0}{1}{1}

\end{tikzpicture}
\end{center}
\begin{center}
\begin{tikzpicture}[scale=0.6]

\myRect{fill=white}{0}{0}{1}{1}
\myRect{fill=white}{1}{0}{1}{1}
\myRect{fill=white}{2}{0}{1}{1}
\myRect{fill=gray!40}{3}{0}{1}{1}
\myRect{fill=gray!40}{4}{0}{1}{1}
\myRect{fill=gray!40}{5}{0}{1}{1}
\myRect{fill=white}{6}{0}{1}{1}

\end{tikzpicture}
\end{center}
\begin{center}
\begin{tikzpicture}[scale=0.6]

\myRect{fill=white}{0}{0}{1}{1}
\myRect{fill=white}{1}{0}{1}{1}
\myRect{fill=white}{2}{0}{1}{1}
\myRect{fill=white}{3}{0}{1}{1}
\myRect{fill=gray!40}{4}{0}{1}{1}
\myRect{fill=gray!40}{5}{0}{1}{1}
\myRect{fill=gray!40}{6}{0}{1}{1}

\end{tikzpicture}
\end{center}
\caption{
    The support of any vector in $\myset^7_3$ is included in one of the shadow regions.
    }
\label{draw2}
\end{figure}

%% file: source/lemmas/Description_epsnet.tex
Let $N_{\eps}(\myset^m_{s_1})$ and $N_{\eps}(\myset^{d-m}_{s_2})$ be proper $\varepsilon$-nets of $\myset^m_{s_1}$ and $\myset^{d-m}_{s_2}$, respectively, w.r.t. Euclidean distance. Let us explicitly demonstrate how we construct them; namely, we construct $N_{\eps}(\myset^m_{s_1})$, while $N_{\eps}(\myset^{d-m}_{s_2})$  can be constructed similarly. Consider a proper $\epsilon$-net of $\Sph^{s_1-1}$ w.r.t. Euclidean distance and denote it as $N_{\eps}(\Sph^{s_1-1})$. Take
    \begin{equation}
        \begin{aligned}
            N_{\eps}(\myset^m_{s_1}) = \bigg\{ [0_{k}^\T, v^\T, 0_{m-k-s_1}^\T]^\T\,\big|\; v\in N_{\eps}(\Sph^{s_1-1}),\; k \in \{ 0, \ldots, m-s_1\} \bigg\}.
        \label{eps-net construction}
        \end{aligned}
    \end{equation}
By Definition~\ref{myset} it is trivial to see that $N_{\eps}(\myset^m_{s_1})$ is indeed an $\varepsilon$-net of $\myset^m_{s_1}$.

    Consider all possible pairs $(\Gamma_1 v, \Gamma_2 w)$ such that $v \in N_{\eps}(\myset^m_{s_1}), w \in N_{\eps}(\myset^{d-m}_{s_2})$. Enumerate them $\{(v_j, w_j)\}_{j=1}^p$ with $p = p(\eps, d, m, s_1, s_2) =  |N_{\eps}(\myset^m_{s_1})| \cdot |N_{\eps}(\myset^{d-m}_{s_2})|$. Note that the constructed $\eps$-net is different for different $\Pp$ and $\Gamma$.

The following lemma provides standard approximation of infinite supremum by finite maximum over the $\eps$-net.

\begin{lemma}[Discretization] \label{L:discretization}
    Let $\Pp$, $\Gamma$, $s_1$, $s_2$ be as in Definition~\ref{D:norm}. For any symmetric matrix $A\in\R^{d\times d}$ satisfying 
    \begin{equation}
        \begin{aligned}
        A = \Pp A (\Id_d-\Pp) + (\Id_d-\Pp) A \Pp
        \nonumber
        \end{aligned}
    \end{equation}
    the following bounds hold:
    \begin{equation}
        \begin{aligned}
            \max\limits_{j\in[p]}  v_j^\T A w_j
        	\leq
        	\| A \|_{(\Pp, \Gamma,s_1,s_2)} \leq \frac{1}{1-2\eps} \max\limits_{j\in[p]} v_j^\T A w_j,
        \nonumber
    \end{aligned}
    \end{equation}
       
\end{lemma}

The size of the $\eps$-net can be bounded according to the following lemma.
\begin{lemma}[Covering number] \label{L:epsnet_size}
    The following bound holds:
    \begin{equation}
        \begin{aligned}
            \log(p(\eps, d, m, s_1, s_2)) \leq (s_1+s_2) \log\left(\frac{3}{ \eps}\right) + 2\log(d).
        \nonumber
        \end{aligned}
    \end{equation}
\end{lemma}
In our proofs, we will use $\eps = 1/n$. This fixes $p$ to be
    \begin{equation}
            \begin{aligned}
            	p \leq \exp\left( (s_1 + s_2)\log(3n) + 2\log(d)\right),
                \nonumber
            \end{aligned}
        \end{equation}

%% file: source/lemmas/Lemmas_CCK.tex
We will be using Gaussian approximation, anti-concentration and comparison results for maximum of a random vector.
Before we state the specific results from \cite{GAR+Max,Comparison}, let us introduce the framework from these papers. Suppose we have a collection of $n$ independent zero-mean random vectors in $\R^p$:
    \begin{equation}
        \begin{aligned}
        	x_i = \{ x_{ij} \}_{j=1}^p,\;\;\; i \in [n].
            \nonumber
        \end{aligned}
    \end{equation}
Let $y_i \sim \mathcal{N}_p(0, \Cov(x_i)) , i\in[n]$ be a collection of Gaussian vectors in $\R^p$ with the same covariances as these of $x_i$'s.
Denote
    \begin{equation}
        \begin{aligned}
        	\Eb\left[\;\cdot\;\right] \eqdef \frac{1}{n} \sum\limits_{i=1}^n \E[\;\cdot\;], \;\;\text{ e.g. }\;\;\Eb\left[x_{ij}^2\right] = \frac{1}{n} \sum\limits_{i=1}^n \E[x_{ij}^2].
            \nonumber
        \end{aligned}
    \end{equation}
    Based on that, introduce
    \begin{equation}
        \begin{aligned}
        	M_k \eqdef \max\limits_{j\in[p]}\left(\Eb\left[|x_{ij}|^k\right]\right)^{1/k}.
            \nonumber
        \end{aligned}
    \end{equation}
    Finally, define $u_x(\gamma)$ as the smallest $u$ such that
    \begin{equation}
        \begin{aligned}
        	 \Prob\left[ |x_{ij}| \leq u\cdot\left(\Eb\left[|x_{ij}|^2\right]\right)^{1/2} \;\;\forall i\in[n]\;\; \forall j\in[p]\right] \geq 1 - \gamma,
            \nonumber
        \end{aligned}
    \end{equation}
    and define $u_y(\gamma)$ similarly for the Gaussian counterpart $y_{ij}$, and denote $u(\gamma) = u_x(\gamma) \lor u_y(\gamma)$. Now we are ready to state the results. The first one is Gaussian approximation for maxima of sum of random vectors.
\begin{lemma}[\cite{GAR+Max}, Theorem 2.2: Main result 1, Gaussian approximation] \label{L:CCK_GA}
        Suppose that there are some constants $0 < c_1 < C_1$ such that $c_1 \leq \Eb\left[x_{ij}^2\right] \leq C_1$ for all $j \in [p]$. Then for every $\gamma \in (0, 1)$,
        \begin{equation}
            \begin{aligned}
            	&\sup\limits_{z\in\R}\left| \Prob\left[ \max\limits_{j\in[p]} \frac{1}{\sqrt{n}}\sum\limits_{i=1}^n x_{ij} \leq z  \right] - \Prob\left[ \max\limits_{j\in[p]} \frac{1}{\sqrt{n}}\sum\limits_{i=1}^n y_{ij} \leq z  \right] \right| \leq \\
            	&\qquad \qquad \leq C\left\{ n^{-1/8} \left(M_3^{3/4}\lor M_4^{1/2}\right) \left(\log(pn/\gamma)\right)^{7/8} + n^{-1/2}\left(\log(pn/\gamma)\right)^{3/2}u(\gamma) + \gamma\right\},
                \nonumber
            \end{aligned}
        \end{equation}
        where $C > 0$ is a constant that depends on $c_1$ and $C_1$ only.
    \end{lemma}
\noindent The following is the anti-concentration result from \cite{Comparison}.
    \begin{lemma}[\cite{Comparison}, Corollary 1: Anti-concentration] \label{L:CCK_AC}
        Let $(Z_1, \ldots, Z_p)^{\top}$ be a centered Gaussian random vector in $\R^p$ with $\sigma_j^2 \eqdef \E\left[Z_j^2\right] > 0$ for all $j \in [p]$. Let $\underline{\sigma} \eqdef \min\limits_{j\in [p]} \sigma_{j}$ and $\overline{\sigma} \eqdef \max\limits_{j\in [p]} \sigma_j$. Then for every $\epsilon > 0$,
        \begin{equation}
            \begin{aligned}
            	\sup\limits_{z\in\R} \Prob \left[ \left| \max\limits_{j\in[p]} Z_j - z \right| \leq \epsilon\right] \leq C\epsilon \sqrt{1 \lor \log(p/\epsilon)},
                \nonumber
            \end{aligned}
        \end{equation}
        where $C > 0$ depends only on $\underline{\sigma}$ and $\overline{\sigma}$.
    \end{lemma}
\noindent The following is the comparison result from \cite{Comparison}. 
\begin{lemma}[\cite{Comparison}, Theorem 2: Comparison of distributions] \label{L:CCK_GC}
    Let $Z = (Z_1, \ldots, Z_p)^{\top}$ and $Y = (Y_1, \ldots, Y_p)^{\top}$ be centered Gaussian random vectors in $\R^p$ with covariances $\{ \sigma_{jk}^{Z} \}_{j,k=1}^p$ and $\{ \sigma_{jk}^{Y} \}_{j,k=1}^p$, respectively. Define 
    \begin{equation}
        \begin{aligned}
        	\Delta \eqdef \max\limits_{j,k\in[p]} |\sigma_{jk}^{Z} - \sigma_{jk}^{Y}| \text{ and } a_p \eqdef \E\left[ \max\limits_{j\in [p]} (Y_j/\sigma_{jj}^{Y}) \right].
            \nonumber
        \end{aligned}
    \end{equation}
    Suppose that $p \geq 2$ and $\sigma_{jj}^{Y} > 0$ for all $j \in [p]$. Then
    \begin{equation}
        \begin{aligned}
        	 &\sup\limits_{z \in \R} \left| \Prob\left[ \max\limits_{j\in[p]} Z_j \leq z \right] - \Prob\left[ \max\limits_{j\in[p]} Y_j \leq z \right]\right| \leq\\
        	 &\qquad \qquad \leq C\Delta^{1/3} \left\{ 1 
        	 \lor a_p^2 \lor \log(1/\Delta) \right\}^{1/3} \; (\log p)^{1/3},
            \nonumber
        \end{aligned}
    \end{equation}
    where $C > 0$ depends only on $\min\limits_{j\in [p]} \sigma_{jj}^{Y}$ and $\max\limits_{j\in [p]} \sigma_{jj}^{Y}$. Moreover, in the worst case, \\$a_p \leq \sqrt{2 \log p}$, so that
    \begin{equation}
        \begin{aligned}
        	 &\sup\limits_{z \in \R} \left| \Prob\left[ \max\limits_{j\in[p]} Z_j \leq z \right] - \Prob\left[ \max\limits_{j\in[p]} Y_j \leq z \right]\right| 
        	 \leq
        	 C^{\prime}\Delta^{1/3} \left\{ 1  \lor \log(p/\Delta) \right\}^{2/3} ,
            \nonumber
        \end{aligned}
    \end{equation}
     where as before $C^{\prime} > 0$ depends only on $\min\limits_{j\in [p]} \sigma_{jj}^{Y}$ and $\max\limits_{j\in [p]} \sigma_{jj}^{Y}$.
    \end{lemma}
    Now we are equipped to proceed with the proof of the main results.

%% file: source/proofs/Proof_GA1.tex
    \subsubsection{Approximation by finite maximum}
    Throughout the proof we work with $\| \cdot \|_{(\Pt, \Gamma^\circ, s_1, s_2)}$-norm, and apply all lemmas from previous subsection with $\Pp = \Pt$ and $\Gamma = \Gamma^\circ$.
    
    In accordance with Lemma~\ref{L:KL_gen} (applied with $\Sp = \St$, $\Ppp = \Pt$), we start by working with the linear part $\sqrt{n}\,L_\J(\Se-\St)$ of $\sqrt{n} (\Pe - \Pt)$. Let $\{(v_j, w_j)\}_{j=1}^p$ be the $\eps$-net constructed in subsection \ref{eps-net-arg} for $\Pt$ and $\Gamma^\circ$ with $\eps = 1/n$. Observe that $L_\J(\Se - \St)$ satisfies the condition of Lemma~\ref{L:discretization} with $\Pp = \Pt$, hence, for
    \begin{equation}
        \begin{aligned}
        	L_{disc} &\eqdef 
        	\max\limits_{j\in[p]} \, v_j^{\top} \left(\sum\limits_{r\in\J} \sum\limits_{s\notin \J}\frac{ \Ptr(\Se-\St)\Pts + \Pts(\Se-\St)\Ptr}{\mu_r - \mu_s}\right) w_j
        	\\&= \max\limits_{j\in[p]} \, \sum\limits_{r\in\J} \sum\limits_{s\notin \J}\frac{v_j^{\top} \Ptr(\Se-\St)\Pts w_j}{\mu_r - \mu_s}
            \nonumber
        \end{aligned}
    \end{equation}
    we have
    \begin{equation}
        \begin{aligned}
        	L_{disc} \leq \| L_\J(\Se-\St) \|_{(\Pt,\Gamma^\circ, s_1, s_2)} \leq \frac{1}{1-2\eps} L_{disc}.
            \label{Eq: disc_bounds}
        \end{aligned}
    \end{equation}
    Now let us represent $L_{disc}$ in a different way. Introduce for $i\in [n], j \in [p]$
    \begin{equation}
        \begin{aligned}
        	x_{ij} \eqdef 
        	\sum\limits_{r\in\J} \sum\limits_{s\notin \J}\frac{v_j^{\top} \Ptr(X_i X_i^{\top}-\St)\Pts w_j}{\mu_r - \mu_s} 
        	=
        	\sum\limits_{r\in\J} \sum\limits_{s\notin \J}\frac{v_j^{\top} \Ptr X_i X_i^{\top}\Pts w_j}{\mu_r - \mu_s}.
            \nonumber
        \end{aligned}
    \end{equation}
    Therefore,
    \begin{equation}
        \begin{aligned}
        	 \sqrt{n} L_{disc}
        	 =
        	 \max\limits_{j\in[p]} \frac{1}{\sqrt{n}}\sum\limits_{i=1}^n x_{ij}.
            \nonumber
        \end{aligned}
    \end{equation}
    We can arrange this random variables as i.i.d. centered random vectors in $\R^p$:
    \begin{equation}
        \begin{aligned}
        x_i \eqdef \{ x_{ij} \}_{j=1}^p.
        \nonumber
        \end{aligned}
    \end{equation}
    Lemma~\ref{L:CCK_GA} suggests that the distribution of $\max\limits_{j\in[p]} \frac{1}{\sqrt{n}}\sum\limits_{i=1}^n x_{ij}$ can be approximated by the distribution of its Gaussian analogue $\max\limits_{j\in[p]} \frac{1}{\sqrt{n}}\sum\limits_{i=1}^n y_{ij}$, where $y_i \eqdef \{ y_{ij} \}_{j=1}^p$ are i.i.d. centered Gaussian random vectors with the same covariance structure as $x_i$'s. In other terms, introducing 
    \begin{equation}
        \begin{aligned}
        	Y \eqdef \frac{1}{\sqrt{n}} \sum\limits_{i=1}^n y_i \;\sim \mathcal{N}_p\left(0, \;\frac{1}{n}\sum\limits_{i=1}^n \Cov(x_i) \right) \sim
        	\mathcal{N}_p\left(0, \Cov(x_1) \right),
            \nonumber
        \end{aligned}
    \end{equation}
    we would like to use the distribution of $\max\limits_{j\in[p]} Y_j$ as the approximation for the distribution of $\sqrt{n}\,L_{disc}$, consequently for the distribution of $\sqrt{n}\,\| L_\J(\Se-\St)\|_{(\Pt, \Gamma^\circ, s_1, s_2)}$, and eventually for the distribution of $\sqrt{n}\,\| \Pe - \Pt\|_{(\Pt, \Gamma^\circ, s_1, s_2)}$.
    
    \subsubsection{Verifying the conditions} \label{S:verify}
    In order to apply Lemma~\ref{L:CCK_GA} to our situation, we need to verify the conditions and compute some quantities. \\ \\
    \underline{Computing the covariance}: Let us start with computing the covariance of $x_i$'s (and hence, the covariance of the Gaussian $y_i$'s and $Y$). First, using Assumption \ref{A: independence}, we compute for $i\in[n], j\in [p], r,r^\prime\in\J, s,s^\prime\notin\J$
    \begin{equation}
        \begin{aligned}
        	&\E\left[ v_j^{\top} \Ptr X_i X_i^{\top} \Pts w_j \cdot v_k^{\top} \Ptrp X_i X_i^{\top} \Ptsp w_k  \right] \\
        	&\qquad=
        	\E\left[ v_j^{\top} \Ptr (\Pt X_i) \cdot (\Ptc X_i)^{\top} \Pts w_j \cdot v_k^{\top} \Ptrp (\Pt X_i) \cdot (\Ptc X_i)^{\top} \Ptsp w_k  \right] \\
        	&\qquad = \E\left[ v_j^{\top} \Ptr X_i X_i^{\top} \Ptrp v_k\right] \cdot \E\left[ w_j^{\top} \Pts X_i X_i^{\top} \Ptsp w_k  \right]
        	\\&\qquad= \delta_{r,r^{\prime}} \, \delta_{s,s^\prime} \, \mu_r \mu_s \cdot (v_j^{\top} \Ptr v_k) \cdot (w_j^{\top} \Pts w_k),
            \nonumber
        \end{aligned}
    \end{equation}
    where $\delta_{\cdot,\cdot}$ is the Kronecker delta.
    Further,
    \begin{equation}
        \begin{aligned}
        	\E[x_{ij}x_{ik}] &= 
        	\sum\limits_{r\in\J}\sum\limits_{s\notin\J} \frac{\mu_r \mu_s}{(\mu_r - \mu_s)^2}\, \, (v_j^{\top} \Ptr v_k) \cdot (w_j^{\top} \Pts w_k).
            \nonumber
        \end{aligned}
    \end{equation}
    So,
    \begin{equation}
        \begin{aligned}
        	&\Cov(x_i) = \{ \sigma^{\Su}_{jk} \}_{j,k=1}^p, \;\text{ where}\\
        	&\sigma^{\Su}_{jk} = \sum\limits_{r\in\J}\sum\limits_{s\notin\J} \frac{\mu_r \mu_s}{(\mu_r - \mu_s)^2}\, \, (v_j^{\top} \Ptr v_k) \cdot (w_j^{\top} \Pts w_k).
            \nonumber
        \end{aligned}
    \end{equation}

    Then, let us show the existence of $c_1$ and $C_1$ required in Lemma~\ref{L:CCK_GA}, bound $M_3$ and $M_4$ and estimate $u(\gamma)$. \\ \\
    \underline{Showing the existence of $c_1$ and $C_1$}: to do that, write
    \begin{equation}
        \begin{aligned}
        	\Eb\left[ x_{ij}^2 \right] =
        	\frac{1}{n} \sum\limits_{i=1}^n \E[x_{ij}^2] \stackrel{i.i.d.}{=}
        	\E[x_{1j}^2] = \sigma^{\Su}_{jj} = \sum\limits_{r\in\J}\sum\limits_{s\notin\J} \frac{\mu_r \mu_s}{(\mu_r - \mu_s)^2}\, \, (v_j^{\top} \Ptr v_j) \cdot (w_j^{\top} \Pts w_j).
            \nonumber
        \end{aligned}
    \end{equation}
    Note that
    \begin{equation}
        \begin{aligned}
        	\sigma^{\Su}_{jj} &\leq \max\limits_{r\in\J, s\notin\J} \frac{\mu_r \mu_s}{(\mu_r - \mu_s)^2} \sum\limits_{r\in\J}\sum\limits_{s\notin\J} \, \, (v_j^{\top} \Ptr v_j) \cdot (w_j^{\top} \Pts w_j) = \chigh^2 \cdot (v_j^{\top} \Pt v_j) \cdot (w_j^{\top} \Ptc w_j) = \chigh^2,\\
        	\sigma^{\Su}_{jj} &\geq \min\limits_{r\in\J, s\notin\J} \frac{\mu_r \mu_s}{(\mu_r - \mu_s)^2} \sum\limits_{r\in\J}\sum\limits_{s\notin\J} \, \, (v_j^{\top} \Ptr v_j) \cdot (w_j^{\top} \Pts w_j) = \clow^2 \cdot (v_j^{\top} \Pt v_j) \cdot (w_j^{\top} \Ptc w_j) =\clow^2.
            \nonumber
        \end{aligned}
    \end{equation}
    This implies that there exist $c_1 = \clow^2 > 0$ and $C_1 = \chigh^2 > 0$ satisfying the condition
    \begin{equation}
        \begin{aligned}
        	c_1 \leq \Eb\left[ x_{ij}^2 \right] \leq C_1
            \nonumber
        \end{aligned}
    \end{equation}
    for all $j \in [p]$. 
    \\ \\
    \underline{Upperbounding $M_3$ and $M_4$}: 
    we have
    \begin{equation}
        \begin{aligned}
        	 |x_{ij}| &= \left| \sum\limits_{r\in\J} \sum\limits_{s\notin \J}\frac{v_j^{\top} \Ptr X_i X_i^{\top}\Pts w_j}{\mu_r - \mu_s} \right| \leq 
        	 \\&\leq \max\limits_{r\in\J, s\notin\J}  \frac{\sqrt{\mu_r \mu_s}}{|\mu_r - \mu_s|}  \cdot  \sum\limits_{r\in\J} \left| v_j^{\top} \Ptr \St^{-1/2} X_i \right| \cdot \sum\limits_{s\notin \J} \left| w_j^{\top}\Pts \St^{-1/2}X_i \right|.
            \nonumber
        \end{aligned}
    \end{equation}
    Let us deal with $\sum\limits_{r\in\J} \left| v_j^{\top} \Ptr \St^{-1/2} X_i \right|$.
    Represent
    \begin{equation}
        \begin{aligned}
        	 \left| v_j^{\top} \Ptr \St^{-1/2}X_i \right|
        	 =
        	 \overline{v}_{j,r}^{\top} \St^{-1/2} X_i,
            \nonumber
        \end{aligned}
    \end{equation}
    where $\overline{v}_{j,r}$ is either $\Ptr v_j$ or $(-\Ptr v_j)$ depending on the sign of $v_j^{\top} \Ptr \St^{-1/2}X_i$. Note that $\{ \overline{v}_{j,r} \}_{r\in\J}$ are orthogonal. Define $\overline{v}_j \eqdef \sum\limits_{r\in\J} \overline{v}_{j,r}$ with the squared norm 
    \begin{equation}
        \begin{aligned}
        	 \| \overline{v}_j \|^2 
        	 = \sum\limits_{r\in\J} \|\overline{v}_{j,r}\|^2
        	 = \sum\limits_{r\in\J} \| \Ptr v_j\|^2 = v_j^{\T} \Pt v_j = 1,
            \nonumber
        \end{aligned}
    \end{equation}
    where the first equality is due to orthogonality. Hence,
    \begin{equation}
        \begin{aligned}
        	 \sum\limits_{r\in\J} \left| v_j^{\top} \Ptr \St^{-1/2}X_i \right| &=
        	 \sum\limits_{r\in\J}  \overline{v}_{j,r}^{\top} \St^{-1/2}X_i =
        	 \overline{v}_{j}^{\top} \St^{-1/2}X_i.
            \nonumber
        \end{aligned}
    \end{equation}
    Similarly, 
    \begin{equation}
        \begin{aligned}
        	 \sum\limits_{s\notin\J} \left| w_j^{\top} \Pts \St^{-1/2}X_i \right| &=
        	 \sum\limits_{s\notin\J}  \overline{w}_{j,s}^{\top} \St^{-1/2}X_i =
        	 \overline{w}_{j}^{\top} \St^{-1/2}X_i
            \nonumber
        \end{aligned}
    \end{equation}
    with some $\| \overline{w}_{j}\| = 1$.
    Thus,
    \begin{equation}
        \begin{aligned}
        	 |x_{ij}| \leq \chigh \cdot\overline{v}_{j}^{\top} \St^{-1/2}X_i \cdot \overline{w}_{j}^{\top} \St^{-1/2}X_i
            \label{Eq:bound_x_ij}
        \end{aligned}
    \end{equation}
    with $\| \overline{v}_{j}\| = 1$, $\| \overline{w}_{j}\| = 1$ and $\overline{v}_{j}, \overline{w}_{j}$ are orthogonal. 
    Therefore,
    \begin{equation}
        \begin{aligned}
        	\left(\Eb\left[ x_{ij}^4\right] \right)^{1/4} &= \left(\frac{1}{n}\sum\limits_{i=1}^n \E\left[ x_{ij}^4\right] \right)^{1/4} \stackrel{i.i.d.}{=} \E\left[ x_{1j}^4\right]^{1/4}
        	\leq
            \chigh\cdot \E\left[ (\overline{v}_j^{\top} \St^{-1/2}X_1)^4 (\overline{w}_j^{\top}\St^{-1/2} X_1)^4\right]^{1/4} \\
        	&=
        	\chigh \cdot \E\left[ (\overline{v}_j^{\top} \St^{-1/2}X_1)^8 \right]^{1/8} \E\left[(\overline{w}_j^{\top} \St^{-1/2} X_1)^8\right]^{1/8}
        	\lesssim 8^{2/\beta}\,\chigh \;\;\;\text{for all } j\in[p],
            \nonumber
        \end{aligned}
    \end{equation}
     where we used Cauchy-Schwarz inequality and the moment bound for sub-Weibull distributions, see \cite{Weibull}, p.7, together with Assumption~\ref{A: tails}.
     So, $M_4 \lesssim 8^{2/\beta} \,\chigh$. By Jensen's inequality, $M_3 \leq M_4$. 
     \\ \\
     \underline{Estimating $u_x(\gamma)$, $u_y(\gamma)$ and $u(\gamma)$}: using $|x_{ij}| \leq \chigh \cdot | \overline{v}_j^{\top} \St^{-1/2} X_i | \cdot | \overline{w}_j^{\top} \St^{-1/2} X_i |$, for arbitrary $u > 0$ we write
     \begin{equation}
            \begin{aligned}
            	&\Prob\left[ |x_{ij}| > u\cdot \left( \Eb\left[x_{ij}^2\right] \right)^{1/2} \;\;\forall i\in [n] \;\forall j\in [p] \right] 
            	\\
            	&\qquad\leq\Prob\left[\,\chigh\cdot| \overline{v}_j^{\top} \St^{-1/2}X_i | \cdot | \overline{w}_j^{\top} \St^{-1/2}X_i | > u\cdot \clow \;\;\;\forall i\in [n] \; \forall j\in [p] \right]\\ 
            	&\qquad=\Prob\left[\,| \overline{v}_j^{\top} \St^{-1/2}X_i | \cdot | \overline{w}_j^{\top} \St^{-1/2}X_i | > u\cdot \clow/\chigh \;\;\;\forall i\in [n] \; \forall j\in [p] \right] \\
            	&\qquad
            	\stackrel{\substack{union\\bound}}{\leq}
            	\sum\limits_{i=1}^n \sum\limits_{j=1}^p \Prob\left[\,| \overline{v}_j^{\top} \St^{-1/2} X_i | \cdot | \overline{w}_j^{\top} \St^{-1/2} X_i | > u/\cond \right] 
            	\leq 
            	n\,p\cdot 2\exp\left( - (u/\cond c^2)^{\beta/2}  \right),
                \nonumber
            \end{aligned}
        \end{equation}
     with the last inequality is due to Assumption \ref{A: tails} and Lemma \ref{L:KC_Prod}. This, by definition of $u_x(\gamma)$, implies
     \begin{equation}
            \begin{aligned}
            	u_x(\gamma) \leq \cond c^2\left( \log(2pn) + \log(1/\gamma) \right)^{2/\beta}.
                \nonumber
            \end{aligned}
        \end{equation}
    At the same time, $y_{ij}$ is centered Gaussian random variable with the variance $\E\left[y_{ij}^2\right] = \Eb\left[y_{ij}^2\right]$, so $y_{ij}/\left(\Eb\left[y_{ij}^2\right]\right)^{1/2} \sim \mathcal{N}(0,1)$. Thus,
    \begin{equation}
            \begin{aligned}
            	&\Prob\left[ |y_{ij}| > u\cdot \left( \Eb\left[y_{ij}^2\right] \right)^{1/2} \;\;\forall i\in [n] \; \forall j\in [p] \right] 
            	=
            	\Prob\left[ |y_{ij}/\left( \Eb\left[y_{ij}^2\right] \right)^{1/2}| > u \;\;\forall i\in [n] \;\forall j\in [p] \right]  \\
            	&\qquad\stackrel{\substack{union\\bound}}{\leq}
            	\sum\limits_{i=1}^n \sum\limits_{j=1}^p \Prob\left[ |y_{ij}/\left( \Eb\left[y_{ij}^2\right] \right)^{1/2}| > u \right]  
            	=
            	n\,p\cdot\Prob\left[ |\mathcal{N}(0,1)| > u \right] 
            	\leq n\,p\cdot 2\,\exp\left( -u^2/2 \right),
                \nonumber
            \end{aligned}
        \end{equation}
    which yields by definition of $u_y(\gamma)$
    \begin{equation}
            \begin{aligned}
                u_y(\gamma) \leq \sqrt{2(\log(2pn) + \log(1/\gamma))}.
                \nonumber
            \end{aligned}
        \end{equation}
    For our purposes we can take $\gamma = 1/n$.
    Therefore, 
    \begin{equation}
            \begin{aligned}
                u(\gamma) = u_x(\gamma) \lor u_y(\gamma) \lesssim \cond c^2 \left( \log(2pn^2) \right)^{2/\beta}.
                \nonumber
            \end{aligned}
        \end{equation}
    \subsubsection{Applying Gaussian approximation and anti-concentration} \label{S:GA_AC}
    To catch the dependence on $\clow$ and $\chigh$ more carefully, we apply Lemma~\ref{L:CCK_GA} not to $x_{ij}$ and $y_{ij}$, but rather to $x^\prime_{ij} := x_{ij}/\clow$ and $y^\prime_{ij} := y_{ij}/\clow$. Then, the conditions verified above translate into
    \begin{equation}
        \begin{aligned}
            &1 \leq \Eb\left[ {x_{ij}^\prime}^2 \right] \leq \frac{\chigh^2}{\clow^2} = \cond^2,\\
            &M_k^\prime = \frac{M_k}{\clow} \lesssim 8^{2/\beta}\cond\;\;\;\text{ for }k=3,4,\\
            &u^\prime(\gamma) = u(\gamma) \lesssim \cond c^2 \left( \log(2pn^2) \right)^{2/\beta}.
        \nonumber
        \end{aligned}
    \end{equation}
   Obviously, passing from $x_{ij}, y_{ij}$ to $x^\prime_{ij}, y^\prime_{ij}$ does not change Kolmogorov distance, so we proved
        \begin{equation}
            \begin{aligned}
            	&\sup\limits_{z\in\R} \left| \Prob\left[ \sqrt{n}L_{disc} \leq z \right] - \Prob\left[ \max\limits_{j\in[p]} Y_j \leq z \right] \right| =
            	\sup\limits_{z\in\R} \left| \Prob\left[ \sqrt{n}L_{disc}/\clow \leq z \right] - \Prob\left[ \max\limits_{j\in[p]} Y_j/\clow \leq z \right] \right|\\
            	&\qquad\leq C(1, \cond^2)
            	\left\{
            	8^{3/(2\beta)} \cond^{3/4} \left( \frac{( \log(pn^2))^7}{n} \right)^{1/8} 
            	+ \cond c^2 \left( \frac{(\log(2pn^2))^{3+4/\beta}}{n} \right)^{1/2} 
            	+ \frac{1}{n}
            	\right\}\\
            	&\qquad\leq C_\cond
            	\left\{
            	8^{3/(2\beta)} \left( \frac{( \log(pn^2))^7}{n} \right)^{1/8} 
            	+ c^2 \left( \frac{(\log(2pn^2))^{3+4/\beta}}{n} \right)^{1/2}
            	\right\}.
                \label{Eq: disc}
            \end{aligned}
        \end{equation}
    Crucial observation here is that the obtained bound depends on $\clow, \chigh$ only through $\cond = \chigh/\clow$.
    \\ \\
    \underline{Getting back from finite maximum to supremum of infinite-state process:}
    now we want to derive the same result, but for $\sqrt{n}\,\| L_\J(\Se - \St) \|_{(\Pt, \Gamma^\circ,s_1,s_2)}$ rather than $\sqrt{n} L^{disc}$. To do that, we clearly need to use \eqref{Eq: disc_bounds} and \eqref{Eq: disc}.
    Let's bound
     \begin{equation}
        \begin{aligned}
        	\Diamond &\eqdef \sup\limits_{z\in\R}\left| \Prob\left[ \max\limits_{j\in[p]} Y_j \leq z \right] - \Prob\left[ \frac{1}{1-2\eps}\max\limits_{j\in[p]} Y_j \leq z \right]\right|
        	\\&= \sup\limits_{z\in\R}
        	\Prob\left[ \max\limits_{j\in[p]} Y_j \in \left[(1-2\eps)z, \; z \right] \right]
        	= \Prob\left[ \max\limits_{j\in[p]} (Y_j/\clow) \in \left[(1-2\eps)z, \; z \right] \right],
            \nonumber
        \end{aligned}
    \end{equation}
    where we again pass from $\max\limits_{j\in[p]} Y_j$ to $\max\limits_{j\in[p]} (Y_j/\clow)$.
    First, notice that since each $(Y_j/\clow)$ is Gaussian with variance at most $\chigh^2/\clow^2$, then all $(Y_j/\clow)$ are sub-Gaussian with parameter $\cond^2$. Then, e.g. by Lemma 5.2 (Maximal tail inequality for sub-Gaussuan random variables) from \cite{Lectures}, we have for all $\delta \in (0,\;1)$
    \begin{equation}
        \begin{aligned}
        	\Prob\left[ \max\limits_{j\in[p]} (Y_j/\clow) \leq \cond \sqrt{\log(p) + \log(1/\delta)} \right] \geq 1-\delta.
            \nonumber
        \end{aligned}
    \end{equation}
    Thus, taking $\delta = 1/n$ and assuming $\eps \leq 1/4$, for $ z \geq 2\cond \sqrt{\log(pn)} \geq \frac{\cond \sqrt{\log(p) + \log(1/\delta)}}{1-2\eps}$ we have
    \begin{equation}
        \begin{aligned}
            \Prob\left[ \max\limits_{j\in[p]} (Y_j/\clow) \leq (1-2\eps)z \right] \geq 1-\frac{1}{n},
            \nonumber
        \end{aligned}
    \end{equation}
    which implies $\Diamond \leq 1/n$. On the other hand, for $ z \leq 2\cond \sqrt{\log(pn)}$ it is better to apply the anti-concentration for Gaussian random vector.
    Applied to our setting, Lemma~\ref{L:CCK_AC} implies
    \begin{equation}
            \begin{aligned}
            	& \Prob\left[ \max\limits_{j\in[p]} (Y_j/\clow) \in \left[(1-2\eps)z, \; z \right] \right] 
            	\leq 
            	C\eps z \sqrt{1 \lor \log(p/\eps z)} 
            	\leq C\eps z\sqrt{\log(ep/\eps z)}\\
            	&\qquad\leq 
            	2C\cond\eps \sqrt{\log(pn)\cdot \log\left(\frac{ep}{2\cond\eps\sqrt{\log(pn)}}\right)}
            	\leq
            	2C\cond\eps \sqrt{\log(pn)\cdot \log\left(\frac{p}{\cond\eps}\right)},
                \nonumber
            \end{aligned}
        \end{equation}
    where $C$ depends only on $\min\limits_{j\in [p]} (\sigma_{jj}^{\Su}/\clow^2)$ and $\max\limits_{j\in [p]} (\sigma_{jj}^{\Su}/\clow^2)$, but effectively on $\clow^2/\clow^2=1$ and $\chigh^2/\clow^2=\cond^2$. Here we used also that $\eps z\sqrt{\log(ep/\eps z)}$ is increasing in $z$ together with assumption $2\cond\eps \sqrt{\log(pn)} \leq ep$ (which anyway should be fulfilled, otherwise our results makes no sense). Combining the bounds on $\Diamond$ for two different regimes of $z$ and recalling $\eps = 1/n$, we get
    \begin{equation}
            \begin{aligned}
            	\Diamond \leq  2C_\cond\cond\eps \sqrt{\log(pn)\cdot \log\left(\frac{p}{\cond\eps}\right)} \;\lor\; \frac{1}{n}
            	\leq \frac{2C_\cond\cond\log(pn)}{n} + \frac{1}{n} \leq \frac{(2C_\cond+1)\cond\log(pn)}{n}.
                \nonumber
            \end{aligned}
        \end{equation}
    This bound, together with \eqref{Eq: disc} and bounds \eqref{Eq: disc_bounds}, yields
        \begin{equation}
            \begin{aligned}
            	&\sup\limits_{z\in\R} \left| \Prob\left[  \sqrt{n}\,\|L_\J(\Se-\St)\|_{(\Pt,\Gamma^\circ,s_1,s_2)} \leq z \right] - \Prob\left[ \max\limits_{j\in[p]} Y_j \leq z \right] \right| 
            	\leq \\
            	&\qquad\leq C_\cond
            	\left\{
            	8^{3/(2\beta)} \left( \frac{( \log(pn^2))^7}{n} \right)^{1/8} 
            	+ c^2 \left( \frac{(\log(2pn^2))^{3+4/\beta}}{n} \right)^{1/2} 
            	\right\} + \Diamond .
                \nonumber
            \end{aligned}
        \end{equation}
    Adjusting the dependence on $\cond$ in $C_\cond$ makes $\Diamond$ negligible compared to the current error term. We obtained
    \begin{equation}
            \begin{aligned}
            	&\sup\limits_{z\in\R} \left| \Prob\left[  \sqrt{n}\,\|L_\J(\Se-\St)\|_{(\Pt,\Gamma^\circ,s_1,s_2)} \leq z \right] - \Prob\left[ \max\limits_{j\in[p]} Y_j \leq z \right] \right| 
            	\leq \\
            	&\qquad\leq C_\cond 
            	\left\{
            	8^{3/(2\beta)} \left( \frac{( \log(pn^2))^7}{n} \right)^{1/8} 
            	+ c^2 \left( \frac{(\log(2pn^2))^{3+4/\beta}}{n} \right)^{1/2}
            	\right\}.
                \label{Eq:linear_final}
            \end{aligned}
        \end{equation}
    
    \subsubsection{From linear part to projectors} \label{S:nonlin}
    We have for all $z\in\R$ with $\delta = \sqrt{n} C\scc^2\rr_\J(\St)^{3/2}$ ($C$ from \eqref{LA_bound} of Lemma~\ref{L:KL_gen}) the following:
    \begin{equation}
        \begin{aligned}
            &\Prob\left[ \sqrt{n}\,\|\Pe-\Pt\|_{(\Pt,\Gamma^\circ,s_1,s_2)} \geq z\right]
            - \Prob\left[ \max\limits_{j\in[p]} Y_j \geq z \right] \\
            &\qquad\leq
            \Prob\left[ \sqrt{n}\,\|\Pe-\Pt\|_{(\Pt,\Gamma^\circ,s_1,s_2)} - \sqrt{n}\,\| L_\J(\Se-\St)\|_{(\Pt,\Gamma^\circ,s_1,s_2)} \geq \delta\right] + \\
            &\qquad\qquad+\left\{\Prob\left[ \sqrt{n}\,\| L_\J(\Se-\St)\|_{(\Pt,\Gamma^\circ,s_1,s_2)} \geq z-\delta\right]
            - \Prob\left[ \max\limits_{j\in[p]} Y_j \geq z-\delta \right] \right\} +\\
            &\qquad\qquad+
            \left\{ \Prob\left[ \max\limits_{j\in[p]} Y_j \geq z-\delta \right]
            - \Prob\left[ \max\limits_{j\in[p]} Y_j \geq z \right] \right\}.
        \nonumber
        \end{aligned}
    \end{equation}
    The first term is bounded by $2/n$ due to Lemma~\ref{L:KL_gen} (applied with $\Sp = \St$, $\Ppp = \Pt$) combined with Lemma~\ref{L:KC_SCconc} and Assumption~\ref{A:additional} (i).
    Next, the second term is bounded according to \eqref{Eq:linear_final}. For the third term we apply Lemma~\ref{L:CCK_AC} again and obtain
    \begin{equation}
        \begin{aligned}
            &\Prob\left[ \max\limits_{j\in[p]} Y_j \geq z-\delta \right]
            - \Prob\left[ \max\limits_{j\in[p]} Y_j \geq z \right] 
            = \Prob\left[ \max\limits_{j\in[p]} (Y_j/\clow) \geq z/\clow-\delta/\clow \right]
            - \Prob\left[ \max\limits_{j\in[p]} (Y_j/\clow) \geq z/\clow \right] 
            \\&\qquad\leq
            C(1, \cond^2) \frac{\delta}{\clow} \sqrt{1 \lor \log\left(\frac{p}{\delta/\clow}\right)}
            \leq C_\cond \frac{\delta}{\clow} \sqrt{\log\left(\frac{ep\clow}{\delta}\right)}.
        \nonumber
        \end{aligned}
    \end{equation}
    The opposite inequality for 
    \begin{equation}
        \begin{aligned}
            &\Prob\left[ \max\limits_{j\in[p]} Y_j \geq z \right] - \Prob\left[ \sqrt{n}\,\|\Pe-\Pt\|_{(\Pt, \Gamma^\circ, s_1, s_2)} \geq z\right]
        \nonumber
        \end{aligned}
    \end{equation}
    can be obtained in a similar way. Putting all the bounds together, we obtain the desired approximation.
    
    As to the spectral norm test statistic $\Qtsp$, since by Lemma~\ref{P:Properties1} (ii)
    \begin{equation}
        \begin{aligned}
    \| L_\J(\Se-\St) \| = \| L_\J(\Se-\St) \|_{(\Pt, \Gamma^\circ, m, d-m)},
        \nonumber
        \end{aligned}
    \end{equation}
    and trivially
    \begin{equation}
        \begin{aligned}
            &\left| \|\Pe-\Pt\| - \| L_\J(\Se-\St)\| \right| \leq \|R_\J(\Se-\St)\|,
        \nonumber
        \end{aligned}
    \end{equation}
    the same proof applies and yields the desired bound with $s_1 = m$ and $s_2 = d-m$.

%% file: source/proofs/Proof_B1.tex
Proof is quite similar to the proof of Theorem \ref{Th:GA1}, so we skip the technical details and focus only on the key parts that are different.

As in the proof of Theorem~\ref{Th:GA1}, we start with concentration written in specific form, but this time we are interested in concentration of $\Sb$ around $\Se$ conditionally on $\data$.
\begin{lemma} \label{L:BootConc}
    With probability $1-1/n$ it holds
    \begin{equation}
        \begin{aligned}
        	\Prob\left( \max\limits_{s,t\in[q]}\frac{\| \Pts (\Sb - \Se) \Ptt \|_{\Fr}}{\sqrt{m_s \mu_s m_t \mu_t}} \geq \bcc \,\Bigg|\, \data\right) \leq \frac{1}{n},
            \nonumber
        \end{aligned}
    \end{equation}
    where
    \begin{equation}
        \begin{aligned}
        	\bcc \eqdef C c^2 \frac{ \left( \log(n) + \log(2d^2) \right)^{\frac{2}{\beta} + \frac{1}{2}}}{\sqrt{n}}.
            \nonumber
        \end{aligned}
    \end{equation}
\end{lemma}
Due to Lemma~\ref{L:KL_gen} (this time applied with $\Sp = \Sb$, $\Ppp = \Pb$) combined with Lemma~\ref{L:BootConc} and Assumption~\ref{A:additional},
    \begin{equation}
        \begin{aligned}
    	\Prob\Big[ &\left| \sqrt{n}\| \Pb - \Pe \|_{(\Pt, \Gamma^\circ, s_1, s_2)} - \sqrt{n}\| L_\J(\Sb-\Se) \|_{(\Pt,\Gamma^\circ,s_1,s_2)} \right|
    	\geq
    	\\
    	&\qquad \geq
    	\sqrt{n} C(\scc+\bcc)^2 \rr_\J(\St)^{3/2} \,\Big|\, \data \Big] \leq \frac{1}{n}
        \label{LA_boot}
        \end{aligned}
    \end{equation}
with probability $1-1/n$.
Therefore, it again makes sense to work with the linear part $\sqrt{n}\,L_\J(\Sb-\Se)$ of $\sqrt{n}\,(\Pb-\Pe)$.
Introduce
\begin{equation}
        \begin{aligned}
        	L_{disc}^B &\eqdef \max\limits_{j\in[p]} \, v_j^{\top}\left( \sum\limits_{r\in\J} \sum\limits_{s\notin \J}\frac{ \Ptr(\Sb-\Se)\Pts + \Pts(\Sb-\Se)\Ptr}{\mu_r - \mu_s} \right) w_j
        	\\& = \max\limits_{j\in[p]} \,\sum\limits_{r\in\J} \sum\limits_{s\notin \J}\frac{ v_j^{\top}\Ptr(\Sb-\Se)\Pts w_j }{\mu_r - \mu_s}.
            \nonumber
        \end{aligned}
    \end{equation}
We again apply discretization step: by Lemma \ref{L:discretization},
    \begin{equation}
        \begin{aligned}
        	 L_{disc}^B \leq \| L_\J(\Sb-\Se)\|_{(\Pt,\Gamma^\circ,s_1,s_2)}\leq \frac{1}{1-2\eps} \, L_{disc}^B.
            \label{disc_bounds}
        \end{aligned}
    \end{equation}
Now let us represent $\sqrt{n} \, L_{disc}^B$ in a different way. Introduce for $i\in [n], j \in [p]$
    \begin{equation}
        \begin{aligned}
        	x_{ij}^B \eqdef 
        	\sum\limits_{r\in\J} \sum\limits_{s\notin \J}\frac{v_j^{\top} \Ptr(\eta_i X_i X_i^{\top}-X_i X_i^{\top})\Pts w_j}{\mu_r - \mu_s} 
        	=
        	(\eta_i - 1)\cdot x_{ij},
            \nonumber
        \end{aligned}
    \end{equation}
    where the random variables $x_{ij},\;i\in[n],\;j\in[p]$ are from the proof of Theorem \ref{Th:GA1}.
    Therefore,
    \begin{equation}
        \begin{aligned}
        	 \sqrt{n} L_{disc}^B
        	 =
        	 \max\limits_{j\in[p]} \frac{1}{\sqrt{n}}\sum\limits_{i=1}^n x_{ij}^B.
            \nonumber
        \end{aligned}
    \end{equation}
    Similarly to $x_i = \{ x_{ij} \}_{j=1}^p,\;i\in[n]$, we can arrange these random variables as i.i.d. $p$-dimensional centered random vectors 
    \begin{equation}
        \begin{aligned}x_i^B \eqdef \left\{ x_{ij}^B \right\}_{j=1}^p.
        \nonumber
        \end{aligned}
    \end{equation}
    Note that conditionally on the data $\data$ these vectors are automatically Gaussian with the covariance
    \begin{equation}
        \begin{aligned}
        	 \Cov(x_i^B\,|\,\data) = x_i x_i^{\top},\;\;\;i = 1,\ldots,n.
            \nonumber
        \end{aligned}
    \end{equation}
    Hence, conditionally on $\data$, the random vector
    \begin{equation}
        \begin{aligned}
        	 Y^B \eqdef \frac{1}{\sqrt{n}}\sum\limits_{i=1}^n x_{i}^B
            \nonumber
        \end{aligned}
    \end{equation}
    is Gaussian with covariance
    \begin{equation}
        \begin{aligned}
        	 \Cov(Y^B\,|\,\data) = \frac{1}{n}\sum\limits_{i=1}^n x_i x_i^{\top}.
            \nonumber
        \end{aligned}
    \end{equation}
    So, the Gaussian approximation step is not needed, unlike in the proof of Theorem \ref{Th:GA1}.
    To proceed, we need to show that $\Cov(Y^B\,|\,\data)$ is close to $\Cov(Y)$ with $Y$ from the proof of Theorem~\ref{Th:GA1}. Recall that
    \begin{equation}
            \begin{aligned}
            \Cov(Y) = \Cov(x_1) = \E\left[ \Cov(Y^B\,|\,\data) \right].
            \nonumber
            \end{aligned}
        \end{equation}
    The next lemma takes place. 
    \begin{lemma} \label{L:CovBound}
        With probability $1-1/n$ it holds
        \begin{equation}
            \begin{aligned}
            	 \left\| \Cov(Y^B\,|\,\data) - \Cov(Y) \right\|_{\mymax}\leq \clow^2\bcd,
                \nonumber
            \end{aligned}
        \end{equation}
        where
        \begin{equation}
            \begin{aligned}
            	\bcd \eqdef 
            	C_{\beta} \,c^4 \, \cond^2 \left( \sqrt{\frac{\log(pn)}{n}} + \frac{(\log(n)^{2/\beta} ( \log(pn))^{4/\beta}}{n}\right).
                \nonumber
            \end{aligned}
        \end{equation}
        Moreover, by Assumption~\ref{A:additionalB}, $\bcd \leq 1/2$.
    \end{lemma}

    \noindent Define $\Omega$ to be the event from Lemma \ref{L:CovBound}, $\Prob[\Omega] \geq 1-1/n$. On $\Omega$ for all $j\in [p]$
    \begin{equation}
        \begin{aligned}
        	 |\Var(Y_j^B) - \Var(Y_j)| = \left| \left[\Cov(Y^B\,|\,\data) - \Cov(Y) \right]_{j,j} \right|\leq \clow^2\bcd \leq \frac{\clow^2}{2}.
            \nonumber
        \end{aligned}
    \end{equation}
    So, since $\sqrt{n} L^B_{disc} = \max\limits_{j\in[p]} Y^B_j$ and 
    \begin{equation}
        \begin{aligned}
        	 &\Var(Y_j^B) \leq \Var(Y_j) + \frac{\clow^2}{2} \leq \chigh^2 + \frac{\clow^2}{2} \leq 2\chigh^2,\\
        	 &\Var(Y_j^B) \geq \Var(Y_j) - \frac{\clow^2}{2} \geq \clow^2 - \frac{\clow^2}{2} = \frac{\clow^2}{2},
            \nonumber
        \end{aligned}
    \end{equation}
    the same approach as in subsection~\ref{S:GA_AC}, applied conditionally on $\data$, together with bounds \eqref{disc_bounds} implies
    \begin{equation}
            \begin{aligned}
            	&\sup\limits_{z\in\R} \left| \Prob\left[  \sqrt{n}\,\|L_\J(\Sb-\Se)\|_{(\Pt,\Gamma^\circ,s_1,s_2)} \leq z\,|\,\data \right] - \Prob\left[ \max\limits_{j\in[p]} Y_j^B \leq z \right] \right| 
            	\leq \\
            	&\qquad\leq \frac{C_\cond \log(pn)}{n}
                \nonumber
            \end{aligned}
        \end{equation}
    on $\Omega$.
    
    Next, we deal with nonlinearity similarly to subsection~\ref{S:nonlin}, this time taking $\delta = \sqrt{n} C(\scc+\bcc)^2 \rr_\J(\St)^{3/2}$ and applying Lemma~\ref{L:KL_gen} with $\Sp = \Sb$, $\widetilde{\psi}_n=\bcc$. Omitting the details, we obtain
    \begin{equation}
            \begin{aligned}
            	&\sup\limits_{z\in\R} \left| \Prob\left[  \sqrt{n}\,\|\Pb-\Pe\|_{(\Pt, \Gamma^\circ,s_1,s_2)} \leq z\,\big|\,\data \right] - \Prob\left[ \max\limits_{j\in[p]} Y_j^B \leq z  \,\big|\,\data \right] \right| 
            	\leq \\
            	&\qquad\leq
            	C_\cond \left( \frac{\log(pn)}{n} + \frac{\delta}{\clow}\sqrt{\log\left(\frac{ep}{\delta/\clow}\right)} \right)
            	+\frac{1}{n}
                \nonumber
            \end{aligned}
        \end{equation}
    on $\Omega$.
    
    The only thing left is to compare the distributions of
    \begin{equation}
            \begin{aligned}
            \left(\max\limits_{j\in[p]} Y_j^B \,\big|\,\data\right) \;\;\text{and}\;\;\; \max\limits_{j\in[p]} Y_j.
            \nonumber
            \end{aligned}
        \end{equation}
    On $\Omega$ we showed
    \begin{equation}
        \begin{aligned}
        	 \max\limits_{j,k\in[p]} \left|[\Cov(Y^B/\clow\,|\,\data) - \Cov(Y/\clow)]_{jk}\right|
        	  \leq \bcd.
            \nonumber
        \end{aligned}
    \end{equation}
    Applying Lemma \ref{L:CCK_GC} results in
    \begin{equation}
            \begin{aligned}
                &\sup\limits_{z\in\R} \left| \Prob\left[  \max\limits_{j\in[p]} Y_j^B \leq z\,|\,\data \right] - \Prob\left[ \max\limits_{j\in[p]} Y_j \leq z \right] \right| 
            	\leq \\
            	&\qquad\leq
            	\sup\limits_{z\in\R} \left| \Prob\left[  \max\limits_{j\in[p]} (Y_j^B/\clow) \leq z\,|\,\data \right] - \Prob\left[ \max\limits_{j\in[p]} (Y_j/\clow) \leq z \right] \right| 
            	\leq \\
            	&\qquad\leq
                C_\cond \bcd^{1/3} \left(\log(ep/\bcd)\right)^{2/3}
                \nonumber
            \end{aligned}
        \end{equation}
    on $\Omega$.
    Putting all the bounds together with Theorem~\ref{Th:GA1}, we obtain the desired.

%% file: source/proofs/Proof_F1.tex
The proof is quite similar to Theorem~\ref{Th:B1}, so we follow the same steps.
First we formulate the concentration result, this time for $\Sf$.
\begin{lemma} \label{L:BayesConc}
    With probability $1-1/n$ it holds
    \begin{equation}
        \begin{aligned}
        	\Prob\left( \max\limits_{s,t\in[g]}\frac{\| \Pts (\Sf - \Se) \Ptt \|_{\Fr}}{\sqrt{m_s \mu_s m_t \mu_t}} \geq \fcc \,\Bigg|\, \data\right) \leq \frac{1}{n},
            \nonumber
        \end{aligned}
    \end{equation}
    where
    \begin{equation}
        \begin{aligned}
        	\fcc \eqdef C c^2 \frac{ (\log(n) + \log(2d^2))^{\frac{2}{\beta}+\frac{1}{2}} }{\sqrt{n}}.
            \nonumber
        \end{aligned}
    \end{equation}
\end{lemma}

Due to Lemma~\ref{L:KL_gen} (applied with $\Sp = \Sf$, $\Ppp = \Pf$) combined with Lemma~\ref{L:BayesConc} and Assumption~\ref{A:additional},
    \begin{equation}
        \begin{aligned}
    	\Prob\Big[ &\left| \sqrt{n}\| \Pf - \Pe \|_{(\Pt, \Gamma^\circ,s_1,s_2)} - \sqrt{n}\| L_\J(\Sf-\Se) \|_{(\Pt,\Gamma^\circ,s_1,s_2)} \right|
    	\geq
    	\\
    	&\qquad \geq
    	\sqrt{n} C(\scc+\fcc)^2 \rr_\J(\St)^{3/2} \,\Big|\, \data \Big] \leq \frac{1}{n}
        \label{LA_boot}
        \end{aligned}
    \end{equation}
with probability $1-1/n$.
We again elaborate on the linear term $L_\J(\Sf-\St)$.
Define its discretized version
\begin{equation}
        \begin{aligned}
        	L_{disc}^F &\eqdef \max\limits_{j\in[p]} \, v_j^{\top}\left(\sum\limits_{r\in\J} \sum\limits_{s\notin \J}\frac{ \Ptr(\Sf-\Se)\Pts + \Pts(\Sf-\Se)\Ptr}{\mu_r - \mu_s} \right)w_j
        	\\&= \max\limits_{j\in[p]} \, \sum\limits_{r\in\J} \sum\limits_{s\notin \J}\frac{v_j^{\top} \Ptr(\Sf-\Se)\Pts w_j}{\mu_r - \mu_s},
            \nonumber
        \end{aligned}
    \end{equation}
and by Lemma \ref{L:discretization} obtain the bounds
    \begin{equation}
        \begin{aligned}
        	 L_{disc}^F \leq \| L_\J(\Sf-\Se)\|_{(\Pt,\Gamma^\circ,s_1,s_2)}\leq \frac{1}{1-2\eps}\, L_{disc}^F.
            \nonumber
        \end{aligned}
    \end{equation}
Introducing for $i\in [n], j \in [p]$
    \begin{equation}
        \begin{aligned}
        	x_{ij}^F \eqdef 
        	\sum\limits_{r\in\J} \sum\limits_{s\notin \J}\frac{v_j^{\top} \Ptr( Z_i Z_i^{\top}-\Se)\Pts w_j}{\mu_r - \mu_s},
            \nonumber
        \end{aligned}
    \end{equation}
we represent
    \begin{equation}
        \begin{aligned}
        	 \sqrt{n} L_{disc}^F
        	 =
        	 \max\limits_{j\in[p]} \frac{1}{\sqrt{n}}\sum\limits_{i=1}^n x_{ij}^F.
            \nonumber
        \end{aligned}
    \end{equation}
    The $p$-dimensional centered random vectors 
    \begin{equation}
        \begin{aligned}
            x_i^F \eqdef \left\{ x_{ij}^F \right\}_{j=1}^p
        \nonumber
        \end{aligned}
    \end{equation}
    are i.i.d., and we need to compute their covariance conditionally on $\data$.
    For any fixed indices $i\in[n]$, $j,k\in[p]$, $r,r^{\prime}\in\J$, $s,s^{\prime}\notin\J$, similarly to Subsection~\ref{S:verify},  ``Computing the covariance'' part, we have
    \begin{equation}
        \begin{aligned}
        	&\E\left[ v_j^{\top} \Ptr Z_i Z_i^{\top} \Pts w_j \cdot v_k^{\top} \Ptrp Z_i Z_i^{\top} \Ptsp w_k \,|\,\data \right] = \\
        	&\qquad = (v_j^{\top} \Ptr \Se \Ptrp v_k) \cdot (w_j^{\top} \Pts \Se \Ptsp w_k),
            \nonumber
        \end{aligned}
    \end{equation}
    where we take into account that for Gaussian $Z$ it holds that $\Pt Z$ and $\Ptc Z$ are independent (which implies Assumption~\ref{A: independence} for $Z$). Thus,
    \begin{equation}
        \begin{aligned}
        	&\E\left[ v_j^{\top} \Ptr (Z_i Z_i^{\top} - \Se) \Pts w_j \cdot v_k^{\top} \Ptrp (Z_i Z_i^{\top} - \Se) \Ptsp w_k \,|\,\data \right] = \\
        	&\qquad = (v_j^{\top} \Ptr \Se \Ptrp v_k) \cdot (w_j^{\top} \Pts \Se \Ptsp w_k) - (v_j^{\top} \Ptr \Se \Pts w_j) \cdot (v_k^{\top} \Ptrp \Se \Ptsp w_k),
            \nonumber
        \end{aligned}
    \end{equation}
    and 
    \begin{equation}
        \begin{aligned}
        	&\left[\Cov(x^F_i\,|\,\data) \right]_{j,k} = \E\left[ x^F_{ij} x^F_{ik} \,|\,\data \right] = \\
        	&\qquad = \sum\limits_{\substack{r\in\J\\r^{\prime} \in \J}} \sum\limits_{\substack{s\notin\J\\s^{\prime} \notin \J}} \frac{(v_j^{\top} \Ptr \Se \Ptrp v_k) \cdot (w_j^{\top} \Pts \Se \Ptsp w_k) - (v_j^{\top} \Ptr \Se \Pts w_j) \cdot (v_k^{\top} \Ptrp \Se \Ptsp w_k)}{(\mu_r-\mu_s)(\mu_{r^{\prime}} - \mu_{s^{\prime}})}.
            \nonumber
        \end{aligned}
    \end{equation}
    We again have to use Lemma~\ref{L:CCK_GA} to pass from $x^F_i$'s to their Gaussian counterparts $y^F_i$'s with the same covariance. Introduce
    \begin{equation}
        \begin{aligned}
          Y^F \eqdef \frac{1}{\sqrt{n}} \sum\limits_{i=1}^n y_i^F \sim \mathcal{N}_p\left(0, \Cov(x_1^F\,|\,\data) \right)
            \nonumber
        \end{aligned}
    \end{equation}
    to approximate the distribution of interest by $\max\limits_{j\in[p]} Y^F_j$.
    Recall that in the proof of Theorem~\ref{Th:GA1} we had for $Y$ and $x_i$'s
    \begin{equation}
        \begin{aligned}
          \left[\Cov(Y) \right]_{j,k} = \left[\Cov(x_1) \right]_{j,k} =  \sum\limits_{r\in\J}\sum\limits_{s\notin\J} \frac{\mu_r \mu_s}{(\mu_r - \mu_s)^2}\, \, (v_j^{\top} \Ptr v_k) \cdot (w_j^{\top} \Pts w_k),
            \nonumber
        \end{aligned}
    \end{equation}
    which alternatively can be expressed as
    \begin{equation}
        \begin{aligned}
        	&\left[\Cov(Y) \right]_{j,k} =  \left[\Cov(x_1) \right]_{j,k} = \\
        	&\qquad = \sum\limits_{\substack{r\in\J\\r^{\prime} \in \J}} \sum\limits_{\substack{s\notin\J\\s^{\prime} \notin \J}} \frac{(v_j^{\top} \Ptr \St \Ptrp v_k) \cdot (w_j^{\top} \Pts \St \Ptsp w_k) - (v_j^{\top} \Ptr \St \Pts w_j) \cdot (v_k^{\top} \Ptrp \St \Ptsp w_k)}{(\mu_r-\mu_s)(\mu_{r^{\prime}} - \mu_{s^{\prime}})}.
            \nonumber
        \end{aligned}
    \end{equation}
    Observe that the difference between this expression and the expression for \begin{equation}
        \begin{aligned}
        \left[\Cov(Y^F\,|\,\data) \right]_{j,k} = \left[\Cov(x^F_1\,|\,\data) \right]_{j,k}
        \nonumber
        \end{aligned}
    \end{equation}
    above is that the true covariance $\St$ replaces the empirical one $\Se$. In the next lemma we bound the maximal element-wise absolute difference between this two covariances.
\begin{lemma} \label{L:FCovBound}
    With probability $1-1/n$ it holds
    \begin{equation}
        \begin{aligned}
        	\left\| \Cov(Y^F\,|\,\data) - \Cov(Y) \right\|_{\mymax} \leq \clow^2 \fcd,
            \nonumber
        \end{aligned}
    \end{equation}
    where
    \begin{equation}
        \begin{aligned}
        	\fcd \eqdef
        	|\J|\, C_\beta c^2 \cond^2 \left( \sqrt{\frac{\log(np) + \log(|\J|)}{n}} + \frac{(\log(n))^{1/\beta} (\log(np) + \log(|\J|))^{2/\beta}}{n} \right).
            \label{Def:fcd}
        \end{aligned}
    \end{equation}
    Moreover, by Assumption~\ref{A:additionalF}, $\fcd \leq 1/2$.
\end{lemma}
\noindent Let $\Omega$ be the event from Lemma~\ref{L:FCovBound}, $\Prob[\Omega] \geq 1-1/n$. On this event
\begin{equation}
        \begin{aligned}
        	 &\Var(Y_j^F) \leq \Var(Y_j) + \frac{\clow^2}{2} \leq \chigh^2 + \frac{\clow^2}{2} \leq 2\chigh^2,\\
        	 &\Var(Y_j^F) \geq \Var(Y_j) - \frac{\clow^2}{2} \geq \clow^2 - \frac{\clow^2}{2} = \frac{\clow^2}{2},
            \nonumber
        \end{aligned}
    \end{equation}
and all the arguments from the proof of Theorem~\ref{Th:GA1} work with slightly shifted variances and $\scc^2$ replaced by $(\scc + \fcc)^2$. Nonlinearity is treated similarly to subsection~\ref{S:nonlin}, this time taking $\delta = \sqrt{n} C(\scc+\fcc)^2 \rr_\J(\St)^{3/2}$ and applying Lemma~\ref{L:KL_gen} with $\Sp = \Sf$, $\widetilde{\psi}_n=\fcc$. So,
        \begin{equation}
            \begin{aligned}
            	&\sup\limits_{z\in\R} \left| \Prob\left[ \sqrt{n}\| \Pf-\Pe \|_{(\Pt, \Gamma^\circ,s_1,s_2)} \leq z \,|\,\data\right] - \Prob\left[ \max\limits_{j\in[p]} Y_j^F \leq z\,|\,\data \right] \right| \leq
            	\\&\qquad\leq
            	C_\cond \left\{\Diamond^{GA} + \zeta\left[\sqrt{n}(\scc+\fcc)^2\rr_\J(\St)^{3/2}/\clow\right] \right\}
            \nonumber
            \end{aligned}
        \end{equation}
on $\Omega$, in addition to already established result
        \begin{equation}
            \begin{aligned}
            	&\sup\limits_{z\in\R} \left| \Prob\left[ \sqrt{n}\| \Pe-\Pt \|_{(\Pt,\Gamma^\circ,s_1,s_2)} \leq z \right] - \Prob\left[ \max\limits_{j\in[p]} Y_j \leq z \right] \right| \leq \\&\qquad\leq C_\cond \left\{\Diamond^{GA} + \zeta\left[\sqrt{n}\scc^2\rr_\J(\St)^{3/2}/\clow\right] \right\}.
            \nonumber
            \end{aligned}
        \end{equation}
The only thing left is to apply Gaussian comparison Lemma~\ref{L:CCK_GC} to
\begin{equation}
    \begin{aligned}
    \max\limits_{j\in[p]} (Y_j/\clow) \;\;\text{ and }\;\;(\max\limits_{j\in[p]} (Y_j^F/\clow)\,|\,\data)
    \nonumber
    \end{aligned}
\end{equation}
with $\Delta = \fcd$ from Lemma~\ref{L:FCovBound}.

%% file: source/proofs/Proof_GA2.tex
We start with the following lemma, which is modification of Lemma~\ref{L:KL_gen}. Not only does it allow to get rid of the remainder term, but at the same time replaces $\| \cdot \|_{(\Pbar, \overline{\Gamma},s_1,s_2)}$ by some \\$\|\cdot \|_{(\Pstar, \Gamma^*,s_1,s_2)}$ with deterministic $\Pstar$. We provide a simplified version to be used in this proof; a more general one, suitable for proofs of Theorem~\ref{Th:B2} and Theorem~\ref{Th:F2}, can be established similarly to Lemma~\ref{L:KL_gen}. 
\begin{lemma} \label{L:2approx}
    Let under $H_0^{(2)}$  be $\Pstar \eqdef \Pta = \Ptb$.
    For any $\overline{\Gamma}$, there exists $\Gamma^* = [\Gamma^*_1 \;\Gamma^*_2] \in \R^{d\times d}$ with $\Gamma^*_1\in\R^{d\times m} , \Gamma^*_2 \in \R^{d\times(d-m)}$ satisfying
\begin{equation} 
    \begin{aligned}
        &\Gamma^*_1 {\Gamma^*_1}^\T = \Pstar, \;\;{\Gamma^*_1}^\T \Gamma^*_1 = \Id_{m},\\
        &\Gamma^*_2 {\Gamma^*_2}^\T = \Id_d - \Pstar, \;\; {\Gamma^*_2}^\T \Gamma_2^* = \Id_{d-m},
    \nonumber 
    \end{aligned} 
\end{equation}
such that the following holds:
    \begin{equation}
        \begin{aligned}
        	&\Prob\left[\left| \| \Pea-\Peb\|_{(\Pbar, \overline{\Gamma},s_1,s_2)} - \| L_a(\Se_a-\St_a) - L_b(\Se_b-\St_b)\|_{(\Pstar,\Gamma^*,s_1,s_2)} \right| \leq
        	C \sccab^2 \left(\rre_{a,b} + \rr_{a,b}^{3/2}\right) \,\Big|\,\overline{\Gamma} \right] \geq
        	\\&\qquad \geq 1-\frac{1}{n_a}-\frac{1}{n_b}
        \nonumber
        \end{aligned}
    \end{equation}
    with probability $1-1/n_a - 1/n_b$.
\end{lemma}

Now we start elaborating on $ \| L_a(\Se_a-\St_a) - L_b(\Se_b-\St_b)\|_{(\Pstar,\Gamma^*,s_1,s_2)}$ in a similar fashion as the proof of Theorem~\ref{Th:GA1}. The only difference is that in the rest of the proof, all the probabilities,  expectations and variances are conditional on $\overline{\Gamma}$.
Let this time $\{(v_j, w_j)_{j=1}^p\}$ enumerate all pairs $(\Gamma_1^* v, \Gamma_2^*w)$ for $v\in N_{\eps}(\myset^m_{s_1})$, $w\in N_{\eps}(\myset^{d-m}_{s_2})$. We again take $\eps = 1/n$, which fixes $p$ to be
    \begin{equation}
        \begin{aligned}
        	p \leq \exp\left( (s_1+s_2)\log(3n) + \log(2d) \right).
        \nonumber
        \end{aligned}
    \end{equation}
Note that both $L_a(\Se_a-\St_a)$ and $L_b(\Se_b-\St_b)$ satisfy
\begin{equation}
        \begin{aligned}
        	&L_a(\Se_a-\St_a) = \Pstar L_a(\Se_a-\St_a) (\Id_d - \Pstar) + (\Id_d - \Pstar) L_a(\Se_a-\St_a) \Pstar,\\
        	&L_b(\Se_b-\St_b) = \Pstar L_b(\Se_b-\St_b) (\Id_d - \Pstar) + (\Id_d - \Pstar) L_b(\Se_b-\St_b) \Pstar,
        \nonumber
        \end{aligned}
    \end{equation}
so Lemma~\ref{L:discretization} yield
\begin{equation}
        \begin{aligned}
        	L_{disc}^{a,b} \leq \| L_a(\Se_a-\St_a) - L_b(\Se_b-\St_b)\|_{(\Pstar,\Gamma^*,s_1,s_2)} \leq \frac{1}{1-2\eps} L_{disc}^{a,b},
        \label{bounds_ab}
        \end{aligned}
    \end{equation}
where $L_{disc}^{a,b} = \max\limits_{j\in[p]} \;v_j^\T (L_a(\Se_a-\St_a) - L_b(\Se_b-\St_b)) w_j$. The same quantity can be expressed as a sum.
For all $j\in[p]$ introduce $x_{ij}^a$ for $i\in[n_a]$ and $x_{ij}^b$ for $i\in[n_b]$, which are analogs of $x_{ij}$. Then
\begin{equation}
        \begin{aligned}
        	L_{disc}^{a,b}  = \max\limits_{j\in[p]} \left( \frac{1}{n_a}\sum\limits_{i=1}^{n_a} x_{ij}^a - \frac{1}{n_b}\sum\limits_{i=1}^{n_b} x_{ij}^b\right).
        \nonumber
        \end{aligned}
    \end{equation}
Gaussian counterpart of $\sqrt{n_an_b/(n_a+n_b)}\, L_{disc}^{a,b}$ is given by $\max\limits_{j\in[p]} Y^{a,b}_j$, where
\begin{equation}
        \begin{aligned}
        	Y^{a,b}  = \sqrt{\frac{n_a n_b}{n_a+n_b}} \left( \frac{1}{n_a}\sum\limits_{i=1}^{n_a} y_{i}^a - \frac{1}{n_b}\sum\limits_{i=1}^{n_b} y_{i}^b\right)
        \nonumber
        \end{aligned}
    \end{equation}
with $y_i^a\,|\,\overline{\Gamma} \sim \mathcal{N}\left(0, \Cov(x_i^a\,|\,\overline{\Gamma})\right)$ for all $i\in[n_a]$ and $y_i^b\,|\,\overline{\Gamma} \sim \mathcal{N}\left(0, \Cov(x_i^b\,|\,\overline{\Gamma})\right)$ for all $i\in[n_b]$.

The conditions of Lemma~\ref{L:CCK_GA}, verified in Subsection~\ref{S:verify}, can be treated here likewise.
Note that similar to $\Var(x_{ij})$, we can lower and upper bound
\begin{equation}
        \begin{aligned}
        	&\clow_{\J_a}(\St_a)^2 \leq \Var(x_{ij}^a\,|\,\overline{\Gamma}) \leq \chigh_{\J_a}(\St_a)^2,\\
        	&\clow_{\J_b}(\St_b)^2 \leq \Var(x_{ij}^b\,|\,\overline{\Gamma}) \leq \chigh_{\J_b}(\St_b)^2.
        \nonumber
        \end{aligned}
    \end{equation}
Furthermore, direct computation shows
\begin{equation}
        \begin{aligned}
        	 \Var(Y_j^{a,b}\,|\,\overline{\Gamma}) = \frac{n_b\Var(x_{1j}^a\,|\,\overline{\Gamma}) + n_a\Var(x_{1j}^b\,|\,\overline{\Gamma})}{n_a+n_b},
        \nonumber
        \end{aligned}
    \end{equation}
implying
\begin{equation}
        \begin{aligned}
        	&\clow_{\J_a}(\St_a)^2 \land \clow_{\J_b}(\St_b)^2 \leq \Var(Y_j^{a,b}\,|\,\overline{\Gamma}) \leq \chigh_{\J_a}(\St_a)^2 \lor \chigh_{\J_b}(\St_b)^2.
        \nonumber
        \end{aligned}
    \end{equation}
So, the existence of $c_1, C_1 > 0$ lower- and upperbounding the variance is established.
Upper bound on $M_3, M_4$ can be obtained as well, and it will be $8^{2/\beta} \cdot\left( \chigh_{\J_a}(\St_a) \lor \chigh_{\J_b}(\St_b) \right)$ instead of $8^{2/\beta} \chigh_{\J}(\St)$. Upper bound on $u(\gamma)$ again follows likewise and becomes (for $\gamma = 1/n_a + 1/n_b$)
\begin{equation}
        \begin{aligned}
        	u(\gamma) \lesssim \condab \,c^2 \left(\log\left(2p(n_a+n_b)^2\right)\right)^{2/\beta}.
        \nonumber
        \end{aligned}
    \end{equation}
Lemma~\ref{L:CCK_GA} then yields almost surely
        \begin{equation}
            \begin{aligned}
            	&\sup\limits_{z\in\R} \left| \Prob\left[  \sqrt{\frac{n_an_b}{n_a+n_b}}\,L_{disc}^{a,b} \leq z\,\Big|\,\overline{\Gamma} \right] - \Prob\left[ \max\limits_{j\in[p]} Y_j^{a,b} \leq z\,\Big|\,\overline{\Gamma} \right] \right| 
            	\leq \\
            	&\qquad\leq C_{\condab}
            	\Bigg\{
            	8^{3/(2\beta)} \left( \frac{( \log(p(n_a+n_b)^2))^7}{n_a+n_b} \right)^{1/8} 
            	+ c^2 \left( \frac{(\log(2p(n_a+n_b)^2))^{3+4/\beta}}{n_a+n_b} \right)^{1/2} +\\
            	&\qquad\qquad\qquad
            	+ \frac{1}{n_a+n_b}
            	\Bigg\},
                \nonumber
            \end{aligned}
        \end{equation}
and by similar to Subsection~\ref{S:GA_AC} reasoning, bounds \eqref{bounds_ab} allow to pass from discretized version to infinite-state supremum, omitting the negligible additional error term:
\begin{equation}
            \begin{aligned}
            	&\sup\limits_{z\in\R} \left| \Prob\left[  \sqrt{\frac{n_an_b}{n_a+n_b}}\,\|L_a(\Se_a-\St_a) - L_a(\Se_b-\St_b)\|_{(\Pstar, \Gamma^*)} \leq z \,\Big|\,\overline{\Gamma}\right] - \Prob\left[ \max\limits_{j\in[p]} Y_j^{a,b} \leq z \,\Big|\,\overline{\Gamma}\right] \right| 
            	\leq \\
            	&\qquad\leq C_{\condab}
            	\Bigg\{
            	8^{3/(2\beta)} \left( \frac{( \log(p(n_a+n_b)^2))^7}{n_a+n_b} \right)^{1/8} 
            	+ c^2 \left( \frac{(\log(2p(n_a+n_b)^2))^{3+4/\beta}}{n_a+n_b} \right)^{1/2} +
            	\\&\qquad\qquad\qquad
            	+ \frac{1}{n_a+n_b}
            	\Bigg\}.
                \label{L_disc_ab_bound}
            \end{aligned}
        \end{equation}

The last step is to use Lemma~\ref{L:2approx} to finalize the result for projectors.
As in Subsection~\ref{S:nonlin}, we have with $\delta = \sqrt{\frac{n_an_b}{n_a+n_b}}C\sccab^2 \left(\rre_{a,b} + \rr_{a,b}^{3/2}\right)$
\begin{equation}
        \begin{aligned}
            &\Prob\left[ \sqrt{\frac{n_an_b}{n_a+n_b}}\,\|\Pea-\Peb\|_{(\Pbar, \overline{\Gamma},s_1,s_2)} \geq z\,\Big|\,\overline{\Gamma}\right]
            - \Prob\left[ \max\limits_{j\in[p]} Y_j^{a,b} \geq z \,\Big|\,\overline{\Gamma}\right] \\
            &\qquad\leq
            \Prob\left[ \sqrt{\frac{n_an_b}{n_a+n_b}}\,\left( \|\Pea-\Peb\|_{(\Pbar, \overline{\Gamma}, s_1,s_2)} - \| L_a(\Se_a-\St_a) - L_b(\Se_b-\St_b)\|_{(\Pstar,\Gamma^*,s_1,s_2)}  \right) \geq \delta\,\Big|\,\overline{\Gamma}\right] + \\
            &\qquad\qquad+\Big\{\Prob\left[ \sqrt{\frac{n_an_b}{n_a+n_b}}\,\| L_a(\Se_a-\St_a) - L_b(\Se_b-\St_b)\|_{(\Pstar,\Gamma^*,s_1,s_2)} \geq z-\delta\,\Big|\,\overline{\Gamma}\right]
            -
            \\&\qquad\qquad\qquad- \Prob\left[ \max\limits_{j\in[p]} Y_j^{a,b} \geq z-\delta \,\Big|\,\overline{\Gamma}\right] \Big\} +\\
            &\qquad\qquad+
            \left\{ \Prob\left[ \max\limits_{j\in[p]} Y_j^{a,b} \geq z-\delta \,\Big|\,\overline{\Gamma}\right]
            - \Prob\left[ \max\limits_{j\in[p]} Y_j^{a,b} \geq z \,\Big|\,\overline{\Gamma}\right] \right\}.
        \nonumber
        \end{aligned}
    \end{equation}
By Lemma~\ref{L:2approx}, the first term is at most $1/n_a + 1/n_b$ with probability $1-1/n_a-1/n_b$.
The second term is bounded by \eqref{L_disc_ab_bound}. For the third term we apply  Lemma~\ref{L:CCK_AC} and get
\begin{equation}
        \begin{aligned}
            \Prob\left[ \max\limits_{j\in[p]} Y_j^{a,b} \geq z-\delta \,\Big|\,\overline{\Gamma}\right]
            - \Prob\left[ \max\limits_{j\in[p]} Y_j^{a,b} \geq z \,\Big|\,\overline{\Gamma}\right]
            \leq
            C_{\condab} \frac{\delta}{\clow_{\J_a}(\St_b) \land \clow_{\J_b}(\St_a)} \left( \log\left(\frac{ep}{\delta}\right)\right)^{1/2}.
        \nonumber
        \end{aligned}
    \end{equation}
The opposite inequality can be obtained similarly. This concludes the proof.

%% file: source/proofs/Proof_BF2.tex
The proofs repeat proofs of Theorem~\ref{Th:B1} and Theorem~\ref{Th:F1}, respectively, with Theorem~\ref{Th:GA2} used in place of Theorem~\ref{Th:GA1}.

%% file: source/proofs/Proof_C12.tex
To prove Corollary~\ref{Corollary1}, it is enough to repeat the proof of Corollary 2.3 of \cite{Silin_1}, applied with our Theorem~\ref{Th:B1} and Theorem~\ref{Th:F1}.  Corollary~\ref{Corollary2} is slightly trickier, since we condition on $\overline{\Gamma}$. However, still the proof of Corollary 2.3 of \cite{Silin_1}, applied with, for example, Theorem~\ref{Th:B2}, yields
\begin{equation}
            \begin{aligned}
                \sup\limits_{\alpha\in(0; 1)} \left| \Prob\left[\Qtab > \gamma_\alpha^B \,\big|\, \overline{\Gamma}\right] - \alpha\right| \leq \Diamond_B + \frac{1}{n_a} + \frac{1}{n_b}
            \nonumber
            \end{aligned}
        \end{equation}
with probability $1-1/n_a-1/n_b$. Integrating $\overline{\Gamma}$ out, we obtain the desired.

%% file: source/proofs/Proof_Power1.tex
Let us demonstrate the proof of Theorem~\ref{Th:Power1} first.
The proofs for (i) and (ii) are identical, so let us only focus on (i).

Recall $\gamma_B^{(1)}(\alpha)$ from Corollary~\ref{Corollary1}. By triangle inequality and the assumption of the theorem,

\begin{equation}
    \begin{aligned}
        &\Prob\left[ \sqrt{n} \| \Pe - \Ph\|_{(\Ph, \Gamma^\circ, s_1, s_2)} \geq \gamma_B^{(1)}(\alpha) \right]
        \geq \\
        &\qquad \geq 
        \Prob\left[ \sqrt{n} \| \Pt - \Ph\|_{(\Ph, \Gamma^\circ, s_1, s_2)} - \sqrt{n} \| \Pe - \Pt\|_{(\Ph, \Gamma^\circ, s_1, s_2)} \geq \gamma_B^{(1)}(\alpha) \right]\\
        &\qquad \geq 
        \Prob\left[ \lambda_n - \sqrt{n} \| \Pe - \Pt\|_{(\Ph, \Gamma^\circ, s_1, s_2)} \geq \gamma_B^{(1)}(\alpha) \right]\\
        &\qquad =
        \Prob\left[  \sqrt{n} \| \Pe - \Pt\|_{(\Ph, \Gamma^\circ, s_1, s_2)} \leq  \lambda_n - \gamma_B^{(1)}(\alpha) \right].
    \nonumber
    \end{aligned}
\end{equation}
Proposition~\ref{P:Properties0} (ii) implies
\begin{equation}
    \begin{aligned}
        &\Prob\left[ \sqrt{n} \| \Pe - \Ph\|_{(\Ph, \Gamma^\circ, s_1, s_2)} \geq \gamma_B^{(1)}(\alpha) \right]\geq\\
        &\qquad \geq
        \Prob\left[  \sqrt{n} \| \Pe - \Pt\| \leq  \lambda_n - \widetilde{\gamma}_B^{(1)}(\alpha) \right],
    \nonumber
    \end{aligned}
\end{equation}
where $\widetilde{\gamma}_B^{(1)}(\alpha)$ is $\alpha$-quantile of $\sqrt{n}\| \Pb - \Pe \|$.
In the proof of Lemma~\ref{L:2approx} we will show that with probability $1-1/n$ we have the bound
\begin{equation}
    \begin{aligned}
       \sqrt{n}\| \Pe - \Pt \| \leq \sqrt{n} C \left( \scc \rre^{1/2} + \scc^2\rr^{3/2}\right),
    \nonumber
    \end{aligned}
\end{equation}
and similarly 
\begin{equation}
    \begin{aligned}
       \sqrt{n}\| \Pb - \Pe \| \leq \sqrt{n} C \left( \bcc \rre^{1/2} + \bcc^2\rr^{3/2}\right),
    \nonumber
    \end{aligned}
\end{equation}
which means that $\widetilde{\gamma}_B^{(1)}(\alpha)$ is at most of the same order treating $\alpha$ as constant.
Denoting for shortness $\Phi_n = \sqrt{n} C\left((\scc+\bcc)\rr^{1/2} + (\scc+\bcc)^2\rre^{3/2} \right)$, if $\CONST \geq 2C$,  this ensures that 
\begin{equation}
    \begin{aligned}
        \Prob\left[  \sqrt{n} \| \Pe - \Pt\| \leq  \lambda_n - \widetilde{\gamma}_B^{(1)}(\alpha) \right] =
        \Prob\left[  \sqrt{n} \| \Pe - \Pt\| \leq  \Phi_n \left( \frac{\lambda_n}{\Phi_n} - \frac{\widetilde{\gamma}_B^{(1)}(\alpha)}{\Phi_n} \right) \right]  \to 1
    \nonumber
    \end{aligned}
\end{equation}
as $n\to\infty$, since $\liminf\limits_{n\to\infty} \lambda_n/\Phi_n \geq \CONST/C \geq 2$ by condition \eqref{Cond: lambda} and $\widetilde{\gamma}_B^{(1)}(\alpha)/\Phi_n \leq 1$. This concludes the proof.

The proof of Theorem~\ref{Th:Power2} repeats the proof above, with the only difference that we will also need to apply the inequality
\begin{equation}
    \begin{aligned}
        \| \Pta - \Ptb \|_{(\Pbar,\Gamma,s_1,s_2)}
        \geq
        \frac{1}{2}\sqrt{\frac{s_1s_2}{m(d-m)}} \| \Pta-\Ptb\| \geq \lambda_{n_a,n_b} \sqrt{\frac{n_a+n_b}{n_an_b}}.
    \nonumber
    \end{aligned}
\end{equation}

%% file: source/AuxLiterature.tex
\subsection{Results from \cite{Weibull}}

\begin{proposition}[\cite{Weibull}, Proposition S.3.2] \label{L:KC_Prod}
    If $W_i$, $i \in [k]$ are (possibly dependent) random variables satisfying $\| W_i \|_{\psi_{\alpha_i}} < \infty$ for some $\alpha_i > 0$, then
    \begin{equation}
        \begin{aligned}
	    \left\| \prod\limits_{i=1}^k W_i \right \|_{\psi_{\beta}} \leq \prod\limits_{i=1}^k \| W_i\|_{\psi_{\alpha_i}}
	    \;\;\;
	    \text{where }\frac{1}{\beta} \eqdef \sum\limits_{i=1}^k \frac{1}{\alpha_i}.
        \nonumber
        \end{aligned}
    \end{equation}
\end{proposition}

\begin{theorem}[\cite{Weibull}, Theorem 4.1] \label{Th:KC_SCmax}
    Let $X_1, \ldots, X_n$ be independent random vectors in $\R^p$ satisfying
    \begin{equation}
        \begin{aligned}
	    \max\limits_{i\in[n], j\in[p]} \| X_{i}(j) \|_{\psi_\beta} \leq K_{n,p} < \infty \;\;\text{ for some }\; 0 < \beta \leq 2.
        \nonumber
        \end{aligned}
    \end{equation}
    Fix $n, p \geq 1$. Then for any $t \geq 0$, with probability at least $1-3e^{-z}$,
    \begin{equation}
        \begin{aligned}
	    &\max\limits_{j,k\in[p]}
	    \left| \frac{1}{n}\sum\limits_{i=1}^n X_{i}(j) X_{i}(k) - \E[X_{i}(j) X_{i}(k)] \right| \leq \\
	    &\qquad \leq
	    7 A_{n,p} \sqrt{\frac{z+2\log(p)}{n}} + \frac{C_{\beta} K_{n,p}^2(\log(2n))^{1/\beta} (z+2\log(p))^{2/\beta}}{n},
        \nonumber
        \end{aligned}
    \end{equation}
    where $C_\beta>0$ is a constant depending only on $\beta$, and $A_{n,p}^2$ is given by
    \begin{equation}
        \begin{aligned}
	    A_{n,p} \eqdef \max\limits_{j,k\in[p]} \frac{1}{n} \sum\limits_{i=1}^n \Var(X_i(j) X_i(k)).
        \nonumber
        \end{aligned}
    \end{equation}
\end{theorem}
\begin{remark}
    Remark 4.1 from \cite{Weibull} claims $A_{n,p} \leq C_\beta K_{n,p}^2$, so in every application of Theorem~\ref{Th:KC_SCmax} we use this small fact without further notice.
\end{remark}

\subsection{Results from \cite{Wahl}}

Recent paper \cite{Wahl} considers infinite-dimensional Hilbert space $\mathcal{H}$ and two covariance operators $\St$, $\Se$ with perturbation $\Pert \eqdef \Se-\St$. The notations for eigenvalues (and distinct eigenvalues), eigenvectors and projectors are similar to ours. Furthermore, the following intuitive notations are employed:
\begin{equation}
        \begin{aligned}
    	\Tr_{\geq r_0}(\St) \eqdef \sum\limits_{r \geq r_0} m_r \mu_r
        \nonumber
        \end{aligned}
    \end{equation}
and 
\begin{equation}
        \begin{aligned}
    	\Ptrr \eqdef \sum\limits_{r \geq r_0} \Ptr.
        \nonumber
        \end{aligned}
    \end{equation}
They also define the resolvent
    \begin{equation}
        \begin{aligned}
    	\Resr = \sum\limits_{s\neq r} \frac{1}{\mu_s-\mu_r} \Pts.
        \nonumber
        \end{aligned}
    \end{equation}
Now we are ready to state relative perturbation bounds for eigenvalues and projectors.
\begin{theorem} [\cite{Wahl}, Theorem 3]
    Let $r \geq 1$. Consider $r_0 \geq 1$ such that $\mu_{r_0} \leq \mu_r/2$. Let $x > 0$ be such that for all $s, t < r_0$,
    \begin{equation}
        \begin{aligned}
    	\frac{\| \Pts \Pert \Ptt \|_{\Fr}}{\sqrt{m_s \mu_s m_t \mu_t}},
    	\frac{\| \Pts \Pert \Ptrr \|_{\Fr}}{\sqrt{m_s \mu_s \Tr_{\geq r_0}(\St)}},
    	\frac{\| \Ptrr \Pert \Ptrr \|_{\Fr}}{\Tr_{\geq r_0}(\St)} \leq x.
        \label{(2.8)}
        \end{aligned}
    \end{equation}
    Suppose that
    \begin{equation}
        \begin{aligned}
    	\relr_r(\St) \leq 1/(6x).
        \label{(2.9)}
        \end{aligned}
    \end{equation}
    Then we have
    \begin{equation}
        \begin{aligned}
    	\frac{1}{m_r\mu_r} \sum\limits_{k=1}^{m_r} \left| \lambda_k(\Per(\Se-\mu_r\Id)\Per) - \lambda_k(\Ptr \Pert \Ptr) \right| \leq Cx^2 \relr_r(\St),
        \label{(2.10)}
        \end{aligned}
    \end{equation}
    where $\lambda_k(\cdot)$ denotes the $k$-th largest eigenvalue. In particular, if $j$ is the smallest integer such that $j\in\mathcal{I}_r$, then
    \begin{equation}
        \begin{aligned}
    	\frac{1}{m_r\mu_r} |\widehat{\lambda}_j - \mu_r - \lambda_1(\Ptr \Pert \Ptr)| \leq Cx^2 \relr_r(\St). 
        \nonumber
        \end{aligned}
    \end{equation}
\end{theorem}  

\begin{theorem}[\cite{Wahl}, Theorem 4] \label{Th:JW4}
    Let $r \geq 1$. Consider $r_0 \geq 1$ such that $\mu_{r_0} \leq \mu_r/2$. Let $x$ be such that \eqref{(2.8)} holds. Moreover, suppose that Condition~\eqref{(2.9)} holds. Then we have
    \begin{equation}
        \begin{aligned}
    	\| \Per - \Ptr - \Resr\Pert \Ptr - \Ptr \Pert \Resr\|_{\Fr} \leq Cx^2 \relr_r(\St) \sqrt{\sum\limits_{s\neq r} \frac{m_r \mu_r m_s \mu_s}{(\mu_r-\mu_s^2)}}
        \label{(2.11)}
        \end{aligned}
    \end{equation}
    and 
    \begin{equation}
        \begin{aligned}
    	\left|\| \Per - \Ptr\|_{\Fr}^2 - 2\|\Resr\Pert \Ptr\|_{\Fr}^2 \right| \leq Cx^3 \relr_r(\St) \sum\limits_{s\neq r} \frac{m_r \mu_r m_s \mu_s}{(\mu_r-\mu_s^2)}.
        \label{(2.12)}
        \end{aligned}
    \end{equation}
\end{theorem}

%% file: source/proofs/Proof_TechLemmas.tex
\begin{proof}[Proof of Proposition~\ref{P:Properties0}]
\mbox{}\\
\textit{(i)} 
Homogeneity is trivial. 
Triangle inequality follows directly from triangle inequality for spectral norm combined with triangle inequality for maximum. The only property to check is that $\| A \|_{(\Pp,\Gamma, s_1, s_2)} = 0$ implies $A = 0$. Indeed, as we will see in $(ii)$, $\| A \|_{(\Pp, \Gamma, s_1, s_2)} = 0$ implies $\| A \| = 0$, and thus $A=0$, since spectral norm is a norm.
\\
\noindent \textit{(ii)} 
To prove the desired bounds, it is more convenient to use the representation from Lemma~\ref{P:Properties1} (i) (which is proved independently slightly later).
The upper bound is trivially implied by the inequalities
\begin{equation} 
        \begin{aligned}
            &\| \Gamma_1^\T A \Gamma_1 \| \leq \| A \|,\;\;\;\;\| \Gamma_2^\T A \Gamma_2 \| \leq \| A \|,\\
            &\sup\limits_{\substack{v \in \myset^m_{s_1} \\ w \in \myset^{d-m}_{s_2}}} v^\T \Gamma_1^\T A \,\Gamma_2 w \leq
            \sup\limits_{\substack{v \in \Sph^{m-1} \\ w \in \Sph^{d-m-1}}} v^\T \Gamma_1^\T A \,\Gamma_2 w  \leq 
            \sup\limits_{\substack{v \in \Sph^{d-1} \\ w \in \Sph^{d-1}}} v^\T A w = \|A\|.
        \nonumber 
        \end{aligned} 
    \end{equation}
Let us prove the lower bound.
Let $\widetilde{u}$ be the eigenvector corresponding to largest absolute eigenvalue of $A$. Define $\widetilde{v} = \Gamma_1^\T\widetilde{u} \in\R^m$ and $\widetilde{w} = \Gamma_2^\T\widetilde{u}\in\R^{d-m}$, so that $\widetilde{u} = \Gamma_1 \widetilde{v} + \Gamma_2 \widetilde{w}$. Note that $\|\widetilde{v}\| \leq 1$ and $\|\widetilde{w}\| \leq 1$. Then
    \begin{equation} 
        \begin{aligned}
            \| A \| &= |\widetilde{u}^\T A \widetilde{u}| 
            = \left| (\Gamma_1 \widetilde{v} + \Gamma_2 \widetilde{w})^\T A (\Gamma_1 \widetilde{v} + \Gamma_2 \widetilde{w}) \right|\\
            &\leq |\widetilde{v}^\T \Gamma_1^\T A \Gamma_1 \widetilde{v}| + |\widetilde{w}^\T \Gamma_2^\T A \Gamma_2 \widetilde{w}|
            + 2 |\widetilde{v}^\T \Gamma_1^\T A \Gamma_2 \widetilde{w}|\\
            &\leq \| \Gamma_1^\T A \Gamma_1 \| + \| \Gamma_2^\T A \Gamma_2 \|
            + 2 |\widetilde{v}^\T \Gamma_1^\T A \Gamma_2 \widetilde{w}|.
        \nonumber 
        \end{aligned} 
    \end{equation}
To bound $|\widetilde{v}^\T \Gamma_1^\T A \Gamma_2 \widetilde{w}|$ in terms of $\sup\limits_{\substack{v\in\myset^m_{s_1}\\ w\in\myset^{d-m}_{s_2}}} v^\T \Gamma_1^\T A \Gamma_2 w$, let us decompose
\begin{equation} 
        \begin{aligned}
            &\widetilde{v} = \sum\limits_{k=1}^{\lceil m/s_1\rceil} \widetilde{v}^{(k)},\\
            &\widetilde{w} = \sum\limits_{l=1}^{\lceil (d-m)/s_2\rceil} \widetilde{w}^{(l)},
        \nonumber 
        \end{aligned} 
    \end{equation}
with
    \begin{equation} 
        \begin{aligned}
            &\supp(\widetilde{v}^{(k)}) \subseteq \{ (k-1)s_1+1,\ldots, ks_1 \} \;\;\;\text{ for all }\;k\in\lceil m/s_1\rceil,\\
            &\supp(\widetilde{w}^{(l)}) \subseteq \{ (l-1)s_2+1,\ldots, ls_2 \} \;\;\;\text{ for all }\;l\in\lceil (d-m)/s_2\rceil,
        \nonumber 
        \end{aligned} 
    \end{equation}
where $\supp(\cdot)$ denotes support of a vector.  Therefore, 
\begin{equation} 
        \begin{aligned}
            &|\widetilde{v}^\T \Gamma_1^\T A \Gamma_2 \widetilde{w}| =
            \left|\left( \sum\limits_{k=1}^{\lceil m/s_1\rceil} \widetilde{v}^{(k)} \right)^\T \Gamma_1^\T A \Gamma_2 \left( \sum\limits_{l=1}^{\lceil (d-m)/s_2\rceil} \widetilde{w}^{(l)} \right)\right|
             \\&\qquad \leq
            \sum\limits_{k=1}^{\lceil m/s_1\rceil} \sum\limits_{l=1}^{\lceil (d-m)/s_2\rceil} \left| (\widetilde{v}^{(k)})^\T \Gamma_1^\T A \Gamma_2  \widetilde{w}^{(l)} \right|
            \\&\qquad=
            \sum\limits_{k=1}^{\lceil m/s_1\rceil} \sum\limits_{l=1}^{\lceil (d-m)/s_2\rceil} \left| \frac{(\widetilde{v}^{(k)})^\T}{\|\widetilde{v}^{(k)}\|} \Gamma_1^\T A \Gamma_2  \frac{\widetilde{w}^{(l)}}{\|\widetilde{w}^{(l)}\|} \right|\cdot \| \widetilde{v}^{(k)}\| \|\widetilde{w}^{(l)}\|
            \\
            &\qquad\leq  \sup\limits_{\substack{v\in\myset^m_{s_1}\\w\in\myset^{d-m}_{s_2}}} v^\T \Gamma_1^\T A \Gamma_2 w
            \cdot  \sum\limits_{k=1}^{\lceil m/s_1\rceil}
            \| \widetilde{v}^{(k)}\| \cdot \sum\limits_{l=1}^{\lceil (d-m)/s_2\rceil}
             \|\widetilde{w}^{(l)}\|
             \\&\qquad
             \leq
             \sqrt{\left\lceil \frac{m}{s_1} \right\rceil \cdot \left\lceil \frac{d-m}{s_2} \right\rceil} \cdot \sup\limits_{\substack{v\in\myset^m_{s_1}\\w\in\myset^{d-m}_{s_2}}} v^\T \Gamma_1^\T A \Gamma_2 w.
        \nonumber 
        \end{aligned} 
    \end{equation}
Here we used 
\begin{equation} 
        \begin{aligned}
        \sum\limits_{k=1}^{\lceil m/s_1\rceil} \| \widetilde{v}^{(k)} \| \leq  \sqrt{\left\lceil \frac{m}{s_1} \right\rceil} \sum\limits_{k=1}^{\lceil \frac{m}{s_1}\rceil}\| \widetilde{v}^{(k)} \|^2 = 
         \sqrt{\left\lceil \frac{m}{s_1} \right\rceil} \| \widetilde{v} \|^2 \leq \sqrt{\left\lceil \frac{m}{s_1} \right\rceil},
        \nonumber 
        \end{aligned} 
    \end{equation}
since $\{ \widetilde{v}^{(k)} \}_{k=1}^{\lceil m/s_1\rceil}$ are orthogonal, and similarly $
        \sum\limits_{l=1}^{\lceil (d-m)/s_2\rceil} \| \widetilde{w}^{(l)} \| \leq \sqrt{\left\lceil (d-m)/s_2 \right\rceil}$.

Hence,
\begin{equation} 
        \begin{aligned}
            \|A\| &\leq
            \| \Gamma_1^\T A \Gamma_1 \| + \| \Gamma_2^\T A \Gamma_2 \| + 2|\widetilde{v}^\T \Gamma_1^\T A \Gamma_2 \widetilde{w}| 
            \\& \leq \| \Gamma_1^\T A \Gamma_1 \| + \| \Gamma_2^\T A \Gamma_2 \| + 2\sqrt{\left\lceil \frac{m}{s_1} \right\rceil \cdot \left\lceil \frac{d-m}{s_2} \right\rceil} \cdot \sup\limits_{\substack{v\in\myset^m_{s_1}\\w\in\myset^{d-m}_{s_2}}}
            v^\T \Gamma_1^\T A \Gamma_2 w
            \\& \leq 2 \sqrt{\left\lceil \frac{m}{s_1} \right\rceil \cdot \left\lceil \frac{d-m}{s_2} \right\rceil} \cdot \| A\|_{(\Pp, \Gamma, s_1,s_2)}.
        \nonumber 
        \end{aligned} 
    \end{equation}
\end{proof}

\begin{proof}[Proof of Lemma~\ref{L:KC_SCconc}]
    Fix $s, t \in [q]$. Expanding squared Frobenius norm over the basis of eigenvectors $\{ u_j \}_{j=1}^d$, we have
    \begin{equation} 
        \begin{aligned}
            &\| \Pts (\Se-\St) \Ptt \|_{\Fr}^2 
            = \sum\limits_{j,k=1}^d \left( u_j^\T \Pts(\Se-\St)\Ptt u_k \right)^2
            =\sum\limits_{j\in\mathcal{I}_s}\sum\limits_{k\in\mathcal{I}_t} \left( u_j^\T (\Se-\St) u_k \right)^2\\
            &\hspace{1cm} \leq m_s m_t \max\limits_{\substack{j\in \mathcal{I}_s\\ k\in\mathcal{I}_t}} \left( u_j^\T (\Se-\St) u_k \right)^2
            = m_s \mu_s m_t \mu_t \max\limits_{\substack{j\in \mathcal{I}_s\\ k\in\mathcal{I}_t}} \left( u_j^\T (\St^{-1/2}\Se\St^{-1/2}-\Id_d) u_k \right)^2.
        \nonumber 
        \end{aligned} 
    \end{equation}
    Hence,
    \begin{equation} 
        \begin{aligned}
            &\max\limits_{s,t\in[q]}\frac{\| \Pts (\Se-\St) \Ptt \|_{\Fr}}{\sqrt{m_s \mu_s m_t \mu_t}}
            \leq  \max\limits_{s,t\in[q]} \max\limits_{\substack{j\in \mathcal{I}_s\\ k\in\mathcal{I}_t}} \left| u_j^\T (\St^{-1/2}\Se\St^{-1/2}-\Id_d) u_k \right|\\
            &\qquad= \max\limits_{j,k\in[d]} \left| u_j^\T (\St^{-1/2}\Se\St^{-1/2}-\Id_d) u_k \right|
            = \max\limits_{j,k\in[d]} \left| \left[ U^\T\St^{-1/2}\Se\St^{-1/2}U-\Id_d \right]_{j,k} \right|\\
            &\qquad= \|U^{\T} \St^{-1/2}\Se\St^{-1/2}U-\Id_d \|_{\mymax},
        \nonumber 
        \end{aligned} 
    \end{equation}
    which is maximum absolute elementwise norm of the difference between sample and true covariance matrices of random vectors $\{ U^\T \St^{-1/2} X_i \}_{i=1}^n$, where columns of $U$ are eigenvectors $\{ u_j \}_{j=1}^d$.
    This fits the framework of Theorem~\ref{Th:KC_SCmax}.
    The joint Orlicz norm of these vectors is 
    \begin{equation} 
        \begin{aligned}
            & \| U^\T\St^{-1/2} X_i \|_{J,\phi_\beta} = \| \St^{-1/2} X_i \|_{J,\phi_\beta} \leq c < \infty, \;\;\;i\in[n],
        \nonumber 
        \end{aligned} 
    \end{equation}
    due to Assumption~\ref{A: tails}.
    Therefore, Theorem~\ref{Th:KC_SCmax} applied with $U^\T\St^{-1/2} X_i$ instead of $X_i$, $d$ instead of $p$, $K_{n,p} = c$ and $z=\log(3n)$ implies
    \begin{equation} 
        \begin{aligned}
            &\|U^\T\St^{-1/2}\Se\St^{-1/2}U-\Id_d \|_{\mymax} \leq \\
            &\hspace{1cm}\leq
            C_\beta c^2 \left( \sqrt{\frac{\log(3n)+2\log(d)}{n}} + \frac{(\log(2n))^{1/\beta} (\log(3n)+2\log(d))^{2/\beta}}{n}\right)
        \nonumber 
        \end{aligned} 
    \end{equation}
    with probability $1-1/n$.
    
\end{proof}

\begin{proof}[Proof of Lemma \ref{L:KL}]
    Theorem~\ref{Th:JW4} is stated for infinite-dimensional Hilbert space $\mathcal{H}$, so we can take $\mathcal{H}$ to be some space, in which $\R^d$ is embedded. Consider covariance operator $\St_\mathcal{H}$ that acts on an element of $\mathcal{H}$ is the same way as $\St$ acts on the first $d$ components of this element. Similarly, $\Sp_\mathcal{H}$ is a counterpart of $\Sp$. Operator $\St_\mathcal{H}$ has $q+1$ distinct eigenvalues: the first $q$ are all the eigenvalues of $\St$, specifically $\mu_1, \ldots, \mu_q$, and the last one is $\mu_{q+1} = 0$. The corresponding projectors for the first $q$ eigenvalues coincide with $\Pp_1, \ldots, \Pp_q$ and the last projector is $\Pp_{q+1} = \Id - \sum\limits_{r\in[q]}\Ptr$ (here $\Id$ is identity operator in $\mathcal{H}$).
    
    Now we apply Theorem~\ref{Th:JW4} for every $r\in\J$ with $r_0 = q+1$.
    Let us verify the conditions.
    Note that $\mu_{r_0} = 0 \leq \mu_r/2$.
    The first inequality of Condition~\eqref{(2.8)} is satisfied automatically by the specific choice of $x$. A bit tricky things are happening to the second and third inequalities of Condition~\eqref{(2.8)}. Observe that
    \begin{equation}
        \begin{aligned}
    	\| \Pts (\Sp_\mathcal{H}-\St_\mathcal{H}) \Ptrr \|_{\Fr} = 0,\;\;
    	\| \Ptrr (\Sp_\mathcal{H}-\St_\mathcal{H}) \Ptrr \|_{\Fr} = 0,\;\;
    	\Tr_{\geq r_0}(\St_\mathcal{H}) = 0,
        \nonumber
        \end{aligned}
    \end{equation}
    so the second and third inequalities of Condition~\eqref{(2.8)} become $0/0 \leq x$, which doesn't allow us to apply this result rigorously. However, from the analysis of the proof of Theorem~4 of \cite{Wahl} it is clear that these inequalities can be replaced by
    \begin{equation}
        \begin{aligned}
    	&\| \Pts (\Sp_\mathcal{H}-\St_\mathcal{H})\Ptrr \|_{\Fr} \leq x\cdot\sqrt{m_s \mu_s \Tr_{\geq r_0}(\St_\mathcal{H})},\\
    	&\| \Ptrr (\Sp_\mathcal{H}-\St_\mathcal{H}) \Ptrr \|_{\Fr} \leq x\cdot\Tr_{\geq r_0}(\St_\mathcal{H}),
        \nonumber
        \end{aligned}
    \end{equation}
    and all the derivation stays true (division by $\Tr_{\geq r_0}(\St_\mathcal{H})$ actually never appears in the proof). In our situation, these inequalities reduce to $0 \leq x\cdot 0$, which holds true.
    Finally, Condition~\eqref{(2.9)} is fulfilled due to Condition~\eqref{A:JW}.
    
    Thus, we obtain the following: for all $r\in\J$
    \begin{equation}
        \begin{aligned}
    	\| \Ppr - \Ptr - \Resr(\Sp-\St) \Ptr - \Ptr (\Sp-\St) \Resr\|_{\Fr} \leq Cx^2 \relr_r(\St) \sqrt{\sum\limits_{s\neq r} \frac{m_r\mu_r m_s \mu_s}{(\mu_r-\mu_s)^2}},
        \nonumber
        \end{aligned}
    \end{equation}
    which, by triangle inequality leads to
    \begin{equation}
        \begin{aligned}
    	&\left\| \sum\limits_{r\in\J} \Ppr - \sum\limits_{r\in\J} \Ptr - \sum\limits_{r\in\J} \left(\Resr(\Sp-\St) \Ptr + \Ptr (\Sp-\St) \Resr\right)\right\|_{\Fr} \leq 
    	\\&\qquad \leq Cx^2  \sum\limits_{r\in\J} \left( \relr_r(\St) \sqrt{\sum\limits_{s\neq r} \frac{m_r\mu_r m_s \mu_s}{(\mu_r-\mu_s)^2}} \right).
        \nonumber
        \end{aligned}
    \end{equation}
    The only thing left to note is
    \begin{equation}
        \begin{aligned}
    	\sum\limits_{r\in\J} \left(\Resr(\Sp-\St) \Ptr + \Ptr (\Sp-\St) \Resr\right)
    	& = &
    	\sum\limits_{r\in\J} \sum\limits_{s\notin\J} 
    			\frac{\Ptr (\Sp-\St) \Pts + \Pts (\Sp-\St) \Ptr}{\mu_r - \mu_s},
        \nonumber
        \end{aligned}
    \end{equation}
    which can be seen from inserting the resolvents and observing that the terms of the type $\frac{\Ptr(\Sp-\St)\Ptrp}{\mu_r - \mu_{r^{\prime}}}$, $r,r^{\prime} \in \J, r\neq r^{\prime}$ cancel out since they appear exactly twice in the sum with different signs.
\end{proof}

\begin{proof}[Proof of Lemma~\ref{L:KL_gen}]
Denote
\begin{equation}
        \begin{aligned}
        	x = \max\limits_{s,t\in[q]} \frac{\| \Pts (\Se-\St) \Ptt \|_{\Fr}}{\sqrt{m_s\mu_s m_t\mu_t}},\;\;\;
        	\widetilde{x} = \max\limits_{s,t\in[q]} \frac{\| \Pts (\Sp-\Se) \Ptt \|_{\Fr}}{\sqrt{m_s\mu_s m_t\mu_t}}.
        \nonumber
        \end{aligned}
    \end{equation}
    Let $\Omega$ be the event on which $x \leq \scc$ and $\widetilde{\Omega} = \widetilde{\Omega}(\data)$ be the event on which $(\widetilde{x} \leq \widetilde{\psi}_n\,|\,\data)$. By Lemma~\ref{L:KC_SCconc}, $\Prob[\Omega] \geq 1-1/n$. By Condition~\eqref{prob_assumption} $\Prob[\widetilde{\Omega} \,|\, \data] \geq 1-1/n$ on some event $\Omega^\prime$ with $\Prob[\Omega^\prime] \geq 1-1/n$. By union bound, $\Prob\left[\Omega \cap \Omega^\prime\right] \geq 1-2/n$. 
    
As in Lemma~\ref{L:KL}, we decompose
    \begin{equation}
        \begin{aligned}
            \Ppp - \Pe &
            = (\Ppp - \Pt) - (\Pe - \Pt)
            \\&= L_\J(\Sp-\St) + R_\J(\Sp-\St)
            - L_\J(\Se-\St) - R_\J(\Se-\St)
            \\&= L_\J(\Sp-\Se) + R_\J(\Sp-\St) - R_\J(\Se-\St).
        \nonumber
        \end{aligned}
    \end{equation}
Then, by Proposition~\ref{P:Properties0} (i), (ii)
\begin{equation}
        \begin{aligned}
        	&\left| \| \Ppp-\Pe\|_{(\Pt,\Gamma^\circ,s_1, s_2)} - \| L_\J(\Sp-\Se) \|_{(\Pt, \Gamma^\circ, s_1, s_2)} \right| 
        	\leq
        	\\&\qquad\leq
        	\| R_\J(\Sp-\St) \|_{(\Pt,\Gamma^\circ,s_1, s_2)} + \| R_\J(\Se-\St) \|_{(\Pt, \Gamma^\circ, s_1, s_2)}
        	\\&\qquad\leq2\| R_\J(\Sp-\St) \| + 2\| R_\J(\Se-\St) \|.
        \nonumber
        \end{aligned}
    \end{equation}
Further, on $\Omega \cap \Omega^\prime$ with $\delta = 4C(\scc + \widetilde{\psi}_n)^2 \rr_\J(\St)^{3/2}$ (with $C$ from Lemma~\ref{L:KL}) we have
    \begin{equation}
        \begin{aligned}
        	&\Prob\left[ 2\| R_\J(\Sp-\St) \| + 2\| R_\J(\Se-\St) \| > \delta \,\big|\,\data \right] =
        	\\&\qquad =\Prob\left[ 2\| R_\J(\Sp-\St) \| + 2\| R_\J(\Se-\St) \| > \delta \,\big|\,\data; \widetilde{x} > \widetilde{\psi}_n \right] \cdot \Prob\left[ \widetilde{x} > \widetilde{\psi}_n\,|\,\data\right] +
        	\\&\qquad\qquad \Prob\left[ 2\| R_\J(\Sp-\St) \| + 2\| R_\J(\Se-\St) \| > \delta \,\big|\,\data; \widetilde{x} \leq  \widetilde{\psi}_n \right] \cdot \Prob\left[ \widetilde{x} \leq \widetilde{\psi}_n\,|\,\data\right]
        	\\&\qquad \leq 1\cdot\frac{1}{n} + \Prob\left[ 2\| R_\J(\Sp-\St) \| + 2\| R_\J(\Se-\St) \| > \delta \,\big|\,\data; \widetilde{x} \leq  \widetilde{\psi}_n \right]\cdot 1.
        \nonumber
        \end{aligned}
    \end{equation}
So far we have used only that we are on $\Omega^\prime$. Since we are also on $\Omega$, $x\leq \scc$ implies $x\max\limits_{r\in\J} \relr_r(\St) \leq 1/12$, yielding that Condition~\eqref{A:JW} is fulfilled. Thus, by Lemma~\ref{L:KL}
    \begin{equation}
        \begin{aligned}
        	&\| R_\J(\Se-\St) \| \leq Cx^2\rr_\J(\St)^{3/2}.
        \nonumber
        \end{aligned}
    \end{equation}
Similarly, when $\widetilde{x} \leq \widetilde{\psi}_n$, Condition~\eqref{A:JW} is satisfied for $(x+\widetilde{x})$, and Lemma~\ref{L:KL} claims
    \begin{equation}
        \begin{aligned}
        	&\| R_\J(\Sp-\St) \| \leq C(x+\widetilde{x})^2\rr_\J(\St)^{3/2}.
        \nonumber
        \end{aligned}
    \end{equation}
Therefore, on $\Omega$
    \begin{equation}
        \begin{aligned}
            \Prob\left[ 2\| R_\J(\Sp-\St) \| + 2\| R_\J(\Se-\St) \| > \delta \,\big|\,\data; \widetilde{x} \leq  \widetilde{\psi}_n \right] = 0,
        \nonumber
        \end{aligned}
    \end{equation}
yielding the desired.
\end{proof}

\begin{proof}[Proof of Lemma~\ref{P:Properties1}]
\mbox{}\\
\noindent \textit{(i)} The first two terms coincide with the original definition, so we have to make sure that the third one coincides as well. From the definitions of $\myset^m_{s_1}$ and $\myset^{d-m}_{s_2}$, we have
\begin{equation} 
        \begin{aligned}
            &\sup\limits_{\substack{v \in \myset^m_{s_1} \\ w \in \myset^{d-m}_{s_2}}} v^\T \Gamma_1^\T A \,\Gamma_2 w 
            =
            \max\limits_{\substack{k \in \{0,\ldots,m-s_1\} \\ l \in \{0,\ldots,d-m-s_2\}}} 
            \sup\limits_{\substack{v\in\Sph^{s_1-1}\\w\in\Sph^{s_2-1}}} 
            \begin{bmatrix} 0_k^\T, v^\T, 0_{m-k-s_1}^\T \end{bmatrix} \Gamma_1^\T A \Gamma_2 \begin{bmatrix} 0_l \\ w\\ 0_{d-m-l-s_2} \end{bmatrix}\\
            &\qquad = 
            \max\limits_{\substack{k \in \{0,\ldots,m-s_1\} \\ l \in \{0,\ldots,d-m-s_2\}}} 
            \sup\limits_{\substack{v\in\Sph^{s_1-1}\\w\in\Sph^{s_2-1}}} 
            v^\T [\Gamma_1^\T A \Gamma_2]_{[k+1:k+s_1], [l+1:l+s_2]}\, w\\
            &\qquad =
             \max\limits_{\substack{k \in \{0,\ldots,m-s_1\} \\ l \in \{0,\ldots,d-m-s_2\}}} 
             \left\| [\Gamma_1^\T A \Gamma_2]_{[(k+1):(k+s_1)], [(l+1):(l+s_2)]} \right\|,
        \nonumber 
        \end{aligned} 
    \end{equation}
as desired.
\\
\noindent \textit{(ii)}
Note that $\Pp \Gamma_1 = \Gamma_1$, $\Pp \Gamma_2 = \Oo_{d\times(d-m)}$, $(\Id_d - \Pp) \Gamma_1 = \Oo_{d\times m}$, $(\Id_d - \Pp) \Gamma_2 = \Gamma_2$. Hence, plugging $A = \Pp A (\Id_d - \Pp) + (\Id_d - \Pp)A\Pp$ into the definition, we notice that the first two terms $\| \Gamma_1^\T A \Gamma_1 \|/2$ and $\| \Gamma_2^\T A \Gamma_2 \|/2$ disappear and $\| A\|_{(\Pp, \Gamma, s_1, s_2)}$ is expressed by the third term only.

Let us take $s_1 = m$ and $s_2 = d-m$. We can represent spectral norm as
    \begin{equation} 
        \begin{aligned}
            \| A \| &=
            \sup\limits_{u\in\Sph^{d-1}} |u^\T A u|
            = \sup\limits_{u\in\Sph^{d-1}} |u^\T \left[ \Pp A (\Id_d - \Pp) + (\Id_d - \Pp) A \Pp \right] u|
            \\ &= 2 \sup\limits_{u\in\Sph^{d-1}} |u^\T \Pp A (\Id_d - \Pp)  u|.
        \nonumber 
        \end{aligned} 
    \end{equation}
Note that
    \begin{equation} 
        \begin{aligned}
            2 \sup\limits_{u\in\Sph^{d-1}} |u^\T \Pp A (\Id_d - \Pp)  u| 
            =
            \sup\limits_{\substack{v\in\Sph^{m-1} \\ w\in\Sph^{d-m-1} }} v^\T \Gamma_1^\T A \Gamma_2 w.
        \label{intermediate_spectral} 
        \end{aligned} 
    \end{equation}
Indeed, for any $u\in\Sph^{d-1}$ we can take 
\begin{equation} 
        \begin{aligned}
        v = \pm \Gamma_1^\T \Pp u/\| \Pp u\|  \;\text{ and }\; 
        w = \pm \Gamma_2^\T (\Id_d-\Pp) u/\| (\Id_d-\Pp) u\|
        \nonumber 
        \end{aligned} 
    \end{equation}
     (it is straightforward to check $v\in \Sph^{m-1}$ and $w\in\Sph^{d-m-1}$) and obtain
\begin{equation} 
        \begin{aligned}
            2 |u^\T \Pp A (\Id_d - \Pp)  u|
            = 2 |v^\T \Gamma_1^\T A \Gamma_2 w| \cdot \| \Pp u\| \| (\Id_d-\Pp) u\| \leq |v^\T \Gamma_1^\T A \Gamma_2 w|.
        \nonumber 
        \end{aligned} 
    \end{equation}
Conversely, for any $v\in\Sph^{m-1}$ and $w\in\Sph^{d-m-1}$ we can take $u = (\Gamma_1 v + \Gamma_2 w)/\sqrt{2}$ (again, easy to see that $u\in\Sph^{d-1}$) and obtain
\begin{equation} 
        \begin{aligned}
            2 u^\T \Pp A (\Id_d - \Pp)  u 
            =
            v^\T \Gamma_1^\T A \Gamma_2 w.
        \nonumber 
        \end{aligned} 
    \end{equation}
This proves \eqref{intermediate_spectral}, and, consequently,
\begin{equation} 
        \begin{aligned}
            \sup\limits_{\substack{v\in\myset^{m}_m \\ w\in\myset^{d-m}_{d-m} }} v^\T \Gamma_1^\T A \Gamma_2 w = \| A \|.
        \nonumber
        \end{aligned} 
    \end{equation}
    
If we take $s_1 = 1$ and $s_2 = 1$, then 
\begin{equation} 
        \begin{aligned}
            &\sup\limits_{\substack{v \in \myset^m_{1} \\ w \in \myset^{d-m}_{1}}} v^\T \Gamma_1^\T A \,\Gamma_2 w 
            =
            \max\limits_{j,k\in[d]} \left|e_j^\T \Gamma_1^\T A \,\Gamma_2 e_k\right| =
            \max\limits_{j,k\in[d]} \left|\left[ \Gamma_1^\T A \,\Gamma_2 \right]_{j,k}\right| = \| \Gamma_1^\T A \,\Gamma_2 \|_{\mymax},
        \nonumber 
        \end{aligned} 
    \end{equation}
    where $\{ e_j \}_{j=1}^d$ are standard basis vectors in $\R^d$.
\end{proof}

\begin{proof}[Proof of Lemma~\ref{L:discretization}]
For any $v \in \myset^m_{s_1}$ denote the closest to $v$ vector of $N_{\eps}(\myset^m_{s_1})$ as $\pi(v)$, that is, $\| v - \pi(v) \| \leq \eps$. Similarly, for any $w \in \myset^{d-m}_{s_2}$ denote the closest to $w$ vector of $N_{\eps}(\myset^{d-m}_{s_2})$ as $\rho(w)$, that is, $\| w - \rho(w) \| \leq \eps$. The construction \eqref{eps-net construction} allows without loss of generality assume $(v - \pi(v)) \in \myset^m_{s_1}$ and $(w - \rho(w)) \in \myset^{d-m}_{s_2}$.

    By Lemma~\ref{P:Properties1} (ii),
    \begin{equation} 
        \begin{aligned}
            \| A \|_{(\Pp, \Gamma, s_1, s_2)} = \sup\limits_{\substack{v \in \myset^m_{s_1}\\ w \in \myset^{d-m}_{s_2}}}  v^\T \Gamma_1^\T  A \Gamma_2 w .
        \nonumber 
        \end{aligned} 
    \end{equation}
    We have the following standard chain of equalities and inequalities:
    \begin{equation}
        \begin{aligned}
        	&  \| A \|_{(\Pp, \Gamma, s_1, s_2)} 
        	=
        	\sup\limits_{\substack{v \in \myset^m_{s_1}\\ w \in \myset^{d-m}_{s_2}}} v^{\T} \Gamma_1^\T A \Gamma_2 w\\
        	&\qquad=
        	\sup\limits_{\substack{v \in \myset^m_{s_1}\\ w \in \myset^{d-m}_{s_2}}} \left[ v^{\T} \Gamma_1^\T A \Gamma_2 w - \pi(v)^{\T} \Gamma_1^\T A \Gamma_2 \rho(w) + \pi(v)^{\T} \Gamma_1^\T A \Gamma_2 \rho(w) \right]\\
        	&\qquad\leq
        	\sup\limits_{\substack{v \in \myset^m_{s_1}\\ w \in \myset^{d-m}_{s_2}}} \left[ v^{\T} \Gamma_1^\T A \Gamma_2 w - \pi(v)^{\T} \Gamma_1^\T A \Gamma_2 \rho(w)\right] + \sup\limits_{\substack{v \in \myset^m_{s_1}\\ w \in \myset^{d-m}_{s_2}}} \pi(v)^{\T} \Gamma_1^\T A \Gamma_2 \rho(w)\\
        	&\qquad=
        	\sup\limits_{\substack{v \in \myset^m_{s_1}\\ w \in \myset^{d-m}_{s_2}}} \left\{ (v-\pi(v))^{\T} \Gamma_1^\T A \Gamma_2 w + \pi(v)^{\T} \Gamma_1^\T A \Gamma_2 (w-\rho(w)) \right\}
        	+ \\
        	&\qquad\qquad+\; \max\limits_{(v, w) \in N_{\eps}(\myset^{m}_{s_1}) \times N_{\eps}(\myset^{d-m}_{s_2})} v^{\T} \Gamma_1^\T A \Gamma_2 w
        	\\
        	&\qquad\leq
        	\eps\cdot\sup\limits_{\substack{v \in \myset^m_{s_1}\\ w \in \myset^{d-m}_{s_2}}} \left\{ \frac{(v-\pi(v))^{\T}}{\| v-\pi(v) \|_2} \Gamma_1^\T A \Gamma_2 w + \pi(v)^{\T} \Gamma_1^\T A \Gamma_2 \frac{(w-\rho(w))}{\|w-\rho(w)\|_2} \right\}
        	+  \max\limits_{j\in[p]} v_j^{\T} A w_j
        	\\
        	&\qquad\leq
        	2\eps\cdot\sup\limits_{\substack{v \in \myset^{m}_{s_1}\\ w \in \myset^{d-m}_{s_2}}} v^{\T} \Gamma_1^\T A \Gamma_2 w + \max\limits_{j\in[p]}\; v_j^{\T} A w_j
        	=
        	2\eps\cdot \| A \|_{(\Pp, \Gamma, s_1,s_2)} + \max\limits_{j\in[p]}\; v_j^{\T} A w_j.
            \nonumber
        \end{aligned}
    \end{equation}
    Therefore, we obtain
    \begin{equation}
        \begin{aligned}
        	\max\limits_{j\in[p]} \; v_j^{\T} A w_j
        	\leq
        	 \; \| A \|_{(\Pp,\Gamma, s_1, s_2)} \leq \frac{1}{1-2\eps} \max\limits_{j\in[p]} \; v_j^{\T} A w_j.
            \nonumber
        \end{aligned}
    \end{equation}

\end{proof}

\begin{proof}[Proof of Lemma \ref{L:epsnet_size}]
    It is a well-known fact that
    \begin{equation}
        \begin{aligned}
        	N_\eps(\Sph^{s_1-1}) \leq \left( \frac{3}{\eps} \right)^{s_1},\;\;\;
        	N_\eps(\Sph^{s_2-1}) \leq \left( \frac{3}{\eps} \right)^{s_2},
            \nonumber
        \end{aligned}
    \end{equation}
    e.g. see Lemma 5.13 of \cite{Lectures}. From the construction \eqref{eps-net construction} follows
    \begin{equation}
        \begin{aligned}
        	&N_\eps(\myset^{m}_{s_1}) \leq (m-s_1+1)\cdot N_\eps(\Sph^{s_1-1}) \leq (m-s_1+1)\cdot\left( \frac{3}{\eps} \right)^{s_1},\\
        	&N_{\eps}(\myset^{d-m}_{s_2}) \leq (d-m-s_2+1)\cdot N_\eps(\Sph^{s_2-1}) \leq (d-m-s_2+1)\cdot\left( \frac{3}{\eps} \right)^{s_2}.
            \nonumber
        \end{aligned}
    \end{equation}
    Taking logarithm of $p(\eps, d, m, s_1, s_2) = |N_\eps(\myset^{m}_{s_1})| \cdot |N_{\eps}(\myset^{d-m}_{s_2})|$, we get the desired bound.
\end{proof}

\begin{proof}[Proof of Lemma \ref{L:BootConc}]
    We start with following the proof of Lemma~\ref{L:KC_SCconc}.
    Fix $s, t \in [q]$. Expanding squared Frobenius norm over the basis of eigenvectors $\{ u_j \}_{j=1}^d$ and using the definition of $\Sb$, we have
    \begin{equation} 
        \begin{aligned}
            &\| \Pts (\Sb-\Se) \Ptt \|_{\Fr}^2 
            = \sum\limits_{j,k=1}^d \left( u_j^\T \Pts(\Sb-\Se)\Ptt u_k \right)^2
            =\sum\limits_{j\in\mathcal{I}_s}\sum\limits_{k\in\mathcal{I}_t} \left( u_j^\T (\Sb-\Se) u_k \right)^2\\
            &\hspace{1cm} \leq m_s m_t \max\limits_{\substack{j\in \mathcal{I}_s\\ k\in\mathcal{I}_t}} \left( u_j^\T (\Sb-\Se) u_k \right)^2
            =
            m_s m_t \max\limits_{\substack{j\in \mathcal{I}_s\\ k\in\mathcal{I}_t}} \left( u_j^\T \frac{1}{n}\sum\limits_{i=1}^n (\eta_i-1)X_i X_i^\T u_k \right)^2\\
            &\hspace{1cm}= m_s \mu_s m_t \mu_t \max\limits_{\substack{j\in \mathcal{I}_s\\ k\in\mathcal{I}_t}} \left( u_j^\T \frac{1}{n}\sum\limits_{i=1}^n (\eta_i-1) (\St^{-1/2} X_i) (\St^{-1/2} X_i)^\T u_k \right)^2.
        \nonumber 
        \end{aligned} 
    \end{equation}
    Hence,
    \begin{equation} 
        \begin{aligned}
            &\max\limits_{s,t\in[q]}\frac{\| \Pts (\Sb-\Se) \Ptt \|_{\Fr}}{\sqrt{m_s \mu_s m_t \mu_t}}
            \leq \max\limits_{s,t\in[q]}\max\limits_{\substack{j\in \mathcal{I}_s\\ k\in\mathcal{I}_t}} \left| u_j^\T \frac{1}{n}\sum\limits_{i=1}^n (\eta_i-1) (\St^{-1/2} X_i) (\St^{-1/2} X_i)^\T u_k \right|\\
            &\hspace{1cm}= \max\limits_{j,k\in[d]}  \left| u_j^\T \frac{1}{n}\sum\limits_{i=1}^n (\eta_i-1) (\St^{-1/2} X_i) (\St^{-1/2} X_i)^\T u_k \right|\\
            &\hspace{1cm}= \max\limits_{j,k\in[d]}  \left|\left[ \frac{1}{n}\sum\limits_{i=1}^n (\eta_i-1) (U^\T\St^{-1/2} X_i) (U^\T\St^{-1/2} X_i)^\T \right]_{j,k} \right|\\
            &\hspace{1cm}= \max\limits_{j,k\in[d]}  \left| \frac{1}{n}\sum\limits_{i=1}^n (\eta_i-1) (U^\T\St^{-1/2} X_i)_j (U^\T\St^{-1/2} X_i)_k \right|.
        \nonumber 
        \end{aligned} 
    \end{equation}
    For arbitrary $j,k\in[d]$, since $\eta_i \stackrel{i.i.d.}{\sim} \mathcal{N}(1,1)$, conditionally on $\data$ we have
    \begin{equation} 
        \begin{aligned}
            &\frac{1}{n}\sum\limits_{i=1}^n (\eta_i-1) (U^\T\St^{-1/2} X_i)_j (U^\T\St^{-1/2} X_i)_k \sim
            \\&\qquad\qquad \sim \mathcal{N}\left( 0, \frac{1}{n^2}\sum\limits_{i=1}^n (U^\T\St^{-1/2} X_i)_j^2 (U^\T\St^{-1/2} X_i)_k^2\right).
        \nonumber 
        \end{aligned} 
    \end{equation}
    Consider the event $\Omega$ defined as
    \begin{equation} 
        \begin{aligned}
            \left\{ \max\limits_{j,k\in[d]} \frac{1}{n}\sum\limits_{i=1}^n (U^\T\St^{-1/2} X_i)_j^2 (U^\T\St^{-1/2} X_i)_k^2 \leq \sigma^2 \right\} 
        \nonumber 
        \end{aligned} 
    \end{equation}
    with $\sigma \eqdef c^2\left(\log(n) + \log(2d^2) \right)^{2/\beta}$.
    
    Let us verify that $\Prob(\Omega) \geq 1-1/n$.
    By Proposition~\ref{L:KC_Prod} and Assumption~\ref{A: tails},
    \begin{equation} 
        \begin{aligned}
            &\left\| \frac{1}{n}\sum\limits_{i=1}^n (U^\T\St^{-1/2} X_i)_j^2\, (U^\T\St^{-1/2} X_i)_k^2 \right\|_{\psi_{\beta/4}} 
            \leq 
            \left\| (U^\T\St^{-1/2} X_1)_j^2\, (U^\T\St^{-1/2} X_1)_k^2 \right\|_{\psi_{\beta/4}}\\
            &\hspace{1cm} \leq
            \left\| (U^\T\St^{-1/2} X_1)_j \right\|_{\psi_\beta}^2 \left\| (U^\T\St^{-1/2} X_1)_k \right\|_{\psi_\beta}^2
            \leq c^4 < \infty,
        \nonumber 
        \end{aligned} 
    \end{equation}
    yielding
    \begin{equation} 
        \begin{aligned}
            &\Prob\left(  \frac{1}{n}\sum\limits_{i=1}^n (U^\T\St^{-1/2} X_i)_j^2\, (U^\T\St^{-1/2} X_i)_k^2  \geq \sigma^2 \right)
            \leq 
            2\exp\left(-(\sigma/c^2)^{\beta/2} \right),
        \nonumber 
        \end{aligned} 
    \end{equation}
    By union bound,
    \begin{equation} 
        \begin{aligned}
            &\Prob\left( \max\limits_{j,k\in[d]}  \frac{1}{n}\sum\limits_{i=1}^n (U^\T\St^{-1/2} X_i)_j^2\, (U^\T\St^{-1/2} X_i)_k^2 \geq \sigma^2 \right)
            \leq 
            2d^2\exp\left(-(\sigma/c^2)^{\beta/2}\right),
        \nonumber 
        \end{aligned} 
    \end{equation}
    which, plugging $\sigma$ in and using definition of $\Omega$ can be rewritten as $\Prob(\Omega^{c}) \leq 1/n$.
    
    Further, on $\Omega$ we have for all $j,k\in[d]$ Gaussian tail inequality
    \begin{equation} 
        \begin{aligned}
            \Prob\left( \left| \frac{1}{n}\sum\limits_{i=1}^n (\eta_i-1) (U^\T\St^{-1/2} X_i)_j (U^\T\St^{-1/2} X_i)_k \right| \geq z \;\Bigg|\; \data\right) \leq 2e^{-nz^2/\sigma^2},
        \nonumber 
        \end{aligned} 
    \end{equation}
    thus, by union bound
    \begin{equation} 
        \begin{aligned}
            \Prob\left( \max\limits_{j,k\in[d]} \left| \frac{1}{n}\sum\limits_{i=1}^n (\eta_i-1) (U^\T\St^{-1/2} X_i)_j (U^\T\St^{-1/2} X_i)_k \right| \geq z \;\Bigg|\; \data\right) \leq 2d^2 e^{-nz^2/\sigma^2}.
        \nonumber 
        \end{aligned} 
    \end{equation}
    We conclude the proof by taking $z = \sqrt{\frac{\sigma^2(\log(n) + \log(2d^2))}{n}}$.
\end{proof}

\begin{proof}[Proof of Lemma \ref{L:CovBound}]
    Fix arbitrary $i \in [n], j\in [p]$. Let us bound $\| x_{ij} \|_{\psi_{\beta/2}}$. From \eqref{Eq:bound_x_ij}
    \begin{equation}
        \begin{aligned}
        	 |x_{ij}| \leq \chigh\cdot \overline{v}_{j}^{\T} \St^{-1/2} X_i \cdot \overline{w}_{j}^{\T} \St^{-1/2}X_i,
            \nonumber
        \end{aligned}
    \end{equation}
    where $\overline{v}_{j}, \overline{w}_{j} \in S^{d-1}$.
    Hence,
    \begin{equation}
        \begin{aligned}
        	\| x_{ij} \|_{\psi_{\beta/2}} &\leq
        	\chigh\, \| \overline{v}_{j}^{\T} \St^{-1/2} X_i \|_{\psi_\beta}
        	\| \overline{w}_{j}^{\T} \St^{-1/2} X_i \|_{\psi_\beta}
        	\leq \chigh c^2,
            \nonumber
        \end{aligned}
    \end{equation}
    where we used Proposition~\ref{L:KC_Prod} and Assumption~\ref{A: tails}.
    Now the claim follows from Theorem~\ref{Th:KC_SCmax} with $X_i = x_i$, $\beta/2$ taken as $\beta$, $K_{n,p} = \chigh c^2$ and $z = \log(3n)$.
\end{proof}

\begin{proof}[Proof of Lemma~\ref{L:BayesConc}]
    The idea is similar to the proof of Lemma~\ref{L:BootConc}.
    Fix $s, t \in [q]$. Expanding squared Frobenius norm over the basis of eigenvectors $\{ u_j \}_{j=1}^d$ and using the definition of $\Sf$, we have
    \begin{equation} 
        \begin{aligned}
            &\| \Pts (\Sf-\Se) \Ptt \|_{\Fr}^2 
            = \sum\limits_{j,k=1}^d \left( u_j^\T \Pts(\Sf-\Se)\Ptt u_k \right)^2
            =\sum\limits_{j\in\mathcal{I}_s}\sum\limits_{k\in\mathcal{I}_t} \left( u_j^\T (\Sf-\Se) u_k \right)^2\\
            &\hspace{1cm} \leq m_s m_t \max\limits_{\substack{j\in \mathcal{I}_s\\ k\in\mathcal{I}_t}} \left( u_j^\T (\Sf-\Se) u_k \right)^2 = m_s \mu_s m_t \mu_t \max\limits_{\substack{j\in \mathcal{I}_s\\ k\in\mathcal{I}_t}} \left( u_j^\T \St^{-1/2} (\Sf-\Se)\St^{-1/2} u_k \right)^2.
        \nonumber 
        \end{aligned} 
    \end{equation}
    Hence,
    \begin{equation} 
        \begin{aligned}
            \max\limits_{s,t\in[q]}\frac{\| \Pts (\Sf-\Se) \Ptt \|_{\Fr}}{\sqrt{m_s \mu_s m_t \mu_t}}
            &\leq \max\limits_{s,t\in[q]}\max\limits_{\substack{j\in \mathcal{I}_s\\ k\in\mathcal{I}_t}} \left| u_j^\T \St^{-1/2} (\Sf-\Se)\St^{-1/2} u_k \right|\\
            &= \max\limits_{j,k\in[d]} \left| u_j^\T \St^{-1/2} (\Sf-\Se)\St^{-1/2} u_k \right|.
        \nonumber 
        \end{aligned} 
    \end{equation}
    For arbitrary $j,k\in[d]$, by the definition of $\Sf$,  we have
    \begin{equation} 
        \begin{aligned}
           & u_j^\T \St^{-1/2} (\Sf-\Se)\St^{-1/2} u_k =
            \\&\qquad=
            \frac{1}{n}\sum\limits_{i=1}^n \left\{ (u_j^\T \St^{-1/2} Z_i) (u_k^\T \St^{-1/2} Z_i) - \E\left[ (u_j^\T \St^{-1/2} Z_i) (u_k^\T \St^{-1/2} Z_i) \right] \right\},
        \nonumber 
        \end{aligned} 
    \end{equation}
    where, conditionally on $\data$,
    \begin{equation} 
        \begin{aligned}
            u_j^\T \St^{-1/2} Z_i \sim \mathcal{N}(0,u_j^\T \St^{-1/2}\Se \St^{-1/2} u_j),\\
            u_k^\T \St^{-1/2} Z_i \sim \mathcal{N}(0,u_k^\T \St^{-1/2}\Se \St^{-1/2} u_k),
        \nonumber 
        \end{aligned} 
    \end{equation}
    since $Z_i \stackrel{i.i.d.}{\sim} \mathcal{N}_d(0, \Se)$ given $\data$.
    Consider the event $\Omega$ defined as
    \begin{equation} 
        \begin{aligned}
            \left\{ \max\limits_{j\in[d]} \,u_j^\T \St^{-1/2}\Se \St^{-1/2} u_j \leq \sigma^2 \right\} 
        \nonumber 
        \end{aligned} 
    \end{equation}
    with $\sigma \eqdef c\left((\log(n) + \log(2d)) \right)^{1/\beta}$.
    
    Let us verify that $\Prob(\Omega) \geq 1-1/n$.
    Fix $j\in[d]$.
    By Proposition~\ref{L:KC_Prod} and Assumption~\ref{A: tails},
    \begin{equation} 
        \begin{aligned}
            &\left\| u_j^\T \St^{-1/2}\Se \St^{-1/2} u_j \right\|_{\psi_{\beta/2}} 
            =
            \left\| \frac{1}{n}\sum\limits_{i=1}^n (u_j^\T \St^{-1/2} X_i)^2 \right\|_{\psi_{\beta/2}}\\
            &\hspace{1cm} \leq
            \left\| (u_j^\T \St^{-1/2} X_1)^2 \right\|_{\psi_{\beta/2}}
            \leq \left\| u_j^\T \St^{-1/2} X_1 \right\|_{\psi_\beta}^2
            \leq c^2 < \infty,
        \nonumber 
        \end{aligned} 
    \end{equation}
    yielding
    \begin{equation} 
        \begin{aligned}
            &\Prob\left(  u_j^\T \St^{-1/2}\Se \St^{-1/2} u_j \geq \sigma^2 \right)
            \leq 
            2\exp\left(-(\sigma/c)^\beta\right).
        \nonumber 
        \end{aligned} 
    \end{equation}
    By union bound,
    \begin{equation} 
        \begin{aligned}
            &\Prob\left( \max\limits_{j\in[d]} \,u_j^\T \St^{-1/2}\Se \St^{-1/2} u_j \geq \sigma^2 \right)
            \leq 
            2d\exp\left(-(\sigma/c)^\beta\right),
        \nonumber 
        \end{aligned} 
    \end{equation}
    which, plugging $\sigma$ in and using definition of $\Omega$, can be rewritten as $\Prob(\Omega^{c}) \leq 1/n$.
    
    Further, on $\Omega$ for arbitrary $j\in[d]$, $(u_j^\T \St^{-1/2} Z_i)$ is $\sigma^2$-subgaussian, and 
    $\| u_j^\T \St^{-1/2} Z_i \|_{\psi_2} \leq C\sigma$. Hence,
    for arbitrary $j,k\in[d]$
    \begin{equation} 
        \begin{aligned}
            \| (u_j^\T \St^{-1/2} Z_i)(u_k^\T \St^{-1/2} Z_i) \|_{\psi_1} \leq \| u_j^\T \St^{-1/2} Z_i \|_{\psi_2}\| u_k^\T \St^{-1/2} Z_i \|_{\psi_2}
            \leq C\sigma^2
        \nonumber 
        \end{aligned} 
    \end{equation}
    due to Proposition~\ref{L:KC_Prod}. So, $ (u_j^\T \St^{-1/2} Z_i)(u_k^\T \St^{-1/2} Z_i)$ is subexponential and by, for instance, Exercise~2.7.10 of \cite{Vershynin} the centered version is also subexponential (but with different multiplicative constant factor in the Orlicz norm):
    \begin{equation} 
        \begin{aligned}
            \left\| (u_j^\T \St^{-1/2} Z_i)(u_k^\T \St^{-1/2} Z_i) - \E\left[(u_j^\T \St^{-1/2} Z_i)(u_k^\T \St^{-1/2} Z_i) \right] \right\|_{\psi_1} \
            \leq C\sigma^2.
        \nonumber 
        \end{aligned} 
    \end{equation}
    Further, by Bernstein inequality (e.g. Corollary 2.8.3 of \cite{Vershynin}) on $\Omega$
    \begin{equation} 
        \begin{aligned}
            &\Prob\left( \left| \frac{1}{n}\sum\limits_{i=1}^n \left\{ (u_j^\T \St^{-1/2} Z_i) (u_k^\T \St^{-1/2} Z_i) - \E\left[ (u_j^\T \St^{-1/2} Z_i) (u_k^\T \St^{-1/2} Z_i) \right] \right\} \right| \geq z \,\Bigg|\,\data\right) \leq \\
            &\hspace{1cm}\leq 2\exp\left( -C n \left\{\frac{z^2}{\sigma^4} \land \frac{z}{\sigma^2} \right\} \right),
        \nonumber 
        \end{aligned} 
    \end{equation}
    and by union bound
    \begin{equation} 
        \begin{aligned}
            &\Prob\left( \max\limits_{j,k\in[d]}\left| u_j^\T \St^{-1/2} (\Sf-\Se)\St^{-1/2} u_k \right| \geq z \,\Big|\,\data\right) \leq  2d^2\exp\left( -C n \left\{\frac{z^2}{\sigma^4} \land \frac{z}{\sigma^2} \right\} \right).
        \nonumber 
        \end{aligned} 
    \end{equation}
    We conclude the proof by taking $z = C \sigma^2 \sqrt{\frac{\log(n) + \log(2d^2)}{n}}$ (assuming $\log(d)/n < 1$).
\end{proof}

\begin{proof}[Proof of Lemma~\ref{L:FCovBound}]
    Fix $j,k \in [p]$. We have to bound
    \begin{equation}
        \begin{aligned}
        	&\left[\Cov(Y^F\,|\,\data) - \Cov(Y) \right]_{j,k} = \\
        	&\qquad = 
        	 \sum\limits_{\substack{r\in\J\\r^{\prime} \in \J}} \sum\limits_{\substack{s\notin\J\\s^{\prime} \notin \J}} \frac{(v_j^{\T} \Ptr \Se \Ptrp v_k) \cdot (w_j^{\T} \Pts \Se \Ptsp w_k) - (v_j^{\T} \Ptr \Se \Pts w_j) \cdot (v_k^{\T} \Ptrp \Se \Ptsp w_k)}{(\mu_r-\mu_s)(\mu_{r^{\prime}} - \mu_{s^{\prime}})}
        	-\\
        	&\qquad\qquad -
        	\sum\limits_{\substack{r\in\J\\r^{\prime} \in \J}} \sum\limits_{\substack{s\notin\J\\s^{\prime} \notin \J}} \frac{(v_j^{\T} \Ptr \St \Ptrp v_k) \cdot (w_j^{\T} \Pts \St \Ptsp w_k) - (v_j^{\T} \Ptr \St \Pts w_j) \cdot (v_k^{\T} \Ptrp \St \Ptsp w_k)}{(\mu_r-\mu_s)(\mu_{r^{\prime}} - \mu_{s^{\prime}})}.
            \nonumber
        \end{aligned}
    \end{equation}
    To simplify the expression, define auxiliary matrices
    \begin{equation}
        \begin{aligned}
        	&B_j = \sum\limits_{r\in \J}\sum\limits_{s\notin\J} \frac{\Ptr v_j w_j^\T \Pts}{\mu_r - \mu_s}  \in \R^{d\times d},\\
        	&B_k = \sum\limits_{r\in \J}\sum\limits_{s\notin\J} \frac{\Ptr v_k w_k^\T \Pts}{\mu_r - \mu_s}  \in \R^{d\times d}.
            \nonumber
        \end{aligned}
    \end{equation}
    Using cyclic property of the trace, we obtain
    \begin{equation}
        \begin{aligned}
        	&\sum\limits_{\substack{r\in\J\\r^{\prime} \in \J}} \sum\limits_{\substack{s\notin\J\\s^{\prime} \notin \J}} \frac{(v_j^{\T} \Ptr \Se \Ptrp v_k) \cdot (w_j^{\T} \Pts \Se \Ptsp w_k)}{(\mu_r-\mu_s)(\mu_{r^{\prime}} - \mu_{s^{\prime}})}
        	=
        	\sum\limits_{\substack{r\in\J\\r^{\prime} \in \J}} \sum\limits_{\substack{s\notin\J\\s^{\prime} \notin \J}} \frac{v_j^{\T} \Ptr \Se \Ptrp v_k w_k^{\T} \Ptsp \Se \Pts w_j}{(\mu_r-\mu_s)(\mu_{r^{\prime}} - \mu_{s^{\prime}})}
        	\\&\qquad=
        	\Tr\left[ \sum\limits_{\substack{r\in\J\\r^{\prime} \in \J}} \sum\limits_{\substack{s\notin\J\\s^{\prime} \notin \J}} \frac{v_j^{\T} \Ptr \Se \Ptrp v_k w_k^{\T} \Ptsp \Se \Pts w_j}{(\mu_r-\mu_s)(\mu_{r^{\prime}} - \mu_{s^{\prime}})} \right]
        	\\&\qquad=
        	\Tr\left[ \sum\limits_{\substack{r\in\J\\r^{\prime} \in \J}} \sum\limits_{\substack{s\notin\J\\s^{\prime} \notin \J}} \frac{\Se \Ptrp v_k w_k^{\T} \Ptsp \Se \Pts w_j v_j^{\T} \Ptr }{(\mu_r-\mu_s)(\mu_{r^{\prime}} - \mu_{s^{\prime}})} \right]
        	\\&\qquad=
        	\Tr\left[ \Se \left( \sum\limits_{r^\prime\in \J}\sum\limits_{s^\prime\notin\J} \frac{\Ptrp v_k w_k^\T \Ptsp}{\mu_{r^\prime} - \mu_{s^\prime}} \right) \Se \left( \sum\limits_{r\in \J}\sum\limits_{s\notin\J} \frac{\Pts w_j v_j^\T \Ptr}{\mu_r - \mu_s} \right) \right] = \Tr\left[ \Se B_k \Se B_j^\T \right].
            \nonumber
        \end{aligned}
    \end{equation}
    Similarly,
    \begin{equation}
        \begin{aligned}
        	&
        	\sum\limits_{\substack{r\in\J\\r^{\prime} \in \J}} \sum\limits_{\substack{s\notin\J\\s^{\prime} \notin \J}} \frac{(v_j^{\T} \Ptr \St \Ptrp v_k) \cdot (w_j^{\T} \Pts \St \Ptsp w_k)}{(\mu_r-\mu_s)(\mu_{r^{\prime}} - \mu_{s^{\prime}})}
        	= \Tr\left[ \St B_k \St B_j^\T\right].
            \nonumber
        \end{aligned}
    \end{equation}
    In a slightly different fashion, we derive
    \begin{equation}
        \begin{aligned}
        	&\sum\limits_{\substack{r\in\J\\r^{\prime} \in \J}} \sum\limits_{\substack{s\notin\J\\s^{\prime} \notin \J}} \frac{(v_j^{\T} \Ptr \Se \Pts w_j) \cdot (v_k^{\T} \Ptrp \Se \Ptsp w_k)}{(\mu_r-\mu_s)(\mu_{r^{\prime}} - \mu_{s^{\prime}})}
        	\\&\qquad= \sum\limits_{r\in\J} \sum\limits_{s\notin\J} \frac{v_j^{\T} \Ptr \Se \Pts w_j}{\mu_r-\mu_s} \cdot
        	\sum\limits_{r^{\prime} \in \J} \sum\limits_{s^{\prime} \notin \J} \frac{v_k^{\T} \Ptrp \Se \Ptsp w_k}{\mu_{r^{\prime}} - \mu_{s^{\prime}}}
        	\\&\qquad=
        	\Tr\left[\sum\limits_{r\in\J} \sum\limits_{s\notin\J} \frac{v_j^{\T} \Ptr \Se \Pts w_j}{\mu_r-\mu_s}\right]\cdot
        	\Tr\left[\sum\limits_{r^{\prime} \in \J} \sum\limits_{s^{\prime} \notin \J} \frac{v_k^{\T} \Ptrp \Se \Ptsp w_k}{\mu_{r^{\prime}} - \mu_{s^{\prime}}}\right]
        	\\&\qquad=
        	\Tr\left[\sum\limits_{r\in\J} \sum\limits_{s\notin\J} \frac{\Se \Pts w_j v_j^{\T} \Ptr }{\mu_r-\mu_s}\right]\cdot
        	\Tr\left[\sum\limits_{r^{\prime} \in \J} \sum\limits_{s^{\prime} \notin \J} \frac{ \Se \Ptsp w_k v_k^{\T} \Ptrp}{\mu_{r^{\prime}} - \mu_{s^{\prime}}}\right]
        	\\&\qquad=
        	\Tr\left[ \Se\left(\sum\limits_{r\in\J} \sum\limits_{s\notin\J} \frac{ \Pts w_j v_j^{\T} \Ptr }{\mu_r-\mu_s}\right)\right]\cdot
        	\Tr\left[ \Se\left(\sum\limits_{r^{\prime} \in \J} \sum\limits_{s^{\prime} \notin \J} \frac{  \Ptsp w_k v_k^{\T} \Ptrp}{\mu_{r^{\prime}} - \mu_{s^{\prime}}}\right)\right]
        	\\&\qquad = \Tr\left[\Se B_j^\T \right] \cdot \Tr\left[ \Se B_k^\T \right].
            \nonumber
        \end{aligned}
    \end{equation}
    Similarly,
    \begin{equation}
        \begin{aligned}
        	&
        	\sum\limits_{\substack{r\in\J\\r^{\prime} \in \J}} \sum\limits_{\substack{s\notin\J\\s^{\prime} \notin \J}} \frac{(v_j^{\T} \Ptr \St \Pts w_j) \cdot (v_k^{\T} \Ptrp \St \Ptsp w_k)}{(\mu_r-\mu_s)(\mu_{r^{\prime}} - \mu_{s^{\prime}})}
        	= \Tr\left[\St B_j^\T \right] \cdot \Tr\left[ \St B_k^\T \right].
            \nonumber
        \end{aligned}
    \end{equation}
    Hence, our expression reduces to
    \begin{equation}
        \begin{aligned}
        	&\left[\Cov(Y^F\,|\,\data) - \Cov(Y) \right]_{j,k} =
        	\\&\qquad=
        	\Tr\left[ \Se B_k \Se B_j^\T \right]
        	- \Tr\left[ \St B_k \St B_j^\T \right]
        	- \Tr\left[\Se B_j^\T \right] \cdot \Tr\left[ \Se B_k^\T \right]
        	+ \Tr\left[\St B_j^\T \right] \cdot \Tr\left[ \St B_k^\T \right].
            \nonumber
        \end{aligned}
    \end{equation}
    Note that actually $\Tr\left[\St B_j^\T \right] = \Tr\left[\St B_k^\T \right] = 0$, so
    \begin{equation}
        \begin{aligned}
        	&\left[\Cov(Y^F\,|\,\data) - \Cov(Y) \right]_{j,k} =
        	\\&\qquad=
        	\Tr\left[ \Se B_k \Se B_j^\T \right]
        	- \Tr\left[ \St B_k \St B_j^\T \right]
        	- \Tr\left[(\Se-\St) B_j^\T \right] \cdot \Tr\left[ (\Se-\St) B_k^\T \right].
            \nonumber
        \end{aligned}
    \end{equation}
    Adding and subtracting $\Tr\left[ \St B_k \Se B_j^\T \right]$, we get
    \begin{equation}
        \begin{aligned}
        	&\left[\Cov(Y^F\,|\,\data) - \Cov(Y) \right]_{j,k} = \\
        	&\qquad =
        	\Tr\left[ (\Se-\St) B_k \Se B_j^\T \right]
        	+ \Tr\left[ \St B_k (\Se-\St) B_j^\T \right] - \Tr\left[ (\Se-\St) B_j^\T \right] \cdot \Tr\left[ (\Se-\St) B_k^\T \right] \\
        	&\qquad =
        	\Tr\left[ (\Se-\St) B_k (\Se-\St) B_j^\T \right]
        	+ \Tr\left[ (\Se-\St) B_k \St B_j^\T \right]
        	+ \Tr\left[ (\Se-\St) B_j^\T \St B_k \right] -\\
        	&\qquad\qquad
        	- \Tr\left[ (\Se-\St) B_j^\T \right] \cdot \Tr\left[ (\Se-\St) B_k^\T \right].
            \nonumber
        \end{aligned}
    \end{equation}
    It is easy to see, that if we define
    \begin{equation}
        \begin{aligned}
        	&\widetilde{B}_j = \sum\limits_{r\in \J}\sum\limits_{s\notin\J} \frac{\sqrt{\mu_r \mu_s}}{\mu_r - \mu_s} \, \Ptr v_j w_j^\T \Pts \in \R^{d\times d},\\
        	&\widetilde{B}_k = \sum\limits_{r\in \J}\sum\limits_{s\notin\J} \frac{\sqrt{\mu_r \mu_s}}{\mu_r - \mu_s} \, \Ptr v_k w_k^\T \Pts \in \R^{d\times d},
            \nonumber
        \end{aligned}
    \end{equation}
    then we can rewrite
    \begin{equation}
        \begin{aligned}
        	&\left[\Cov(Y^F\,|\,\data) - \Cov(Y) \right]_{j,k} = \\
        	&\qquad =
        	\Tr\left[ (\St^{-1/2}\Se\St^{-1/2}-\Id_d) \widetilde{B}_k (\St^{-1/2}\Se\St^{-1/2}-\Id_d) \widetilde{B}_j^\T \right] +
        	\\&\qquad\qquad
        	+ \Tr\left[ (\St^{-1/2}\Se\St^{-1/2}-\Id_d) (\widetilde{B}_k \widetilde{B}_j^\T + \widetilde{B}_j^\T\widetilde{B}_k) \right] -\\
        	&\qquad\qquad
        	- \Tr\left[ (\St^{-1/2}\Se\St^{-1/2}-\Id_d) \widetilde{B}_j^\T \right] \cdot \Tr\left[ (\St^{-1/2}\Se\St^{-1/2}-\Id_d) \widetilde{B}_k^\T \right]
        	=: T_1 + T_2 - T_3.
            \nonumber
        \end{aligned}
    \end{equation}
    It suffices to bound each $|T_1|, |T_2|, |T_3|$ with probability $1-1/(3np^2)$ to get the bound on $\left| \left[\Cov(Y^F\,|\,\data) - \Cov(Y) \right]_{j,k} \right|$ with probability $1-1/(np^2)$. Before we do so, let us bound some important quantities. 
    To slightly simplify some expressions, introduce for all $r\in \J$
    \begin{equation}
        \begin{aligned}
        	&\overline{w}_{j,r} \eqdef \sum\limits_{s\notin\J} \frac{\sqrt{\mu_r \mu_s}}{\mu_r - \mu_s}\Pts w_j,\\
        	&\overline{w}_{k,r} \eqdef \sum\limits_{s\notin\J} \frac{\sqrt{\mu_r \mu_s}}{\mu_r - \mu_s} \Pts w_k,
            \nonumber
        \end{aligned}
    \end{equation}
    so that 
    \begin{equation}
        \begin{aligned}
        	 \widetilde{B}_j =\sum\limits_{r\in\J} \Ptr v_j \overline{w}_{j,r}^\T,\\
        	 \widetilde{B}_k =\sum\limits_{r\in\J} \Ptr v_k \overline{w}_{k,r}^\T.
            \nonumber
        \end{aligned}
    \end{equation}
    Note that $\| \overline{w}_{j,r} \| \leq \chigh$.
    We have
    \begin{equation}
        \begin{aligned}
        	\| B_j \|_* = \| B_j^\T \|_* &= 
        	\left\| \sum\limits_{r\in\J} \Ptr v_j \overline{w}_{j,r}^\T\right\|_*
        	\leq
        	\sum\limits_{r\in\J} \| \Ptr v_j \overline{w}_{j,r}^\T \|_*
        	=
        	\sum\limits_{r\in\J} \| \Ptr v_j \| \, \|\overline{w}_{j,r} \|\\
        	&\leq
        	\chigh \sum\limits_{r\in\J} \| \Ptr v_j \|
        	\leq \sqrt{|\J|} \chigh \sqrt{\sum\limits_{r\in\J} \| \Ptr v_j \|^2} = \sqrt{|\J|}\; \chigh.
            \nonumber
        \end{aligned}
    \end{equation}
    Similarly, $\| B_k \|_* = \| B_k^\T \|_*\leq \sqrt{|\J|}\;\chigh $.
    Moreover,
    \begin{equation}
        \begin{aligned}
        	&\| \widetilde{B}_k \widetilde{B}_j^\T + \widetilde{B}_j^\T \widetilde{B}_k\|_* \leq 2\| B_k \|_* \| B_j^\T \|_* \leq 2 |\J| \,\chigh^2.
            \nonumber
        \end{aligned}
    \end{equation}
    Note also that $\rank(\widetilde{B}_k) \leq |\J|$,  $\rank(\widetilde{B}_j) \leq |\J|$ and $\rank(\widetilde{B}_k \widetilde{B}_j^\T + \widetilde{B}_j^\T\widetilde{B}_k) \leq 2|\J|$.
    Now we want to show, that for any matrix $D \in \R^{d\times d}$ of rank $h\in[d]$ it holds
    \begin{equation}
        \begin{aligned}
        	&\left| \Tr\left[ (\St^{-1/2}\Se\St^{-1/2}-\Id_d) D \right] \right| \leq\\
        	&\qquad\leq\| D \|_* \cdot C_\beta c^2 \left( \sqrt{\frac{\log(n) + \log(h)}{n}} + \frac{(\log(n))^{1/\beta} (\log(n) + \log(h))^{2/\beta}}{n} \right)
            \label{Tr_conc}
        \end{aligned}
    \end{equation}
    with probability $1-1/(3np^2)$. Indeed, let $D = \sum\limits_{l=1}^h \sigma_l(D)\, a_l b_l^\T$ be SVD of $D$ with singular values $\sigma_l(D)$ and left and right singular vectors $a_l, b_l\in\Sph^{d-1},\,l\in[h]$. Then
    \begin{equation}
        \begin{aligned}
        	&\left| \Tr\left[ (\St^{-1/2}\Se\St^{-1/2}-\Id_d) D \right] \right| =
        	\left| \Tr\left[ (\St^{-1/2}\Se\St^{-1/2}-\Id_d) \sum\limits_{l=1}^h \sigma_l(D) a_l b_l^\T \right] \right|
        	\\&\qquad=
        	\left| \sum\limits_{l=1}^h \sigma_l(D)\; b_l^\T (\St^{-1/2}\Se\St^{-1/2}-\Id_d) a_l\right| \leq
        	\max\limits_{l\in[h]} \left|b_l^\T (\St^{-1/2}\Se\St^{-1/2}-\Id_d) a_l\right|\; \sum\limits_{l=1}^h \sigma_l(D)
        	\\&\qquad = \| D \|_* \cdot \max\limits_{l\in[h]} \left|b_l^\T (\St^{-1/2}\Se\St^{-1/2}-\Id_d) a_l\right| .
            \nonumber
        \end{aligned}
    \end{equation}
    Now \eqref{Tr_conc} follows from Theorem~\ref{Th:KC_SCmax} applied with $2h$ instead of $p$, $K_{n,p} = c^{-1/\beta}$, $X_i(2l-1) = a_l^\T \St^{-1/2} X_i$, $X_i(2l) = b_l^\T \St^{-1/2} X_i$, $l\in[h], i\in[n]$ and $z = \log(9np^2)$.
    
    Applying \eqref{Tr_conc} to $|T_2|$ gives with probability $1-1/(3np^2)$
    \begin{equation}
        \begin{aligned}
        	& |T_2| \leq   \chigh^2 \nu_n
            \nonumber
        \end{aligned}
    \end{equation}
    where
    \begin{equation}
        \begin{aligned}
        	&\nu_n \eqdef  C_\beta c^2 \left( \sqrt{\frac{\log(np) + \log(|\J|)}{n}} + \frac{(\log(n))^{1/\beta} (\log(np) + \log(|\J|))^{2/\beta}}{n} \right).
            \nonumber
        \end{aligned}
    \end{equation}
    Similarly, with probability $1-1/(3np^2)$
    \begin{equation}
        \begin{aligned}
        	&|T_3| \leq  \chigh^2 \nu_n^2.
            \nonumber
        \end{aligned}
    \end{equation}
    Finally, $|T_1|$ can be bounded in the same way. Let 
    \begin{equation}
        \begin{aligned}
            &\widetilde{B}_k = \sum\limits_{l=1}^{|\J|} \sigma_l(B_k)\, a_l^{(k)} {b_l^{(k)}}^\T,\\
        	&\widetilde{B}_j^\T = \sum\limits_{l=1}^{|\J|} \sigma_l(B_j^\T)\, a_l^{(j)} {b_l^{(j)}}^\T
            \nonumber
        \end{aligned}
    \end{equation}
    be SVD of $\widetilde{B}_k$ and $\widetilde{B}_j^\T$. Then
    \begin{equation}
        \begin{aligned}
        	|T_1| &= 
        	\left| \Tr\left[ (\St^{-1/2}\Se\St^{-1/2}-\Id_d)\widetilde{B}_k (\St^{-1/2}\Se\St^{-1/2}-\Id_d)\widetilde{B}_j^\T \right]\right|\\
        	&=
        	\Bigg| \Tr\Bigg[ (\St^{-1/2}\Se\St^{-1/2}-\Id_d) \left( \sum\limits_{l=1}^{|\J|} \sigma_l(\widetilde{B}_k)\, a_l^{(k)} {b_l^{(k)}}^\T \right) \times\\
        	&\qquad\times (\St^{-1/2}\Se\St^{-1/2}-\Id_d) \left(\sum\limits_{l=1}^{|\J|} \sigma_l(\widetilde{B}_j^\T)\, a_l^{(j)} {b_l^{(j)}}^\T \right) \Bigg]\Bigg|\\
        	&\leq \sum\limits_{l_1=1}^{|\J|} \sum\limits_{l_2=1}^{|\J|}
        	\sigma_{l_1}(\widetilde{B}_k) \sigma_{l_2}(\widetilde{B}_j^\T) \cdot \left|{b_{l_2}^{(j)}}^\T (\St^{-1/2}\Se\St^{-1/2}-\Id_d)a_{l_1}^{(k)}\right| \times \\
        	&\qquad\qquad \times\left|{b_{l_1}^{(k)}}^\T (\St^{-1/2}\Se\St^{-1/2}-\Id_d)a_{l_2}^{(j)}\right|\\
        	&\leq 
        	\| \widetilde{B}_k \|_* \| \widetilde{B}_j^\T \|_* \max\limits_{l_1, l_2 \in [|\J|]} \left|{b_{l_2}^{(j)}}^\T (\St^{-1/2}\Se\St^{-1/2}-\Id_d)a_{l_1}^{(k)}\right|  \times
        	\\&\qquad\qquad\times \left|{b_{l_1}^{(k)}}^\T (\St^{-1/2}\Se\St^{-1/2}-\Id_d)a_{l_2}^{(j)}\right|.
            \nonumber
        \end{aligned}
    \end{equation}
    By yet another application of Theorem~\ref{Th:KC_SCmax}
    \begin{equation}
        \begin{aligned}
        	|T_1| &\leq |\J|\chigh^2 \nu_n^2
            \nonumber
        \end{aligned}
    \end{equation}
    with probability $1-1/(3np^2)$. Putting all the bounds together, we derive
    \begin{equation}
        \begin{aligned}
        	\left|\left[\Cov(Y^F\,|\,\data) - \Cov(Y) \right]_{j,k}\right| &\leq |\J| \chigh^2 (\nu_n + \nu_n^2)
            \nonumber
        \end{aligned}
    \end{equation}
    with probability $1-1/(np^2)$. Union bound concludes the proof:
    \begin{equation}
        \begin{aligned}
        	\max\limits_{j,k\in[p]} \left|\left[\Cov(Y^F\,|\,\data) - \Cov(Y) \right]_{j,k}\right| &\leq |\J| \chigh^2 (\nu_n + \nu_n^2)
            \nonumber
        \end{aligned}
    \end{equation}
    with probability $1-1/n$.
\end{proof}

\begin{proof}[Proof of Lemma~\ref{L:2approx}]
We first construct proper $\Gamma^*$. Recall that $\overline{\Gamma}_1$ and $\overline{\Gamma}_2$ which satisfy $\overline{\Gamma}_1\overline{\Gamma}_1^\T = \Pbar$,
$\overline{\Gamma}_1^\T\overline{\Gamma}_1 = \Id_m$,
$\overline{\Gamma}_2\overline{\Gamma}_2^\T = \Id_d - \Pbar$ and
$\overline{\Gamma}_2^\T\overline{\Gamma}_2 = \Id_{d-m}$ are fixed at the beginning of our procedure.
At the same time, by Davis-Kahan theorem (e.g. Theorem 2 from \cite{DavisKahan}), there exist ${\Gamma}_1^*$ and ${\Gamma}_2^*$ such that ${\Gamma}_1^*{\Gamma^*_1}^\T = \Pstar$,
${\Gamma^*_1}^\T{\Gamma}_1^* = \Id_m$,
${\Gamma}_2^*{\Gamma_2^*}^\T = \Id_d - \Pstar$ and
${\Gamma_2^*}^\T{\Gamma}_2^* = \Id_{d-m}$, and 
\begin{equation}
        \begin{aligned}
        	&\| \Gamma_1^* - \overline{\Gamma}_1 \|_{\Fr} \leq 2^{3/2} \| \Pbar - \Pstar \|_{\Fr},\;\;\;
        	\| \Gamma_2^* - \overline{\Gamma}_2 \|_{\Fr} \leq 2^{3/2} \| \Pbar - \Pstar \|_{\Fr}.
        \nonumber
        \end{aligned}
    \end{equation}
    
 Then, as in Lemma~\ref{L:KL}, denoting the linear parts of $\Pea - \Pta$ and $\Peb-\Ptb$ as  $L_a(\Se_a-\St_a)$ and $L_b(\Se_b-\St_b)$ and the remainder terms as $R_a(\Se_a-\St_a)$ and $R_b(\Se_b-\St_b)$, we decompose
    \begin{equation}
        \begin{aligned}
            \Pea - \Peb &
            = (\Pea - \Pstar) - (\Peb - \Pstar)
            = (\Pea - \Pta) - (\Peb - \Ptb)
            \\&= L_a(\Se_a-\St_a) + R_a(\Se_a-\St_a)
            - L_b(\Se_b-\St_b) - R_b(\Se_b-\St_b)
            .
        \nonumber
        \end{aligned}
    \end{equation}
 We first state some auxiliary bounds. \\
 \underline{Bounds on the linear parts:} Let $\widehat{x}^a$, $\widehat{x}^b$ $\overline{x}^a$ and $\overline{x}^b$  be the same quantities as $x$ in Lemma~\ref{L:KL}, but now for $(\Se_a-\St_a)$, $(\Se_b-\St_b)$, $(\So_a-\St_a)$ and $(\So_b-\St_b)$, respectively. Since we can bound
\begin{equation}
        \begin{aligned}
        	& 
        	\| L_\J(\Se-\St) \|_{\Fr}^2 =\\
        	&\qquad
        	= \left\| \sum\limits_{r\in\J} \sum\limits_{s\notin\J} \frac{\Ptr(\Se-\St)\Pts + \Pts(\Se-\St)\Ptr}{\mu_r-\mu_s} \right\|_{\Fr}^2 = 2\sum\limits_{r\in\J} \sum\limits_{s\notin\J} \frac{\| \Ptr(\Se-\St)\Pts\|_{\Fr}^2}{(\mu_r-\mu_s)^2}\\
        	&\qquad
        	\leq 2\left( \sum\limits_{r\in\J} \sum\limits_{s\notin\J} \frac{m_r\mu_r m_s \mu_s}{(\mu_r-\mu_s)^2} \right) \left( \max\limits_{r\in\J, s\notin\J} \frac{\| \Ptr (\Se-\St) \Pts\|_{\Fr}}{\sqrt{m_r\mu_r m_s\mu_s}}\right)^2
        	\leq 2 \rre_\J(\St) x^2,
        \nonumber
        \end{aligned}
    \end{equation}
    similar bounds apply to $L_a(\Se_a-\St_a)$, $L_b(\Se_b-\St_b)$, $L_a(\So_a-\St_a)$ and $L_b(\So_b-\St_b)$ and it holds
    \begin{equation}
        \begin{aligned}
        	&\|L_a(\Se_a-\St_a)\|_{\Fr} \leq \sqrt{2} \rre_{\J_a}(\St_a)^{1/2} \widehat{x}^a,\\
        	&\|L_b(\Se_b-\St_b)\|_{\Fr} \leq \sqrt{2} \rre_{\J_b}(\St_b)^{1/2} \widehat{x}^b,\\
        	&\|L_a(\So_a-\St_a)\|_{\Fr} \leq \sqrt{2} \rre_{\J_a}(\St_a)^{1/2} \overline{x}^a,\\
        	&\|L_b(\So_b-\St_b)\|_{\Fr} \leq \sqrt{2} \rre_{\J_b}(\St_b)^{1/2} \overline{x}^b.
        \label{L_bound}
        \end{aligned}
    \end{equation}
\\
\underline{Main part:} 
Denote for shortness for the rest of the proof
    \begin{equation}
        \begin{aligned}
        	A \eqdef L_a(\Se_a-\St_a) -
        	L_b(\Se_b-\St_b).
        \nonumber
        \end{aligned}
    \end{equation}
Let us bound $\left| \| \Pea-\Peb\|_{(\Pbar, \overline{\Gamma}, s_1,s_2)} - \| A\|_{(\Pstar, \Gamma^*,s_1,s_2)} \right|$. First, if $\widehat{x}^a \leq \scca$ and $\widehat{x}^b \leq \sccb$, then
    \begin{equation}
        \begin{aligned}
        	&\left| \| \Pea-\Peb\|_{(\Pbar, \overline{\Gamma},s_1,s_2)} - \| A\|_{(\Pbar, \overline{\Gamma},s_1,s_2)} \right| 
        	\leq
        	\| R_a(\Se_a-\St_a) \|_{(\Pbar,\overline{\Gamma},s_1,s_2)} + \| R_b(\Se_b-\St_b) \|_{(\Pbar,\overline{\Gamma},s_1,s_2)}
        	\\&\qquad\leq2\| R_a(\Se_a-\St_a) \| +  2\| R_b(\Se_b-\St_b) \|
        	\leq 2C(\widehat{x}^a)^2\rr_{\J_a}(\St_a)^{3/2} + 2C(\widehat{x}^b)^2\rr_{\J_b}(\St_b)^{3/2}
        	\\&\qquad\leq 2C\sccab^2 \rr_{a,b}^{3/2},
        \nonumber
        \end{aligned}
    \end{equation}
where we used Proposition~\ref{P:Properties0} (i), (ii) and Lemma~\ref{L:KL} (Condition~\ref{A:JW} is fulfilled by Assumption~\ref{A:additional} (i)).
Next, we bound $\left| \| A \|_{(\Pbar,\overline{\Gamma},s_1,s_2)} - \| A \|_{(\Pstar, \Gamma^*,s_1,s_2)}\right|$.
By Definition~\ref{D:norm}, it is clear that
    \begin{equation}
        \begin{aligned}
        	\left| \|A\|_{(\Pbar, \overline{\Gamma},s_1,s_2)} - \| A \|_{(\Pstar, \Gamma^*,s_1,s_2)} \right|
        	&\leq
        	\left| \|\overline{\Gamma}_1^\T A \overline{\Gamma}_1 \| - \| {\Gamma_1^*}^\T A \Gamma_1^* \| \right|
        	+
        	\left| \|\overline{\Gamma}_2^\T A \overline{\Gamma}_2 \| - \| {\Gamma_2^*}^\T A \Gamma_2^* \| \right| +
        	\\&\qquad +\left| \sup\limits_{\substack{v\in\myset^m_{s_1} \\w\in\myset^{d-m}_{s_2} }} v^\T \overline{\Gamma}_1^\T A \overline{\Gamma}_2 w- \sup\limits_{\substack{v\in\myset^m_{s_1} \\w\in\myset^{d-m}_{s_2} }} v^\T {\Gamma_1^*}^\T A {\Gamma_2^*} w \right|.
        \nonumber
        \end{aligned}
    \end{equation}
Each of the three terms can be bounded similarly, so let us bound just the first one. Adding and subtracting the mixed term $\|\overline{\Gamma}_1^\T A \Gamma^*_1 \|$, we obtain
\begin{equation}
        \begin{aligned}
        	&\left| \|\overline{\Gamma}_1^\T A \overline{\Gamma}_1 \| - \| {\Gamma_1^*}^\T A \Gamma_1^* \| \right| \leq
        	\left| \|\overline{\Gamma}_1^\T A \overline{\Gamma}_1 \| - \| \overline{\Gamma}_1^\T A \Gamma_1^* \| \right|
        	+ \left| \|\overline{\Gamma}_1^\T A \Gamma^*_1 \| - \| {\Gamma_1^*}^\T A \Gamma_1^* \| \right|
        	\\&\qquad\leq
        	\| \overline{\Gamma}_1^\T A (\overline{\Gamma}_1 - \Gamma^*_1) \| + \| (\overline{\Gamma}_1 - \Gamma^*_1)^\T A \Gamma^*_1 \|
        	\leq
        	 2\|A\| \| \overline{\Gamma}_1 - \Gamma^*_1 \| 
        	 \leq 2^{5/2} \|A\|\| \Pbar - \Pstar \|_{\Fr}.
        \nonumber
        \end{aligned}
    \end{equation}
    So, if $\widehat{x}^a \leq \scca$ and $\widehat{x}^b \leq \sccb$, then
    \begin{equation}
        \begin{aligned}
        	&\left| \|A\|_{(\Pbar, \overline{\Gamma},s_1,s_2)} - \| A \|_{(\Pstar, \Gamma^*,s_1,s_2)} \right|
        	 \leq 3\cdot2^{5/2} \|A\| \| \Pbar - \Pstar \|_{\Fr}
        	 \\&\qquad
        	 \leq 3\cdot2^{5/2} \left( \|L_a(\Se_a-\St_a)\| + \|L_b(\Se_b-\St_b)\|_{\Fr} \right)\| \Pbar - \Pstar \|_{\Fr}
        	 \\&\qquad \leq 3\cdot2^{5/2} \left( \widehat{x}^a \rre_{\J_a}(\St_a)^{1/2} + \widehat{x}^b \rre_{\J_b}(\St_b)^{1/2}\right) \cdot \| \Pbar - \Pstar \|_{\Fr}
        	 \\&\qquad \leq 3\cdot2^{5/2} \sccab \rre_{a,b}^{1/2}\cdot \| \Pbar - \Pstar \|_{\Fr}.
        \nonumber
        \end{aligned}
    \end{equation}
    Hence, if $\widehat{x}^a \leq \scca$ and $\widehat{x}^b \leq \sccb$, then
    \begin{equation}
        \begin{aligned}
        	&\left| \| \Pea-\Peb\|_{(\Pbar, \overline{\Gamma},s_1,s_2)} - \| A\|_{(\Pstar, \Gamma^*,s_1,s_2)} \right| 
        	\leq 2C\sccab^2 \rr_{a,b}^{3/2} + 3\cdot2^{5/2} \sccab \rre_{a,b}^{1/2}\cdot \| \Pbar - \Pstar \|_{\Fr},
        \nonumber
        \end{aligned}
    \end{equation}
    which implies, introducing event $\Omega = \{\widehat{x}^a \leq \scca, \widehat{x}^b \leq \sccb \}$ with $\Prob[\Omega] \geq 1-1/n_a-1/n_b$ (by Lemma~\ref{L:KC_SCconc} and union bound),
    \begin{equation}
        \begin{aligned}
        	&\Prob\left[\left| \| \Pea-\Peb\|_{(\Pbar, \overline{\Gamma},s_1,s_2)} - \| A\|_{(\Pstar, \Gamma^*,s_1,s_2)} \right| 
        	> C\left(\sccab^2 \rr_{a,b}^{3/2} + \sccab \rre_{a,b}^{1/2}\cdot \| \Pbar - \Pstar \|_{\Fr}\right) \,\Big|\, \overline{\Gamma} \right]=
        	\\&\hspace{0.5cm}=
        	\Prob\left[\left| \| \Pea-\Peb\|_{(\Pbar, \overline{\Gamma},s_1,s_2)} - \| A\|_{(\Pstar, \Gamma^*,s_1,s_2)} \right| 
        	> C\left( \sccab^2 \rr_{a,b}^{3/2} + \sccab \rre_{a,b}^{1/2}\cdot \| \Pbar - \Pstar \|_{\Fr} \right) \,\Big|\, \overline{\Gamma}; \Omega \right] \times
        	\\&\qquad\qquad \times\Prob[\Omega] +
        	\\&\qquad+
        	\Prob\left[\left| \| \Pea-\Peb\|_{(\Pbar, \overline{\Gamma},s_1,s_2)} - \| A\|_{(\Pstar, \Gamma^*,s_1,s_2)} \right| 
        	> C\left( \sccab^2 \rr_{a,b}^{3/2} + \sccab \rre_{a,b}^{1/2}\cdot \| \Pbar - \Pstar \|_{\Fr} \right) \,\Big|\, \overline{\Gamma}; \Omega^c \right] \times
        	\\&\qquad\qquad \times \Prob[\Omega^c]
        	\\&\hspace{0.5cm} \leq0\cdot 1 + 1\cdot \left(\frac{1}{n_a} + \frac{1}{n_b}\right) = \frac{1}{n_a} + \frac{1}{n_b}.
        \nonumber
        \end{aligned}
    \end{equation}
    Now it is left to bound $\| \Pbar - \Pstar \|_{\Fr}$ with high probability.
        
    By definition of $\Pbar$ given by \eqref{Def:Pbar}
    \begin{equation}
        \begin{aligned}
        	\| \Pbar - \Poa \|_{\Fr}^2 + \| \Pbar - \Pob \|_{\Fr}^2
        	\leq
        	\| \Pstar - \Poa \|_{\Fr}^2 + \| \Pstar - \Pob \|_{\Fr}^2
        	=
        	\| \Poa - \Pta \|_{\Fr}^2 + \| \Pob - \Ptb \|_{\Fr}^2.
        \nonumber
        \end{aligned}
    \end{equation}
Therefore,
    \begin{equation}
        \begin{aligned}
        	&\| \Pbar - \Pstar \|_{\Fr}^2 
        	=
        	\frac{1}{2} \left( \| \Pbar - \Pta \|_{\Fr}^2 + \| \Pbar - \Ptb \|_{\Fr}^2 \right)
        	\\&\qquad\leq
        	 \| \Pbar - \Poa \|_{\Fr}^2 + \| \Pbar - \Pob \|_{\Fr}^2 + \| \Poa - \Pta \|_{\Fr}^2 + \| \Pob - \Ptb \|_{\Fr}^2 
        	\\&\qquad\leq 2\left( \| \Poa - \Pta \|_{\Fr}^2 + \| \Pob - \Ptb \|_{\Fr}^2 \right).
        \nonumber
        \end{aligned}
    \end{equation}
Hence, if $\overline{x}^a \leq \scca$ and $\overline{x}^b \leq \sccb$, then
    \begin{equation}
        \begin{aligned}
        	&\| \Pbar - \Pstar \|_{\Fr}
        	\leq \sqrt{2} \left( \| \Poa - \Pta \|_{\Fr} + \| \Pob - \Ptb \|_{\Fr} \right)
        	\\&\qquad\leq \sqrt{2} \left( \|L_a(\So_a-\St_a)\|_{\Fr} + \| R_a(\So_a-\St_a) \|_{\Fr} + \|L_b(\So_b-\St_b)\|_{\Fr} + \| R_b(\So_b-\St_b) \|_{\Fr}\right)
        	\\&\qquad\leq C \left( \sccab \left[ \rre_{\J_a}(\St_a)^{1/2} + \rre_{\J_b}(\St_b)^{1/2}\right] + \sccab^2 \left[ \rr_{\J_a}(\St_a)^{3/2} + \rr_{\J_b}(\St_b)^{3/2} \right]\right)
            \\&\qquad= C\left( \sccab \rre_{a,b}^{1/2} + \sccab^2 \rr_{a,b}^{3/2}\right),
        \nonumber
        \end{aligned}
    \end{equation}
where we used bounds \eqref{L_bound} and Lemma~\ref{L:KL} (again, Condition~\ref{A:JW} is fulfilled by Assumption~\ref{A:additional} (i)).
Since, probability of the event $\overline{x}^a \leq \scca$ and $\overline{x}^b \leq \sccb$ is at least $1-1/n_a-1/n_b$ (again by Lemma~\ref{L:KC_SCconc} and union bound), we conclude (adjusting the constants and using technical assumption to simplify the bound):
\begin{equation}
        \begin{aligned}
        	&\Prob\left[\left| \| \Pea-\Peb\|_{(\Pbar, \overline{\Gamma},s_1,s_2)} - \| A\|_{(\Pstar, \Gamma^*,s_1,s_2)} \right| 
        	> C\sccab^2\left( \rr_{a,b}^{3/2} + \rre_{a,b} \right) \,\Big|\, \overline{\Gamma} \right]\leq
        	\frac{1}{n_a} + \frac{1}{n_b}
        \nonumber
        \end{aligned}
    \end{equation}
    with probability $1-1/n_a-1/n_b$.
\end{proof}

\begin{proof}[Proof of Proposition~\ref{FM1}]
    Let us first prove that the spectral projector onto the sum of $m$ eigenspaces corresponding to non-zero eigenvalues of $\BB\Cov[\ff_1]\BB^\T$ is given by $\BB(\BB^\T \BB)^{-1}\BB^\T$. Consider the eigendecomposition of $(\BB^\T \BB)^{1/2} \Cov[\ff_1] (\BB^\T \BB)^{1/2}$:
    \begin{equation}
        \begin{aligned}
        	(\BB^\T \BB)^{1/2} \Cov[\ff_1] (\BB^\T \BB)^{1/2} = Q D Q^\T,
        \nonumber
        \end{aligned}
    \end{equation}
    where $Q\in\R^{m\times m}$ is orthogonal and $D\in\R^{m\times m}$ is diagonal. Take $\HH = (\BB^\T \BB)^{-1/2} Q$.
    Then
    \begin{equation}
        \begin{aligned}
            (\BB\HH)^\T \BB\HH = \HH^\T \BB^\T \BB \HH = Q^\T (\BB^\T \BB)^{-1/2} \BB^\T \BB (\BB^\T \BB)^{-1/2} Q = Q^\T Q = \Id_m,
        \nonumber
        \end{aligned}
    \end{equation}
    i.e. columns of $\BB\HH$ are orthogonal and have unit length. Also,
    \begin{equation}
        \begin{aligned}
            \HH^{-1} \Cov[\ff_1] (\HH^{-1})^\T =  Q^\T (\BB^\T \BB)^{1/2} \Cov[\ff_1] (\BB^\T \BB)^{1/2} Q = Q^\T Q D Q^\T Q = D,
        \nonumber
        \end{aligned}
    \end{equation}
    i.e. diagonal. Therefore,
     \begin{equation}
        \begin{aligned}
            (\BB\HH) \left[ \HH^{-1} \Cov[\ff_1] (\HH^{-1})^\T \right] (\BB\HH)^\T
        \nonumber
        \end{aligned}
    \end{equation}
    is a valid eigendecomposition of $\BB\Cov[\ff_1] \BB^\T$. The spectral projector of interest is then exactly
    \begin{equation}
        \begin{aligned}
            (\BB\HH) (\BB\HH)^\T = \BB \HH \HH^\T \BB^\T = \BB (\BB^\T \BB)^{-1/2} QQ^\T (\BB^\T \BB)^{-1/2} \BB^\T = \BB(\BB^\T \BB)^{-1} \BB^\T.
        \nonumber
        \end{aligned}
    \end{equation}
    
    Now we just apply Davis-Kahan theorem to $\St$ and $\BB \Cov[\ff_1] \BB^\T$ to get
    \begin{equation}
        \begin{aligned}
            \| \Pt - \BB(\BB^\T \BB)^{-1} \BB^\T\| \lesssim \frac{\| \Cov[\ixi_1] \|}{\mu_m - \mu_{m+1}} = O\left( \frac{1}{d} \right),
        \nonumber
        \end{aligned}
    \end{equation}
    where we used Assumption~\ref{FM_assumption}.
\end{proof}

\begin{proof}[Proof of Proposition~\ref{FM2}]
    The condition $\Cov[\ixi_1] \BB = 0_{d\times m}$ implies that any projector of $\Cov[\ixi_1]$ corresponding to non-zero eigenvalue is orthogonal to $\BB(\BB^\T \BB)^{-1} \BB^\T$. This means that the projectors of $\St$ and $\BB\Cov[\ff_1]\BB^\T$ onto the first $m$ eigenspaces coincide.
\end{proof}